\documentclass[reqno]{amsart}

\usepackage{graphicx}
 \usepackage{bbm} 
  \usepackage{bm}
\usepackage{amsfonts,}
\usepackage{amsbsy}
\usepackage{amssymb}
 \usepackage[]{amsmath}
\usepackage{graphicx}
\usepackage{oldgerm}
\usepackage{yhmath}
\usepackage{wrapfig}
 \usepackage{dsfont}
 \usepackage{cite}
 \usepackage{mathabx} 
\usepackage{color}
\usepackage[dvipsnames]{xcolor}

   \usepackage[bookmarks=false,colorlinks,linkcolor=blue,citecolor=blue]{hyperref}

\theoremstyle{plain}
\newtheorem{theorem}{Theorem}
\newtheorem{theorem*}{Theorem}
\newtheorem{lemma}{Lemma}
\newtheorem{proposition}{Proposition}
\newtheorem{corollary}{Corollary}

\theoremstyle{definition}

\theoremstyle{plain}
\newtheorem{remark}{Remark}
\newtheorem*{remark*}{Remark}
\numberwithin{equation}{section}

 \newcommand{\lb}{\left|}
 \newcommand{\rb}{\right|}
 \newcommand{\lc}{\left[}
 \newcommand{\rc}{\right]}
 \newcommand{\lp}{\left(}
 \newcommand{\rp}{\right)}
 \newcommand{\la}{\left\{}
 \newcommand{\ra}{\right\}}
 \newcommand{\lpt}{\left.}
  \newcommand{\rpt}{\right.}
 \newcommand{\e}{\varepsilon}
  
  \newcommand{\lime}{\lim_{\e\to 0}} 
 \newcommand{\dto}{\rightharpoonup \!\!\!\!\!\!\!\!\rightharpoonup}
 \newcommand{\mdto}{{\,\buildrel m_\e  \over \dto\,}}

\newcommand{\bfCap}{{\rm \bf Cap }}
\newcommand{\Capa}{{\rm Cap }}
\newcommand{\capsca}{{\rm cap}}
  
\newcommand{\bfnabla}{\pmb{  \nabla}} 

\newcommand{\bfdiv}{{\rm \bf div}} 
\newcommand{\dist}{{\rm dist}}
\newcommand{\diam}{{\rm diam}} 
\newcommand{\spt}{{\rm supp\,}}

\newcommand{\intb}{{\int\!\!\!\!\!\!-}}
\newcommand{\ov}{\overline}

 \newcommand{\bfa}{\pmb{a}} 
 \newcommand{\bfb}{\pmb{b}}

 \newcommand{\bfe}{\pmb{e}}  
 \newcommand{\bff}{\pmb{f}}  
 \newcommand{\bfg}{\pmb{g}}  
 \newcommand{\bfh}{\pmb{h}}  
 \newcommand{\bfi}{\pmb{i}}

 \newcommand{\bfl}{\pmb{l}} 
    
 \newcommand{\bfn}{\pmb{n}}  
   
 \newcommand{\bfp}{\pmb{p}}  
 \newcommand{\bfq}{\pmb{q}}  
   
 \newcommand{\bfs}{\pmb{s}}  
 \newcommand{\bft}{\pmb{t}}  
 \newcommand{\bfu}{\pmb{u}}  
 \newcommand{\bfv}{\pmb{v}}  
 \newcommand{\bfw}{\pmb{w}}  
 \newcommand{\bfx}{\pmb{x}}   
 \newcommand{\bfy}{\pmb{y}}  
 \newcommand{\bfz}{\pmb{z}}  

\newcommand{\bfA}{\pmb{A}} 
\newcommand{\bfB}{\pmb{B}} 
\newcommand{\bfC}{\pmb{C}}

\newcommand{\bfI}{\pmb{I}}

 \newcommand{\bfM}{\pmb{M}}

\newcommand{\B}{{\mathcal B}}
\newcommand{\C}{{\mathcal C}}
\newcommand{\D}{{\mathcal D}}

\newcommand{\F}{{\mathcal F}}

\newcommand{\calH}{{\mathcal H}}
\newcommand{\I}{{\mathcal I}}

\newcommand{\K}{{\mathcal K}}
\newcommand{\calL}{{\mathcal L}}
\newcommand{\M}{{\mathcal M}}

\newcommand{\calP}{{\mathcal P}}

\newcommand{\R}{{\mathcal R}} 

\newcommand{\T}{{\mathcal T}}

\newcommand{\W}{{\mathcal W}}

\newcommand{\NN}{\mathbb{N}}

\newcommand{\RR}{\mathbb{R}}  
 \newcommand{\SSym}{\mathbb{S}}

\newcommand{\ZZ}{\mathbb{Z}}

\newcommand{\bfalpha}{{\pmb{ \alpha}}}

\newcommand{\bfeta}{{\pmb{  \eta}}}

\newcommand{\bfnu}{{\pmb{ \nu}}} 
\newcommand{\bfxi}{{\pmb{ \xi}}}

\newcommand{\bfphi}{{\pmb{  \phi}}} 
\newcommand{\bfchi}{{\pmb{  \chi}}} 
\newcommand{\bfpsi}{{\pmb{  \psi}}}

\newcommand{\bfvarphi}{{\pmb{ \varphi}}}

\newcommand{\bfXi}{\pmb{ \Xi}}

\newcommand{\bfUpsilon}{\pmb{ \Upsilon}}

\newcommand{\bfPsi}{\pmb{  \Psi}}

\newcommand{\gots}{{{\mathfrak{s}}}}

\newcommand{\gotv}{{{\mathfrak{v}}}}
\newcommand{\gotw}{{{\mathfrak{w}}}}

\newcommand{\bfgotp}{{{\textswab{\pmb p}}}}

\newcommand{\bfgotr}{{{\textswab{\pmb r}}}}

 \newcommand{\bfgotu}{{{\textswab{\pmb u}}}}
\newcommand{\bfgotv}{{{\textswab{\pmb v}}}}
\newcommand{\bfgotw}{{{\textswab{\pmb w}}}}

\newcommand{\gotN}{{{\mathfrak{N}}}}
\newcommand{\gotO}{{{\mathfrak{O}}}}

\newcommand{\bfgotD}{{\bm{\textgoth{D}}}}

\newcommand{\bfgotI}{{\bm{\textgoth{I}}}}
\newcommand{\bfgotJ}{{\bm{\textgoth{J}}}}

\newcommand{\bfgotL}{{\bm{\textgoth{L}}}}

\newcommand{\bfgotW}{{\bm{\textgoth{W}}}}

\begin{document}

\bibliographystyle{plain}

\title[]{Homogenization of  Norton-Hoff  fibered composites with high contrast
}


\author{Michel Bellieud}
\address{LMGC (Laboratoire de M\'ecanique et de G\'enie Civil de Montpellier)
UMR-CNRS 5508, Universit\'e Montpellier II, Case courier
048, Place Eug\`ene Bataillon, 34095 Montpellier Cedex 5, France}
\email{michel.bellieud@umontpellier.fr}

\date{}


\begin{abstract}
 We study the steady creep flow of a  
perfectly viscoplastic 
solid  reinforced by 
fibers with high viscosity contrast. 
Our study unveils  new effects  
related to  anisotropy and conditioned by    the Norton exponent. 
\end{abstract}

\maketitle

\tableofcontents 

\section{Introduction}
\label{intro}

Given   a    bounded Lipschitz   cylindrical    domain  $\Omega:=\Omega'\times(0,L) $  of $\RR^3$ and   two   strictly convex   functions  $f$,  $g$  satisfying  a  growth condition  of order $p\in (1,+\infty)$, we consider the  problem 
     \begin{equation}
 (\calP_\e):   \la
\begin{aligned}
& \inf_{\bfu_\e\in W^{1,p}_{b} (\Omega;\RR^3) } F_\e(\bfu_\e)-\int_\Omega \bff. \bfu_\e dx , 
 \\&
  F_\e(\bfu_\e):= \int_{\Omega\setminus {T_{r_\e}}}  f(\bfe(\bfu_\e)) dx +   
l_\e
    \int_{{T_{r_\e}}}  g(\bfe(\bfu_\e)) dx , 
\\&
  \bfe(\bfu_\e)=\frac{1}{2} (\bfnabla \bfu_\e +\bfnabla^T\bfu_\e),   \quad \bff\in L^{p'}(\Omega;\RR^3), 
 \quad   \frac{1}{p}+\frac{1}{p'}=1 ,
  \\& W^{1,p}_{b} (\Omega;\RR^3):= \la \bfpsi\in W^{1,p} (\Omega;\RR^3), \ \bfpsi=0 \ \text{ on } \Omega'\times\{0\} \ra ,
 \end{aligned} 
\right.
 \label{Pe}
  \end{equation}
where $T_{r_\e}$  is 
a  distribution  of    disjoint    cylinders     of very small volume fraction
and $l_\e$ 
 a large parameter.
  The solution $\bfu_\e$   to \eqref{Pe}  represents  
    the   
Eulerian velocity field in a 
Norton-Hoff material  
of Norton exponent $\frac{1}{p-1}$  undergoing  a steady 
creep flow
 under the influence of  a density $\bff$ of
applied body forces  \cite{Ho,LeCh,No}.
Composites  comprising   a    small volume fraction   of fibers with strong  properties    have been studied in various   contexts
\cite{BeJCA,BeBo,BeLiOr,BoFe,BrianeStokes,CaDi,FeKh,Kh,LeCh,MaKh,PiSe}.
They are characterized  by    the   interaction    of   concentration phenomena     
in  the fibers and in  a small region  of space surrounding them.
The equilibrium problem  in linear elasticity,  a   special case  of  Problem \ref{Pe} corresponding    to positive definite quadratic forms $f$ and $g$,
  has  been    investigated in the isotropic case for fibers of circular cross-sections 
in  \cite{BeGr}.
 
  Our study unveils  new effects  
related to  anisotropy and conditioned by      the  growth parameter   $p$ and the shape  of the cross-sections of the fibers,  described by some bounded  connected  open subset $S$ of $\RR^2$. 
The distinctive  feature of Problem \ref{Pe} lies in  the interaction of  concentrations of stress in the fibers and  rate of deformation  in their close  outer  neighborhood.
For ease of exposition, we assume that the fibers are $\e$-periodically distributed. Our  analysis goes through in  the non-periodic case   provided the fibers are well separated. 
We show that the asymptotic  behavior  of  $\bfu_\e$   depends on the order of magnitude of  $l_\e$,  characterized 
  in terms of the size $r_\e$ of the cross-sections of the fibers
     by  the parameters 
\begin{equation}
\begin{aligned}
 & k:= \lime   k_\e    \in [0,+\infty];\quad  \kappa:=\lime r_\e^p   k_\e   \in [0,+\infty];
 \quad  k_\e:=  l_\e\frac{r_\e^2|S|}{\e^2},
 \end{aligned} 
\label{kkappa}
   \end{equation}
and    on the effective  $p$-capacity   of their  cross-sections in $\Omega'$,  represented
by
      \begin{equation}
\begin{aligned}
& \gamma^{(p)}:=  \lim_{\e\to0} \gamma_\e^{(p)}(r_\e)\in [0,+\infty];
\quad   \gamma_\e^{(p)}(r)  :=  \begin{cases}
   \frac{r^{2-p}}{\e^2}   &  \text{ if } \  p\not = 2, \quad 
\\   \frac{1}{\e^2|\log r|}  &   \text{ if } \  p  =   2.
\end{cases}
\end{aligned}
\label{defgammae}
\end{equation} 
When  $0<\kappa<+\infty$, we establish  that  the      fibers  locally behave like   rigid bodies  rotating around their principal axes (parallel to $\bfe_3$)
with  an  angular  velocity $\delta$,
   the axes moving at the velocity  $\bfv$.
   The field $\bfv$ approximates as $\e\to0$ to the   local  average  $\bfv_\e$  of  $\bfu_\e$   over  the cross-sections of the  fibers, and $\delta$
to  their  mean angular velocity $\delta_\e$. 
We   show that  $\frac{v_{\e3}}{r_\e}$ converges to some function $w$  and   demonstrate  that
   the  effective contribution of the fibers 
 is described by a functional of $\bfv$, $\delta$, and $w$
    of 
   the  form
   \begin{equation}
\begin{aligned}
\Phi_{fibers}(\bfv,0,\delta,w)= \int_\Omega g^{hom}\lp \frac{\partial^2v_1}{\partial^2x_3},\frac{\partial^2v_2}{\partial^2x_3},
\frac{\partial w}{\partial x_3}, \frac{\partial \delta}{\partial x_3}\rp dx.
 \end{aligned} 
\label{Phistiffer}
  \end{equation}
If  $0<k<+\infty$, we prove   that  the fibers  possibly display a larger  angular velocity,   characterized 
by a function  $\theta$ approximating to  $\theta_\e:=r_\e \delta_\e $. Tangential and normal velocities 
 are  then   of the same order at the surface of the fibers. 
We show  that 
the  concentration of stress  in  the fibers leads  to   a contribution of  the type
\begin{equation}
\begin{aligned}
\Phi_{fibers}(\bfv,\theta)= \int_\Omega g^{hom}\lp \frac{\partial v_3}{\partial x_3}, \frac{\partial \theta}{\partial x_3}\rp dx.
 \end{aligned} 
\label{Phistiff}
  \end{equation}
  In the other cases, we establish that  $\Phi_{fibers}$ vanishes. Setting
           \begin{equation}
   \begin{aligned}
&  \bfv^{tuple}:= (  \bfv, \theta)\quad   \text{if }  \kappa=0;
\qquad 
\bfv^{tuple}:= ( \bfv, \theta, w,\delta) \quad   \text{if } 0<\kappa\leq +\infty,
\end{aligned} 
 \label{tuple}
  \end{equation}
 we demonstrate  that 
the  limit problem associated with \eqref{Pe} takes the form 
 \begin{equation}
\begin{aligned}
&\inf_{(\bfu,\bfv^{tuple})\in W^{1,p}_b(\Omega;\RR^3)\times \D} \Phi(\bfu,\bfv^{tuple})-\int_\Omega \bff\cdot\bfu,
\\&\Phi(\bfu,\bfv^{tuple})=\int_\Omega f(\bfe(\bfu)) dx +  \Phi_{cap}(\bfv-\bfu,\theta) +\ \Phi_{fibers}(\bfv^{tuple}),
 \end{aligned} 
\label{Phihomisot}
  \end{equation}
  for some suitable domain $\D$.

The second term of $\Phi$ in \eqref{Phihomisot} stems from    a concentration 
 of rate of deformations  
   in  the close  outer neighborhood of the fibers, resulting from the   interaction  between the matrix and the fibers. 
   We prove  that 
  \begin{equation}
\begin{aligned}
\Phi_{cap}(\bfv-\bfu,\theta) = \int_\Omega 
c^f(\bfv-\bfu, \theta) dx, 
 \end{aligned} 
\label{cf}
  \end{equation}
for some   convex function    $c^f$
satisfying,  if $p\not=2$, the growth condition 
 \begin{equation}
\begin{aligned}
&c\gamma^{(p)}(|\bfv-\bfu|^p + |\theta|^p)  \le  c^f(\bfv-\bfu,\theta)\le C \gamma^{(p)}(|\bfv-\bfu|^p + |\theta|^p),
\ \end{aligned}
\label{estimcfp<2}
    \end{equation}
 with  $c,C>0$. 
A 
  phenomenon 
  related to  the  Stokes'   paradox
induces a different behavior when $p=2$:  we   show  that  then,
     \begin{equation} 
\begin{aligned}
& c^f(\bfv-\bfu,\theta)=+\infty  \qquad \qquad & &  \hbox{if }  \theta\not=0,
\\&c\gamma^{(2)} |\bfv-\bfu|^2  \le  c^f(\bfv-\bfu,0)\le C\gamma^{(2)}|\bfv-\bfu|^2.
 \end{aligned}
\label{estimcfp=2}
    \end{equation}  
This  means  that  the interaction between the  matrix and the fibers   precludes  large  angular velocities 
   of the fibers.
If  $\gamma^{(p)}=+\infty$, 
 which always holds when $p>2$, 
  the limit problem is simply obtained by setting   $\bfu=\bfv$ 
  and $\theta=0$  in \eqref{Phihomisot}.

We turn to a more detailed description of the mappings $g^{hom}$ and $c^f$, using  standard notations recalled in the next section. 
We establish  that the  functions  
 $ g^{hom}$ in \eqref{Phistiffer} 
  and  \eqref{Phistiff} 
are, up to a multiplicative constant,  the infima with respect
   to $\bfq$ over $W^{1,p}(S;\RR^3)$ of   
    \begin{equation}
\begin{aligned}
 \int_{S} g^{0,p}  \lp   \bfe_{y}( \bfq)+
\frac{2}{\diam S}\frac{\partial \delta}{\partial x_3}\right.
(-y_2 &\bfe_1\odot\bfe_3  +y_1\bfe_2\odot\bfe_3)
\\&\lpt+\lp  \frac{\partial w}{\partial x_3}  - \sum_{\alpha=1}^2  \frac{\partial^2 v_\alpha}{\partial x_3^3}  y_\alpha \rp\bfe_3\otimes\bfe_3 \rp dy ,
 \end{aligned} 
\label{ghom2intro}
  \end{equation}
 and
 \begin{equation}
\begin{aligned}
   \int_{S} g^{0,p}  \lp   \bfe_{y}( \bfq)+
\frac{2}{\diam S}\frac{\partial \theta}{\partial x_3}
(-y_2 \bfe_1\odot\bfe_3 +y_1\bfe_2\odot\bfe_3)
+\frac{\partial v_3}{\partial x_3} \bfe_3\otimes\bfe_3 \rp dy,
 \end{aligned} 
\label{ghom1intro}
  \end{equation}
respectively,   where $g^{0,p}$  denotes a  $p$-homogeneous approximation of $g$ near $0$.
On the other hand, assuming without loss of generality that $\int_S \bfy dy =0$, we  prove that
 for every   $p\in (1,+\infty)$,
\begin{equation}
\begin{aligned}
&c^f(\bfv-\bfu,\theta)
 = 
\lim_{\e\to0} \frac{1}{\e^2} \capsca^f(\bfv-\bfu,\theta\bfe_3; r_\e S,R_\e D),
 \end{aligned}
\label{defcfintro}
    \end{equation}
for some $R_\e\gg r_\e$, where  $D$ is the unit ball of $\RR^2$. The mapping  $\capsca^f$ is   the  variant    of  the notion of capacity introduced in \cite{BeArma} (see also \cite{Vi}) defined,  for  any   $(\bfa,\bfalpha)\in (\RR^3)^2$ and any couple  $(U,V)$  of  open subsets   of $\RR^2$
with    $U$  connected,  bounded, $\int_U \bfy dy=0$, and 
    $\ov U \subset V$,   by
  \begin{equation}
\hskip-0,1cm\begin{aligned}
&\capsca^f(\bfa,\bfalpha;U,V)  := \inf \calP^f(\bfa,\bfalpha;U,V), 
\\&
 \calP^f(\bfa,\bfalpha;U,V): \inf_{\bfpsi\in \W^p(\bfa,\bfalpha;U,V)  }   \int_V   f(  \bfe_y(\bfpsi)) dy,  
 \\& \W^p(\bfa,\bfalpha;U,V)  := \Big\{\bfpsi\in W^{1,p}_0 (V;\RR^3)  ,  
     \bfpsi (y)=  \bfa+  \tfrac{2}{\diam U}\bfalpha\wedge \bfy  \hbox{   in }  U\Big\}.
 \end{aligned}
\label{infW}
    \end{equation}
The extended real  $c^f(\bfa,\zeta)$  can be seen as a capacity density
  approximating   to the sum of the images of the sections of the fibers under 
$\capsca^f(\bfa,\zeta\bfe_3; .,\Omega')$ per unit surface.
Our   investigation into 
 the properties of   $\capsca^f$ leads to    the   formula 
\begin{equation}
\begin{aligned}
&c^f(\bfv-\bfu,\theta)
 = 
  \gamma^{(p)} 
\capsca^{f^{\infty,p}}(\bfv-\bfu,\theta \bfe_3; S,\RR^2) \quad
  \hbox{ if } \ p<2,
 \end{aligned}
\label{cfp<2}
    \end{equation}
    where $f^{\infty,p}$ is a $p$-homogeneous approximation of $f$ at $\infty$.   If 
   $p=2$,  we  show  that   $c^f$ is independent of the shape of the cross-sections of the fibers.
   
   An interesting  feature  of our results lies in the dependence of the limit problem  on  the effective rescaled angular velocity  $\theta$.
A non-vanishing  $\theta$  may only    arise  when $\gamma^{(p)}<+\infty$, $p<2$, and $0<k<+\infty$. 
This explains why this dependence was  not detected   in  \cite{BeGr}.
It  arose in another context in \cite{BeSiam}.
One can guess, from \eqref{infW} and \eqref{cfp<2}, that $\theta$ is conditioned 
 by  the shape of the  cross-sections of  the fibers:  
the matrix  is more likely to induce a rotating motion in fibers of pear-shaped cross-sections  than circular ones. 
Formulae  \eqref{ghom1intro}, \eqref{infW} and \eqref{cfp<2}    suggest  that 
such  rotations   can also   be    brought about   by  the  anisotropy of either   matrix or fibers. 
Besides, we  prove  that 
$\theta$  can    be influenced   by
 large twisting  body forces applied on the fibers. 
Another distinctive aspect  of our work lies in 
the   dependence  of $\Phi_{fibers}$ on the effective angular velocity $\delta$ and   microscopic longitudinal velocity  $w$ when $0<\kappa<+\infty$. This  was not perceived  in  the setting of linear  isotropic elasticity
because, as  we show,    these functions vanish  when   the material constituting the fibers is linear  isotropic, whatever the shape of their cross-sections.
By contrast,  we prove that $\delta$   does not vanish,  
in general,  for  anisotropic fibers. The same is likely to hold for $w$. 
We  also  establish that, like for $\theta$,     large twisting  body forces  applied on the fibers may have an effect  on  $\delta$ and $w$ and, in particular,
 lead to a non-vanishing couple $(w,\delta)$   for isotropic fibers. 
Notice that,
unlike $\theta$,
the matrix exerts no influence on   $(\delta,w)$.

Taking advantage of the above study, we  next  examine  the problem 
  \begin{equation}
 \begin{aligned}
 &
  \inf_{\bfu_\e\in W^{1,p}_b(\Omega;\RR^3) } F_\e^{soft}(\bfu_\e) 
  -\int_\Omega \bff \cdot \bfu dx,
    \\&F_\e^{soft}(\bfu_\e) :=
   \e^p \int_{\Omega\setminus  {T_{\e}}}   f(\bfe(\bfu_\e)) dx +   
l_\e
    \int_{{T_{\e}}}  g(\bfe(\bfu_\e)) dx,
 \end{aligned} 
  \label{Pesoft}
  \end{equation}
when $T_{\e}$  is  an $\e$-periodic
   distribution  of    disjoint    cylinders     of  volume fraction of order $1$.
High-contrast homogenization problems of this type   have  been  and are  intensively investigated in many contexts
  \cite{Al,ArDoHo,AvGrMiRo,BaKaSm,BoFe2,BoBoFe,Ch,ChCo,ChLi,Co,KaSm,Pa,Sm,Zh1,Zh2,ZhPa}.
Problem \ref{Pesoft}
has   been   studied
 in detail   in \cite{BeSiam}  in the  setting of 
 linear  isotropic   elasticity  for fibers of circular cross-sections, correcting   results  previously obtained  in \cite{BeBo2} where the influence of $\theta$ had failed to be taken into account.
 We show that the effective problem  associated  to \eqref{Pesoft}  takes the form 
\begin{equation}
\begin{aligned}
&\inf_{(\bfu,\bfv^{tuple})\in L^p(\Omega;\RR^3)\times \D} \int_\Omega   \Phi^{soft}(\bfu,\bfv^{tuple})-\int_\Omega \bff\cdot\bfu dx,
 \\& \Phi^{soft}(\bfu,\bfv^{tuple})=\int_\Omega c^f_{soft}(\bfv-\bfu, \theta)dx +    \Phi_{fibers}(\bfv^{tuple}).
 \end{aligned} 
\nonumber
  \end{equation}
The second term of $\Phi^{soft}$ is  common  with  
  \eqref{Phihomisot}. The first one,  which   unlike \eqref{cf}  emanates    from large  rates of deformation arising  in the entire matrix,
 is    given by
\begin{equation}
\begin{aligned}
&
c^f_{soft}(\bfa,\zeta) :=  \inf_{\bfpsi\in \W(\bfa,\zeta)}
   \int_{Y\setminus S}\hskip-0,5cm  f^{\infty,p}(\bfe_y(\bfpsi)) dx   , \qquad Y:=[-1/2, 1/2[^2, 
\\&  \W(\bfa,\zeta):= \la 
 \bfpsi \in  W^{1,p}_\sharp(Y;\RR^3)\lb \ \begin{aligned} & \int_Y \bfpsi  dy=0, 
\   
\\& \bfpsi (y)=\bfa  +\frac{2}{\diam S}\zeta \bfe_3 \wedge \bfy \  \hbox{ in } \  S 
\end{aligned}\rpt
 \ra,
 \end{aligned} 
\label{defcfsoft}
  \end{equation}
  where $W^{1,p}_\sharp(Y;\RR^3)$ denotes the set of $Y$-periodic members of $W^{1,p}_{\rm loc}(\RR^2;\RR^3)$.

The paper is organized as follows: our main results are presented  in Section  \ref{secapplhom} in the periodic case.
In Section \ref{secvariants}, 
 we discuss their extension to a  non-periodic or random  setting
and  
  the case of large applied body forces.  
  Section \ref{secpreliminaries} is devoted to the asymptotic analysis   of the  sequence of the solutions to \eqref{Pe} 
and of    some auxiliary  sequences characterizing the behavior of the fibers.  Section \ref{secpropcap} 
comprises    a  detailed study of the mapping 
$\capsca^f$ on which  our proofs  crucially rely.
The demonstrations  of our main results, based on the $\Gamma$-convergence method \cite{Da},  are situated in Section \ref{secproof}.
The appendix comprises two technical lemmas relating  to the lower bound and   the proof of  the convergence  \eqref{defcfintro} in the case $p=2$.
Our results were partially  announced in \cite{BeVisco}.
\section{Notations}
\label{secnot}
In this paper,  $\{ \bfe_1, \bfe_2, \bfe_3\}$ stands for the canonical basis of   $\RR^3$. 
Points in $\RR^3$   and real-valued functions  are represented by symbols beginning with a  lightface minuscule  (example $x, i, \det \bfA...$), vectors and vector-valued functions   by symbols
beginning  with a boldface minuscule (examples:     $\bfi$,   $\bfu$,
$\bff$, $\bfg$,
$\bfdiv  \bfPsi $,...),  
matrices and matrix-valued functions   by symbols beginning with a boldface majuscule
with the following exceptions:
$\bfnabla\bfu$ (velocity gradient), $\bfe(\bfu)$ rate of strain tensor). The symbol $\bfI_n$ represents the $n\times n$ identity matrix. 
We denote by $u_i$ or
$(\bfu)_i$  the  components of a vector $\bfu$ and   by $A_{ij}$ or $(\bfA)_{ij}$ those of a matrix $\bfA$ (that is $\bfu=
\sum_{i=1}^3 u_i \bfe_i =\sum_{i=1}^3 (\bfu)_i
\bfe_i$; $\bfA=\sum_{i,j=1}^3 A_{ij} \bfe_i\otimes\bfe_j =\sum_{i,j=1}^3 (\bfA)_{ij} \bfe_i\otimes\bfe_j $). We do not employ  the usual repeated index
convention for summation.    We denote 
 by $\bfA\!:\!\bfB=\sum_{i,j=1}^3 A_{ij}B_{ij}$ the   inner product of two matrices, 
 by $ \e_{ijk} $  the three-dimensional alternator, by 
$\bfu\wedge \bfv=\sum_{i,j,k=1}^3\e_{ijk}u_jv_k\bfe_i$ the exterior product in $\RR^3$,   by $\SSym^n$
   the set of all real symmetric matrices of order
$n$,  by    $\sharp A$   the cardinality of a finite  set $A$, by $B$ (resp. $D$)   the open  unit   ball of $\RR^3$ (resp. $\RR^2$), by $\spt f$ the support of a function $f$, 
 by $\R$  the space of rigid motions in dimension $2$ or $3$.  
For any two vectors $\bfa$, $\bfb$ in $\RR^3$, we set $\bfa\odot \bfb:=\frac{1}{2}(\bfa\otimes \bfb+\bfb\otimes\bfa)$. For any two symmetric matrices $\bfA$, $\bfB$, we write $\bfA\leq \bfB$ if $\bfB-\bfA$ is semi-definite positive.
 Given  $x\in \RR^3$, we write $x=(x', x_3)$ for $x'= (x_1,x_2)$. 
    For any two  weakly differentiable fields  $\bfpsi: \RR^2\to \RR^3$  and    $\bfu: \RR^3\to \RR^3$, we set   
\begin{equation}
\hskip-0,1cm \begin{aligned} 
  & \bfe_y(\bfpsi):= \hskip-0,1cm \sum_{\alpha,\beta=1}^2 \frac{1}{2} \lp  \frac{\partial \psi_\alpha }{\partial y_\beta}+  \frac{\partial \psi_\beta }{\partial y_\alpha }\rp \bfe_\alpha \hskip-0,1cm \otimes\bfe_\beta
  + \sum_{\alpha=1}^2  \frac{\partial \psi_3}{\partial y_\alpha } \bfe_\alpha \odot\bfe_3,
\\&\bfe_{x'}(\bfu):= \hskip-0,1cm 
\sum_{\alpha, \beta=1}^2 \frac{1}{2} \lp  \frac{\partial u_\alpha  }{\partial x_\beta}+ \frac{\partial u_\beta }{\partial  x_\alpha }\rp \bfe_\alpha \hskip-0,1cm \otimes\bfe_\beta
  + \sum_{\alpha=1}^2  \frac{\partial u_3}{\partial x_\alpha } \bfe_\alpha \odot\bfe_3,
  \\& \bfpsi':= \psi_1\bfe_1+\psi_2\bfe_2,  \quad \bfu':= u_1\bfe_1+u_2\bfe_2.
\end{aligned} 
\label{defexprim}  
   \end{equation} 
 Given a topological space $X$,  the symbol $\B(X)$ represents  the $\sigma$-algebra  of the Borel subsets of $X$. 
For any    Radon measure $\nu$ on $X$ and any Banach space $E$, we denote by 
  $L^p_{\nu}(X;E)$   the set of $E$-valued Borel fields $\bfpsi$ on $X$ such that $\int_X |\bfpsi|^p_E d\nu<+\infty$.
 The letter
$C$  (resp. the symbol  $C(a)$ if a dependence on some variable $a$ is indicated) stands for different   positive  constants  whose precise values may vary. 
 %
%
\section{Main results}\label{secapplhom}
    \subsection{Problem \ref{Pe}}\label{secresPe}
\begin{figure}[!h]
  \centering
 \includegraphics[height=6cm]{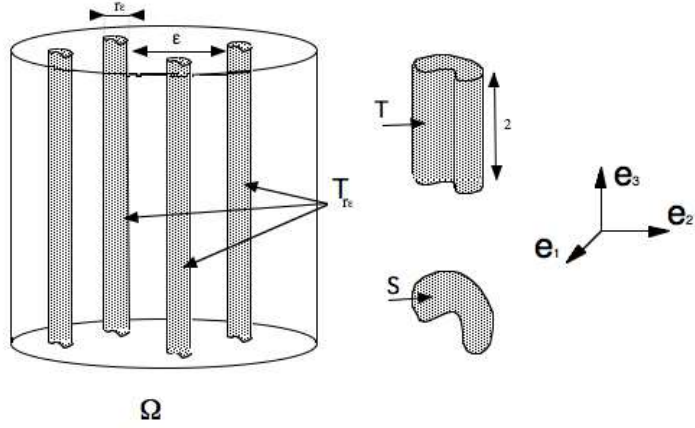}  
  \caption{ The fibered structure }\label{fig1} 
    \end{figure}

  \noindent Let  $\Omega:= \Omega'\times (0,L)$ be   a  bounded Lipschitz domain of $\RR^3$,   $S$  a  bounded Lipschitz domain    of $\RR^2$ verifying 
   \begin{equation}
 D \subset S, \qquad  \bfy_S:= \intb_S \bfy dy=0,  \label{DsubsetS}
   \end{equation}
and  $(r_\e)_{\e>0}$  a sequence of  real numbers   such that
  \begin{equation}
\begin{aligned}
 &   0< r_\e\ll \e \qquad \hbox{ as } \quad \e\to0. 
 \end{aligned} 
\label{repetit}
   \end{equation}  
\noindent    
We consider the $\e$-periodic distribution of fibers defined by   (see fig. \ref{fig1})
\begin{equation}
\begin{aligned}
& T_{r_\e}:= S_{r_\e}\times(0,L)= \bigcup_{i\in I_\e} T_{r_\e}^{i}, \qquad   T_{r_\e}^{i}:= S_{r_\e}^{i} \times (0,L),
\\&  I_\e:= \{ i\in \ZZ^2, \quad Y_\e^{i}\subset  \Omega'\},
\quad    Y^{i}_\e:= \e i+ \e Y, \quad   Y:= \lc -\frac{1}{2}, \frac{1}{2}\rp^2,  
  \end{aligned} 
\label{defTre}
\end{equation}
where for any subset $A$ of $\RR^2$ and  any   $b_\e>0$,  
  \begin{equation}
\begin{aligned}
  & A_{b_\e}:= \bigcup_{i\in I_\e} A_{b_\e}^{i},\quad & & A_{b_\e}^{i}:=  \e i+ b_\e A.  
  \end{aligned} 
\label{defAbi}
\end{equation}
\hskip-0,12cm  
We are concerned with the homogenization of Problem     \ref{Pe} when  $f, g: \SSym^3\to \RR$   
     are    strictly convex  and  satisfy  
\begin{equation}
\begin{aligned}
& c|\bfM|^p\leq f(\bfM) \leq  C |\bfM|^p , \qquad  \ c|\bfM|^p\leq g(\bfM) \leq  C |\bfM|^p \quad  \forall \bfM\in \SSym^3, 
   \end{aligned} 
\label{growthp}
\end{equation}
for some $p\in (1,+\infty)$ and some  positive constants $c,C$.
  We prove  that the limit problem depends on the parameters $k$,    $\kappa$     and  $\gamma^{(p)}$ defined  by \eqref{kkappa}, \eqref{defgammae}.  
We    focuse on the case 
   \begin{equation}
\begin{aligned}
 & k>0,  \qquad \quad \gamma^{(p)}>0. 
 \end{aligned} 
\label{kgamma>0}
   \end{equation}
 Setting  
 \begin{equation}
 f^{\infty, p}(\bfM):= \limsup_{t\to+\infty} \frac{f(t\bfM)}{t^p}; \quad  g^{0,p}(\bfM):= \liminf_{t\to0^+} \frac{g(t\bfM)}{t^p} \qquad \forall \bfM\in \SSym^3,
 \label{deffinfty}
 \end{equation}
we  suppose  that
 %
\begin{equation} 
\begin{aligned}
&
\lb  f(\bfM)- f^{\infty, p}(\bfM)\rb \leq C  (1+\vert \bfM\vert^{p-\varsigma}) & & \forall \bfM\in\SSym^3  \quad \hskip1cm \hbox{ if } \ 1<p\le2, 
\\
 &  \lb \frac{g(t\bfM)}{t^p} -g^{0,p}(\bfM)\rb\leq \varpi(t) |\bfM|^p& &   \forall \bfM\in \SSym^3, \ \forall t\in (0,1)   \quad  \hbox{ if } \kappa>0,
 \end{aligned} 
\label{finftyfg0g}
   \end{equation}  
  for some $\varsigma $, $\varpi$ verifying 
\begin{equation}
\varsigma\in (0,p), \quad   \varpi\in L^\infty(0,1), \quad \lim_{t\to0}\varpi(t)=0.
\label{condvarpi}
   \end{equation}
  Under these hypotheses,   we show that   the solution  $\bfu_\e$   to  (\ref{Pe})
    converges   in a sense defined below  to $\bfu$,   the   sequence $\lp \bfv_\e ,   \theta_\e   \rp$  given 
 by  
    \begin{equation}
    \hskip-0,2cm \begin{aligned}
    &  \bfv_\e (x ) :=\sum_{i\in I_\e} \left (
 {\int\!\!\!\!\!\!-}_{ S_{r_\e}^{i}}
\bfu_\e (s, x_3  )   d\calH^{2}(s) \right)  \mathds{1}_{Y_\e^{i}}(x' ), \quad \tau:= \frac{2\, {\int\!\!\!\!-}_S |y|^2 dy}{\diam S}, 
\\&    \theta_\e (x  ) \hskip-0,1cm: =  \hskip-0,1cm
  \sum_{i\in I_\e} 
  \frac{1}{\tau}  {\int\!\!\!\!\!\!-}_{{ \hskip-0,1cmS_{r_\e}^{i}}}  \hskip-0,2cm
 \lp  -  \frac{y_{\e2}}{r_\e } u_{\e1}+ \frac{y_{\e1}}{r_\e }u_{\e2} \rp\hskip-0,1cm(s,x_3) d\calH^2(s  )    \mathds{1}_{Y_\e^{i}}(x') ,
\end{aligned}
\label{defvethetae}
  \end{equation}
    where 
   \begin{equation}
   \begin{aligned}
   &\bfy_\e( x'):=  \sum_{i\in I_\e}  
  \lp (x_1-\e i_1)\bfe_1+ (x_2-\e i_2)\bfe_2  \rp \mathds{1}_{Y_\e^{i}} (x'),
 \end{aligned} 
\label{defye}
  \end{equation}
 \noindent converges, up to a subsequence   to 
 $(\bfv,\theta)$, and,   if $\kappa>0$,    $\lp  \frac{1}{r_\e}v_{\e3}, \frac{1}{r_\e}  \theta _\e   \rp$,  converges up to a subsequence,   to  $(w,\delta)$, 
  where 
    the couple $(\bfu, \bfv^{tuple})$ defined by  \eqref{tuple}  is  a  solution to  
 %
     \begin{equation}
\begin{aligned}
&(\calP^{hom}):  \inf_{(\bfu, \bfv^{tuple})  \in W^{1,p}_b (\Omega;\RR^3)\times\D }  \Phi (\bfu, \bfv^{tuple})-\int_\Omega \bff.\bfu dx,
\end{aligned} 
\label{Phom}    \end{equation}
and   $\bfu$ is   the unique solution to
\begin{equation}
\begin{aligned}
&   \inf_{\bfu \in W^{1,p}_b (\Omega;\RR^3) }  F(\bfu)-\int_\Omega \bff.\bfu dx,
\qquad   \qquad F(\bfu):=\inf_{\bfv^{tuple}\in\D}\Phi(\bfu,\bfv^{tuple}).
\end{aligned} 
\label{defF}    \end{equation}
  The functional $\Phi$  is defined by 
%
 %
     \begin{equation}
    \label{defPhi}
   \begin{aligned}
&\Phi (\bfu, \bfv^{tuple}):= \int_\Omega f(\bfe(\bfu))  
+ c^f(\bfv-\bfu, \theta)  
 + g^{hom}(\bfgotD\,\bfv^{tuple}) dx, 
 \end{aligned}
\end{equation}
where,   if $\kappa=0$,  
          \begin{equation}
\begin{aligned}
&\bfgotD\,\bfv^{tuple}:=\lp  \frac{\partial v_3}{\partial x_3},  \frac{\partial \theta}{\partial x_3}\rp,
\\&g^{hom}\lp a,  \beta\rp:=
  \inf_{\bfq\in  W^{1,p}(S;\RR^3)}   k \intb_{S} g  \lp  \bfgotL\lp \bfq,  a,  \beta \rp \rp dy  ,
\\& \bfgotL(\bfq,a,\beta):=  \bfe_{y}( \bfq)+
\frac{2}{\diam S}\beta
(-y_2 \bfe_1\odot\bfe_3 +y_1\bfe_2\odot\bfe_3)
+a \bfe_3\otimes\bfe_3, 
\end{aligned} 
\label{defghom1}   \end{equation}
 if $\kappa>0$,  
      \begin{equation}
\begin{aligned}
&\bfgotD\,\bfv^{tuple}
:= \lp\frac{\partial^2  v_1}{\partial x_3^2},   \frac{\partial^2  v_2}{\partial x_3^2} , \frac{\partial w}{\partial x_3}, \frac{\partial \delta}{\partial x_3}\rp,
\\& g^{hom}(\zeta_1,\zeta_2,a,\beta):=  \! \!\inf_{\bfq\in W^{1,p}(S;\RR^3)}   \kappa  \intb_{S} g^{0,p}  \lp\bfgotJ(\bfq, \zeta_1,\zeta_2,a,\beta) \rp dy    ,
\\&\bfgotJ (\bfq, \zeta_1,\zeta_2,a,\beta):=  \bfgotL(\bfq,a,\beta) -\lp \sum_{\alpha=1}^2  \zeta_\alpha  y_\alpha \rp\bfe_3\otimes\bfe_3, 
\end{aligned} 
\label{defghom2}   \end{equation}
and   $c^f$ is   defined   by \eqref{defcfintro} in terms of $\capsca^f$ 
given   by  \eqref{infW}, and of any sequence $(R_\e)$ verifying
      \begin{equation}
 \begin{aligned}
&0 < r_\e\ll R_\e \ll \e\ll 1,  \quad  
  R_\e\ll r_\e^{2-p},  \quad   
  R_\e\ll\frac{1}{\sqrt{|\log r_\e|}} \  \hbox{ if } \ p=2.
\end{aligned}
 \label{Re}
   \end{equation}
We show that
 $c^f$   is well defined, convex, independent of $(R_\e)$ and,  if $p\ge2$, 
 of $S$,  and  satisfies
\begin{equation}
\begin{aligned}
  & \lpt \begin{aligned} &c^f(\bfa, \zeta)   =  
   \gamma^{(p)}    \capsca^{f^{\infty,p}}( \bfa, \zeta\bfe_3 ;S,\RR^2),  
\\&c\gamma^{(p)}(|\bfa|^p + |\zeta|^p)  \le c^f(\bfa, \zeta)\le C\gamma^{(p)}(|\bfa|^p + |\zeta|^p) 
\end{aligned}\ra & & \hbox{if } 1<p<2,
\\&
  \lpt  \begin{aligned} &  c\gamma^{(2)} |\bfa|^2  \le  c^f(\bfa,0)\le C\gamma^{(2)}|\bfa|^2                           \\&  c^f(\bfa, \zeta)= +\infty  \qquad \hbox{ if } \zeta\not=0,           
  \end{aligned}\ra & & \hbox{if }\  p=2, 
  \\&\ \,
 c^f(0,0) =    0, \quad 
    c^f(\bfa, \zeta) = +\infty  \ \   \hbox{if } \  \ (\bfa, \zeta) \not=(0,0),
 & &\hbox{if }\  p>2.
\end{aligned}
 \label{cf=}
\end{equation} 
The set   $\D$   in \eqref{Phom} 
is   the Banach space given  by  
\begin{equation}
   \begin{aligned}
&  \D := 
\Big\{ ( \bfv, \theta)  \in   L^p(\Omega;\RR^3) \times L^p( \Omega'; W^{1,p} (0,L)),
\\&\qquad v_3  \in  L^p( \Omega'; W^{1,p} (0,L)),   v_3 = \theta = 0  \text{ on }    \Omega'\times\{0\} \Big\}\ 
& & \hskip-3,2cm \text{ if } 0<k<+\infty,
\\&\D := 
\la ( \bfv, \theta)   \in   L^p(\Omega;\RR^3) \times L^p(\Omega), \  v_3= \theta=0\ra \  & & \hskip-3,2cm\text{ if }  (k,\kappa)\hskip-0,1cm=\hskip-0,1cm(+\infty ,0),
\\&\D :=   
\Big\{ \hskip-0,1cm( \bfv, \theta,w,\delta) \hskip-0,1cm \in \hskip-0,1cm L^p (\Omega';W^{2,p} (0, L;\RR^3) ) \times   L^p(\Omega';W^{1,p} (0, L ) )^3,    v_3 \!=\! \theta\! = 0,  
\\& 
\qquad  w\!=\!\delta\!=\! v_\alpha =\!\frac{\partial v_\alpha}{\partial x_3}\!=0   \text{ on } \Omega'\times\{0\} \ (\alpha\in \{1,2\}) 
\Big\}\qquad 
 & & \hskip-3,2cm\text{ if } 0\! <\! \kappa\! <\! +\infty, 
\\&\D :=    \la 0 \ra\hskip3cm & & \hskip-3,2cm\text{ if }  k= \kappa=+\infty.  
  \end{aligned} 
 \label{defD}
  \end{equation}
   %
\begin{theorem} \label{th}  Assume \eqref{repetit},  \eqref{defTre}, \eqref{growthp},     \eqref{kgamma>0}, \eqref{finftyfg0g},  then  the solution $ \bfu_\e $  to  \eqref{Pe}
 weakly  converges  in  $W_b^{1,p} (\Omega; \RR^3)$ to $\bfu$,   the couple 
   $(\bfv_\e, \theta_\e)$  defined by \eqref{defvethetae}   weakly  converges  in  $L^p( \Omega;\RR^3 )\times L^p(\Omega)$, up to a subsequence,   to 
$(\bfv,\theta)$, and,   if $\kappa>0$,   $(\frac{1}{r_\e}  v_{\e3},\frac{1}{r_\e}\theta_\e)$  weakly  converges 
in $(L^p( \Omega))^2$, up to a subsequence,  to $(w,\delta)$, where   $(\bfu,\bfv^{tuple})$  given  by \eqref{tuple}  
 is  a solution   to \eqref{Phom} and $\bfu$ is the unique solution to \eqref{defF}.
    \end{theorem}

\begin{remark}\label{remgeneral} (i) 
The functional $F$  in \eqref{defF}  is  the $\Gamma$-limit of the sequence $(F_\e)$  in the weak topology of $W_b^{1,p}(\Omega;\RR^3)$ (see Remark \ref{remFGamma}). 
\\
(ii)  If  $\gamma^{(p)}=0$ or  $k=0$, 
  the solution  $\bfu_\e$   to  (\ref{Pe}) weakly  converges  in  $W_b^{1,p} (\Omega; \RR^3)$  to the unique solution to:\ 
$   \inf_{ \bfu  \in W^{1,p}_b (\Omega;\RR^3)  } \int_\Omega f(\bfe(\bfu)) dx -\int_\Omega \bff\cdot\bfu dx$.
\\ (iii) If  $f^{\infty,p}$ and   $g^{0,p}$  are   strictly convex, which   the strict convexities of $f$ and $g$ do not ensure,
and $p\not=2$, then $\Phi$ 
  is strictly convex, the solution 
  to \eqref{Phom} unique, and  all  convergences stated in Theorem \ref{th} hold for the whole sequences.  
\\(iv) {\it Linear isotropic elasticity.}  If    $f=f_{\lambda_0,\mu_0}$ and $g=f_{\lambda_1,\mu_1}$  for some $\lambda_0,\lambda_1\ge 0$, $\mu_0,\mu_1>0$, where
    \begin{equation}
\label{fisot}
\begin{aligned}
f_{\lambda,\mu}(\bfM): =\frac{\lambda}{2}({\rm tr}\bfM)^2+\mu\bfM:\bfM \quad \forall \bfM\in \SSym^3,
\end{aligned}
\end{equation}
we prove  (see  Remark \ref{remMf2}) that    
   \begin{equation}\begin{aligned} 
 & c^{f_{\lambda_0,\mu_0} }(\bfa,0)= \gamma^{(2)} \bfa\cdot 2 \mu_0\pi  \lp 2\frac{\lambda_0+2\mu_0}{\lambda_0+3\mu_0}(\bfe_1\otimes\bfe_1+\bfe_2\otimes\bfe_2)+\bfe_3\otimes\bfe_3 \rp\bfa,
 \end{aligned}
\label{explicit2}
  \end{equation} 
and deduce from 
   \eqref{defghom1} and \eqref{defghom2} 
that, 
if $\kappa=0$,   
          \begin{equation}
\begin{aligned}
& \Phi_{fiber}(\bfv^{tuple})=k\mu_1  \int_\Omega \frac{3\lambda_1+2\mu_1}{2(\lambda_1+\mu_1)}   \lb \frac{\partial v_3}{\partial x_3} \rb^2 
+  \frac{2   m}{(\diam S)^2}    \lb \frac{\partial \theta}{\partial x_3} \rb^2dx , 
\end{aligned} 
\label{Phiisotk} \end{equation}
  and, if $\kappa>0$, 
      \begin{equation}
\begin{aligned}
&\Phi_{fiber}(\bfv^{tuple})=
 \kappa  \mu_1 \int_\Omega   \frac{3\lambda_1+2\mu_1}{2(\lambda_1+\mu_1)}
\Bigg[
  \sum_{\alpha=1}^2   \lp\intb_S y_\alpha^2 dy \rp \lb \frac{\partial^2 v_\alpha}{\partial x_3^2}\rb^2
\\&\qquad\qquad  - 2 \lp\intb_S y_1y_2dy \rp \frac{\partial^2 v_1}{\partial x_3^2} \frac{\partial^2 v_2}{\partial x_3^2}
 + \lb \frac{\partial w}{\partial x_3}\rb^2
\Bigg] +  \frac{2  m}{(\diam S)^2}    \lb \frac{\partial \delta}{\partial x_3}\rb^2 dx,
\end{aligned} 
\label{Phiisotkappa}   \end{equation}
 where  $m:= \inf_{\varphi\in H^1(S)} \intb_S \lp\tfrac{\partial\varphi}{\partial y_1}-y_2\rp^2+ \lp\tfrac{\partial\varphi}{\partial y_2}+y_1\rp^2 dy$.
This infimum being attained, it results from  the Schwarz theorem that  $m>0$.   Noticing that, by  \eqref{Phiisotkappa},  the  infimum \eqref{Phom} is achieved by 
$(w,\delta)=(0,0)$   and that, by \eqref{cf=}, $\theta=0$, we 
  recover the 
 formulae
obtained   for  $S=D$   in  \cite[Remark 2.2]{BeGr}.
\\   (v)     The following    example  shows that  $\delta$   
doesn't vanish,  in general, when anisotropic fibers are considered: assume that  $0<\kappa<+\infty$ and   
 $$
  g(\bfM):= \sum_{\alpha, \beta=1}^2 M_{\alpha\beta}M_{\alpha\beta}+ M_{13}^2+M_{33}^2+M_{13}M_{33}+M_{23}^2.
  $$
Then, a solution to the minimization problem  \eqref{defghom2} is  $\bfq:=( \zeta_1\varphi^{(1)}+\zeta_2\varphi^{(2)}+ a \varphi^{(3)} +\tilde \beta\varphi^{(4)})\bfe_3$, where  $\tilde\beta:= \tfrac{1}{\diam S}\beta$ and 
  $\varphi^{(1)},..,\varphi^{(4)}$  are   solutions    in $H^1(S)/\RR$
to  the  Neumann  problems ($\bfn$:  outer unit normal to $\partial S$)
  \begin{equation}
    \begin{aligned}
 & (\calP^{(1)}): \Delta\varphi^{(1)}=1\quad \hbox{ in }  S, \quad & &\bfnabla \varphi^{(1)}\cdot \bfn= y_1 n_1 \quad \quad & &\hbox{ on }  \partial S, 
 \\ & (\calP^{(2)}): \Delta\varphi^{(2)}=0\quad \hbox{ in }  S, \quad & &\bfnabla \varphi^{(2)}\cdot \bfn= y_2 n_2 \quad \quad & &\hbox{ on }  \partial S, 
 \\ & (\calP^{(3)}): \Delta\varphi^{(3)}=0\quad \hbox{ in }  S, \quad & &\bfnabla \varphi^{(3)}\cdot \bfn= - n_1 \quad \quad & &\hbox{ on }  \partial S, 
 \\ & (\calP^{(4)}): \Delta\varphi^{(4)}=0\quad \hbox{ in }  S, \quad & &\bfnabla \varphi^{(4)}\cdot \bfn= -2\begin{pmatrix} -y_2\\y_1\end{pmatrix}\cdot \bfn   & &\hbox{ on }  \partial S,
 \end{aligned}
\nonumber
\end{equation}
respectively.
We deduce that  $g^{hom} (\zeta_1,\zeta_2,a,\beta)$ $  = \begin{pmatrix} \zeta_1\\ \zeta_2\\ a\\ \tilde\beta\end{pmatrix} \cdot \bfC \begin{pmatrix} \zeta_1\\ \zeta_2\\ a\\ \tilde\beta\end{pmatrix}$,
where 
\begin{equation}
    \begin{aligned}
 &\bfC\in \SSym^4, \quad C_{33}= \kappa, \quad C_{44}= \kappa\intb_S|y|^2 dy, \quad 2 C_{34}=\kappa \intb_S  \tfrac{1}{2}\tfrac{\partial \varphi^{(4)}}{\partial y_1}dy
  \\&2C_{\alpha3}=\kappa \intb_S \tfrac{1}{2}\tfrac{\partial \varphi^{(\alpha)}}{\partial y_1} dy,
\quad 2 C_{\alpha4}
 =\kappa\intb_S -y_\alpha \lp \tfrac{1}{2}\tfrac{\partial \varphi^{(4)}}{\partial y_1}-y_2 \rp dy \quad \forall \alpha\in  \{1,2\}.
 \end{aligned}
 \nonumber
\end{equation}
If    $(\bfu,\bfv^{tuple})$ is a solution to \eqref{Phom}, then    $(w,\delta)$ is a solution to 
  $\inf_{(w,\delta)} \Phi(\bfu, \bfv, 0,w,\delta)$. 
 We infer 
$$ \begin{pmatrix} C_{33}&C_{34}\\C_{34} & C_{44}\end{pmatrix}  \begin{pmatrix} w \\ \tfrac{1}{\diam S}\delta \end{pmatrix}  = 
\sum_{\alpha=1}^2 
 \begin{pmatrix} -C_{\alpha3}\frac{\partial v_\alpha }{\partial x_3}
 \\ -C_{\alpha4}\frac{\partial v_\alpha }{\partial x_3} \end{pmatrix}.$$
If $S=D$,  we obtain  $\varphi^{(4)}=0$,
   $C_{14}=C_{34}=0$,  $C_{44}=4 C_{24}>0$, and deduce
   $\delta= -\tfrac{\diam S}{2}  \tfrac{\partial v_2}{\partial x_3}$.
\end{remark}

%
  %
    \subsection{Problem \ref{Pesoft}}\label{secresPesoft}    Assuming     
 \begin{equation}
 \begin{aligned}
 &
r_\e=  \e, \quad \ov S \subset \ring Y,
 \end{aligned} 
  \label{resoft}
  \end{equation}
  and 
 \begin{equation}
 \begin{aligned}
 &
|\bfvarphi|_{L^p(\Omega;\RR^3)} \le C F^{soft}_\e(\bfvarphi) \quad \forall \bfvarphi\in W^{1,p}_b(\Omega;\RR^3),
 \end{aligned} 
  \label{condsoft}
  \end{equation}
we  prove  that the     effective problem  associated to \eqref{Pesoft}  is
\begin{equation}
\begin{aligned}
&\inf_{(\bfu,\bfv^{tuple})\in L^p(\Omega;\RR^3) \times \D}  \quad  \Phi^{soft}(\bfu,\bfv^{tuple})-\int_\Omega \bff\cdot\bfu,
  \\&  \Phi^{soft}(\bfu,\bfv^{tuple})=\int_\Omega c^f_{soft}(\bfv-\bfu, \theta)dx +   \int_\Omega  g^{hom}(\bfgotD\,\bfv^{tuple}) dx,
 \end{aligned} 
 \label{Phomsoft}
  \end{equation}
  where $c^f_{soft}$,   $g^{hom}$ and $\D$  are defined by \eqref{defcfsoft},  \eqref{defghom1}, \eqref{defghom2} and  \eqref{defD}, respectively.
\begin{theorem} \label{thsoft} Assume \eqref{resoft} and \eqref{condsoft}, then 
  the solution   $\bfu_\e$  to \eqref{Pesoft}    
 weakly  converges in $L^p(\Omega;\RR^3)$,   up to a subsequence,  to $\bfu$,   and the convergences 
of   $(\bfv_\e, \theta_\e)$ and $( \frac{1}{\e}  v_{\e3},\frac{1}{\e}\theta_\e)$ stated in Theorem \ref{th}  hold,
where   $(\bfu,\bfv^{tuple})$
 is a  solution   to \eqref{Phomsoft}.  
\end{theorem}

 \begin{remark}\label{remcomp} (i)  Theorems \ref{th} and  \ref{thsoft} can be generalized to multiphase media comprising non-intersecting families of parallel    fibers distributed in different directions  (for Theorem \ref{thsoft}, case $p=2$,  see    \cite[Sec. 4]{BeSiam}).
%

\noindent (ii)  The condition \eqref{condsoft} 
is  guaranteed  when   $\kappa>0$  or when  suitable multiphase media  described in    \cite[Sec. 5]{BeSiam}  are considered,  for instance 
when $k>0$ and  the fibers are distributed in at least three independent directions.  

\noindent (iii) Most of Remark \ref{remgeneral}   applies to  Problem \ref{Pesoft}.
 
      \end{remark} 
      
\section{Variants}\label{secvariants} 
%
 \subsection{Non-periodic case}
 \label{secnonper}
 The study of  Problem  \ref{Pe}
  can be extended to a non-periodic setting:
 we then  parametrize 
it   by  
 the size $r$ of the cross-sections of the fibers.
The collection of fibers   is   described  in terms of   the  image $G_r$ of the set of their  principal axes   under  the  orthogonal projection  onto $\Omega'$,
by setting  
\begin{equation}\label{TG}
\begin{aligned}
T_{r}(G_r):= \bigcup_{t\in G_r}   (t +r S) \times     (0,L).
\end{aligned}
\end{equation}
We assume  that $(G_r)_{r>0}$  is a family of   finite  subsets of $\Omega'$ satisfying,    as $r\to0$,
 \begin{equation}
\begin{aligned}
&  
r \ll \inf_{t_1,t_2\in G_r, t_1\not=t_2} |t_1-t_2| \ll 1;
\qquad  r\ll  \dist(G_r, \partial \Omega') .
 \end{aligned}
\label{Omegae}
\end{equation}
We focuse on the case  $p\in (1,2]$. 
To take advantage of the notations introduced in  the periodic case, we fix   an arbitrary positive real number $\gamma^{(p)}$ and introduce the small parameter $\e_r$   defined by 
\begin{equation}
\begin{aligned}
 & \e_r:= \frac{r^{\frac{2-p}{2}}}{\lp\gamma^{(p)}\rp^{1/2}}
 \quad \hbox{if }  \ p\in (1,2);  \qquad 
  \e_r:= \frac{1}{\sqrt{\gamma^{(2)}|\log{r}|}} \quad \hbox{if }  \ p=2.
  \end{aligned}
\label{gammaenonper}
\end{equation}
By \eqref{defgammae} we have $\gamma_{\e_r}^{(p)}(r)=\gamma^{(p)}$ for every $r>0$.
Given a sequence of positive real numbers  $(l_{\e_r})_{r>0}$, we   consider  the sequence of problems     $(\calP_{r}(G_{r}))_{r>0}$  formally   deduced from  \eqref{Pe} by substituting $T_{r}(G_{r})$ for $T_{r_\e}$ and $l_{\e_r}$ for $l_\e$.
The function  defined by 
 \begin{equation}	\label{defne}
\begin{aligned}
  	& n_{G_{r}}(x') :=  \sum_{ i\in I_{\e_r}}\sharp   \lp  G_{r}  \cap  Y_{\e_r}^{i}  \rp\  \mathds{1}_{ Y_{\e_r}^{i} }(x') \qquad \forall x'\in \Omega',
	\end{aligned}
\end{equation}
  where  $I_{\e_r}$ and $Y_{\e_r}^{i}$ are  given by   \eqref{defTre}, 
  locally approximates to  the number of fibers crossing a square of size $\e_r$ in the cross-section of $\Omega$.
We suppose that   $(n_{G_{r}})_{r>0}$  is bounded in $L^\infty(\Omega')$ and 
 weakly$^\star$ converges to some   $n$  as $r\to0$. Setting   $\mu:= n\calL^2_{\lfloor \ov \Omega'}$, 
by combining  the argument of  the proof of Theorem \ref{th} with  the one developed in \cite{BeLiOr},
  one can prove that  
   the effective problem      takes the form 
     \begin{equation}
\begin{aligned}
&(\calP^{hom, \, \mu}):
  \inf_{\bfu \in W^{1,p}_b (\Omega) } F^\mu (\bfu)-\int_\Omega \bff.\bfu dx; \quad 
   F^\mu (\bfu):= \hskip-0,3cm  \inf_{\bfv^{tuple}\in   \D^{\mu} } \Phi^\mu (\bfu, \bfv^{tuple}),
\\
 &\Phi^\mu (\bfu, \bfv^{tuple})=  \int_\Omega  f(\bfe(\bfu)) dx + \int_\Omega   c^f(\bfv-\bfu, \theta)  +   g^{hom}(\bfgotD\,\bfv^{tuple})
 d\mu\otimes\calL^1\hskip-0,15cm,
 \end{aligned}
 \label{Phomn} 
\end{equation}
\noindent  where   $c^f$  and  $g^{hom}$  are  given by  \eqref{defghom1},  \eqref{defghom2}  and  \eqref{cf=} in terms  
of  $\gamma^{(p)}$ introduced above  and  $k$, $\kappa$ defined by 
substituting $\e_r$ for $\e$ in \eqref{kkappa}.  The domain 
   $\D^{\mu}$   is   deduced   from \eqref{defD}   by   substituting  
 $L^p_{\mu\otimes\calL^1}(\Omega;\RR^3)$ for $L^p(\Omega;\RR^3)$,    the spaces 
 $L^p_\mu( \Omega'; E)$  for $L^p( \Omega'; E)$ for any Banach space $E$, and by replacing the  homogeneous Dirichlet conditions on 
$\Omega'\times\{0\}$ by  homogeneous Dirichlet conditions  $\mu\otimes\delta_{\{0\}}$-a.e.  on 
$\Omega'\times\{0\}$.  
\begin{theorem} \label{thnonper}  
Under the assumptions stated above, 
  the solution   $\bfu_r$  to $(\calP_{r}(G_{r}))$
   weakly converges  
in $W^{1,p}_b(\Omega;\RR^3)$ as $r\to0$  to   the unique solution $\bfu$   to \eqref{Phomn}. 
   \end{theorem}

\begin{remark}\label{remnonper} 
(i)  Under the same assumptions, 
one can show that the measure $\tilde m_r$,  defined  by  substituting $T_{r}(G_r)$ for $S_{r_\e}\times (0,L)$ and $(r,\e_r)$ for $(r_\e, \e)$ 
 in 
\eqref{defme}, weakly$^\star$ converges  to $\mu\otimes\calL^1$ in $\M(\ov\Omega)$, 
and  the 
sequences  $(\tilde\bfv_r  \tilde m_r)$,  $(\tilde\theta_r   \tilde m_r)$ and, if $\kappa>0$,  $\lp \tfrac{1}{\e_r} \tilde v_{r3}  \tilde m_r \rp$ and  $\lp \tfrac{1}{\e_r} \tilde \theta_r  \tilde m_r \rp$,
where
  \begin{equation}
    \hskip-0,2cm \begin{aligned}
    &  \tilde\bfv_r (x ) :=\sum_{t\in G_{r}} \left (
 {\int\!\!\!\!\!\!-}_{{ t+r S}}
\bfu_r (s, x_3  )   d\calH^{2}(s) \right)  \mathds{1}_{ t+r  S}(x' ),
\\&   \tilde \theta_{r} (x  ): = 
  \sum_{t\in G_r} 
  \frac{1}{\tau}  {\int\!\!\!\!\!\!-}_{t+r S}  \hskip-0,2cm
 \lp  -  \frac{\tilde y_{r2}}{r } u_{r1}+ \frac{\tilde y_{r1}}{r } u_{r2} \rp\hskip-0,1cm(s,x_3) d\calH^2(s  )    \mathds{1}_{t+r S}(x')  , 
\\
   &\tilde\bfy_r( x'):=  \sum_{t\in G_r}  
  \lp (x_1-t_1)\bfe_1+ (x_2-t_2)\bfe_2  \rp \mathds{1}_{t+rS} (x'),
 \end{aligned} 
\nonumber
  \end{equation}
weakly$^\star$ converge up to a subsequence  in $\M(\ov\Omega;\RR^3)$ and $\M(\ov \Omega)$ to $\bfv \mu\otimes \calL^1$,  $\theta\mu\otimes \calL^1$,   $w \mu\otimes \calL^1$,  and $\delta \mu\otimes \calL^1$, respectively, 
where $\bfv^{tuple}$, defined by \eqref{tuple}, is  a  solution to the minimization problem defining  $F^\mu(\bfu)$ in \eqref{Phomn}.  
\\
(ii)  Theorem \ref{thnonper}    holds, in some cases,  when  $(n_{G_r}\calL^2_{\lfloor\Omega'})$ is   
 bounded in $\M(\ov\Omega')$ and  
  weakly$^\star$ converges   to some  singular 
measure $\mu$ which vanishes on all Borel  subset  of  $\ov\Omega'$  of  $p$-capacity zero. For instance, 
assume  that   
$\Omega'=(0,L)^2$  and set 
$$
G_{r}:=\la (L/2,\e_r^2i), \ i\in \ZZ\ra \cap\  \{L/2\}\times (\e_r^2, L-\e_r^2).
$$
The set $T_{r}(G_{r})$ defined by  \eqref{TG} represents a  family of fibers whose principal axes are $\e_r^2$-periodically distributed  on the surface
$\Sigma:= \{L/2\}\times (0,L)^2$. The sequence $(n_{G_{r}}\calL^2_{\lfloor{\Omega'}})$ defined by \eqref{defne}  is bounded in $\M(\ov\Omega')$  
 and weakly$^\star$ converges  to  $\mu:=\calH^1_{\lfloor \{L/2\}\times (0,L)}$.
One can show, by
 combining  the argument of the proof of Theorem \ref{th}  with that   developed in  \cite{BeGeKr},  
 that   Theorem \ref{thnonper}  holds, as well as 
the convergences stated    in (i).
\\
\noindent (iii) {\it Dirichlet problems in varying domains.}  Under the assumptions of Theorems \ref{th} or  \ref{thnonper}  (or  Remark \ref{remnonper} (ii)),  the solution to 
     \begin{equation}
\begin{aligned}
& \inf_{\bfu_r\in W^{1,p}_{b} (\Omega;\RR^3) } \la \int_\Omega f(\bfe(\bfu_r)) dx -\int_\Omega \bff. \bfu_r dx ,   \quad \bfu_r= 0  \hbox{ in } \ T_{r}(G_{r})\ra,
 \end{aligned} 
 \label{PeDir}
  \end{equation}
weakly converges in $W^{1,p}_{b} (\Omega;\RR^3)$   to the unique  solution to   
   \begin{equation}
\begin{aligned}
&  \inf_{\bfu \in W^{1,p}_b (\Omega;\RR^3)} \int_\Omega   f(\bfe(\bfu)) dx +  \int_\Omega   c^f(\bfu, 0) d\mu\otimes\calL^1 -\int_\Omega \bff.\bfu dx,
\end{aligned} 
\label{PhomDir}    \end{equation}
 deduced  from  $(\calP^{hom, \, \mu})$
   by substituting $0$ for $\bfv^{tuple}$. 
   The second term is  analogous  to the so-called   "strange term" 
  \cite{CiMu,MaKh}.
If  $(n_{G_r})_{r>0}$ is only assumed to be bounded in $L^1(\Omega')$ and to weakly$^\star$  converge to some arbitrary 
measure $\mu$,  we expect the limit problem  associated with \eqref{PeDir} to depend,  not on $\mu$, but, up to a subsequence, on some  
 measure 
$\mu_0$  
defined
through  a variant of the    $\gamma$-convergence  introduced in \cite{DaMo}.
Many references on this subject can be found in   
 \cite{DaMu}.
In particular, compactness results ensuring the existence of a $\gamma$-converging subsequence have been established in various contexts.
Under \eqref{Omegae}, the limit  problem associated with \eqref{PeDir} is  likely to be  deduced from \eqref{PhomDir}, formally,  by substituting such a measure   $\mu_0$ for $\mu$. 
This   suggests that the limit problem associated with $(\calP_{r}(G_{r}))$ might  possibly   be $(\calP^{hom, \, \mu_0})$.

\noindent (iv) The assumption \eqref{gammaenonper}   on the choice  of  $\e_r$ was  omitted in \cite{BeLiOr}.
In order the results stated in   \cite{BeLiOr} to  be  correct, one should  either assume  in \cite[(17)]{BeLiOr} that $\gamma^{(p)}$  is finite (hence $1<p\le2$), or  in \cite[(5)]{BeLiOr}  that $n_\e$ is bounded from below by a positive constant, otherwise  the  asymptotic  behavior of  the sequence $(\gamma_{\e}(r)n_\e)$, whose  knowledge  is necessary   to obtain  for instance \cite[(126)]{BeLiOr}, would be undetermined.
  \end{remark}
  
%
  \subsection{Random case}
  \label{secrandom}
Assuming $1<p\le2$, fixing    $d>0$, we set 
 \begin{equation}
\label{defgotO}
   \begin{aligned}
  & \gotO:=\Big\{ \omega
 \in 2^{\RR^2 }, \ \forall (\omega_1,\omega_2)\in \omega^2, \  \omega_1\not= \omega_2 \ \Rightarrow \ 
\vert \omega_1-\omega_2\vert \geq d \Big \}.
\end{aligned}
\end{equation}
\noindent 
One can check that  the  finite metric     $d_\gotO: (\omega,\omega')\in \gotO^2\to   \inf \{1, d_\calH(\omega,\omega')\}$, where $d_\calH$ stands for  the Hausdorff distance,   turns   $\gotO$ into a complete metric space.
Denoting  by   $\B(\gotO)$   the associated  Borel $\sigma$-algebra on $\gotO$, we consider   a probability measure  $P$    on  $(\gotO, \B(\gotO))$ satisfying
 \begin{equation}
\nonumber
P(A+z )= P(A) \quad \forall \  z\in \ZZ^2, \quad \forall \,A\in \B(\gotO).
\end{equation}
We  introduce   the $\sigma$-algebra $\F$    and the random variable $n_0$  on $\gotO$ defined by    
\begin{equation}
\nonumber
\F:= \la A\in \B(\gotO), \  A+z =A \quad \forall z\in \ZZ^2\ra, \qquad n_0(\omega):= \sharp \lp \omega\cap \left[ -\frac{1}{2}, \frac{1}{2}\right[^2\rp.
\end{equation}
We denote by  ${E}_{{P}}^{\F} n_0$  
    the   conditional expectation    of  $n_0$   given $\F$ w.r.t. ${P}$.
We fix a positive real number $\gamma^{(p)}$ and  set 
 \begin{equation}
\nonumber
   \begin{aligned}
 G_{r} (\omega):= \e_r \omega \cap \{x\in \Omega', \ \dist(x,\partial \Omega')\ge \e_r\} \qquad \forall r>0,\  \forall \omega\in \gotO,
\end{aligned}
\end{equation}
where $\e_r$ is given by  \eqref{gammaenonper}.
The following statement is proved in    \cite[Th. 2.4.2]{BeHDR}:
 \begin{theorem}\label{thrandom}
Under the assumptions stated above,   there exists   a  sequence $(r_k)_{k\in\NN}$  converging to $0$ and  a $P$-negligible set  $\gotN\in \B(\gotO)$  such that,  for each $\omega\in \gotO\setminus \gotN$,  the sequence  $(n_{G_{r_k}(\omega)}) $  defined by    (\ref{defne})  weakly$^\star$ converges in $L^\infty(\Omega')$  to 
  the constant function ${E}_{{P}}^{\F} n_0(\omega)$.
    \end{theorem}
    
    \noindent
The assumptions \eqref{gammaenonper} and   \eqref{defgotO} ensure that   $ G_{r} (\omega)$ satisfies \eqref{Omegae}.
By combining     Theorems \ref{thnonper}  and \ref{thrandom},  we obtain: 
 \begin{corollary}\label{corrandom}
There exists a  sequence $(\!r_k\!)_{k\in\NN}$  converging to $0$ and  a $P\!$-negligible set  $\gotN\in \B(\gotO)$  such that, 
for all $\omega\in \gotO\setminus \gotN$, the solution  $\bfu_{r_k}(\omega)$ 
to  the problem
 $(\calP_{r_k}(G_{r_k}(\omega)))$  considered   in Section \ref{secnonper}
  weakly  converges in  $W^{1,p}_b(\Omega;\RR^3)$
 to the unique solution $\bfu(\omega)$  to  
 $ ( \calP^{hom, \ \mu(\omega) })$ defined by  \eqref{Phomn},  where 
 $\mu(\omega):= {E}_{{P}}^{\F} n_0(\omega)\calL^2_{\lfloor \Omega'}$.
  \end{corollary}
  
%
 \subsection{Large  applied   body forces}
  \label{seclargebodyforces} 
We consider the problems
  \begin{equation}
\begin{aligned}
 &
   \inf_{\bfu_\e\in W^{1,p}_b(\Omega;\RR^3) } F_\e(\bfu_\e)-\int_\Omega \bff_\e \cdot \bfu_\e dx,
   \end{aligned} 
 \label{Pefe} 
   \end{equation}
  \begin{equation}
\begin{aligned}
 &
   \inf_{\bfu_\e\in W^{1,p}_b(\Omega;\RR^3) } F^{soft}_\e(\bfu_\e)-\int_\Omega \bff_\e \cdot \bfu_\e dx,
   \end{aligned} 
 \label{Pesoftfe} 
   \end{equation}
%
where, given $\bff\in L^{p'}(\Omega;\RR^3)$,   $ \bfg_0  \in C(\ov { \Omega\times S};\RR^3) $ and   $a_0,   \beta_0 \in C(\ov { \Omega\times S})$, 
$\bff_\e$ is defined in terms of  $y_\e$  given  by  \eqref{defye},  by 
   \begin{equation}
\hskip-0,2cm \begin{aligned}
 &     \bff_\e(x)= \bff (x)+ \frac{|\Omega|}{|{T_{r_\e}}|} \bfg_0 \lp x,\frac{y_\e(x')}{r_\e}\rp \mathds{1}_{{T_{r_\e}}}(x) \hskip 1,6cm & &  \hbox{ if } \ \kappa=0, 
 \\
&      \bff_\e(x) \hskip-0,05cm :=\hskip-0,05cm\bff(x)+ \frac{|\Omega|}{|{T_{r_\e}}|} \bfg_0 \lp\hskip-0,05cm x,\frac {y_\e(x')}{r_\e}\rp\hskip-0,05cm\mathds{1}_{{T_{r_\e}}}(x) + \frac{1}{r_\e}  \frac{|\Omega|}{|{T_{r_\e}}|}  \bfh_\e(x) 
\ & &\hbox{ if } \ 0<\kappa<+\infty,  
\end{aligned} 
\label{fe2}
   \end{equation}
with
     \begin{equation}
\hskip-0,2cm \begin{aligned}
     \bfh_\e(x)\hskip-0,1cm:=\hskip-0,1cm\sum_{i\in I_\e}\hskip-0,1cm \begin{pmatrix}\hskip-0,1cm -  \frac {y_{\e2}(x')}{r_\e}\beta_0\lp x,\frac {y_{\e}(x')}{r_\e}\rp+{\int\!\!\!\!-}_{S_{r_\e}^{i} }\frac {y_{\e2}(s)}{r_\e}
 \beta_0\lp s,x_3,\frac {y_{\e}(s)}{r_\e} \rp ds 
\\\hskip-0,1cm \frac {y_{\e1}(x')}{r_\e} \beta_0\lp x,\frac {y_{\e}(x')}{r_\e}\rp  -{\int\!\!\!\!-}_{S_{r_\e}^{i} }\frac {y_{\e1}(s)}{r_\e}
 \beta_0\lp s,x_3,\frac {y_{\e}(s)}{r_\e} \rp ds 
 \\\hskip-0,1cm a_0\lp x,\frac {y_{\e}(x')}{r_\e}\rp  \end{pmatrix} \hskip-0,05cm  \mathds{1}_{{T_{r_\e}^{i}}}\hskip-0,1cm(x).
\end{aligned} 
\nonumber
   \end{equation}
 We establish that  the  limit problems associated  to 
\eqref{Pefe} and \eqref{Pesoftfe}      are
     \begin{equation}
\begin{aligned}
 \inf_{(\bfu,\bfv^{tuple})\in W^{1,p}_b(\Omega;\RR^3)\times \D} \Phi (\bfu, \bfv^{tuple})-\int_\Omega\bff\cdot\bfu dx -L( \bfv^{tuple}), \end{aligned} 
 \label{PhomL}
    \end{equation} 
     \begin{equation}
\begin{aligned}
 \inf_{(\bfu,\bfv^{tuple})\in L^p(\Omega;\RR^3)\times \D} \Phi^{soft} (\bfu, \bfv^{tuple})-\int_\Omega\bff\cdot\bfu dx -L( \bfv^{tuple}), \end{aligned} 
 \label{PhomsoftL}
    \end{equation} 
 \noindent   respectively,  where,  if $\kappa=0$, 
     \begin{equation}
\hskip-0,2cm  \begin{aligned}
&        L(\bfv^{tuple})=\int_\Omega \lp \intb_S  \hskip-0,1cm\bfg_0(x,y)  dy   \rp   \hskip-0,1cm \cdot \hskip-0,05cm \bfv    
       + \lp \intb_S \bfg_0(x,y)  \hskip-0,05cm \cdot \hskip-0,05cm  \frac{2}{\diam S} \bfe_3\wedge   \bfy  dy \!  \rp  \hskip-0,05cm \theta   dx,  
   \end{aligned}\label{L1}
     \end{equation}  
     and,  if $0<\kappa<+\infty$, 
 \begin{equation}
\hskip-0,2cm  \begin{aligned}
&          L(\bfv^{tuple})\hskip-0,1cm=  \hskip-0,1cm  \int_\Omega       \hskip-0,05cm \lp \intb_S \hskip-0,1cm\bfg_0(x,y)  dy   \rp 
   \hskip-0,1cm \cdot \bfv      
  +
 \tau  \lp \intb_S \hskip-0,1cm\beta_0(x,y) dy\rp  \delta(x)  
     \\&  \hskip2cm  + \lp \intb_S  \lp w(x)- \sum_{\alpha=1}^2\frac{\partial v_\alpha}{\partial x_3}(x) y _\alpha \rp a_0(x,y)  dy\rp dx,
     \end{aligned} \label{L2}  \end{equation}  
being $\tau$ given by \eqref{defvethetae}. 
%
 \begin{theorem}\label{thL} 
 The   statements deduced  from Theorems \ref{th} and \ref{thsoft}
   by substituting 
  \eqref{Pefe},  \eqref{PhomL} for 
    \eqref{Pe}, \eqref{Phom}   and \eqref{Pesoftfe}, \eqref{PhomsoftL} for 
    \eqref{Pesoft}, \eqref{Phomsoft}   hold.
  \end{theorem}
  
   \begin{remark}\label{remdeltaw} 
When  $\kappa>0$ and   $\Phi_{fibers}$ is given by \eqref{Phiisotkappa}, 
if   $(\bfu,\bfv^{tuple})$ is a solution to  \eqref{PhomL}  or \eqref{PhomsoftL}, then
          \begin{equation}
\begin{aligned}
& \delta(x)=\frac{\tau (\diam S)^2}{  \kappa  \mu_1 m  }\lp \intb_S  \beta_0(x',L,y) dy\rp \lp -\frac{x_3^2}{2}+L\rp,
\\&  w(x)= \frac{l+1}{ \kappa  \mu_1(3l+2)}\lp \intb_S  a_0(x',L,y)  dy\rp \lp -\frac{x_3^2}{2}+L\rp.
\ \end{aligned} 
 \label{deltaw=}
    \end{equation} 
       \end{remark}

 \section{Preliminary results and a priori estimates}\label{secpreliminaries}
  %
 This section is devoted to the study of the asymptotic behaviors   of a   sequence $(\bfu_\e)$ satisfying 
\begin{equation}\begin{aligned}
&
\sup_{\e>0} F_\e(\bfu_\e) <+\infty, 
\end{aligned}
  \label{supFeuefini}
  \end{equation}
  and    of  the auxiliary    sequences $(\bfv_\e)$,  $(\theta_\e)$, $\lp\frac{v_{\e3}}{r_\e}\rp$,   $\lp\frac{\theta_\e}{r_\e}\rp$
 defined by  \eqref{defvethetae}.  Our main results, stated in  Section \ref{secapriori},
 will be deduced from   a series   of  inequalities proved   in Section \ref{seckey}.
   \subsection{Auxiliary sequences}\label{secaux}
\begin{figure}[!h]
  \centering
 \includegraphics[height=5cm]{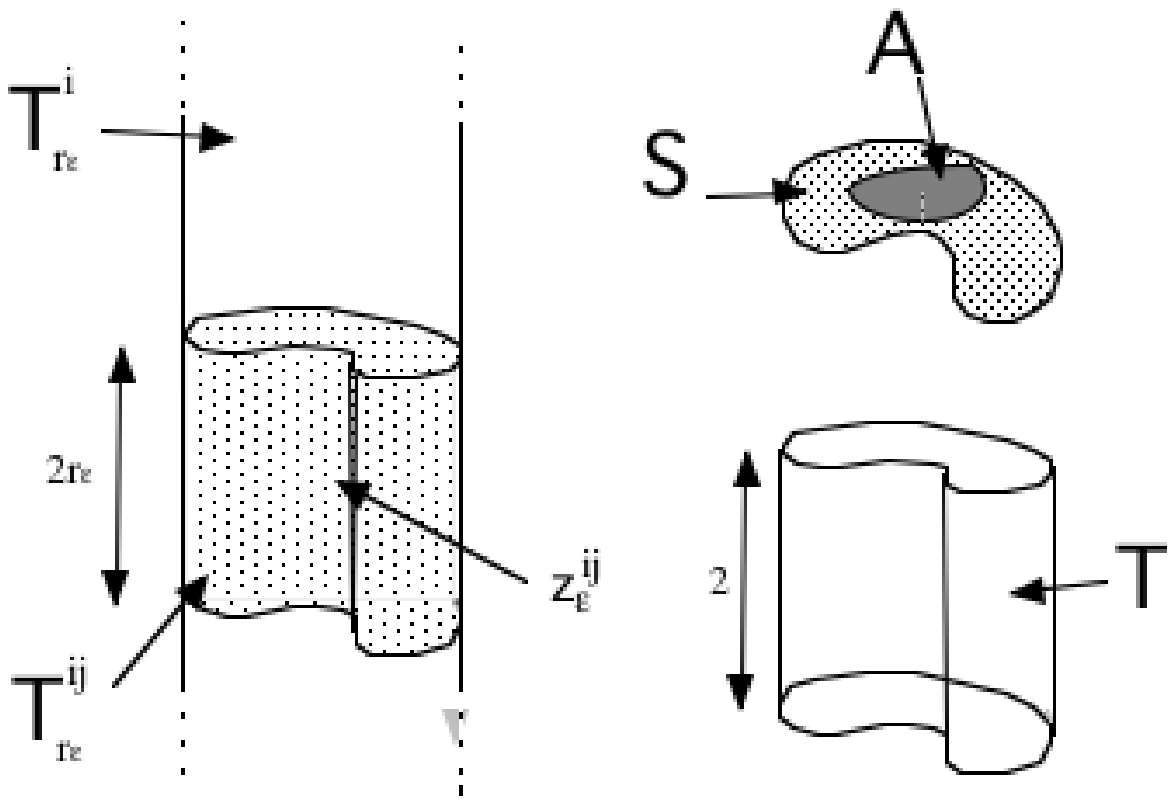}  
  \caption{ The set  $T_{r_\e}^{ij}$.}
  \label{fig2}
  \end{figure}

   The proof of Theorem \ref{th} rests on
  a   choice of  sequences of auxiliary   fields 
    \begin{itemize}
  \item   weakly differentiable w.r.t. $x_3$,
  \item   weakly relatively compact in $L^p$, 
  \item  locally characterizing a rigid motion 
   approximating  to the velocity  in the fibers.
  \end{itemize}
\noindent  The selection of surface integrals in   \eqref{defvethetae}  is  required  to ensure  the first condition, and  is a  hindrance  to  the third.
To  circumvent this  difficulty, we introduce another  couple of auxiliary sequences 
constructed from volume  instead of surface  integrals,
 satisfying the second and third  conditions    and having  the same cluster points as $\bfv_\e$, $\theta_\e$.
For convenience, we  rewrite 
 \eqref{defvethetae} as follows 
    \begin{equation}\begin{aligned}
&\bfv_\e  = \ov\bfgotv_\e^{S} (\bfu_\e), \qquad \theta_\e  = \ov\theta_\e^{S}(\bfu_\e),
\end{aligned} 
\label{ve=gotve}
   \end{equation}
where, for any  Lipschitz  domain $A$ of $\RR^2$  verifying
  \begin{equation}\begin{aligned}
&D\subset A\subset S, \qquad  \bfy_A=0,
\end{aligned} 
\label{hypA}
   \end{equation}
  $ \ov\bfgotv_\e^{A}$ and $\ov\theta_\e^{A}$
are the linear  operators     defined on $L^p(\Omega;\RR^3)$    by     (see  \eqref{defTre},  \eqref{defAbi},  \eqref{defye}):
    \begin{equation}\begin{aligned}
   &   \ov \bfgotv^{A}_\e (\bfvarphi)(x ) :=\sum_{i\in I_\e} \left (
 {\int\!\!\!\!\!\!-}_{{ A_{r_\e}^{i}}}
\bfvarphi (s, x_3  )   d\calH^{2}(s) \right)  \mathds{1}_{Y_\e^{i}}(x' ),
\\&    \ov\theta^{A}_\e (\bfvarphi)(x  ): = 
  \sum_{i\in I_\e} 
 \frac{ \diam S}{2  {\int\!\!\!\!-}_A |\bfy|^2 dy }   {\int\!\!\!\!\!\!-}_{{  A_{r_\e}^{i}}}  
\hskip-0,3cm  \lp  -  \frac{y_{\e2}}{r_\e }  \varphi_1+ \frac{y_{\e1}}{r_\e } \varphi_2 \rp(s,x_3) d\calH^2(s  ) \    \mathds{1}_{Y_\e^{i}}(x'),
 \\ & \ov\bfgotr_\e^{A}(\bfvarphi)(x):=  \ov\bfgotv_\e^{A}(\bfvarphi)(x)  +\frac{2}{\diam S}  \ov\theta_\e^{A} (\bfvarphi)( x )\bfe_3  \wedge   \frac{\bfy_\e(x')}{r_\e}.
\end{aligned}
\label{defovgotA}
  \end{equation}
The "volumic" couterparts $\bfgotv_\e(\bfu_\e)$ and $\theta_\e(\bfu_\e)$ of \eqref{defvethetae} will be defined by splitting 
 each fiber $T^{i}_{r_\e}$ into small cylinders $T_{r_\e}^{ij} $ of size $r_\e$  given by
 (see     fig. \ref{fig2})
\begin{equation}
  \hskip-0,1cm  \begin{aligned}
  &    T_{r_\e}^{ij} :=   S_{r_\e}^{i} \hskip-0,1cm \times \Delta_{r_\e}^j,   \quad  \Delta_{r_\e}^j:=  r_\e\lp 2j+ \lp0,2\rc\rp  \hskip-0,05cm, \quad L_\e:= \la j\in \ZZ, \  \Delta_{r_\e}^j  \cap (0,L)\not=\emptyset \ra\hskip-0,1cm,
 \\&\bfz_\e( x):= \hskip-0,2cm \sum_{(i,j)\in I_\e\times L_\e } \hskip-0,3cm \mathds{1}_{S_{r_\e}^{i} \times \Delta_{r_\e}^j}(x )
  \lp \bfx-\bfz_\e^{ij} \rp, \  
\quad\bfz_\e^{ij}\hskip-0,1cm = \hskip-0,1cm \e (i_1\bfe_1+i_2\bfe_2)+    (2 j +1)r_\e \bfe_3,
    \end{aligned} 
\label{defze}
  \end{equation}
  and setting for every $\bfvarphi\in L^p(\Omega;\RR^3)$,
   \begin{equation}\begin{aligned}
    & \bfgotv_\e(\bfvarphi)(x):=\hskip-0,3cm \sum_{(i,j)\in {I}_\e\times L_\e} \hskip-0,2cm\left (
 {\int\!\!\!\!\!\!-}_{T_{r_\e}^{ij}}  \widehat\bfvarphi(s )   d\calL^3(s) \right)  \mathds{1}_{Y_\e^{i}\times \Delta_{r_\e}^j }(x );\quad 
\theta_\e(\bfvarphi):= \gotw_{\e3}(\bfvarphi),
\\&
\bfgotw_\e(\bfvarphi)  (x  )\!:=\hskip-0,5cm   \sum_{(i,j)\in {I}_\e\times L_\e}\hskip-0,2cm \left (
\frac{\diam S}{2}\bfM{\int\!\!\!\!\!\!-}_{T_{r_\e}^{ij}} \!\!\! \lp
   \frac{\bfz_\e(s)}{r_\e} \wedge \widehat\bfvarphi(s ) \rp \hskip-0,1cm d\calL^3(s) \hskip-0,1cm\right) \hskip-0,1cm  \mathds{1}_{Y_\e^{i}\times \Delta_{r_\e}^j }(x ),  
    \\&  \bfM:=   \sum_{i,j,k=1}^3 \frac{1}{2} \e_{ijk}^2 \lp\intb_T z_i^2+z_j^2 d\calL^3(z)\rp^{-1}\bfe_k\otimes\bfe_k,   
 \\&  \bfgotr_\e(\bfvarphi) := \bfgotv_\e(\bfvarphi)+ \frac{2}{\diam S}\bfgotw_\e(\bfvarphi)\wedge  \frac{\bfz_\e(x)}{r_\e},
\end{aligned}
\label{defgot}
  \end{equation}
  where $\widehat\bfvarphi$  is an extension of $\bfvarphi$
    (see \eqref{defhatfi})
      bringing    meaning to
 $ {\int\!\!\!\!-}_{T_{r_\e}^{ij^m}} \widehat \bfvarphi$ 
 where  
    \begin{equation}\begin{aligned} 
  j^m_\e:=\max L_\e. 
 \end{aligned}
\label{defjme}
  \end{equation} 
We will prove that the    piecewise   rigid motion  $\bfgotr_\e(\bfu_\e)$   
  approximates to the velocity  in each set $T_{r_\e}^{ij}$ and    asymptotically  behaves 
as  $\ov\bfgotr_\e^S(\bfu_\e)$.
To derive     formula \eqref{defcfintro}, 
we  will employ the  auxiliary field 
 $\ov\bfgotu_\e(\bfu_\e)$
   defined   by  
\begin{equation}\begin{aligned}
&  
\ov\bfgotu_\e(\bfvarphi) (x )  \!:=\!
\sum_{i\in I_\e} \left (
\intb_{D_{R_\e}^{i}\setminus D_{R_\e\!/2}^{i}}\!\!\bfvarphi(s, x_3  ) d\calH^{2}(s) \right )  \mathds{1}_{Y_\e^{i}}(x' )
 \end{aligned}
\label{defovgotue}
  \end{equation}
\noindent  
in terms of a sequence  $(R_\e)$ satisfying \eqref{Re}. The assumption \eqref{Re} 
  ensures  that $(\ov\bfgotu_\e(\bfu_\e))$ approximates to  $\bfu_\e$ in $\Omega$.
      The   selection of integrals over $D_{R_\e}^{i}\setminus D_{R_\e\!/2}^{i}$ in \eqref{defovgotue} 
will yield a crucial estimate at  a technical step.
%
\subsection{Two-scale convergence with respect to $(m_\e)$}\label{sectwoscale}
To particularize   the effective behavior of $\bfu_\e$    in the fibers,  
  we introduce  variant   of  the   two-scale convergence   \cite{Al,CiDaGr,Ng}
  defined in terms of the measures 
    \begin{equation}\begin{aligned}
m_\e:= \frac{\e^2}{r_\e^2|S|}   \mathds{1}_{S_{r_\e} \times (0,L)} \calL^3_{\lfloor \Omega}
\end{aligned}
\label{defme}
  \end{equation}
 which are    supported on the fibers and   weakly$^\star$ converging  in $ \M(\ov\Omega)$ to $\calL^3_{\lfloor \Omega}$.
We say that a  sequence 
 $(\bfvarphi_\e)$ in $L^p(\Omega;\RR^k) $   
 two-scale converges    to   $\bfvarphi_0\in L^p (\Omega\times S;\RR^k) $  with respect to $(m_\e)$ if
 (see \eqref{defye})
\begin{equation}\begin{aligned}
\lim_{\e \to 0}\int  \bfvarphi_\e(x ) \cdot \bfeta\lp x,  \frac{y_\e(x')}{r_\e}\rp & dm_\e(x)   = \frac{1}{|S|} \int_ {\Omega\times S}  \bfvarphi_0(x,y) \cdot \bfeta(x,y)  dx dy , 
\\&\ \forall\,   \bfeta\in 
C(\ov{\Omega\times S};\RR^k)  \quad \hbox{(notation: $ \bfvarphi_\e\mdto  \bfvarphi_0$)}.
\end{aligned}
\label{defmdto}
  \end{equation}
The  main properties  of this convergence are stated in the following  Lemma.
   \begin{lemma} \label{lemtwoscale}  
For every $\bfeta\in
C(\ov{\Omega}; L^\infty(S;\RR^k))$,  
  \begin{equation}\begin{aligned}
 \lim_{\e \to 0}\int   \bfeta\lp x,  \frac{y_\e(x')}{r_\e}\rp  dm_\e(x)  = \frac{1}{|S|} \int_ {\Omega\times S}  \bfeta(x,y) dx d y. 
  \end{aligned}
\label{twoscalecontinuous}
  \end{equation}  
Any  sequence  $(\bfvarphi_\e)$  in $ L^p(\Omega;\RR^k) $ 
such that 
\begin{equation}
\begin{aligned}
 \sup_{\e>0} \int |\bfvarphi_\e|^p  dm_\e \leq C, 
\end{aligned}
\label{hypfie}
  \end{equation}
  \noindent  
  two-scale converges w.r.t. $(m_\e)$,     up to a subsequence, to some 
   $\bfvarphi_0$.
Setting $ \bfvarphi(x) :=  \frac{1}{|S|}\int_S \bfvarphi_0(x,y) dy$, the following implications  holds: 
\begin{eqnarray}
 & &
    \bfvarphi_\e \mdto   \bfvarphi_0 \quad\Longrightarrow  \quad \bfvarphi_\e \bfeta\lp x,  \frac{y_\e(x')}{r_\e}\rp  \mdto \   \bfvarphi_0\bfeta\qquad  \forall \bfeta\in
C(\ov{\Omega\times S};\RR^k),
\label{twoscaleeta}
\\
 & &
    \bfvarphi_\e \mdto   \bfvarphi_0 \quad\Longrightarrow  \quad \bfvarphi_\e m_\e \buildrel \star \over \rightharpoonup \bfvarphi \ \ \hbox{weakly$^\star$ in }\  \M(\ov\Omega;\!\RR^k).
\label{twoscalemeas}
  \end{eqnarray}
\noindent  If $ \bfvarphi_\e \mdto   \bfvarphi_0$, for any convex function   $j:\RR^k \to \RR$, we have 
\begin{equation}
\liminf_{\e\to 0} \int  j(\bfvarphi_\e)  dm_\e 
  \ge\frac{1}{|S|} \int_ {\Omega\times S} \hskip-0,2 cm j(  \bfvarphi_0)  dx d y \ge  \int_\Omega  j( \bfvarphi )  dx.
\label{linfcarre}
  \end{equation}
\end{lemma}

\noindent
\noindent  {\bf Proof.}   The proof of \eqref{twoscalecontinuous}  is similar to that of \cite[Lemma  1.3]{Al}.  Setting $\int \bfeta\cdot d\bfnu_\e:=   \int \bfvarphi_\e(x ) \cdot  \bfeta \lp x, \frac{y_\e(x')}{r_\e}\rp dm_\e(x)$ $\forall \bfeta\in C (\ov{\Omega\times S};\RR^k)$,
the 
two-scale convergence w.r.t. $(m_\e)$ of $(\bfvarphi_\e)$ to $\bfvarphi_0$
is equivalent to 
the weak$^\star$ convergence  in $\M(\ov{\Omega\times S})$ of  $(\nu_\e)$  to $\frac{1}{|S|}  \bfvarphi_0(x,y)\calL^5_{\lfloor\Omega\times S}$.
  By   \eqref{hypfie} and  H\"older's  inequality,  we have   
\begin{equation}
\begin{aligned}
 \lb \int \bfeta\cdot  d\bfnu_\e\rb  
 \leq C     \lp \int  \lb\bfeta \lp x, \frac{y_\e(x')}{r_\e}\rp \rb^{p'}  dm_\e(x)\rp^{\frac{1}{p'} }\!\! \leq C |\bfeta |_{L^\infty(\Omega\times S;\RR^k)}, 
\end{aligned}
\label{psid}
  \end{equation}
thus
   $(\bfnu_\e)$,  bounded
 in $ \M(\ov{\Omega\times S};\RR^k)$, 
   weakly$^\star$  converges   up to a subsequence  to some $\bfnu$. 
Passing to the limit  in   (\ref{psid}), taking 
\eqref{twoscalecontinuous} into account,  we obtain  $
\lb \int \bfeta\cdot  d\bfnu  \rb\leq  C     \lp \int_{\Omega\times S} \lb\bfeta \lp x,  y \rp \rb^{p'} dxdy \rp^{\frac{1}{p'} }$, hence, by 
   the Riesz representation theorem  
$\bfnu=\frac{1}{|S|}  \bfvarphi_0(x,y)\calL^5_{\lfloor\Omega\times S}$ for some $\bfvarphi_0\in L^p(\Omega\times S;\RR^k )$.
The assertion \eqref{twoscaleeta} is straightforward.
By  choosing   in \eqref{defmdto}     test fields  independent of $y$, 
 we obtain \eqref{twoscalemeas}.
   Denoting by $j^\star$ the Fenchel transform of $j$,  
we deduce from   Fenchel's inequality,  \eqref{defmdto} and   (\ref{twoscalecontinuous}), that 
\begin{equation}\begin{aligned}
\liminf_{\e\to 0}    \int  \! j( \bfvarphi_\e) dm_\e(x)  
& \ge \lim_{\e\to 0}   \int   \bfvarphi_\e\cdot \bfeta\lp x, \frac{y_\e(x')}{r_\e}\rp
 \!- j^\star \lp  \bfeta\lp x, \frac{y_\e(x')}{r_\e}\rp \hskip-0,1cm\rp  \hskip-0,1cmdm_\e(x)  
\\ & =\frac{1}{|S|} \int_ {\Omega\times S}   \bfvarphi_0\cdot \bfeta-j^\star(\bfeta)   dx d y \qquad \forall \bfeta\in C (\ov{\Omega\times S};\RR^k).
\end{aligned}
   \nonumber \end{equation}
  \noindent  
\noindent  Noticing that by the convexity of  $j$  we have   $j^{\star\star}=j$, we infer  
 \begin{equation}\begin{aligned}
 \liminf_{\e\to 0}    \int  \! j( \bfvarphi_\e) dm_\e(x)  
 & \ge \sup_{\bfeta\in C (\ov{\Omega\times S};\RR^k)} \frac{1}{|S|} \int_ {\Omega\times S}  \bfvarphi_0\cdot \bfeta-j^\star(\bfeta)   dx d y
\\& = \frac{1}{|S|} \int_ {\Omega\times S} \  \sup_{\bfeta\in \RR^k } \{ \bfvarphi_0\cdot \bfeta-j^\star(\bfeta) \} dx d y
 \\&= \frac{1}{|S|}  \int_ {\Omega\times S}   j^{\star\star}(\bfvarphi_0) dx d y
= \frac{1}{|S|} \int_ {\Omega\times S}  j (\bfvarphi_0) dx d y,
\end{aligned}
\label{lirem}\end{equation}
  the second line being justified in   Remark \ref{remghom}.
The second inequality in \eqref{linfcarre} results from \eqref{twoscalemeas} and  Jensen's inequality.    \qed

    \begin{remark}\label{remtwoscalenonper} 
Combining Lemma \ref{lemtwoscale} with \cite[Lemma  4.2]{BeSIAM2},  one can derive   
 a  non-periodic version of  Lemma \ref{lemtwoscale} adapted for the proof of Theorem \ref{thnonper}.
    \end{remark}
    
%
\subsection{Properties of the auxiliary sequences}\label{secop}
     
     In this section, we compare  diverse types of convergence for the    auxiliary sequences   $(\bfgotp_\e(\bfu_\e))$ introduced in Section \ref{secaux}.  It turns out    that,  except for 
 $\bfgotp_\e\in \{ \bfgotr_\e,\ov\bfgotr_\e^{S}\}$,  the  weak limit  of  $(\bfgotp_\e(\bfu_\e))$  in $L^p$,  
 its   two-scale limits   w.r.t. to 
 $(m_\e)$ and the    weak* limits  of $(\bfgotp_\e(\bfu_\e)m_\e)$ in $\M(\ov\Omega)$  are the same, and  that  $(\ov\bfgotv_\e^S (\bfu_\e))$ and $(\bfgotv_\e(\bfu_\e))$ have 
  the same  weak limits in $L^p$.
We will prove in  Proposition \ref{prop1}  that the same holds 
for $(\ov\theta_\e^S (\bfu_\e)\bfe_3)$ and $(\bfgotw_\e(\bfu_\e))$. This, combined with  \eqref{gotrdto}, shows  that  $(\ov\bfgotr_\e^S (\bfu_\e))$ and  $(\bfgotr_\e (\bfu_\e))$  have  the same two-scale limits w.r.t. $(m_\e)$. 
   
  \begin{proposition}\label{lemmdx}
  \noindent Let $\bfgotp_\e$  be a  linear combination  of 
 $\bfgotw_\e$, $\bfgotv_\e$, $\theta_\e\bfe_3$,  ${\ov\bfgotv}^{A}_\e$, ${\ov \theta}^{A}_\e\bfe_3$, $\ov\bfgotu_\e$  
defined by    \eqref{defovgotA},
  \eqref{defgot},   
  \eqref{defovgotue}, and  $(\bfvarphi_\e)$    a sequence in $L^p(\Omega;\RR^3)$ satisfying \eqref{hypfie}. 
  \\ (i)  For all   $\bfvarphi \in W^{1,p}(\Omega;\RR^3)$,   $\e>0$,  $k\in\{1,2,3\}$, we have 
    \begin{equation}\begin{aligned}
&   \int |\bfgotp_\e (\bfvarphi)|^p dm_\e  =   \int_\Omega|\bfgotp_\e (\bfvarphi)|^p dx,
\end{aligned}
\label{mdx=}
  \end{equation}
  \begin{equation}\begin{aligned}
&  \int |\bfgotp_\e (\bfvarphi)|^p dm_\e   \leq C  \int_\Omega|\bfvarphi|^p dm_\e, 
\quad 
 \int |\ov \gotv^{A}_{\e k} (\bfvarphi)|^p dm_\e  \leq \frac{|S|}{|A|}  \int_\Omega|\varphi_k|^p dm_\e.
\end{aligned}\hskip-0,1cm 
\label{mdx}
  \end{equation}
  \\  (ii) For every  $\bfp\in L^p(\Omega;\RR^3)$,  the following equivalences hold 
 \begin{equation}\begin{aligned}
  \hskip-0,2cm 
 \bfvarphi_\e m_\e \ \buildrel\star\over\rightharpoonup\ \bfp  \hbox{ weakly$^\star$ in }  \M(\ov\Omega;\RR^3)
&  \Longleftrightarrow\  \  \ov\bfgotv_\e^S (\bfvarphi_\e)\! \rightharpoonup\! \bfp  \hbox{ weakly  in } L^p(\Omega;\RR^3)
\\&  \Longleftrightarrow\  \  \bfgotv_\e  (\bfvarphi_\e)\! \rightharpoonup\! \bfp  \hbox{ weakly  in } L^p(\Omega;\RR^3).
\end{aligned}
\label{mdx0}
  \end{equation}
  \noindent (iii) If $(\bfgotp_\e (\bfvarphi_\e))$    two-scale converges to $\bfp_0$ w.r.t. $(m_\e)$, then  
     $\bfp_0(x,y)=\bfp(x)$ a. e. in $\Omega\times S$ for some $\bfp\in L^p(\Omega;\RR^3)$.
Furthermore, for every $\bfp\in L^p(\Omega;\RR^3)$,  
   \begin{equation}\begin{aligned}
  \hskip-0,2cm 
 & \bfgotp_\e (\bfvarphi_\e) \mdto \bfp \quad  \Longleftrightarrow   \quad
\bfgotp_\e (\bfvarphi_\e)m_\e\ \buildrel\star\over\rightharpoonup\ \bfp  \hbox{ weakly$^\star$ in }  \M(\ov\Omega;\RR^3)
\\&  \Longleftrightarrow\     \bfgotp_\e (\bfvarphi_\e)\! \rightharpoonup\! \bfp  \hbox{ weakly  in } L^p(\!\Omega;\RR^3)
  \Longrightarrow     \lp \bfgotp_\e (\bfvarphi_\e)\wedge  \frac{\bfy_\e(x')}{r_\e} \rp  \mdto\    \bfp \wedge  \bfy.
\end{aligned}
\label{mdx2}
  \end{equation}
(iv) If   $\,\bfgotp_\e$  is    a  linear combination  of 
 $\bfgotw_\e$, $\bfgotv_\e$, $\theta_\e\bfe_3$, then 
   \begin{equation}\begin{aligned}
    \bfgotp_\e (\bfvarphi_\e)\! \rightharpoonup\! \bfp  \hbox{ weakly  in } L^p(\!\Omega;\RR^3) \quad
&  \Longrightarrow      \lp\bfgotp_\e (\bfvarphi_\e)\wedge  \frac{\bfz_\e(x)}{r_\e} \rp \mdto   \bfp \wedge  \bfy, 
\end{aligned}
\label{mdx222}
  \end{equation} 
 where $\bfz_\e(x)$ is defined by \eqref{defze}.
 In particular,   the following holds:
\begin{equation}
\hskip-0,2cm \begin{aligned} 
&   \lpt \begin{aligned}&\bfgotv_\e (\bfvarphi_\e)\! \rightharpoonup\! \bfv & &\hbox{weakly  in } L^p(\Omega;\RR^3)
\\&   \bfgotw_\e (\bfvarphi_\e)\! \rightharpoonup\! \bfw    & & \hbox{weakly  in } L^p(\Omega;\RR^3)
\end{aligned} \ra
\Longrightarrow \bfgotr_\e(\bfvarphi_\e) \mdto \bfv+\frac{2}{\diam S}\bfw\wedge \bfy .
\end{aligned}
\label{gotrdto}
  \end{equation}
   \end{proposition}

  \noindent    {\bf Proof.} (i)  Since  $\bfgotp_\e (\bfvarphi)$ is  constant in  each set $Y_\e^{i}\times \{x_3\}$, by  \eqref{defTre}  and  \eqref{defme},
 \begin{equation}\begin{aligned}
 \int |\bfgotp_\e (\bfvarphi)|^p dm_\e  &=\frac{\e^2}{r_\e^2|S|} \int_{T_{r_\e}} |\bfgotp_\e (\bfvarphi)|^p dx
 =\frac{\e^2}{r_\e^2|S|} \sum_{i\in I_\e} \int_0^L dx_3 \int_{S_{r_\e}^{i}}   |\bfgotp_\e (\bfvarphi)|^p dx'
 \\&=\frac{\e^2}{r_\e^2|S|} \sum_{i\in I_\e} \int_0^L dx_3 \frac{|S_{r_\e}^{i}|}{|Y_\e^{i}|} \int_{Y_\e^{i}}   |\bfgotp_\e (\bfvarphi)|^p dx'
 =   \int_{\Omega} |\bfgotp_\e (\bfvarphi)|^p dx.
\end{aligned}
\nonumber
  \end{equation}
  \noindent  We deduce from  \eqref{hypA},   \eqref{defovgotA},  \eqref{defme}  and  Jensen's inequality, that
  for $k\in\{1,2,3\}$,
   \begin{equation}\begin{aligned}
   \int |\ov\gotv^{A}_{\e k} (\bfvarphi)|^p dm_\e&= \frac{\e^2}{ r_\e^2|S|}\sum_{i\in I_\e}
   \int_0^L dx_3 \int_{S_{r_\e}^{i}} \lb \intb_{A_{r_\e}^{i}} \varphi_k(s',x_3) d\calH^2(s')   \rb^p dx'
   \\&\hskip-2cm \leq  \frac{\e^2}{r_\e^2|S| }\sum_{i\in I_\e}
   \int_0^L dx_3 \frac{|S_{r_\e}^{i}|}{|A_{r_\e}^{i}|}  \int_{S_{r_\e}^{i}} | \varphi_k|^p(s',x_3) d\calH^2(s')   
  =  \frac{|S|}{|A|}  \int_\Omega|\varphi_k|^p dm_\e .
\end{aligned}
\nonumber
  \end{equation}
  \noindent  The inequality $ \int |\bfgotp_\e (\bfvarphi)|^p dm_\e  \leq C  \int_\Omega|\bfvarphi|^p dm_\e$ 
  is   obtained in a similar way.  \qed

\noindent (ii) We   fix   $\bfpsi\in \D(\Omega;\RR^3)$   and set (see \eqref{defze}, \eqref{defjme}) 
\begin{equation}\
\bfpsi^{(\e)}(x):= \sum_{(i,j)\in I_\e\times (L_\e\setminus\{j^m_\e\})}\lp \intb_{Y_\e^{i}\times \Delta_\e^j} \bfpsi  d\calL^3\rp \mathds{1}_{Y_\e^{i}\times \Delta_\e^j} (x). 
\label{defpsie}
\end{equation}
As  $\bfpsi$ is   compactly supported in $\Omega$, the following estimate holds: 
\begin{equation}
  |\bfpsi^{(\e)}-\bfpsi|_{L^\infty(\Omega;\RR^3)}<C\e.
\label{estimpsie}
\end{equation} 
Taking   \eqref{defovgotA}, \eqref{defgot}, \eqref{defovgotue}, \eqref{defme}, and the constant nature of   $\bfpsi^{(\e)}$  in each set $Y_\e^{i}\times  \Delta_\e^j$ into account,  elementary computations yield
 \begin{eqnarray}
& &\begin{aligned}
   &  \int   \bfpsi^{(\e)} \hskip-0,1cm\cdot \ov\bfgotv^S_\e  (\bfvarphi_\e) dm_\e 
   = \int \bfpsi^{(\e)} \cdot  \bfvarphi_\e dm_\e,
   \\& \int \bfpsi^{(\e)} \cdot  \bfgotv_\e(\bfvarphi_\e) dm_\e= \int \bfpsi^{(\e)} \cdot  \bfvarphi_\e dm_\e, 
\end{aligned}
\label{psievpsievarphi}
 \\& &\begin{aligned}
  \int   \bfpsi^{(\e)} \cdot \ &  \bfgotp_\e (\bfvarphi_\e)  dm_\e 
  = \int_\Omega   \bfpsi^{(\e)} \cdot  \bfgotp_\e (\bfvarphi_\e) dx.
\end{aligned}
\label{psiepdmpsiepdx}
\end{eqnarray} 
  %
By \eqref{hypfie}, \eqref{mdx} and 
  Lemma \ref{lemtwoscale},
   the sequences $(\bfvarphi_\e m_\e)$, $(\ov\bfgotv_\e^{S}(\bfvarphi_\e) m_\e)$, and $(\bfgotv_\e(\bfvarphi_\e)m_\e)$ weakly$^\star$  converge in  $ \M(\ov\Omega;\RR^3)$,    up to a subsequence,  
  to $\bfp_1$, $\bfp_2$, and $\bfp_3$, respectively, for some $\bfp_i\in L^p(\Omega;\RR^3)$.
By passing to the limit as $\e\to0$ in  \eqref{psievpsievarphi}, taking  \eqref{estimpsie}  into account, 
we infer $\int_\Omega\bfpsi\cdot  \bfp_1 dx=\int_\Omega\bfpsi\cdot  \bfp_2 dx=\int_\Omega\bfpsi\cdot  \bfp_3 dx$ and deduce from  the arbitrariness of $\bfpsi$  that 
$\bfp_1 =\bfp_2 =\bfp_3 $.  \qed

%
\noindent  (iii) Assume that   $\bfgotp_\e (\bfvarphi_\e)\hskip-0,05cm  \mdto\bfp_0$ for some $\bfp_0\in L^p(\Omega\times Y;\RR^3)$. Since 
 $\bfgotp_\e(\bfvarphi_\e) $ is constant    in each set $Y_\e^{i}\times\{x_3\}$,   for any  
   $\bfPsi \in \D(\Omega\times S ;\RR^3\times\RR^3)$  we have   $\displaystyle \int \bfgotp_\e(\bfvarphi_\e)  \cdot\bfdiv_y \bfPsi\Big(
x,\frac {y_\e(x')}{r_\e}\Big)dm_\e  =0$. Passing to the limit as $\e\to 0$,   we get 
$ \displaystyle  \frac{1}{|S|}\int_{\Omega\times S}  \hskip-0,4cm   \bfp_0\cdot\bfdiv_y \bfPsi  dx dy=0$ 
and deduce  from the arbitrary choice of $\bfPsi$ that    $\bfp_0(x,y)=\bfp(x)$ a.e. in $\Omega\times S$, where $\bfp(x):= \intb_S \bfp_0(x,y)dy \in L^p(\Omega;\RR^3)$.  It  follows from \eqref{twoscalemeas} that 
 $\bfgotp_\e (\bfvarphi_\e)m_\e\!\buildrel\star\over\rightharpoonup\! \bfp \calL^3$ weakly$^\star$ in $  \M(\ov\Omega;\RR^3)$. 
  On the other hand, by \eqref{hypfie}, \eqref{mdx=}, and \eqref{mdx},  $(\bfgotp_\e(\bfvarphi_\e))$ is bounded in $L^p(\Omega;\RR^3)$ and weakly converges, up to a subsequence, to some $\tilde\bfp$.
   By passing to the limit as $\e\to0$ in \eqref{psiepdmpsiepdx}, taking \eqref{estimpsie}   into account, we obtain  $\int_\Omega  \bfpsi \cdot \bfp dx= \int_\Omega  \bfpsi \cdot \tilde \bfp dx$ and deduce  $\bfp=\tilde\bfp$.    Conversely, if  $\bfgotp_\e (\bfvarphi_\e)\rightharpoonup \bfp$ weakly   in  $L^p(\Omega;\RR^3)$, then by  \eqref{hypfie},  \eqref{mdx}     and Lemma \ref{lemtwoscale},    $(\bfgotp_\e (\bfvarphi_\e))$ two-scale converges w.r.t. $(m_\e)$, up to a subsequence, to some  
  element of $L^p(\Omega;\RR^3)$ which, by virtue of the last established implications,  necessarily equals $\bfp$. 
 The equivalences stated  in  \eqref{mdx2} are  proved and the implication    is obtained 
by   choosing $\bfy  \wedge  \bfeta\lp x,   y \rp$   as a test function  for the two-scale convergence of $(\bfgotp_\e(\bfvarphi_\e))$ to $\bfp$   w.r.t. $(m_\e)$.

\noindent (iv) 
Assume that $\bfgotp_\e$  is    a  linear combination  of 
 $\bfgotw_\e$, $\bfgotv_\e$, $\theta_\e\bfe_3$ and   $\bfgotp_\e(\bfvarphi_\e)\rightharpoonup \bfp$ weakly  in $L^p(\Omega;\RR^3)$. By \eqref{hypfie} and  \eqref{mdx},   $\int \lb \bfgotp_\e(\bfvarphi_\e)\wedge \frac{\bfz_\e}{r_\e}\rb^p dm_\e 
  \leq C$, hence, by Lemma \ref{lemtwoscale},  $\lp\bfgotp_\e(\bfvarphi_\e)\wedge \frac{\bfz_\e}{r_\e}\rp \mdto \bfq_0$, 
 up to a subsequence, for  some $\bfq_0$.
We fix $ \bfeta\in \D(\Omega\times S;\RR^3)$ and set 
 %
  \begin{equation}\begin{aligned}
&
 \bfeta^{[\e]}(x,y):= \sum_{(i,j)\in I_\e\times (L_\e\setminus \{j^m_\e\})}  \lp\intb_{Y_\e^{i}\times\Delta_\e^{j}}  \bfeta(s,y) d\calH^3(s)\rp \mathds{1}_{Y_\e^{i}\times \Delta_\e^j} (x).
 \end{aligned}
\label{defetae}
  \end{equation}
The   fields  $\bfgotp_\e(\bfvarphi_\e)$  and 
  $\bfeta^{[\e]}\lp x,  \frac{y_\e(x')}{r_\e}\rp$  are  independent of $x_3$  in each     $Y_\e^{i}\times \Delta_\e^j$   and,  by \eqref{defye} and   \eqref{defze},  $\int_{ \Delta_\e^j}  \bfz_\e(x',x_3)  dx_3=\int_{ \Delta_\e^j}  \bfy_\e(x')  dx_3$,
therefore
\begin{equation}\begin{aligned}
   \int_{Y_\e^{i}\times \Delta_\e^j}  &\lp\bfgotp_\e (\bfvarphi_\e)\wedge  \frac{\bfz_\e(x)}{r_\e} \rp \cdot \bfeta^{[\e]}\lp x,  \frac{y_\e(x')}{r_\e}\rp  dm_\e(x)  
 \\&\hskip2cm= \int_{Y_\e^{i}\times \Delta_\e^j}   \lp\bfgotp_\e (\bfvarphi_\e)\wedge  \frac{\bfy_\e(x')}{r_\e} \rp \cdot \bfeta^{[\e]}\lp x,  \frac{y_\e(x')}{r_\e}\rp  dm_\e(x) .   
 \end{aligned}
\nonumber
  \end{equation}
Summing   w.r.t. $(i,j)$ over $I_\e\times (L_\e\setminus \{j^m_\e\})$, passing to the limit as $\e\to0$, taking 
the estimate $|\bfeta-\bfeta^{[\e]}|_{L^\infty(\Omega\times S;\RR^3)}\leq C\e$   into account and noticing that,  by \eqref{mdx2},  $\bfgotp_\e (\bfvarphi_\e)\wedge  \frac{\bfy_\e(x')}{r_\e} \mdto  \bfp\wedge\bfy$, we obtain 
$ \frac{1}{|S|}\int_{\Omega\times S} \bfq_0(x,y)\cdot\bfeta(x,y) dxdy
  = \frac{1}{|S|}\int_{\Omega\times S} ( \bfp\wedge\bfy)\cdot\bfeta(x,y) dxdy     
 $
and deduce     $\bfq_0 =\bfp\wedge\bfy$. The assertion \eqref{mdx222} is proved. 
The assertion \eqref{gotrdto} results from \eqref{defgot}, \eqref{mdx2},  and \eqref{mdx222}.     \qed

%
\subsection{Key inequalities}\label{seckey}
 \noindent  In this section,  we establish a  series of inequalities which, combined with Proposition \ref{lemmdx},
 will yield 
a number of  a priori estimates and convergences for a sequence satisfying \eqref{supFeuefini} and its associated auxiliary sequences. 
%
  %
The following   Korn's inequalities, proved in \cite{GeSu,KoOl}, will be employed at several occurrences:

 %
 \begin{lemma}
 \label{lemKorn}  
 (i) We have, for $N\in\{2,3\}$,
  \begin{equation}
\hskip-0,1cm\begin{aligned}
 \int_{\RR^N} |\bfnabla \bfvarphi|^p dy&\leq  
  C   \int_{\RR^N}  |\bfe(\bfvarphi)|^p   dy \qquad   \forall p\in(1,+\infty),\    \forall \bfvarphi\in  \D (\RR^N;\RR^N). 
\end{aligned}
\label{LpfirstKorn0}
    \end{equation}
(ii) If  $U$ is   a bounded Lipschitz domain of $\RR^N$ and $V$ a subspace of $W^{1,p}(U;\RR^N)$ 
such that $V\cap \R=\{0\}$, where $\R$ is the space of
rigid motions  in $U$, then 
   \begin{equation}
 \hskip-0,1cm\begin{aligned}
\int_U |\bfvarphi|^p dx+  \int_U |\bfnabla\bfvarphi|^p dx\le C \int_U  |\bfe(\bfvarphi)|^p dx
  \quad     \forall \bfvarphi\in V. 
\end{aligned}
\label{LpKorn}
    \end{equation}
 \end{lemma}

\noindent
  We set
 \begin{equation}\begin{aligned}
&   q:= \frac{3p}{3-p}\ & &\hbox{  if } \ p<3; \qquad  & & q\in (p,+\infty) \quad  & &\text{ if } p  \ge3,
 \\& c_\e:=     \frac{1}{\gamma^{(p)}_\e\!(r_\e)} \quad
   & &\text{ if } p\not =2;  \qquad  & & c_\e:=   \e^2 \quad & & \text{ if } p  =2. 
 \end{aligned}
\label{defqce}
  \end{equation}

\begin{proposition} \label{propestim}  \noindent Let $A$  be a bounded Lipschitz domain of $\RR^2$ satisfying \eqref{hypA} and  let $ \alpha\in \{1,2\}$. 
The following inequalities   hold   for every $\bfvarphi\in W^{1,p} (\Omega;\RR^3)$:   %
 \begin{eqnarray}
& &   \hskip-0,5cm  \int \lb\bfvarphi-\bfgotr_\e(\bfvarphi) \rb^p d m_\e 
\leq C   r_\e^p\int|\bfe(\bfvarphi)|^pdm_\e,
\label{estimrigid}
\\
& &   \hskip-0,5cm \int \lb \bfvarphi-{ \bfgotv_\e} (\bfvarphi)  \rb^p dm_\e  \leq C c_\e  \int_{{T_{r_\e}}}   |\bfnabla \bfvarphi|^pdx,\label{estimtild11}
\\
& &   \hskip-0,5cm \int \hskip-0,1cm\lb  { \bfgotw_\e} (\bfvarphi)  \rb^p \hskip-0,1cmdm_\e \leq C
r_\e^p \hskip-0,1cm\int \hskip-0,1cm |\bfe ( \bfvarphi)|^p\hskip-0,05cm dm_\e 
+ C c_\e\hskip-0,1cm  \int_{{T_{r_\e}}}  \hskip-0,2cm |\bfnabla \bfvarphi|^pdx, 
\label{estimgotw}
 \\ 
 & &   \hskip-0,5cm \int_\Omega \hskip-0,1cm|\bfvarphi-\ov\bfgotu_\e(\bfvarphi)|^pdx \leq C\lp  \frac{1}{\gamma^{(p)}_\e\!(R_\e)}+\e^{1-\frac{p}{q}} \rp  \hskip-0,1cm \int_\Omega |\bfnabla \bfvarphi|^pdx,\label{estimtild1}
\\
& &   \hskip-0,5cm
  \int_\Omega \hskip-0,1cm| \bfvarphi -\ov \bfgotv^{A}_\e(\bfvarphi)|^p dx \leq C\lp \frac{1}{\gamma^{(p)}_\e\!(r_\e)} +\e^{1-\frac{p}{q}}\rp \hskip-0,1cm\int_\Omega |\bfnabla \bfvarphi|^pdx,
\label{estimtild1b}
\\
& & \hskip-0,5cm
   \int 
\lb\bfvarphi-\ov\bfgotr_\e^{A}(\bfvarphi) \rb^p dm_\e
\leq C   r_\e^p\int  |\bfe_{x'}(\bfvarphi)|^pdm_\e \nonumber
\\
& & \hskip  5cm \ \ \forall  \bfvarphi\in L^p(0,L;W^{1,p} (\Omega';\RR^3)),
\label{estimrigidprim}
\\
 & &   \hskip-0,6cm
  \begin{aligned}\int  \lb 
(\theta_\e - \ov \theta_\e^{A})(\bfvarphi)  \rb^p + \lb (\gotv_{\e\alpha }-\ov \gotv^{A}_{\e\alpha })(\bfvarphi)  \rb^p   dm_\e \leq C&r_\e^p  \int 
  \lb\bfe(\bfvarphi)
 \rb^p dm_\e,
 \end{aligned}\label{estimovA}
\\& &   \hskip-0,5cm
\int   \lb  (\gotv_{\e3} 
-{\ov\gotv}^{A}_{\e3})(\bfvarphi)\rb^p dm_\e
  \leq C  c_\e  \int_{T_{r_\e}}  \lb \bfnabla\bfvarphi  \rb^p dx.
\label{estimovAb}
\end{eqnarray} 
The following inequalities hold   for every $\bfvarphi\in W_b^{1,p} (\Omega;\RR^3)$:
 \begin{eqnarray} 
& &   \hskip-0,5cm 
    \int    \lb      \theta_\e(\bfvarphi) \rb^p +   \lb      \ov\theta^A_\e(\bfvarphi) \rb^p dm_\e  \leq C    \int  |\bfe(\bfvarphi)|^pdm_\e,   \label{estimDirichlet}
\\
& &   \hskip-0,5cm 
\int  |r_\e\varphi_{1}|^p+|r_\e\varphi_{ 2}|^p+\left|{\varphi_{ 3} }\right|^p 
dm_\e
\leq {C }\hskip-0,1cm\int \left|  \bfe(\bfvarphi)\right|^p dm_\e,
\label{estimDirichletb}
\end{eqnarray} 
\begin{equation}
   \begin{aligned}
 \hskip-0,2cm &   \int  \lb    \ov\gots^{A}(\bfvarphi)   \rb^p  \hskip-0,1cm dm_\e    \leq C \hskip-0,1cm  \int  \lb \bfe (\bfvarphi) \rb^p  dm_\e , \hskip-0,1cm
    \end{aligned} 
\label{votheta3}
  \end{equation}
where
     \begin{equation}\begin{aligned}
      \ov\gots^{A}(\bfvarphi):=   -\frac{y_{\e2}}{r_\e} \lp \varphi_1-\ov\gotv^{A}_{\e1}(\bfvarphi)\rp 
 + \frac{y_{\e1}}{r_\e} \lp \varphi_2-\ov\gotv^{A}_{\e2}(\bfvarphi)\rp. 
          \end{aligned} 
  \label{defovgots} 
    \end{equation}   
 \end{proposition}
 
   \noindent
   %
 {\bf Proof.}  \  {\it Proof of \eqref{estimrigid}, \eqref{estimtild11}, \eqref{estimgotw}.}  
 By \eqref{defze}, we have (see fig. \ref{fig2})
\begin{equation}
   \begin{aligned}
  &    T_{r_\e}^{ij}  = \bfz_\e^{ij} +r_\e T,   \qquad T:= S\times (-1,1].
    \end{aligned} 
\label{defT}
  \end{equation}
 We consider  the linear 
  operators  $\bfgotr, \bfgotv, \bfgotw$  defined on $W^{1,p}(T;\RR^3)$  by 
  %
  \begin{equation}\begin{aligned}
&\bfgotr(\bfvarphi)(z)\!:=  \bfgotv(\bfvarphi)+\! \frac{2}{\diam S} \bfgotw(\bfvarphi)\wedge   \bfz ,  \quad \!
\\&   \bfgotv(\bfvarphi):= \! \intb_T  \bfvarphi (s) d\calL^3(s),
   \qquad \! \bfgotw(\bfvarphi):= \!  \frac{\diam S}{2}\bfM \intb_T   (\bfs \wedge \bfvarphi (s ) ) d\calL^3(s) .
\end{aligned} 
\nonumber
  \end{equation} 
We have  $\bfgotr\circ\bfgotr=\bfgotr$ and  $\bfgotr(\bfvarphi)=\bfvarphi \ \forall \bfvarphi\in \R$,  hence     $V:= \{ \bfvarphi-\bfgotr(\bfvarphi),\ 
  \bfvarphi\in W^{1,p}\lp T;\RR^3\rp \}$  satisfies 
  $V\cap \R=\{0\}$.
 Applying Lemma \ref{lemKorn}, noticing that  $\bfe( \bfvarphi)= \bfe( \bfvarphi -\bfgotr( \bfvarphi))$,   
   we infer
   \begin{equation}\begin{aligned}
&|  \bfvarphi -\bfgotr( \bfvarphi)|_{L^p\lp T;\RR^3\rp}
  \leq C |\bfe( \bfvarphi)|_{L^p\lp T;\RR^3\rp}\quad \forall   \bfvarphi\in W^{1,p}\lp T;\RR^3\rp.
\end{aligned} 
 \label{deuC}
  \end{equation}
By making  appropriate changes of variables,  taking \eqref{defgot} and  \eqref{defT} into account,   we deduce that  for every $(i,j) \in {I}_\e\times L_\e$,
  %
   \begin{equation}
\begin{aligned}
&  \int_{{T^{ij}_{r_\e}}} 
\lb\bfvarphi-\bfgotr_\e(\bfvarphi) \rb^p    dx 
\leq C   r_\e^p  \int_{{T^{ij}_{r_\e}}}   |\bfe(\bfvarphi)|^pdx  \quad \forall   \bfvarphi\in  W^{1,p} \lp T_{r_\e}^{ij};\RR^3\rp.   
\end{aligned}
\label{abovein}
  \end{equation}
 If $ \bfvarphi \in W^{1,p}(\Omega;\RR^3)$, its  extension  
   defined on  $\Omega'\times (L,2L)$ by 
    \begin{equation}
   \begin{aligned}
   &  \widehat\varphi_\alpha (x',x_3):= -3 \varphi_\alpha (x',2L-x_3)  +4\varphi_\alpha \lp x', \frac{3L-x_3}{2}\rp 
\\&  \widehat\varphi_3(x',x_3):= 3 \varphi_3(x',2L-x_3)  -2\varphi_3\lp x', \frac{3L-x_3}{2}\rp ,
\end{aligned} 
\label{defhatfi}
  \end{equation}
satisfies   the  following  estimate    on $\Omega'\times(L,2L)$:
   \begin{equation}\begin{aligned} 
  \lb\bfe(\widehat\bfvarphi)(x)\rb  \leq C \lb\bfe(\bfvarphi)\rb(x',2L-x_3)   + C \lb\bfe(\bfvarphi) \rb\lp x', \frac{3L-x_3}{2}\rp.
  \end{aligned}
\label{choicewidehat}
  \end{equation} 
Observing that  \eqref{defjme},   \eqref{abovein} and \eqref{choicewidehat}  imply  
 \begin{equation}
\int_{{T^{ij_\e^m}_{r_\e}}} 
\lb\widehat\bfvarphi-\bfgotr_\e(\widehat\bfvarphi) \rb^p    dx \leq C   r_\e^p   \int_{{T^{i}_{r_\e}}}  |\bfe(\bfvarphi)|^pdx\quad  \bfvarphi\in W^{1,p}(\Omega;\RR^3),
  \ \forall i \in I_\e,
\label{choicewidehat2}
  \end{equation}
we deduce (see  \eqref{defze})
   \begin{equation}
\begin{aligned}
   \int_{T_{r_\e}^{i}} 
\lb\bfvarphi-\bfgotr_\e(\bfvarphi) \rb^p dx& \leq \sum_{j=1}^{j_\e^m-1}  \int_{T_{r_\e}^{ij}} 
\lb\bfvarphi-\bfgotr_\e(\bfvarphi) \rb^p dx +   \int_{T_{r_\e}^{ij_\e^m}} 
\lb\widehat\bfvarphi-\bfgotr_\e(\widehat\bfvarphi) \rb^p dx
\\& \leq C   r_\e^p\int_{T^{i}_{r_\e}}  |\bfe(\bfvarphi)|^pdx \quad  \qquad \forall \bfvarphi\in W^{1,p}(\Omega;\RR^3), \quad \forall i\in I_\e, 
\end{aligned}
\nonumber
  \end{equation}
 \noindent  
yielding  \eqref{estimrigid}.
Applying  a similar  argument to 
  $W:= \{ \bfvarphi\in W^{1,p}\lp T;\RR^3\rp, \ \bfgotv(\bfvarphi)=0\}$, 
  we obtain 
$\int 
\lb\bfvarphi-\bfgotv_\e(\bfvarphi) \rb^p dm_\e \hskip-0,1cm\leq  \hskip-0,1cmC   r_\e^p\int  |\bfnabla(\bfvarphi)|^pdm_\e$ and, taking    \eqref{defgammae},  \eqref{defme}  and  \eqref{defqce}  into account,   deduce
  \eqref{estimtild11}. 
To prove \eqref{estimgotw}, we start from   the  inequality 
 $|\bfb|^p\leq C \int_T |\bfb\wedge   \bfz |^pdz  \  \forall \bfb\in \RR^3$ (easily proved by contradiction). 
 By  suitable   changes  of variables, we obtain
 \begin{equation}
 |\bfb|^p\leq \frac{C}{r_\e^3} \int_{T_{r_\e}^{ij}} \lb \bfb \wedge  \frac{\bfz_\e(  x)}{r_\e}  \rb^pd x \qquad  \forall \,\bfb\in \RR^3,
  \ \forall (i,j)\in I_\e\times L_\e.
\nonumber
  \end{equation}
Substituting  for $\bfb$  the constant value taken by 
   $\bfgotw_\e(\bfvarphi)$  in each set $T_{r_\e}^{ij}$, 
     we infer  
   \begin{equation}\begin{aligned}
 \int   |\bfgotw_\e (\bfvarphi)  
|^p  dm_\e   
&\leq    \sum_{(i,j)\in I_\e\times L_\e} \int_{T_{r_\e}^{ij}} \lp    \frac{C}{r_\e^3} \int_{T_{r_\e}^{ij}}  \lb \bfgotw_\e (\bfvarphi)   \wedge  \frac{\bfz_\e(  x)}{r_\e}  \rb^p  dx \hskip-0,1cm\rp    dm_\e 
\\&  \leq C    \int   \lb \bfgotw_\e (\bfvarphi)   \wedge  \frac{\bfz_\e(  x)}{r_\e}  \rb^p    dm_\e = C \int  \lb  \bfgotr_\e(\bfvarphi)   - \bfgotv_\e(\bfvarphi)
    \rb^p d m_\e
\\& \leq 
 C\int   \lb \bfvarphi-\bfgotr_\e(\bfvarphi)     \rb^p +  \lb  \bfvarphi-\bfgotv_\e(\bfvarphi)
    \rb^p d m_\e,
\end{aligned}
\nonumber
  \end{equation}
\noindent   which, combined   with   \eqref{estimrigid} and \eqref{estimtild11}, yields \eqref{estimgotw}. \qed

\noindent {\it Proof of \eqref{estimtild1}, \eqref{estimtild1b}}. We put $Y_\e := \bigcup_{i\in I_\e}   Y_\e^{i}$.
By \eqref{defovgotA} and \eqref{defovgotue},  we have 
\begin{equation}\begin{aligned}
&  
\int_\Omega \hskip-0,2cm |\bfvarphi-\ov\bfgotp_\e(\bfvarphi)|^pdx  = \int_{Y_\e\times(0,L)} \hskip-1cm |\bfvarphi-\ov\bfgotp_\e(\bfvarphi)|^pdx + \int_{\Omega'\setminus Y_\e\times(0,L)}\hskip -1cm  |\bfvarphi |^pdx
\quad \forall \ \ov\bfgotp_\e\in \{\ov\bfgotv_\e^{A}, \ov\bfgotu_\e\}.\hskip-0,2cm 
 \end{aligned}\hskip-0,2cm 
\label{estimtild01}
  \end{equation}
  \noindent   By  \eqref{defTre},  $ \Omega' \setminus Y_\e  \subset  \{ x'   \in  \Omega', {\rm dist } (x', \partial \Omega')  \leq  \sqrt 2 \e \}$, thus, since $\Omega'$ is Lipschitz,    $\calL^2(\Omega'   \setminus  Y_\e) \leq  C\e$.  
By     H\" older's   and Poincar\'e's inequalities   and the continuous embedding of $W^{ 1,p}( \Omega; \RR^3)$   into $ L^{ q}( \Omega; \RR^3)$ (see  \eqref{defqce}, 
 \cite[Corollary 9.14]{Br}),  the following inequalities  hold  in $W^{1,p}_b(\Omega;\RR^3)$:
\begin{equation}\begin{aligned}  
\int_{\Omega' \setminus Y_\e\times(0,L)}  |\bfvarphi|^p\,dx 
 \leq |\bfvarphi|_{L^{q}(\Omega;\RR^3)}^p |C\e|^{1-\frac{p}{q}}
\leq C& \e^{1-\frac{p}{q}}   \int_{\Omega}   |\bfnabla \bfvarphi|^p\,dx .
\end{aligned}
\label{estimtild02}
  \end{equation}
\noindent  Let $E$ be a bounded Lipschitz  domain of $\RR^2$   such that $\partial D\subset \ov E$. 
We prove below the existence of   $C>0$ such that,  for  every $\alpha\in (0,1]$ verifying  $\alpha E\subset Y$,
 %
  \begin{equation}\begin{aligned}
& 
\int_{Y} \left| \varphi- \intb_{\alpha E } \varphi\,ds' \right|^pdx'
\leq    C  l(\alpha)  \int_Y |\nabla \varphi |^p dx'\qquad  \forall \varphi \in W^{1,p}(Y), 
\\& l(\alpha) =
\alpha^{p-2}   \text{ if }  1<p<2; \quad  l(\alpha) =
 1+|\log \alpha|   \text{ if }  p=2; \quad  l(\alpha) =
1    \text{ if }  p>2. \end{aligned}
\label{A4}
\end{equation} 
 Let us see how the claim follows from \eqref{A4}: 
let $(d_\e)\subset \RR$ be such that $0<d_\e\ll\e$ and $E_{d_\e}^{i}$ defined by \eqref{defAbi}.
 An elementary change of variables yields
  %
 \begin{equation}\begin{aligned}
& 
\int_{Y_\e^{i}} \left| \varphi- \intb_{E_{d_\e}^{i} } \varphi\,ds' \right|^pdx'
\leq    C \e^p l\lp\frac{d_\e}{\e}\rp \int_{Y_\e^{i}}  |\nabla \varphi |^p dx'\qquad  \forall \varphi \in W^{1,p}(Y_\e^{i}) .
 \end{aligned}
\nonumber
\end{equation} 
Summing w.r.t. $i$   over $I_\e$ and integrating w.r.t. $x_3$   over $(0,L)$, 
choosing successively $(E,d_\e):= (D\setminus D/2, R_\e)$ and $(E,d_\e):= (A, r_\e)$,
noticing that, by \eqref{defgammae} and \eqref{A4},  $ \e^p l\lp\frac{d_\e}{\e}\rp\leq \frac{1}{\gamma_\e^{(p)}(d_\e)}$,  we infer
\begin{equation}\begin{aligned}
&  
\int_{Y_\e\times(0,L)} |\bfvarphi-\ov\bfgotu_\e(\bfvarphi)|^pdx \leq    \frac{C}{\gamma^{(p)}_\e\!(R_\e)}  \int_\Omega |\bfnabla \bfvarphi|^pdx,
\\ &
  \int_{Y_\e\times(0,L)}  | \bfvarphi -\ov \bfgotv^{A}_\e(\bfvarphi)|^p dx \leq   \frac{C}{\gamma^{(p)}_\e\!(r_\e)}   \int_\Omega |\bfnabla \bfvarphi|^pdx.
 \end{aligned}
\nonumber
  \end{equation}
\noindent Taking  \eqref{estimtild01} and \eqref{estimtild02} into account,    \eqref{estimtild1} and \eqref{estimtild1b} are proved.  We turn to the proof of \eqref{A4}. 
  A straightforward variant of   \cite[Lemma  A4]{BeBo} yields
     \begin{equation}\begin{aligned}
&
\int_{Y} \left| \varphi- \intb_{\partial  \alpha D } \varphi\,d\calH^1 \right|^pdx
\leq  C  l(\alpha)  \int_{Y} |\nabla \varphi |^p dx\quad  \forall \varphi \in W^{1,p}(Y).
\end{aligned}
\label{A41}
\end{equation} 
On the other hand,  for all $\varphi\in W^{1,p}(Y)$,  
 \begin{equation}\begin{aligned} 
&  \left| \intb_{ \alpha  E} \varphi\,d\calL^2- \intb_{\partial \alpha  D} \varphi\,d\calH^1 \right|^p 
\leq   \intb_{\alpha   E} \lb \varphi - \intb_{\partial \alpha  D } \varphi\,d\calH^1 \rb^p dx
\\& = \frac{1}{\alpha^2|E|}   \int_{\alpha   E} \lb \varphi - \intb_{\partial \alpha  D } \varphi\,d\calH^1 \rb^p dx \leq C\frac{\alpha^p  }{\alpha^2} \int_{\alpha  E } |\nabla \varphi|^p dx
 \leq   C   l(\alpha)  \int_Y |\nabla \varphi |^p dx,
\end{aligned}
\nonumber
\end{equation} 
which, along  with  \eqref{A41}, yields   \eqref{A4}. \qed

\noindent {\it Proof of \eqref{estimrigidprim}}. We introduce  the linear  operators       defined  on $L^p(S;\RR^2)$ by 
  \begin{equation}\begin{aligned}
&\hskip-0,1cm\bfgotr^{A}(\bfpsi)(y)\hskip-0,1cm := \hskip-0,05cm \bfgotv^{A}(\bfpsi) \hskip-0,05cm+\hskip-0,05cm \frac{2}{ \diam A} \theta^{A}(\bfpsi) \bfe_3\wedge   \bfy;   
\ \quad  \bfgotv^{A}(\bfpsi)\hskip-0,05cm:= \hskip-0,1cm  \intb_{\hskip-0,1cm A} \hskip-0,05cm \bfpsi (y) d\calL^2(y),
\\&\hskip-0,1cm \theta^{A}(\bfpsi):=   \frac{\diam A}{2\intb_A|\bfy|^2 dy }   {\int\!\!\!\!\!\!-}_A
 \lp  -  y_{2}  \psi_1(y)+ y_{1}   \psi_2(y) \rp d\calL^2(y).
\end{aligned} 
\label{OA}
  \end{equation} 
\noindent The  two-dimensional version of the argument used to establish \eqref{deuC} yields
\begin{equation}
\int_S |\bfpsi- \bfgotr^{A}(\bfpsi)|^p d\calL^2\leq   C \int_S |\bfe(\bfpsi)|^p d\calL^2  \qquad \forall \bfpsi\in W^{1,p}(S;\RR^2).  
\label{deuC2}
\end{equation}
\noindent By the Poincar\'e-Wirtinger inequality in $W^{1,p}(S)$, we have 
  \begin{equation}
\int_S \lb \eta -\intb_A \eta d\calL^2 \rb^p d\calL^2 \leq   C  \int_S |\bfnabla \eta|^p d\calL^2  \qquad \forall \eta\in W^{1,p}(S).  
\label{PW}
\end{equation}
 Fixing $\bfvarphi \in W^{1,p}(T;\RR^3)$, we set 
 (see \eqref{defexprim})
  \begin{equation}\hskip-0,2cm\begin{aligned}
  &\ov\bfgotv^{A}(\bfvarphi) (z_3)= \bfgotv^{A}(\bfvarphi'(.,x_3))
    +  \intb_A \varphi_3(y,z_3) d\calH^2(y) \bfe_3
\\&\ov\bfgotr^A(\bfvarphi)(y,z_3):= \bfgotr^{A}(\bfvarphi'(.,x_3))
 + \ov\gotv^{A}_3(\bfvarphi)\bfe_3;  \quad\ov\theta^{A}(\bfvarphi)(z_3)= \theta^{A} (\bfvarphi'(.,x_3))
 .\end{aligned} 
\nonumber
  \end{equation} 
Applying \eqref{deuC2}  
   to $\bfpsi \hskip-0,1cm:=  \bfvarphi'(.,x_3)$,  \eqref{PW} to $\eta:=\varphi_3(.,z_3)$,  integrating w.r.t. $z_3$ over $(-1,1)$,  noticing that  $|\bfe(\bfpsi)| + |\bfnabla\eta|\leq C |\bfe_{x'}(\bfvarphi)|$, 
 we infer 
   \begin{equation}
   \begin{aligned}
 \int_T|\bfvarphi-\ov\bfgotr^{A}(\bfvarphi)|^p d\calL^3& \leq C  \int_{-1}^1 dz_3  \int_S |\bfpsi- \bfgotr^{A}(\bfpsi) |^p   +    \lb \eta -\intb_A \eta d\calL^2 \rb^p d\calL^2 
\\&  \leq C  \int_{-1}^1 dz_3  \int_S |\bfe(\bfpsi) |^p   +    \lb \bfnabla\eta \rb^p d\calL^2 
 \leq   C \int_T |\bfe_{x'}(\bfvarphi)|^p d\calL^3,
 \end{aligned}
\nonumber
\end{equation}
 \noindent    
and deduce, by   suitable changes of variables,   
  \begin{equation}
\begin{aligned}
&  \int_{{T^{ij}_{r_\e}}} 
\lb\bfvarphi-\ov\bfgotr_\e^{A}(\bfvarphi) \rb^p dx 
\leq C   r_\e^p\int_{{T^{ij}_{r_\e}}} |\bfe_{x'}(\bfvarphi)|^pdx   \quad  \forall   \bfvarphi\in  W^{1,p} \lp T_{r_\e}^{ij};\RR^3\rp.
\end{aligned}
\nonumber
  \end{equation}
 %
Noting that,   
by  \eqref{defhatfi}, 
   $\int_{{T^{ij^m_\e}_{r_\e}}} |\bfe_{x'}(\widehat\bfvarphi)|^pdx  
 \leq C\int_{{T^{i}_{r_\e}}} |\bfe_{x'}(\bfvarphi)|^pdx $, 
  we infer   \eqref{estimrigidprim}.
\qed

\noindent {\it Proof of \eqref{estimovA}}.  Fixing    $\bfvarphi\in W^{1,p}(\Omega;\RR^3)$,   applying  \eqref{estimrigidprim} to 
$\bfvarphi'$, 
observing  that $|\bfe_{x'}(\bfvarphi')|\leq |\bfe(\bfvarphi)|$ and    $\ov\bfgotr_\e^{A}(\bfvarphi')=  \ov\bfgotr_\e^{A}(\bfvarphi)'$ (see \eqref{defovgotA}),  we obtain
\begin{equation}\begin{aligned}
&    
 \int  \lb\bfvarphi'-  \ov\bfgotr_\e^{A}(\bfvarphi)' \rb^p   dm_\e   \leq Cr_\e^p  \int  \lb\bfe(\bfvarphi) \rb^p dm_\e.
\end{aligned}
\label{2d}
  \end{equation}
\noindent  Combining this with the following inequality, deduced from   \eqref{defgot} and   (\ref{estimrigid}), 
$$
  \int \lb\bfvarphi'-\bfgotr_\e(\bfvarphi)' \rb^p dm_\e   \leq C \int |\bfvarphi-\bfgotr_\e(\bfvarphi)|^p dm_\e 
 \leq C   r_\e^p\int  |\bfe(\bfvarphi)|^pdm_\e,
 $$
we  infer 
  \begin{equation}\begin{aligned}
&       \int 
\lb(\bfgotr_\e-\ov\bfgotr_\e^{A})(\bfvarphi)' 
 \rb^p dm_\e \leq Cr_\e^p  \int 
 \lb\bfe(\bfvarphi)
 \rb^p dm_\e.
\end{aligned}
\nonumber
  \end{equation}
By   \eqref{mdx}, we have
    \begin{equation} 
     \begin{aligned}
 \int  \lb 
 \ov \theta^{A}_\e \lp (\bfgotr_\e-\ov\bfgotr_\e^{A})(\bfvarphi)' \rp \rb^p + &\lb  \ov\bfgotv^{A}_\e \lp(\bfgotr_\e-\ov\bfgotr_\e^{A})(\bfvarphi)'  \rp \rb^p   dm_\e
  \leq C  \int  \lb 
 (\bfgotr_\e-\ov\bfgotr_\e^{A})(\bfvarphi)'  \rb^p   dm_\e.
 \end{aligned} 
\nonumber
  \end{equation}
By \eqref{defovgotA} and  \eqref{defgot},   $ \ov \theta^{A}_\e \lp  (\bfgotr_\e-\ov\bfgotr_\e^{A})(\bfvarphi)' \rp = 
({\theta_\e}- \ov\theta^{A}_\e)(\bfvarphi)$ and 
$\ov\bfgotv^{A}_\e\lp(\bfgotr_\e-\ov\bfgotr_\e^{A})(\bfvarphi)' \rp= (\bfgotv_\e-\ov\bfgotv^{A}_\e)(\bfvarphi))'$. This, along   with   the above inequalities,  proves \eqref{estimovA}. \qed
%

\noindent{\it Proof of \eqref{estimovAb}.}
Given $f\in W^{1,p}(T)$, we set   $g(t):=-\hskip-0,3cm \int_{A}f( x',t  ) d \calH^2(x')$ and  $T^{A}:= A\times (-1,1)$. For every  $x_3\in (-1,1)^2$, we have
 \begin{equation}\begin{aligned}
&\lb -\hskip-0,3cm \int_{T^{A}}  f  d\calL^3 - \intb_{A}\hskip-0,1cm f(x', x_3) d \calH^2(x') \rb^p = 
 \lb -\hskip-0,3cm \int_{-1}^1 g(t)dt  - g(x_3) \rb^p 
 \\&\qquad  \le  -\hskip-0,3cm \int_{-1}^1 |g(t )-g(x_3)|^p dt
  \le C \int_{-1}^1  \lb \int_{-1}^1  g' (s)ds\rb^p dt 
  \le C \intb_{T^{A}}\lb \frac{\partial f}{\partial x_3}\rb^p d\calL^3.
\end{aligned}
\nonumber
  \end{equation}
This, combined with the  following Poincar\'e-Wirtinger inequality 
 \begin{equation}\begin{aligned}
\lb-\hskip-0,3cm \int_T f  d\calL^3 -  -\hskip-0,3cm \int_{T^{A}}  f  d\calL^3\rb^p & \le \intb_T\lb f-\intb_{T^{A}}fd\calL^3 \rb^p d\calL^3
\le C \int_T |\bfnabla f |^p d\calL^3,
\end{aligned}
\nonumber
  \end{equation}
implies
 \begin{equation}\begin{aligned}
&   
\hskip-0,2cm \int_{T} \lb   \intb_T\hskip-0,2cm f d\calL^3 \hskip-0,08cm-\intb_{A}\hskip-0,1cm f(x', x_3) d \calH^2(x')  
\rb^p d\calL^3 \leq  C\hskip-0,1cm
 \int_{T}\hskip-0,1cm  |\bfnabla f|^pd\calL^3  \quad \forall f\in W^{1,p}(T).\hskip-0,3cm
 \end{aligned}
 \nonumber
  \end{equation}
   By    suitable changes of variables, 
  we infer  (see  \eqref{defovgotA},  \eqref{defgot}) 
 %
   \begin{equation}\begin{aligned}   
&\int_{T_{r_\e}^{ij}}   \lb  \gotv_{\e3}(\bfvarphi) 
-{\ov\gotv}^{A}_{\e3}(\bfvarphi)\rb^p dx 
  = \int_{T_{r_\e}^{ij}}  \lb  {\int\!\!\!\!\!\!-}_{T_{r_\e}^{ij}} 
\varphi_3   d\calL^3
  - {\int\!\!\!\!\!\!-}_{{ A_{r_\e}^{i}}}
\varphi_3 (s, x_3  )   d\calH^{2}(s) \rb^p dx 
\\&\hskip 1 cm  \leq C  r_\e^p \int_{T_{r_\e}^{ij}}  \lb \bfnabla\bfvarphi  \rb^p dx
\qquad \forall   \bfvarphi\in  W^{1,p} \lp T_{r_\e}^{ij};\RR^3\rp, \ \forall (i,j) \in {I}_\e\times L_\e.
 \end{aligned}
\nonumber
  \end{equation}
 Summing   w.r.t. $(i,j)$,
  in view of  \eqref{defme}, \eqref{defqce} and  \eqref{defhatfi},  we 
   obtain    (\ref{estimovAb}).  \qed
   
\noindent {\it Proof of \eqref{estimDirichlet}}. Let us fix   $\bfvarphi \in W^{1,p}(\Omega;\RR^3)$. By (\ref{defovgotA}), for   a. e. $x_3\in (0,L)$,  
   \begin{equation}
   \hskip-0,1cm \begin{aligned} 
 \frac{2{\int\!\!\!\!\!-}_D |\bfy|^2 dy}{ \diam S}   \frac{\partial}{\partial x_3}
 \ov\theta_\e^D (\bfvarphi)(x  )& 
     =\!
   {\int\!\!\!\!\!\!-}_{{  D_{r_\e}^{i}}}  \hskip-0,3cm     - \frac{ y_{\e2}}{r_\e} \frac{\partial  \varphi_{1}}{\partial x_3}  (s,x_3  ) +
    \frac{ y_{\e1}}{ r_\e}  \frac{\partial\varphi_{ 2}}{\partial x_3}  (s,x_3  )
 d\calH^2(s   ).\end{aligned}
\label{theta,3}
  \end{equation}
By \eqref{defye},  $\frac{\partial y_{\e1}}{\partial x_2}=\frac{\partial y_{\e2}}{\partial x_1}=0$ in   $T_{r_\e}^{i}$ 
and 
  $\frac{ \bfy_\e }{r_\e} $ coincides  on   $\partial D_{r_\e}^{i}\times \{x_3\}$ with  the outward normal to $\partial D_{r_\e}^{i}$, hence
  $ {\int\!\!\!\!\!\!-}_{{  D_{r_\e}^{i}}}   
  \hskip-0,2cm   - \frac{ y_{\e2} }{r_\e}  \frac{\partial  \varphi_{3}}{ \partial x_1}  (.,x_3)+
   \frac { y_{\e1} }{ r_\e}  \frac{\partial  \varphi_{3}}{ \partial x_2}  (.,x_3)  d\calH^2 =0$. 
Summing this with  (\ref{theta,3}), we obtain
   \begin{equation}
   \hskip-0,1cm \begin{aligned} 
 \frac{2{\int\!\!\!\!\!-}_D |\bfy|^2 dy}{ \diam S}   \frac{\partial}{\partial x_3}
 \ov\theta_\e^D (\bfvarphi)(x  )& 
     =\!
    2{\int\!\!\!\!\!\!-}_{{  D_{r_\e}^{i}}}    \hskip-0,3cm    - \frac{ y_{\e2}}{r_\e}   e_{13}(\bfvarphi)   (s,x_3  )  +
    \frac{ y_{\e1}}{ r_\e}    e_{23}(\bfvarphi)   (s,x_3  )
 d\calH^2(s  ),
\end{aligned}
\nonumber
  \end{equation}
and  deduce  
$\int 
\lb\frac{\partial}{ \partial x_3}  \ov\theta_\e^D(\bfvarphi) \rb^p dm_\e
\leq C    \int  |\bfe(\bfvarphi)|^pdm_\e
$. 
If   $\bfvarphi\in W^{1,p}_b(\Omega;\RR^3)$,  
  by \eqref{defovgotA} we have  $ \ov\theta_\e^D(\bfvarphi) \in  L^p(\Omega'; W_b^{1,p}(0,L))$  where $W_b^{1,p}(0,L):=\{\eta\in W^{1,p}(0,L), \eta(0)=0\}$, therefore,
by \eqref{defme} and Jensen's inequality,  %
  \begin{equation}\begin{aligned} 
   \int \lb      \ov\theta_\e^D (\bfvarphi)\rb^p dm_\e&= \frac{\e^2}{r_\e^2|S|}  \int_{S_{r_\e}} dx' \int_0^L dx_3   \lb  \int_0^{x_3}    
    \frac{\partial}{ \partial x_3}      \ov\theta_\e^D (\bfvarphi) (x',s_3) ds_3\rb^p 
  \\& \hskip-1cm \leq C  
  \int \lb \frac{\partial}{ \partial x_3}      \ov\theta_\e^D (\bfvarphi)\rb^p dm_\e\leq C    \int  |\bfe(\bfvarphi)|^pdm_\e\quad \forall \bfvarphi\in W^{1,p}_b(\Omega;\RR^3).
  \end{aligned}
 \nonumber
  \end{equation}
    Taking  \eqref{estimovA} into account,   the assertion \eqref{estimDirichlet} is proved.  \qed
   
   \noindent {\it Proof of  \eqref{estimDirichletb}}. 
 The Banach space  $V:=\{ \bfpsi \in W^{1,p}(S\times(0,L);\RR^3), \bfpsi=0  \text{ on } S\times \{0\}\}$
 satisfies  $V\cap \R=\{0\}$.
   By applying (\ref{LpKorn}) 
  to
  $\bfpsi\in V$ defined
   by        $
\psi_\alpha(z_1,z_2,z_3):=\varphi_\alpha(r_\e (z_1-i_1),r_\e (z_2-i_2), z_3)$, $ \alpha  \in \{1,2\}$ and 
$\psi_3(z):=  \frac{1}{r_\e}\varphi _3(r_\e(z_1-i_1),r_\e (z_2-i_2), z_3)$,  and by making a suitable change of variables, we infer
$\int_{T_{r_\e}^{i}}  |r_\e\varphi_{1}|^p+|r_\e\varphi_{ 2}|^p+\left|{\varphi_{ 3} }\right|^p 
dx
\leq {C }\int_{T_{r_\e}^{i}} \left|  \bfe(\bfvarphi)\right|^p dx$, yielding  (\ref{estimDirichletb}). 
  \qed  
  
  \noindent{\it Proof of \eqref{votheta3}}. By
by \eqref{defexprim}, \eqref{defovgotA},  \eqref{defovgots}  and  (\ref{2d}),  we have
  \begin{equation}\begin{aligned}
  \int \lb  \ov\gots^{A}(\bfvarphi)
 -2\frac{ \lb\frac{y_\e}{r_\e}\rb^2}{\diam S} \ov\theta^{A}_\e\rb^pdm_\e 
&= \int
 \lb \frac{\bfy_\e}{r_\e} \wedge \lp   \bfvarphi' - (\ov\bfgotr^{A}_\e(\bfvarphi))'\rp
 \rb^p  dm_\e 
\\ & \leq C r_\e^p \int   \lb \bfe (\bfvarphi)(x  )\rb^p  dm_\e .
\end{aligned}
\nonumber
  \end{equation}
This, combined with  \eqref{estimDirichlet}, implies \eqref{votheta3}.  \qed 

  \subsection{Convergences}\label{secapriori} 
 \noindent  In the next three propositions, we 
establish  a series of  convergences
for a  sequence $(\bfu_\e)$ satisfying \eqref{supFeuefini}, the sequence of  its symmetrized gradients and the  associated auxiliary sequences  defined by \eqref{ve=gotve}, \eqref{defovgotue}, \eqref{Re}, \eqref{defovgots}.
  \begin{proposition} \label{prop1}      Let  $\bfu_\e$ be  a sequence  in $W^{1,p}_b (\Omega;\RR^3)$
 satisfying \eqref{supFeuefini}.
The following  convergences hold, up to a subsequence:
\begin{equation}\begin{aligned}
&\!\begin{aligned}& \bfu_\e   \rightharpoonup \bfu   \  \hbox{weakly  in }    W^{1,p}_b (\Omega;\RR^3),
  & & \ov\bfgotu_\e( \bfu_\e)  \to \bfu   \   \hbox{strongly in }       L^p(\Omega;\RR^3),
 \\&  \ov\bfgotv^S_\e (\bfu_\e)  \rightharpoonup \bfv    \  \hbox{weakly in    }    L^p(\Omega;\RR^3 ), 
 & &\ov\theta^S_\e(\bfu_\e)  \rightharpoonup \theta \   \hbox{weakly  in }  L^p(\Omega),
 \end{aligned}
 \\&  \bfu_\e  \mdto  \bfv+\frac{2}{\diam S}\theta\bfe_3 \wedge \bfy  \hskip2,9cm    \hbox{in accordance with } \eqref{defmdto},    
 \\& \bfu_\e m_\e \buildrel \star\over \rightharpoonup \bfv   \hskip5cm     \hbox{weakly$^\star$ in    }     \M(\ov\Omega;\RR^3 ),
 \\&  \tau^{-1}\lp -\frac{y_{\e2}}{r_\e} u_{\e1} +\frac{y_{\e1}}{r_\e} u_{\e2}\rp m_\e
\buildrel \star \over \rightharpoonup 
  \theta \hskip1cm    \quad \text{weakly$^\star$ in }    \M(\ov\Omega),
      \end{aligned}
\label{cvbasic}
  \end{equation}
    where $\tau$ is defined by \eqref{defvethetae}.
  If  in addition $\kappa>0$, then 
\begin{equation}\begin{aligned}
&  \frac{u_{\e3}}{r_\e} m_\e \buildrel\star\over\rightharpoonup w , \quad  \ 
 \tau^{-1}  \frac{\ov\gots^{S}(\bfu_\e)}{r_\e}  m_\e
 \buildrel \star \over \rightharpoonup  
\delta \quad   \hbox{weakly$^\star$  in }   \M (\ov\Omega), \quad 
\\&  \lp  \frac{\ov\gotv^S_{\e3}(\bfu_\e)}{r_\e} ,  \frac{\ov\theta^S_\e(\bfu_\e)}{r_\e}\rp   \rightharpoonup (w, \delta) \    \hbox{weakly  in }  L^p(\Omega)^2.     \end{aligned}
\label{cvbasickappa>0}
  \end{equation} 
  Besides,  
   \begin{equation}
\begin{aligned}
&   \theta  =0   \   \hbox{ if }  \  p=  2; \qquad 
  \theta  =v_3= 0   \ \hbox{ if }  \ k=+\infty;
\\&  \theta  =0     \hbox{ and  }    \bfv=\bfu       \hbox{ if }    \gamma^{(p)} =+\infty;   \quad \bfv=0 \hbox{ and }  w=\delta=0 \ \hbox{ if }  \  \kappa=+\infty.
\end{aligned}
\label{regthetav}
  \end{equation}
    \end{proposition}
 
   \noindent
 \noindent {\bf Proof.}  {\it Proof of \eqref{cvbasic}.}
By     \eqref{Pe}, \eqref{kkappa} and  \eqref{growthp},
\begin{equation}\hskip-0,2cm
\begin{aligned}
 \eqref{supFeuefini} 
\quad  \Longleftrightarrow \quad 
  \sup_{\e>0}\  |\bfu_\e|_{W_b^{1,p}(\Omega;\RR^3)}+ k_\e  \int|\bfe(\bfu_\e)|^pdm_\e<+\infty.
    \end{aligned}\hskip-0,2cm 
\label{uebounded}
  \end{equation}
Noticing that   $W^{1,p}_b(\Omega;\RR^3)\cap\R=\{0\}$, we infer from  \eqref{LpKorn} that
\begin{equation}\begin{aligned}
 \int_\Omega  |\bfu_\e|^p  dx  
+   \int_\Omega |\bfnabla \bfu_\e|^p  dx
  \leq C \int_\Omega |\bfe( \bfu_\e)|^p  dx 
 + C k_\e\int  |\bfe( \bfu_\e)|^p  dm_\e 
 \leq C.
\end{aligned}
 \label{estimbasic1}
  \end{equation}
Hence  $(\bfu_\e)$ is bounded in $W^{1,p}_b(\Omega;\RR^3)$, thus weakly  converges,
 and strongly     in $L^p(\Omega;\RR^3)$, up to a subsequence, to some $\bfu$. Since $\gamma^{(p)}>0$, by    \eqref{defgammae} and \eqref{Re}, $\lime \gamma_\e^{(p)}(R_\e)=+\infty$, therefore, by 
   \eqref{estimtild1}, 
   $(\ov\bfgotu_\e(\bfu_\e))$  strongly converges  to $\bfu$  in $L^p(\Omega;\RR^3)$.
Observing that, by \eqref{estimtild1b}, 
$(\ov\bfgotv_\e^S(\bfu_\e))$ is bounded in $L^p(\Omega;$ $\RR^3)$, we infer  from  \eqref{mdx=},    
   \eqref{estimtild11}, \eqref{estimovA},  \eqref{estimovAb}, and \eqref{estimbasic1} that 
  \begin{equation}\begin{aligned}
 \hskip-0,2cm   \int  \hskip-0,1cm |\bfu_\e|^p d m_\e&  \leq  C\int \lb\bfu_\e- \bfgotv_\e(\bfu_\e) \rb^p +   \lb \bfgotv_\e(\bfu_\e) -\ov\bfgotv_\e^S(\bfu_\e) \rb^p\hskip-0,1cm+    \lb \ov\bfgotv_\e^S(\bfu_\e) \rb^p d m_\e 
 \\& \leq C \lp c_\e\int_{T_{r_\e}} \hskip-0,1cm |\bfnabla \bfu_\e|^pdx+ r_\e^p\hskip-0,1cm \int\hskip-0,1cm |\bfe(\bfu_\e)|^p dm_\e+ \int_\Omega |\ov\bfgotv_\e^S(\bfu_\e)|^p dx \rp\leq C,
  \end{aligned}
\label{uemebounded}
  \end{equation}
%
and from \eqref{mdx} that
%
\begin{equation}\begin{aligned} 
  \int   \lb\bfgotp_\e(\bfu_\e) \rb^p  d m_\e \leq C \qquad \forall \bfgotp_\e\in \{ \ov\bfgotv^S_\e, \bfgotv_\e,  \bfgotw_\e, 
   \theta_\e, \ov\theta^S_\e\} .
\end{aligned}
\label{estimbasic3}
  \end{equation}
By \eqref{mdx2},  up to a subsequence,  each of the above   sequences $(\bfgotp_\e(\bfu_\e))$
  weakly converges in $L^p(\Omega;\RR^3)$ to some $\bfp$ 
and two-scale converges w.r.t. $(m_\e)$ to the same $\bfp$. 
Taking  \eqref{mdx0}, \eqref{gotrdto},  \eqref{estimovA} and \eqref{estimbasic1}  into account, we infer 
\begin{equation}\begin{aligned} 
& \bfgotv_\e (\bfu_\e)  \rightharpoonup \bfv ,  \quad \ov\bfgotv_\e^S(\bfu_\e) \rightharpoonup \bfv ,   \quad  \bfgotw_\e (\bfu_\e)  \rightharpoonup \bfw    \quad  \hbox{weakly in    }    L^p(\Omega;\RR^3 ), 
 \\& \theta_\e(\bfu_\e)=\gotw_{\e3}(\bfu_\e) \rightharpoonup \theta=w_3, \quad \ov\theta^S_\e(\bfu_\e)\rightharpoonup \theta  \quad  \hbox{weakly in    }    L^p(\Omega),
\\&  \bfgotr_\e(\bfu_\e)\  \mdto \  \bfv+\frac{2}{\diam S}\bfw\wedge \bfy,
\end{aligned}
\label{piticlop}
  \end{equation}
 for some suitable $\bfv, \bfw, \theta$. It  follows   from  \eqref{estimrigid} and \eqref{estimbasic1} that 
  \begin{equation}\begin{aligned}
   & \bfu_\e\  \mdto \  \bfv+\frac{2}{\diam S}\bfw\wedge \bfy,
 \end{aligned}
\label{uemdto}
  \end{equation}
  and then from 
  \eqref{DsubsetS}  and 
  \eqref{twoscalemeas}  that 
  \begin{equation}\begin{aligned}
   & \bfu_\e m_\e \buildrel\star\over\rightharpoonup  
   \bfv \calL^3_{\lfloor \Omega} \quad \hbox{ weakly$^\star$ in } \ \M(\ov\Omega;\RR^3).
 \end{aligned}
\label{uemetov}
  \end{equation}
 Given $\psi\in \D(\Omega\times S)$, by passing  to the limit as $\e\to 0$ in the  formula
     \begin{equation}\begin{aligned}
   &  r_\e \int  \lp\frac{\partial u_{\e1}}{ \partial x_3}+ \frac{\partial u_{\e3}}{ \partial x_1}\rp \psi\lp x,\frac {y_\e(x')}{r_\e}\rp dm_\e 
 = -r_\e \int  u_{\e1}\frac{\partial\psi}{\partial x_3}\lp x,\frac {y_\e(x')}{r_\e}\rp dm_\e  
 \\&\hskip2cm-r_\e \int  u_{\e3}\frac{\partial\psi}{\partial x_1} \lp x,\frac {y_\e(x')}{r_\e}\rp dm_\e  -
  \int_\Omega u_{\e3}\frac{\partial \psi}{\partial y_1} \lp x,\frac {y_\e(x')}{r_\e}\rp dm_\e,
 \end{aligned}
\nonumber
  \end{equation}
\noindent  taking   \eqref{estimbasic1}  and \eqref{uemdto}    into account,  we obtain
     \begin{equation}\begin{aligned}
     0 &=  
 -\frac{1}{|S|}  \int_{\Omega\times S} \lp v_3(x)  +  \frac{2}{\diam S}w_1(x)y_2-\frac{2}{\diam S}w_2(x)y_1\rp  \frac{\partial\psi}{\partial y_1}  (x, y)dxd y
 \\&=-\frac{1}{|S|}  \int_{\Omega\times S}  \frac{2}{\diam S}  w_2(x)   \psi  (x, y)dxd y,
 \end{aligned}
\nonumber
  \end{equation}
 and
deduce   that $w_2=0$. We likewise obtain  $w_1=0$, yielding
     \begin{equation}\begin{aligned}
    \bfw=w_3\bfe_3=\theta\bfe_3.
 \end{aligned}
\label{w=thetae3}
  \end{equation}
    By  \eqref{twoscaleeta}, \eqref{uemdto} and \eqref{w=thetae3}, the sequence $\lp -\frac{y_{\e2}}{r_\e} u_{\e1} +\frac{y_{\e1}}{r_\e} u_{\e2}\rp$ two-scale converges w.r.t. $(m_\e)$ to 
$-y_2 v_1+y_1 v_2 +\frac{2}{\diam S} |y|^2 \theta$, therefore, by \eqref{DsubsetS} and  \eqref{twoscalemeas}, 
$  \lp -\frac{y_{\e2}}{r_\e} u_{\e1} +\frac{y_{\e1}}{r_\e} u_{\e2}\rp m_\e \buildrel\star\over\rightharpoonup
  {\int\!\!\!\!-}_S\lp -y_2 v_1+y_1 v_2 +\frac{2}{\diam S} |y|^2 \theta\rp dy =\tau\theta$ weakly$^\star$  in $  \M(\ov\Omega)$. 
      The  
assertion \eqref{cvbasic} is  proved. \qed

\noindent  {\it Proof of \eqref{cvbasickappa>0}.} By \eqref{kkappa},   \eqref{mdx},  \eqref{estimDirichlet}, \eqref{estimDirichletb}, and \eqref{estimbasic1},  we have, if  $\kappa>0$, 
\begin{equation}\begin{aligned} 
  \int   \lb\frac{u_{\e3}}{r_\e} \rb^p  d m_\e+  \int   \lb\frac{\ov\gotv^S_{\e3}(\bfu_\e)}{r_\e} \rb^p  d m_\e+  \int   \lb\frac{\ov\theta^S_\e(\bfu_\e)}{r_\e} \rb^p  d m_\e  \leq \frac{C}{k_\e r_\e^p}\leq C,
  \end{aligned}
\nonumber
  \end{equation}
 hence by Lemma \ref{lemtwoscale}, \eqref{mdx0},  and \eqref{mdx2}, 
   the following convergences hold
\begin{equation}\begin{aligned}
&  \frac{u_{\e3}}{r_\e} m_\e \buildrel\star\over\rightharpoonup w    \hbox{ weakly$^\star$  in }   \M (\ov\Omega), \quad  \frac{\ov\gotv^S_{\e3}(\bfu_\e)}{r_\e}  \rightharpoonup w     \hbox{ weakly  in }  L^p(\Omega),
\\&   \frac{\ov\theta^S_\e(\bfu_\e)}{r_\e}  \rightharpoonup \delta     \hbox{ weakly  in }  L^p(\Omega), \quad \frac{\ov\theta^S_\e(\bfu_\e)}{r_\e} m_\e\ \buildrel\star\over\rightharpoonup\ \delta  \hbox{ weakly$^\star$ in }  \M(\ov\Omega),
     \end{aligned}
\label{clop}
  \end{equation} 
 up to a subsequence, for some $w$, $\delta$. By  
  \eqref{votheta3} and \eqref{estimbasic1}, we have
     \begin{equation}\begin{aligned}
 &   \int \lb \frac{ \ov\gots^{S}(\bfu_\e)}{r_\e}\rb^p   dm_\e \leq \frac{C}{r_\e^p k_\e}\leq C,
  \end{aligned}
\nonumber
  \end{equation}
hence the sequence $\lp   \frac{ \ov\gots^{S}(\bfu_\e)}{r_\e}  m_\e\rp$  weakly$^\star$ converges in $\M(\ov\Omega)$, up to a subsequence,  to  some 
$g\in L^p(\Omega)$. 
We set
\begin{equation}
\begin{aligned}
  &\ov \zeta_\e (x   ):= \sum_{i\in I_\e} \lp \intb_{S_{r_\e}^{i}}
\zeta(s_1,s_2,x_3  )ds_1ds_2\rp  \mathds{1}_{Y_\e^{i}}(x_1,x_2) 
\quad \forall \zeta\in L^p(\Omega)
.
\end{aligned}
\label {defovzeta}
\end{equation} 
One can  check that 
\begin{equation}
\begin{aligned}
   |\ov  \zeta_\e -\zeta|_{L^\infty(T_{r_\e})}+ \lb \frac{\partial \ov  \zeta_\e}{\partial x_3} -\frac{\partial \zeta}{\partial x_3}\rb_{L^\infty(T_{r_\e})}\hskip-0,3cm  \leq Cr_\e  \qquad \forall \zeta \in C^1(\ov\Omega).
\end{aligned}
\label {zetaovzeta}
\end{equation} 
We fix  $\zeta\in C^1(\ov\Omega)$,  $\alpha,  \beta\in \{1,2\}$. 
By   \eqref{DsubsetS} and \eqref{defye}, $\int_{S_{r_\e}^{i}}\bfy_\e(x')dx'=0$, 
hence, 
since   $ \ov\gotv^S_{\e\beta}(\bfu_\e)  \ov\zeta_\e$  
is   constant   in each  
 $S_{r_\e}^{i}\times\{x_3\}$,   $     \int    \frac{y_{\e\alpha}}{r_\e}  \ov\gotv^S_{\e\beta}(\bfu_\e)  \ov\zeta_\e  dm_\e =0$. Taking  \eqref{defovgotA}
 and \eqref{defovgots}  into account, we  deduce 
  \begin{equation}\begin{aligned}
   \int  
 \frac{ \ov\gots^{S}(\bfu_\e)}{r_\e} & \ov\zeta_\e  dm_\e 
  = 
\tau \int     \frac{1}{r_\e}  \ov\theta^S_\e(\bfu_\e)  \ov\zeta_\e dm_\e,
\end{aligned}
\nonumber
  \end{equation}
  \noindent where $\tau$ is given by \eqref{defvethetae}. 
  By passing to the limit as $\e\to 0$, thanks to   \eqref{clop} and  \eqref{zetaovzeta}, we obtain
  $\int_\Omega g\eta  dx =   \tau \int_\Omega \delta \eta  dx$. 
    The  
assertion \eqref{cvbasickappa>0} is  proved. \qed

   \noindent {\it  Proof of \eqref{regthetav}.}
 By  \eqref{defqce},  \eqref{estimgotw},   \eqref{estimDirichlet},   \eqref{estimDirichletb}, and \eqref{estimbasic1},  we have 
  \begin{equation}\begin{aligned}
& \int    \lb      \ov\theta^S_\e(\bfu_\e) \rb^p 
 + \lb u_{\e3}\rb^p dm_\e \leq \frac{C}{k_\e},
 \\&\int  \lb  { \bfgotw_\e} (\bfvarphi)  \rb^p dm_\e  \leq C
\lp r_\e^p  
+   \frac{1}{\gamma_\e(r_\e)}  \rp & &\hbox{ if } \ p\not=2, 
 \\&\int  \lb  { \bfgotw_\e} (\bfvarphi)  \rb^p dm_\e  \leq C
\lp r_\e^2  
+  \e^2 \rp & &\hbox{ if } \ p =2, \end{aligned}
\nonumber
  \end{equation}
  therefore,  by  \eqref{mdx=},   \eqref{piticlop}  and \eqref{uemetov},  $\theta=v_3=0$ if $k=+\infty$, 
  $\bfw=0$ and $\theta=0$ if 
  $p=2$ or  $\gamma^{(p)}=+\infty$. By \eqref{cvbasic},  the sequence $(\ov\bfgotu_\e(\bfu_\e)-\ov\bfgotv^S_\e(\bfu_\e))$ weakly converges in $L^p(\Omega;\RR^3)$ to $\bfu-\bfv$, therefore, by 
  \eqref{estimtild1}, \eqref{estimtild1b} and \eqref{estimbasic1}, 
\begin{equation}\begin{aligned}
 \int_{\Omega}\hskip-0,1cm
   |  \bfu - \bfv |^pdx 
 &  \leq \liminf_{\e\to 0} \hskip-0,1cm \int_{\Omega}
   \lb \ov\bfgotu_\e(\bfu_\e)-\ov\bfgotv^S_\e(\bfu_\e)
   \rb^p\hskip-0,1cmdx 
\hskip-0,1cm \leq \hskip-0,1cmC \liminf_{\e\to 0}\lp  \frac{1}{  \gamma^{(p)}_\e\!(r_\e) }+\e^{1-\frac{p}{q}}\hskip-0,1cm \rp\hskip-0,1cm,
 \end{aligned}
   \nonumber \end{equation}
\noindent
thus $\bfu=\bfv$ if  $\gamma^{(p)}=+\infty$.  
If $\kappa=+\infty$,   by   \eqref{estimDirichlet},   \eqref{estimDirichletb}, and \eqref{estimbasic1},  we have 
\begin{equation}\begin{aligned}
& \int    \lb     \frac{ \ov\theta^S_\e(\bfu_\e)}{r_\e} \rb^p + \lb u_{\e1}\rb^p+ \lb u_{\e2}\rb^p
 + \lb \frac{u_{\e3}}{r_\e}\rb^p    dm_\e \leq \frac{C}{r_\e^p k_\e}=o(1),
  \end{aligned}
\nonumber
  \end{equation}
  hence,  by \eqref{cvbasic} and \eqref{clop}, 
  $\bfv=0$ and $\delta=w=0$.
   Proposition \ref{prop1}  is proved.
 \qed  

   \noindent  In the following proposition, we   identify  the two-scale limits  w.r.t.  $(m_\e)$  of  the sequence 
  $(\bfe(\bfu_\e))$  
 in terms of  the functions $v_3$ and $\theta$ given by   \eqref{cvbasic}. 
       %
  \begin{proposition} \label{prop2}       Let  $(\bfu_\e)$ be  a sequence  in $W^{1,p}_b (\Omega;\RR^3)$
 satisfying  (\ref{supFeuefini},\ref{cvbasic}).  
 Then
   \begin{equation}
\begin{aligned}
&  \theta,   v_3   \in  L^p( \Omega';W^{1,p} (0,L))); 
\qquad      \theta=v_3=0  \  \text{ on } \  \Omega'\times\{0\},
\end{aligned}
\label{regthetav1}
  \end{equation}
and, up to a subsequence,
  \begin{equation}\begin{aligned}
& \bfe (\bfu_\e) \mdto \bfgotL\lp \bfq, \frac{\partial v_3}{\partial x_3}, \frac{\partial \theta}{\partial x_3}\rp \hskip1,5cm \hbox{in accordance with } \eqref{defmdto},
     \end{aligned}
\label{cvbasice}
  \end{equation}
for some $\bfq\in L^p(\Omega; W^{1,p}(S;\RR^3))$, where $\bfgotL$ is defined by \eqref{defghom1}.%
 \end{proposition}
 
   \noindent
  {\bf Proof.}   By  \eqref{estimbasic1}  the sequence  $(\bfe(\bfu_\e))$  satisfies  \eqref{hypfie}, hence, by 
  Lemma \ref{lemtwoscale},  
 %
  \begin{equation}\begin{aligned}
  \bfe (\bfu_\e)\ \mdto\ \bfXi_0 ,
   \end{aligned}
\label{cvbasic1}
  \end{equation}
    up to a subsequence for some $\bfXi_0\in L^p(\Omega\times S; \SSym^3)$.
We fix  $\psi\in C^{\infty}(\ov \Omega; \D(S))$ such that $\psi=0 $ on $\Omega'\times \{L\}$. By   passing  to the limit as $\e\to 0$ in the equation
$$
\int   \frac{\partial u_{\e 3}}{ \partial x_3} \psi\lp x ,\frac {y_\e(x')}{r_\e} \rp dm_\e=-\int    u_{\e3} \frac{\partial   \psi}{ \partial x_3}\lp x,\frac{y_\e(x')}{r_\e}\rp dm_\e.
$$
taking \eqref{cvbasic} and \eqref{cvbasic1} into account,  we obtain 
$$ 
\int_{\Omega\times S} (\bfXi_0)_{33} \psi dx dy  = - \int_{\Omega\times S}    v_3 \frac{\partial   \psi}{ \partial x_3} dx dy,
$$
 and infer,
by the arbitrariness of $\psi$,  
 \begin{equation}
  \begin{aligned} 
  v_3  \in L^p( \Omega'; W^{1,p} (0,L)), \quad  v_3=0 \ \text{ on } \ \Omega'\times\{0\},  \quad  \frac{\partial v_3}{\partial x_3}= ( \bfXi_0)_{33}  .
\end{aligned}
\label{vXi}
\end{equation}  
We   fix   $\eta\in C^\infty(\ov\Omega)$ verifying  $\eta=0 $ on $\Omega'\times \{L\}$ and  $\bfvarphi $ such that  
 \begin{equation}
\begin{aligned}
\bfvarphi\in \D(S; \RR^3), \qquad 
\varphi_3=0, 
\qquad   \frac{\partial \varphi_1}{ \partial y_1}+ \frac{\partial \varphi_2}{ \partial y_2}=0.
 \end{aligned}
\label{divvarphi}
\end{equation}
By   (\ref{defovzeta})  and  (\ref{divvarphi}),  we have
\begin{equation}
\begin{aligned}
 \int \bfe(\bfu_\e) :  \ov \eta_\e(x)
 \sum_{\alpha=1}^2 2 \varphi_\alpha \lp  \frac{y_\e(x')}{r_\e}\rp & \bfe_\alpha \odot\bfe_3
  dm_\e
 \\&=- \int   \bfu_\e\cdot \frac{\partial \ov \eta_\e}{\partial x_3}(x) \bfvarphi \lp \frac{y_\e(x')}{r_\e}\rp    dm_\e .
 \end{aligned}
\label{ippk}
\end{equation}
 We set
 (see  \eqref{defgot}):
\begin{equation}
\begin{aligned}
& \int   \bfu_\e\cdot \frac{\partial \ov\eta_\e}{\partial x_3}(x) \bfvarphi \lp \frac{y_\e(x')}{r_\e}\rp    dm_\e=K_{1\e}+K_{2\e}+K_{3\e},  
 \\&  K_{1\e}:= \int  \lp \bfu_\e-\bfgotr_\e(\bfu_\e)  \rp  \cdot\frac{\partial \ov\eta_\e}{ \partial x_3}(x) \bfvarphi \lp \frac{y_\e(x')}{r_\e}\rp dm_\e,
  \\&  K_{2\e}:= \int  
 \bfgotv_\e(\bfu_\e) 
  \cdot  \frac{\partial \ov\eta_\e}{ \partial x_3}(x) \bfvarphi \lp \frac{y_\e(x')}{r_\e}\rp dm_\e,
  \\&  K_{3\e}:= \int   \lp   \frac{2}{\diam S}\bfgotw_\e(\bfu_\e) \wedge  \frac{\bfz_\e( x)}{ r_\e} \rp  \cdot\frac{\partial \ov\eta_\e}{ \partial x_3}(x) \bfvarphi \lp \frac{y_\e(x')}{r_\e}\rp dm_\e.
  \end{aligned}
\label{splitK}
\end{equation}
\noindent By  H\"older's inequality,    \eqref{estimrigid}  and \eqref{estimbasic1},  we have  
\begin{equation}
\begin{aligned}
 |K_{1\e}|\hskip-0,1cm&\leq C
 \lp \int  \lb  \bfu_\e-\bfgotr_\e(\bfu_\e)  
 \rb^p dm_\e
\rp^\frac{1}{p}\leq C
 r_\e   .
\end{aligned}
\label{K10}
\end{equation}
 \noindent 
 Let us fix $\bfxi\in\RR^3$.
Suitable changes of variables in the  equations
$ 0    =  -  \hskip-0,1cm \int_{ S} 
  (\xi_1y_1 + \xi_2 y_2)  \lp \frac{\partial \varphi_1}{ \partial y_1}+ \frac{\partial \varphi_2}{ \partial y_2}\rp   (y) dy  = 
\int_{ S} 
 \bfxi\cdot \bfvarphi(y)
 dy $, deduced from \eqref{divvarphi}, yield
\begin{equation}
\begin{aligned}
 \int_{ S_{r_\e}^{i}}\bfxi\cdot \bfvarphi\lp\frac{y_\e(x')}{r_\e}\rp dx'=0\qquad  \forall  \bfxi\in\RR^3, \quad   \forall i\in I_\e.
\end{aligned}
\label{xifi}
\end{equation}
 Since 
$ \frac{\partial \ov\eta_\e}{ \partial x_3} \bfgotv_\e(\bfu_\e)$ is constant   in each  $S_{r_\e}^{i}\times\{x_3\}$,  
we infer  
 \begin{equation}
\begin{aligned}
 K_{2\e} = \frac{\e^2}{r_\e^2|S|}  \sum_{i\in I_\e} \int_0^L dx_3 \int_{S_{r_\e}^{i}}  \frac{\partial \ov\eta_\e}{ \partial x_3} \bfgotv_\e(\bfu_\e) \cdot \bfvarphi \lp \frac{y_\e(x')}{r_\e}\rp dx' =
 0.
 \end{aligned}
\label{K2=0}
\end{equation}
By \eqref{defye}, \eqref{defze},  \eqref{defgot} and \eqref{divvarphi},  the following holds:   $\bfz_\e(x)= \bfy_\e(x')+z_{\e3}(x)\bfe_3$,  $\varphi_3=0$ and  $\gotw_{\e3}(\bfu_\e)=\theta_\e(\bfu_\e)$, therefore
%
 \begin{equation}
\begin{aligned}
 & K_{3\e} =  K_{3\e}^{(1)}+ K_{3\e}^{(2)},
 \\&  K_{3\e}^{(1)}:=\int   \lp   \frac{2}{\diam S}\theta_\e(\bfu_\e)\bfe_3 \wedge  \frac{\bfy_\e( x')}{ r_\e} \rp   \cdot\frac{\partial \ov\eta_\e}{ \partial x_3}(x)  \bfvarphi \lp \frac{y_\e(x')}{r_\e}\rp  dm_\e,
 \\&  K_{3\e}^{(2)}\hskip-0,1cm:=\hskip-0,1cm\frac{2}{\diam S}\hskip-0,1cm \int \hskip-0,1cm\frac{z_{\e3}( x)}{ r_\e} \frac{\partial \ov\eta_\e}{ \partial x_3} \hskip-0,1cm  \lp\hskip-0,1cm -\gotw_{\e1}(\bfu_\e) \varphi_2\lp \frac{y_\e(x')}{r_\e}\rp\hskip-0,1cm
   + \hskip-0,1cm\gotw_{\e2}(\bfu_\e)  \varphi_1 \lp \frac{y_\e(x')}{r_\e}\rp \hskip-0,1cm
   \rp \hskip-0,1cmdm_\e.
 \end{aligned}
\label{K3e=}
\end{equation}
Since  $z_{\e3}$,  $\bfgotw_\e(\bfu_\e)$, and $ \frac{\partial \ov\eta_\e}{ \partial x_3}$   are     constant   in 
 $S_{r_\e}^{i}\times\{x_3\}$,  we deduce from   \eqref{xifi}  that
\begin{equation}
\begin{aligned} 
K_{3\e}^{(2)}   = 0. \end{aligned}
\label{K3e2=0}
\end{equation}
By    \eqref{mdx2} and \eqref{piticlop},   we have $\theta_\e(\bfu_\e)m_\e \ \mdto  \theta\bfe_3 \wedge  \bfy$, hence, by \eqref{zetaovzeta}, 
\begin{equation}
\begin{aligned} 
\lime K_{3\e}^{(1)}   =  \frac{1}{|S|} \int_{\Omega\times S}    \lp   \frac{2}{\diam S}\theta(x)\bfe_3 \wedge  {\bfy } \rp  \cdot\frac{\partial  \eta}{ \partial x_3}(x) \bfvarphi \lp y  \rp dxdy.  \end{aligned}
\label{limK3e1=}
\end{equation}
In view of  \eqref{zetaovzeta}, \eqref{cvbasic1}, \eqref{splitK}, \eqref{K10}, \eqref{K2=0}, \eqref{K3e=}, \eqref{K3e2=0}, and \eqref{limK3e1=}, by passing to the limit as $\e\to0$  in (\ref{ippk}),  we obtain    
\begin{equation}
\begin{aligned}
\frac{1}{|S|} \int_{\Omega\times S} &2((\bfXi_0 )_{13}\bfe_1  +  (\bfXi_0 )_{23}\bfe_2) \cdot   \eta (x)\bfvarphi(y)
  dxdy  
  \\&=-  \frac{1}{|S|} \int_{\Omega\times S}    \lp   \frac{2}{\diam S}\theta(x)\bfe_3 \wedge  {\bfy } \rp  \cdot\frac{\partial  \eta}{ \partial x_3}(x) \bfvarphi \lp y  \rp dxdy. 
 \end{aligned}
\label{Xt}
\end{equation}
Therefore $\lb \int_\Omega \theta\frac{\partial \eta}{\partial x_3} dx \rb \leq C |\eta|_{L^{p'}(\Omega)}$
and, by   the arbitrariness of $\eta $,
\begin{equation}
\begin{aligned}
  \theta \in L^p( \Omega'; W^{1,p} (0,L)), \quad \theta=0 \ \text{ on } \ \Omega'\times\{0\}. 
 \end{aligned}
\label{W1p0}
\end{equation}
By   \eqref{vXi}  and \eqref{W1p0},  \eqref{regthetav1}  is proved. 
By integrating  (\ref{Xt})  by parts, 
 we find
\begin{equation}
\begin{aligned}
 & \int_{ \Omega\times S}  
\begin{pmatrix} (\bfXi_0)_{13}  +\frac{1}{2}\frac{2}{\diam S} \frac{\partial \theta}{\partial x_3}(x)  y_2 \\
  (\bfXi_0)_{23}   - \frac{1}{2} \frac{2}{\diam S}  \frac{\partial \theta}{\partial x_3} (x)  y_1 
\end{pmatrix}\cdot\eta(x) \begin{pmatrix}\varphi_1( y)
\\\,  \varphi_2( y)\
\end{pmatrix} dx d y
=0,
 \end{aligned}
\nonumber
\end{equation}
and, by the arbitrary choice of   $\eta\in \D(\Omega)$ and     $\bfvarphi$ verifying \eqref{divvarphi},   deduce that
 \begin{equation}
\begin{aligned}
 &  
 (\bfXi_0)_{13}  =  \frac{\partial q_3}{ \partial y_1}  -\frac{1}{2}\frac{2}{\diam S} \frac{\partial \theta}{\partial x_3}(x) y_2, 
\quad   (\bfXi_0)_{23}   = \frac{\partial q_3}{ \partial y_2} +  \frac{1}{2} \frac{2}{\diam S}  \frac{\partial\theta}{ \partial x_3} (x) y_1,
  \end{aligned}
\label{ipp52}
\end{equation}
   for some $q_3 \in L^p(\Omega; W^{1,p}(S))$. 
We fix    $\bfPsi\in \D(S; \SSym^3)$ such that  $\Psi_{i3}=0  \ \forall i\in \{1,2,3\}$   and  $\bfdiv\, \bfPsi=0$.  By \eqref{defovzeta}, we have 
 $
 \int   \bfe(\bfu_\e) :\ov\eta_\e(x) \bfPsi\lp\frac {y_\e(x')}{r_\e}\rp dm_\e = 0.
 $
Passing to the limit as $\e\to 0$, in view of  \eqref{zetaovzeta} and  \eqref{cvbasic1}, 
  we obtain 
 $
   \int_{\Omega\times S} \bfXi_0:  \eta (x) \bfPsi(y) dxdy  = 0.
 $
By the arbitrary choice of $\eta$ and $\bfPsi$, we deduce from  
  a generalization of the Donati's theorem    \cite[Th. 1]{Mo} that
\begin{equation}
\begin{aligned}
    (\bfXi_0)_{\alpha\beta} (x,y)   = \frac{1}{2} \lp \frac{\partial q_\alpha}{ \partial y_\beta}+\frac{\partial q_\beta}{\partial y_\alpha } \rp (x,  y)
 \quad \forall \alpha, \beta\in   \{1,2\}    \quad \text{ in } \Omega\times S,
 \end{aligned}
\label{ipp01}
\end{equation}
 for some $q_1, q_2 \in L^p(\Omega; W^{1,p}(S))$.
Combining   \eqref{defghom1},   (\ref{vXi}),  (\ref{ipp52}),  and   (\ref{ipp01}),  we  obtain   $ \bfXi_0    = \bfgotL\lp \bfq, \frac{\partial v_3}{\partial x_3}, \frac{\partial \theta}{\partial x_3} \rp$.    Proposition \ref{prop2}  is proved. \qed

 \noindent  In the next  proposition, we  derive  the two-scale limits  w.r.t. $(m_\e)$  of  the sequences 
 $\lp  \frac{u_{\e 3}}{  r_\e}\rp$ and 
  $( \frac{1}{r_\e} \bfe(\bfu_\e))$   in terms of $w,\delta, v_1, v_2$ given by   \eqref{cvbasic} and  \eqref{cvbasickappa>0}, in the case $\kappa>0$.
   \begin{proposition} \label{prop3}    Let  $(\bfu_\e)\!\subset \!W^{1,p}_b \!(\!\Omega;\RR^3)$
 satisfying \eqref{supFeuefini},  \eqref{cvbasic}, \eqref{cvbasickappa>0}.    If $\kappa\!>\!0$, 
 %
  \begin{equation}\begin{aligned}
& \bfv\in
 L^p(\Omega';W^{2,p} (0,L;\RR^3)); \quad \delta, w\in L^p(\Omega' ; W^{1,p} (0,L)),
\\& v_\alpha=  \frac{\partial  v_\alpha }{\partial x_3}=\delta=w=0  \text{ on }  \Omega'\times\{0\} \quad \forall\alpha\in \{1,2\},
    \end{aligned}\label{regvdeltaw}
  \end{equation}
  and,    up to a subsequence,
%
   %
  \begin{equation}\begin{aligned}
&
   \frac{u_{\e 3}}{  r_\e}\   \mdto  \ w -\frac{\partial v_1}{\partial x_3} y_1-\frac{\partial v_2}{\partial x_3} y_2,  
   \quad \frac{1}{r_\e} \bfe(\bfu_\e) \mdto   \bfgotJ\lp \bfl, \frac{\partial^2 v_1  }{ \partial x_3^2}, \frac{\partial^2 v_2 }{\partial x_3^2},
\frac{\partial w}{ \partial x_3}, \frac{\partial \delta }{\partial x_3}\rp ,
  \end{aligned}\label{cvkappa}
  \end{equation}
 for some $\bfl\in L^p(\Omega; W^{1,p}(S;\RR^3))$, being  $\bfgotI$ defined by \eqref{defghom2}.
    Moreover, 
\begin{equation}\begin{aligned}
 \bfv^{tuple}\in \D \qquad \forall (k,\kappa)\in (0,+\infty]\times  [0,+\infty].
 \end{aligned}
 \label{vtupleinD}
  \end{equation}
   \end{proposition}
 
   \noindent
{\bf Proof.} \quad    If $\kappa>0$,   by \eqref{kkappa}, \eqref{estimDirichletb},  and   \eqref{estimbasic1},  we have  $\int \lb\frac{u_{\e3}}{r_\e}\rb^p+  \lb\frac{1}{r_\e} \bfe(\bfu_\e)\rb^p dm_\e\leq C$. Applying Lemma    \ref{lemtwoscale}, taking  \eqref{twoscalemeas} and \eqref{cvbasickappa>0} into account, we infer  
\begin{equation}\begin{aligned}
 \frac{u_{\e 3}}{  r_\e}\   \mdto  w_0; \qquad 
\frac{1}{r_\e} \bfe(\bfu_\e) \mdto \bfUpsilon_0; \qquad \intb_S w_0(x,y)dy=w(x),
 \end{aligned}
 \label{w0U0}
  \end{equation}
  up to a subsequence,
  for some  $w_0\in L^p(\Omega\times S)$,   $\bfUpsilon_0\in L^p(\Omega\times S; \SSym^3)$. 
We fix  $\alpha\in\{1;2\}$ and  $\psi\in C^{\infty}(\ov\Omega; \D(S))$ such that $\psi=0 $ on $\Omega'\times \{L\}$. By   integration  by parts,  
\begin{equation}
\begin{aligned}
&  \int  2 (\bfe(\bfu_\e))_{\alpha3}  \psi \lp x,\frac {y_\e(x')}{r_\e} \rp dm_\e=
\\& -\int u_{\e3}\frac{\partial\psi}{\partial x_\alpha} \lp x,\frac {y_\e(x')}{r_\e} \rp  
+\frac{u_{\e3}}{r_\e} \frac{\partial\psi}{\partial y_\alpha} \lp x,\frac {y_\e(x')}{r_\e} \rp  
+ u_{\e\alpha }\frac{\partial\psi}{\partial x_3} \lp x,\frac {y_\e(x')}{r_\e} \rp dm_\e.
 \end{aligned}
\nonumber
\end{equation}
\noindent 
Passing to the limit as $\e\to0$,  taking 
  \eqref{w0U0} into account and 
noticing that, by  \eqref{cvbasic} and  \eqref{regthetav},    $\bfu_\e\mdto v_1\bfe_1+ v_2\bfe_2$,  we obtain
\begin{equation}
\begin{aligned}
0= -\int_{\Omega\times S} w_0(x,y)  \frac{\partial\psi}{\partial y_\alpha} (x,y)+ v_\alpha(x)\frac{\partial\psi}{\partial x_3} (x,y)dxdy.
 \end{aligned}
\label{w0va}
\end{equation}
We deduce from   the arbitrariness  of $\psi$ that  $
 \frac{\partial v_\alpha}{\partial x_3}  \in L^p(\Omega)$,   $\frac{\partial w_0}{\partial y_\alpha}\in L^p(\Omega\times S)$,  $ 
 \frac{\partial w_0}{\partial y_\alpha}(x,y) = -\frac{\partial v_\alpha}{\partial x_3}(x) $  a. e.  in  $\Omega\times S$ and $ v_\alpha=0$  a. e.  in $ \Omega'\times \{0\}$. 
Thus   $\bfnabla_y \lp w_0(x,y)+ \frac{\partial v_1}{\partial x_3}(x) y_1+ \frac{\partial v_2}{\partial x_3}(x) y_2\rp=0$ and, by   \eqref{DsubsetS} and 
  \eqref{w0U0},
$$
w_0+ \frac{\partial v_1}{\partial x_3}(x) y_1+ \frac{\partial v_2}{\partial x_3}(x) y_2= 
\intb_S \hskip-0,1cm\lp \hskip-0,1cm w_0(x,y)+ \frac{\partial v_1}{\partial x_3}(x) y_1 \hskip-0,1cm+ \frac{\partial v_2}{\partial x_3}(x) y_2 \hskip-0,1cm\rp dy =w(x).
$$
Passing to the limit as $\e\to0$ in 
$\int\frac{1}{r_\e}\frac{\partial u_{\e3}}{\partial x_3}\psi dm_\e=- \int \frac{u_{\e3}}{r_\e} \frac{\partial \psi}{\partial x_3} dm_\e$, 
we obtain 
$
\int_{\Omega\times S} \Upsilon_{033} \psi dxdy= - \int_{\Omega\times S}  \lp  -\frac{\partial v_1}{\partial x_3} y_1-\frac{\partial v_2}{\partial x_3} y_2 + w\rp \frac{\partial \psi}{\partial x_3} dxdy, $ %
and infer   
\begin{equation}
\begin{aligned}
&v_\alpha   \in 
 L^p(\Omega';W^{2,p} (0,L));\quad  
  w \in  L^p(\Omega' ; W^{1,p} (0,L)); \quad \frac{\partial v_\alpha }{\partial x_3}= w=0  \text{ on }  \Omega'\times\{0\},
 \\&\Upsilon_{033}(x,y)= -\frac{\partial^2 v_1}{\partial x_3^2}(x) y_1-\frac{\partial^2 v_2}{\partial x_3^2} (x)y_2 + \frac{\partial w}{\partial x_3}(x)
 .
 \end{aligned}
\nonumber
\end{equation}
\noindent 
Substituting $\lp\frac{\theta_\e(\bfu_\e)}{r_\e}\rp$ for $\lp {\theta_\e(\bfu_\e)} \rp$,  and $\frac{1}{r_\e} \bfe(\bfu_\e)$ for $  \bfe(\bfu_\e)$ 
in    the   argument employed  to establish  \eqref{W1p0}, \eqref{ipp52},  \eqref{ipp01}, 
  we find 
\begin{equation}
\begin{aligned}
  & (\bfUpsilon_0)_{\alpha\beta} (x,y) 
= \frac{1}{2} \lp \frac{\partial l_\alpha}{\partial y_\beta}+\frac{\partial l_\beta  }{\partial y_\alpha } \rp (x,  y)
    \quad \text{ in } \Omega\times S \quad \forall \alpha, \beta\in   \{1,2\},
\\ & \delta \in L^p( \Omega'; W^{1,p} (0,L)), \qquad \delta=0 \ \text{ on } \ \Omega'\times \{0\}, 
 \\ &  
 (\bfUpsilon_0)_{13}  =  \frac{1}{2}\frac{\partial l_3}{ \partial y_1}  -\frac{1}{2}\frac{2}{\diam S} \frac{\partial \delta}{\partial x_3} y_2, 
\quad
   (\bfUpsilon_0)_{23}   =  \frac{1}{2}\frac{\partial l_3}{ \partial y_2} +  \frac{1}{2} \frac{2}{\diam S}  \frac{\partial\delta}{ \partial x_3}   y_1,
  \end{aligned}
\nonumber
\end{equation}
  for some    $\bfl  \in  L^p(\Omega; W^{1,p}(S;\RR^2))$, and deduce  $ \bfUpsilon_0 \hskip-0,1cm= \hskip-0,1cm\bfgotJ\lp \bfl, \frac{\partial^2 v_1}{\partial x_3^2},  \frac{\partial^2 v_2}{\partial x_3^2}, \frac{\partial w}{\partial x_3}, \frac{\partial \delta}{\partial x_3} \rp$ (see  (\ref{defghom2})). 
The proofs of   \eqref{regvdeltaw} and (\ref{cvkappa})  are completed. 
The assertion \eqref{vtupleinD} follows  from \eqref{defD},  \eqref{regthetav}, \eqref{regthetav1},  and \eqref{regvdeltaw}.
 \qed  
 
  \section{Properties of  $\capsca^f$}\label{secpropcap}

 The proof of Theorem \ref{th} relies in the apriori estimates and convergences established in Section \ref{secpreliminaries} and  an investigation into the properties of the capacity  $\capsca^f$ developed below.
In what follows,  $f: \SSym^3\to  \RR$  denotes    a   convex function, not necessarily strictly so,    satisfying  \eqref{growthp}
  for some $p\in (1,+\infty)$ and some positive constants $c,C$. 
For every  nonempty bounded    Lipschitz  domain      $S\subset \RR^2$ 
verifying  \eqref{DsubsetS}
 and   every  open set   $V\subset \RR^2$ 
such that      $\ov S \subset V$, we consider  the mapping  $\capsca^f(.;S,V) :(\RR^3)^2\to [0,+\infty)$ defined by \eqref{infW}.
The infimum  problem $ \calP^f(\bfa,\bfalpha;S,V )$  may  fail    to be    attained 
when     $V$ is not bounded (see Remark \ref{remnotattained}).   
We prove  below that, if $1<p<2$, $\capsca^f$ is  equivalently  defined  by  a   well posed minimum  problem, namely 
\begin{equation}
\begin{aligned}
& \capsca^f(\bfa,\bfalpha;S,V )  =  \min \calP^f_K(\bfa,\bfalpha;S,V ),
\\& \calP^f_K(\bfa,\bfalpha;S,V ):  \min_{\bfvarphi\in  \K^p(\bfa,\bfalpha;S,V)  }   \int_V   f(  \bfe_y(\bfvarphi)) dy,
  \\&   \K^p(\bfa,\bfalpha;S,V):= \la\bfvarphi\in K_0^p (V;\RR^3)  , \quad  \bfvarphi =  \bfa+ \frac{2}{\diam S}\bfalpha\wedge  \bfy   \hbox{ in }  S \ra,
 \end{aligned}
\label{minK}
    \end{equation}
where  $K_0^p (V ; \RR^3)$  is the closure   of $\D(V;\RR^3)$  in   the reflexive Banach space 
      \begin{equation}\begin{aligned}
      &  K^p (V ; \RR^3) := \la \bfvarphi \in L^{p^{\star}}(V;\RR^3), \ \bfnabla \bfvarphi \in L^p(V; \RR^3\times \RR^2)\ra,
       \\&  |\bfvarphi|_{K^p(V;\RR^3)} :=  \lp   \int_V |  \bfvarphi|^{p^{\star}}  dy\rp^{\frac{1} {p^{\star}}} + \lp   \int_V |\bfnabla \bfvarphi|^p  dy \rp^\frac{1}{p},
 \qquad p^{\star} :=  \frac{2p}{2 - p}.
\end{aligned} \label{K0}
  \end{equation}
By the Poincar\'e inequality,   $K_0^p(V;\RR^3)$ is equal to $W^{1,p}_0(V; \RR^3)$ when  $V$ is bounded in one direction. Otherwise, it
 may be    larger.  
 
   \begin{lemma} \label{reg}
(i)   The functional   $\bfvarphi\to  \int_V   f(  \bfe_y(\bfvarphi)) dy$ is 
 strongly  continuous  and weakly lower semi-continuous  on $W^{1,p}(V;\RR^3)$  and,  if $1<p<2$,  
 on $K^p(V;\RR^3)$. 
\\(ii)  Problem   $\calP^f(\bfa,\bfalpha;S,V )$  defined by \eqref{infW} and,  if $1\!<\!p\!<\!2$,   $\calP^f_K(\bfa,\bfalpha;S,V )$  given   by \eqref{minK},        have minimizing sequences   in $\D(V;\RR^3)$.   
\\
(iii) 
If $1<p<2$, 
$\calP^f_K(\bfa,\bfalpha;S,V )$ has a solution, 
unique if   $f$ is strictly convex. Moreover,  \eqref{minK} holds and,  for all   $ \bfvarphi\in  K_0^p(V;\RR^3)$, 
  \begin{equation}
 \hskip-0,1cm\begin{aligned}
 \hskip-0,1cm  \lp\int_V  \hskip-0,1cm|   \bfvarphi |^{p^{\star}}   dy  \rp^{\frac{p}{p^\star}}  & \hskip-0,2cm 
 \leq \frac{3p}{2-p}
 \int_V  \hskip-0,1cm  | \bfnabla  \bfvarphi |^p  dy \leq   C(p)   \int_V  \hskip-0,1cm |\bfe_y(\bfvarphi )|^p dy,
 \end{aligned}
 \label{minKorn}
    \end{equation} 
for some  $C(p)>0$   independent of $V$. 
 \\      (iv)  If  $V$ is bounded in one direction and  $p\in (1,+\infty)$, the same holds   for   $\calP^f(\bfa,\bfalpha;S,V )$.
Moreover,  for all $ \bfvarphi\in  W^{1,p}_0(V;\RR^3)$,
  \begin{equation}
\hskip-0,1cm\begin{aligned}
\hskip-0,1cm \int_V\hskip-0,1cm |   \bfvarphi |^p  dy    \leq   C(p,\hskip-0,05cmV)  \hskip-0,1cm \int_V \hskip-0,1cm  |\bfnabla\bfvarphi |^p  dy \leq  
  C(p,\hskip-0,05cmV)  \hskip-0,1cm \int_V  \hskip-0,1cm |\bfe_y(\bfvarphi )|^p  dy   .
  \end{aligned}
 \label{minKornVbounded}
    \end{equation}
    \\ (v) There exists    positive constants    $c(p,S,V)$, $C(p,S,V)$  such that
  \begin{equation}
\begin{aligned}
  \capsca^f(\bfa,\bfalpha;S,V) \leq C(p,S,V) (|\bfa|^p+|\bfalpha|^p)\quad \forall  (\bfa,\bfalpha)\in (\RR^3)^2,
  \end{aligned}
 \label{boundp}
    \end{equation}
and, if  $1\!<\!p\!<\!2$ or  $V$ is bounded in one direction, 
    %
     \begin{equation}
\begin{aligned}
  c(p,S,V) (|\bfa|^p+|\bfalpha|^p)\leq   \capsca^f(\bfa,\bfalpha;S,V)  \qquad \forall   (\bfa,\bfalpha)\in (\RR^3)^2.
  \end{aligned}
  \label{boundp2}
    \end{equation}
    \end{lemma}

\noindent
   \begin{remark} \label{remnotattained} We prove  below   (see  Remarks \ref{remMf},  \ref{rem0})  that,  if   $p\ge 2$,  $ \capsca^f(\bfa,0;$ $S,\RR^2)=0$,    therefore     \eqref{boundp2} doesn't hold  for $V=\RR^2$ and  the infimum  
   problem
$ \calP^f(\bfa,0;S,\RR^2)$ is not attained when  $\bfa\not=0$.   
  This     is comparable  to  the  Stokes'   paradox.
       \end{remark} 
     
    \noindent  {\bf Proof.}  (i) By  
       \eqref{growthp}, 
    the  functional $\bfvarphi\to  \int_V   f(  \bfe_y(\bfvarphi)) dy$ is convex and  bounded on the unit ball of  $W^{1,p}(V;\RR^3)$,
 hence 
 strongly  continuous, thus  weakly lower semi-continuous. 
The same     holds with $K^p(V;\RR^3)$ in place of  $W^{1,p}(V;\RR^3)$ if  $1<p<2$.
\\  (ii) Follows  from (i) and  a density argument as in   \cite[Lemma  1]{BeArma}.
 \\
(iii) Given  $\bfvarphi\in \D(V;\RR^3)$,  applying \eqref{LpfirstKorn0}
 for $N=2$  to $\bfvarphi'$ defined by 
   \eqref{defexprim}   seen as an element of $\D(V;\RR^2)$, noticing that  $|\bfe(\bfvarphi')| \le |\bfe_y(\bfvarphi)|$ and 
  $ | \bfnabla \varphi_3| \leq C |\bfe_y(\bfvarphi)|$,
  we deduce
  

%
  \begin{equation}
\hskip-0,2cm\begin{aligned}
 \int_{V} |\bfnabla \bfvarphi|^p dy&\leq 2^{p-1} \int_V | \bfnabla \bfvarphi'|^p+  | \bfnabla \varphi_3|^p dy \leq 
  C(p)   \int_V  |\bfe(\bfvarphi')|^p +|\bfe_y(\bfvarphi )|^p dy 
  \\&  \leq   C(p)  \int_V   |\bfe_y(\bfvarphi )|^p  dy \hskip1 cm  \forall p\in(1,+\infty),\   \forall \bfvarphi\in  \D(V;\RR^3).
\end{aligned}
\label{LpfirstKorn}
    \end{equation}
    By the    Sobolev embedding theorem in $\RR^2$
  \cite[Th. 9.9]{Br}, we have 
 \begin{equation}
\hskip-0,2cm \begin{aligned}
  &\lp\int_{V} |\bfvarphi|^{p^{\star}} dy \rp^{\frac{1}{p^{\star}}} \hskip-0,2cm \leq \frac{3p}{2-p}\lp \int_{V} |\bfnabla \bfvarphi|^p dy\rp^{\frac{1}{p}}
 \quad \forall p\in(1,2),\
   \forall \bfvarphi\!\in\!  \D(V;\RR^3).\hskip-0,1cm 
\end{aligned}
\label{GNS}
   \end{equation}
 The estimate  \eqref{minKorn}  results from      \eqref{LpfirstKorn},  \eqref{GNS}, and  the  density  of $\D(V,\RR^3)$ in $K_0^p(V;\RR^3)$. 
We fix a sequence   $(\bfvarphi_n)_{n\in \NN^{\star}}\subset  \K^p(\bfa,\bfalpha;S,V)$   satisfying  
 $\int_V  f(\bfe_y(\bfvarphi_n)) dy $ $\leq \capsca^f(\bfa,\bfalpha;S,V ) + \frac{1}{n}$  $\forall n\in \NN^{\star}$.
The set  $\K^p(\bfa,\bfalpha;S,V)$   is convex and  strongly closed  in $K^p(V;\RR^3)$, thus weakly closed.
By  \eqref{growthp}
and \eqref{minKorn}, 
  $(\bfvarphi_n)_{n\in \NN^{\star}}$  is bounded 
  in  
   $K^p(V;\RR^3)$,
thus weakly converges, up to a subsequence,
to some
  $\bfvarphi \in \K^p(\bfa,\bfalpha;S,V)$. 
Taking  (i) into account, we  deduce 
that  $ \capsca^f(\bfa,\bfalpha;S,V )  \leq \int_V  f(\bfe_y(\bfvarphi)) dy$ $\leq \liminf_{n\to+\infty}\int_V  f(\bfe_y(\bfvarphi_n)) dy$ $ = \capsca^f(\bfa,\bfalpha; S,V ) $, 
thus $\bfvarphi$  is a solution to
 (\ref{minK}). Its   uniqueness when  $f$ is strictly convex  is straightforward. 
Taking  (ii) into account,  we deduce that  $\capsca^f$
  satisfies \eqref{minK}. 
\\   (iv) The estimate  \eqref{minKornVbounded} follows from  \eqref{LpfirstKorn} and the Poincar\'e inequality.
We conclude 
  by   replacing \eqref{minKorn} by \eqref{minKornVbounded} 
 and $K_0^p(V;\RR^3)$ by $W^{1,p}_0(V;\RR^3)$ in the previous  argument.  
\\ (v)  Given   $\eta\in \D(V)$    such   that $\eta=1$ in $S$, the  field    
$\bfvarphi(y):= \eta(y)\lp \bfa+\tfrac{2}{\diam S} \bfalpha\wedge \bfy \rp$ belongs to $\W^p(\bfa,\bfalpha;S,V)$, thus
 %
\begin{equation}
\begin{aligned}
  \capsca^f(\bfa,\bfalpha;S,V) \hskip-0,1cm &\le\hskip-0,1cm \int_V\hskip-0,1cm f(\bfe_y(\bfvarphi))dy\hskip-0,1cm \leq C\hskip-0,1cm \int_V\hskip-0,1cm |\bfe_y(\bfvarphi)|^p dy
 \hskip-0,1cm\leq C (|\bfa|^p+|\bfalpha|^p)\hskip-0,1cm\int_V\hskip-0,1cm|\bfnabla\eta|^p dy,
     \end{aligned}
 \nonumber
    \end{equation}
yielding  \eqref{boundp}.
Assume that  $1<p<2$  and 
  let  
  $V_1$  be  a bounded Lipschitz domain of $\RR^2$ such that $\ov S\subset V_1\subset V$. 
  Let  $\bfvarphi$  be   a solution to \eqref{minK}.
  By   H\"older's inequality,    $|\bfvarphi|_{L^p(V_1;\RR^3)}\leq   |V_1|^{\frac{1}{2}} |\bfvarphi|_{L^{p^\star}(V_1;\RR^3)}$, hence  
 $|\bfvarphi|_{W^{1,p}(V_1;\RR^3)}\leq  (1+|V_1|^{\frac{1}{2}})  |\bfvarphi|_{K^p(V;\RR^3)}$,
yielding 
  $|\bfvarphi|_{L^p_{\calH^1}(\partial S;\RR^3)}\leq    C(S,V) |\bfvarphi|_{K^p(V;\RR^3)}$  by  the continuity of the trace   from $W^{1,p}(V_1;\RR^3)$ to $L^p_{\calH^1}(\partial S;\RR^3)$.  Since $\bfvarphi \subset  \K^p(\bfa,\bfalpha;S,V)$,
taking   \eqref{growthp} and   \eqref{minKorn} into account, we deduce
\begin{equation}
\begin{aligned}
\int_{\partial S}\lb \bfa+\tfrac{2}{\diam S}\bfalpha\wedge \bfy\rb^p \hskip-0,1cm &d\calH^1(y) 
\leq C(S,V) |\bfvarphi|_{K^p (V ; \RR^3) }^p \leq C(S,V) \hskip-0,1cm\int_V \hskip-0,1cm |\bfe_y(\bfvarphi)|^pdx 
\\& \leq C(S,V) \hskip-0,1cm\int_V \hskip-0,1cm f(\bfe_y(\bfvarphi))dx 
 = C(S,V) \capsca^f (\bfa,\bfalpha;S,V) .     \end{aligned}
 \nonumber
    \end{equation}
     Noticing that $|\bfa|^p+|\bfalpha|^p \leq C(S) 
\int_{\partial S}\lb \bfa+\tfrac{2}{\diam S}\bfalpha\wedge \bfy\rb^p d\calH^1(y)$ for some $C(S)>0$, 
\eqref{boundp2}  is proved for $1<p<2$. 
The proof of the other case     is similar.
 \qed
 
     %
\noindent We establish below some continuity properties  for
  $(\bfa,\bfalpha)\to \capsca^f(\bfa,\bfalpha;S,V)$.
      %
     \begin{lemma}\label{lemconvex}  The mapping 
      \begin{equation}
\begin{aligned}
(\bfa,\bfalpha)\in (\RR^3)^2 \to \capsca^f(\bfa,\bfalpha;S,V)\in \RR
 \end{aligned}
 \label{map}
    \end{equation}
    is convex (resp. strictly convex if  $f$ is strictly convex and either  $1<p<2$ or $V$ is bounded in one direction), hence continuous.  The   functional  
      \begin{equation}
\begin{aligned}
(\bfa,\bfalpha)\in (L^p(\Omega;\RR^3))^2
 \to \int_\Omega \capsca^f(\bfa(x),\bfalpha(x);S,V)dx 
   \end{aligned}
 \label{func}
    \end{equation}
    is convex  (resp. strictly convex under the above stated  additional assumptions) and  strongly continuous, hence weakly lower-semicontinuous.
      \end{lemma}
      
\noindent{\bf Proof.} \ 
Let us fix   $\lambda\in (0,1)$,    $t>0$,  and $((\bfa_1,\bfalpha_1), (\bfa_2,\bfalpha_2))\in (\RR^3\times \RR^3)^2$
such that $(\bfa_1, \bfalpha_1)\not=(\bfa_2,\bfalpha_2)$. 
 By \eqref{infW}, there exists     $(\bfvarphi_1 , \bfvarphi_2)\hskip-0,1cm \in\W^p(\bfa_1,\bfalpha_1;S,V)$ $\times\W^p(\bfa_2,\bfalpha_2;S,V)$    satisfying  
 \begin{equation}
\begin{aligned}
  \int_V f(   \bfe_y(\bfvarphi_k )) dy\leq \capsca^f(\bfa_k,\bfalpha_k;S,V) +t \qquad \forall k\in \{1,2\}.
    \end{aligned}
 \label{lb}
    \end{equation}
Observing  that  $   \lambda \bfvarphi_1+ (1-\lambda)\bfvarphi_2 \in\W^p( \lambda \bfa_1 +(1-\lambda) \bfa_2, \lambda \bfalpha_2+(1-\lambda)\bfalpha_2;S,V)$, 
we deduce from 
   the  convexity of $f$ that
    \begin{equation}
\begin{aligned}
  \capsca^f(  \lambda \bfa_1+(1-&\lambda) \bfa_2 , \lambda \bfalpha_1 +(1-\lambda) \bfalpha_2 ;S,V)
\leq 
 \int_V f(    \bfe_y(  \lambda \bfvarphi_1+ (1-\lambda)\bfvarphi_2)) dy
 \\&\leq  \lambda  \int_V f(    \bfe_y(  \bfvarphi_1 )) dy+(1-\lambda) \int_V f(    \bfe_y( \bfvarphi_2 )) dy
 \\& \leq 
  \lambda  \capsca^f(\bfa_1,\bfalpha_1;S,V) + (1-\lambda)  \capsca^f(\bfa_2; \bfalpha_2;S,V) + 
 t ,
 \end{aligned}
\label{inconvex}
    \end{equation}
hence
 the mapping  \eqref{map} 
 is convex.
If  $f$ is strictly convex and  $1<p<2$,   let 
    $\bfvarphi_k$ now   denote   a solution to 
 $\calP^f_K(\bfa_k,\bfalpha_k;S,V)$   (see Lemma \ref{reg} (iii)). Then  \eqref{lb} and     \eqref{inconvex}  
    hold with $t=0$ and, since  $ \bfe_y(\bfvarphi_1 )\not= \bfe_y(\bfvarphi_2 )$,   the second inequality in  \eqref{inconvex}  is strict, thus the mapping 
\eqref{map}  is strictly convex. 
The same holds for  the functional \eqref{func}.
We obtain  the same conclusion when  $V$ is bounded in one direction. 
  By \eqref{boundp},  
the convex  functional   \eqref{func}  is 
   bounded on the unit ball of    $(L^p(\Omega;\RR^3))^2$,  thus    strongly continuous, so  weakly lower-semicontinuous. \qed

\noindent 
We   state below 
some    monotonicity  properties of $\capsca^f(\bfa,\bfalpha;S,V)$ with respect to $V$, 
 $S$ and  $f$,  whose proofs
are so easy that we omit them.
\begin{lemma}\label{lemdecr} 
(i) Let  $V_1$,  $V_2$ be  two open subsets of $\RR^2$ such that  $\ov S\subset V_1\hskip-0,05cm\subset \hskip-0,05cmV_2$, and 
$S_1$,  $S_2$    two nonempty  bounded   Lipschitz  domains of $\RR^2$ such that  $  \ov S_1\hskip-0,05cm \subset\hskip-0,05cm \ov S_2 \hskip-0,05cm\subset V $
  and $\bfy_{S_1}\hskip-0,05cm=\hskip-0,05cm\bfy_{S_2}\hskip-0,05cm=\hskip-0,05cm0$. Then
\begin{equation}
\begin{aligned}
&\capsca^f(\bfa,\bfalpha;S,V_1) \ge\capsca^f(\bfa,\bfalpha;S,V_2) & &  \forall (\bfa,\bfalpha)\in (\RR^3)^2,
\\& \capsca^f\lp\bfa ,\tfrac{\diam S_1}{\diam S_2} \bfalpha;S_1,V\rp
 \le
  \capsca^f(\bfa,\bfalpha;S_2,V)\quad & &\forall (\bfa,\bfalpha)\in (\RR^3)^2.
\end{aligned}
 \label{croissance}
  \end{equation}
  \noindent (ii) Let $f_1$ and $f_2$ be two convex functions on $\SSym^3$ such that $0\leq f_1(\bfM)\leq f_2(\bfM) \ \forall \bfM\in \SSym^3$. Then,
  for all $(\bfa,\bfalpha,\lambda)\in (\RR^3)^2\times (0,+\infty)$,
    \begin{equation}
   \capsca^{f_1}(\bfa,\bfalpha;S,\!V)\!
 \le\!
  \capsca^{f_2}(\bfa,\bfalpha;S,\!V);\quad  \capsca^{ \lambda f_1}(\bfa,\bfalpha;S,\!V)\!=\!  \lambda  \capsca^{f_1}(\bfa,\bfalpha;S,\!V).   
\label{cf1lecf2}
  \end{equation}
   \noindent (iii) We have 
    \begin{equation}
   \capsca^f(\bfa,\bfalpha;\bft+ S,\bft+ V) =   \capsca^f(\bfa,\bfalpha;S, V) \qquad \forall \bft\in \RR^2. 
\label{cftranslation}
  \end{equation}
\end{lemma}

\noindent
We  now address 
 the  asymptotic  behavior of  $\capsca^f(\bfa,\bfalpha;.,.)$  
  w.r.t.   monotonous sequences of sets. Recall  
      that,  for  any  convex function $h:\SSym^3\to \RR$ satisfying   \eqref{growthp}, the following holds   \cite[Prop. 2.32]{Dac}: 
\begin{equation}
  |h(\bfM)-h(\bfM')| \leq C |\bfM-\bfM'|(1+|\bfM|^{p-1}+|\bfM'|^{p-1}) \qquad \forall \bfM, \bfM' \in \SSym^3.
\label{inconv}
\end{equation}
 %
 \begin{lemma}\label{cont}   (i) Let   $(V_n)_{n\in \NN^{\star}}$ be  an nondecreasing  sequence of open subsets of $\RR^2$ and  $(S_n)_{n\in \NN^{\star}}$  a nonincreasing sequence of bounded   Lipschitz domains of $\RR^2 $ such that  
  \begin{equation}
 \begin{aligned}
&   \bigcup_{n\in \NN^{\star}}  \uparrow  V_n=V,\quad \bigcap_{n\in \NN^{\star}}  \downarrow \ov S_n=\ov S,
\quad \ov S_1  \subset  V_1,\quad  \bfy_{S_n}  =0 
 \quad \forall n\in \NN^{\star}.
\end{aligned}
\label{hypVnSn}
    \end{equation}
Then,
\begin{equation}
\lim_{n\to+\infty} \capsca^f\lp\bfa,\tfrac{\diam S_n}{\diam S}\bfalpha; S_n,V_n\rp=\capsca^f(\bfa,\bfalpha;S,V) \qquad \forall (\bfa,\bfalpha)\in (\RR^3)^2.
 \label{VnV}
  \end{equation}
  \noindent  Moreover, for all  $(\bfa,\bfalpha)\in  (L^p(\Omega ;\RR^3))^2$,  
  \begin{equation}
 \begin{aligned}
\hskip-0,2cm\lim_{ n\to +\infty}  \int_\Omega \hskip-0,1cm \capsca^f\hskip-0,1cm \lp  \bfa(x), \tfrac{\diam S_n}{\diam S}\bfalpha(x) ; S_n,V_n\rp dx
 =\hskip-0,1cm\int_\Omega \hskip-0,1cm\capsca^f (\bfa(x),\bfalpha(x);S,V) dx.
\end{aligned}
\label{estimintcap0}
    \end{equation}
  \noindent (ii)  Assume that  either  $p<2$ or $V$ is bounded in one direction  and  let 
 $(S'_n)_{n\in \NN^{\star}}$ be  an nondecreasing sequence of bounded  Lipschitz domains  such that
  \begin{equation}
 \begin{aligned}
\bfy_{S'_n}=0 \quad  \forall n\in \NN^{\star} \quad \hbox{ and } \quad  \bigcup_{n\in \NN^{\star}}\uparrow S'_n= S. 
\end{aligned}
\label{hypSnprim}
    \end{equation}
Then, 
 \begin{equation}\begin{aligned}
&
\hskip-0,1cm  \lim_{n\to+\infty} \capsca^f\lp\bfa,\tfrac{\diam S'_n}{\diam S}\bfalpha; S'_n,V \rp  = \capsca^f(\bfa,\bfalpha;S,V) \quad \forall (\bfa,\bfalpha)\in (\RR^3)^2 .
\end{aligned}
\label{lbS'}
  \end{equation}
 \end{lemma}

\noindent
\noindent{\bf Proof.} 
(i)  We  fix  $t>0$, $\e>0$,  set  $ S^t:= \{ x\in \RR^2, $ $  \dist ( x, S ) <t  \}$ and choose 
   $\eta_t\in \D(V;[0,1])$  and  $\bfvarphi\in \W^p(\bfa,\bfalpha;S,V)\cap \D(V;\RR^3)$  
    such that   
%
   \begin{equation}
   \begin{aligned}
 &\ov{S^{2t}}\subset V ,\qquad  \spt\eta_t \subset S^{2t}, \qquad \eta_t=1 \hbox{ in } \ S^t, \qquad |\bfnabla\eta_t|_{L^\infty(V)}\leq \frac{C}{t},
 \\&  \capsca^f(\bfa,\bfalpha;S,V) \leq   \int_V  f(  \bfe_y(\bfvarphi)) dy \leq \capsca^f(\bfa,\bfalpha;S,V)+\e.
    \end{aligned}
    \nonumber
    \end{equation}
By  \eqref{hypVnSn},  there exists   $n_0\in \NN^{\star}$  such that
 $S_n \subset S^t$  and  $\spt \bfvarphi\subset V_n$,  $ \forall n\ge n_0$.    We  set $ \bfvarphi_t:= \eta_t\lp \bfa+\tfrac{1}{\diam S}\bfalpha\wedge  \bfy  -\bfvarphi \rp +\bfvarphi$. 
Noticing  that    %
    \begin{equation}
\begin{aligned}
&  \bfvarphi_t\in  \W^p\lp\bfa,\tfrac{\diam S_{n}}{\diam S}\bfalpha;S_{n},V_{n}\rp \cap \D(V;\RR^3) \qquad \forall n\ge n_0, 
 \\&
  |\bfvarphi-\bfvarphi_t|\leq C t  \mathds{1}_{S^{2t}\setminus S}; \quad
 |\bfe_y(\bfvarphi)-\bfe_y(\bfvarphi_t)|\leq C \mathds{1}_{S^{2t}\setminus S}; \quad | \bfe_y(\bfvarphi_t)|\leq C, 
  \end{aligned} 
\nonumber  \end{equation}
we deduce from   \eqref{inconv}   that $ \lb  \int_V\hskip-0,05cm f(\bfe_y(\bfvarphi)) dy-\hskip-0,1cm \int_V\hskip-0,05cm f(\bfe_y(\bfvarphi_t)\!) dy\rb\hskip-0,1cm \leq$ $ C |S^{2t}\setminus S|\leq C t$. Taking \eqref{croissance} and  \eqref{hypVnSn}  into account,  we infer
  \begin{equation}
\begin{aligned}
 \capsca^f(\bfa,\bfalpha;S,V)&\leq \capsca^f \lp \bfa,\tfrac{\diam S_n}{\diam S}\bfalpha;S_n,V_n\rp  \leq 
 \int_V f(\bfe_y(\bfvarphi_t)) dy  
 \\& \leq \int_V f(\bfe_y(\bfvarphi)) dy + Ct 
  \leq  \capsca^f(\bfa,\bfalpha;S,V)+\e+ Ct \quad \forall n\ge n_0,  
  \end{aligned} 
\nonumber   \end{equation}
and deduce  \eqref{VnV}.
If    $(\bfa,\bfalpha)\in (L^p(\Omega;\RR^3))^2$, 
by \eqref{boundp}, \eqref{croissance}  and \eqref{hypVnSn} we have 
 %
  \begin{equation}
\begin{aligned}
 0\leq \capsca^f \lp \bfa(.),\tfrac{\diam S_n}{\diam S}\bfalpha(.);S_n,V_n\rp 
  \leq 
 \capsca^f \lp \bfa(.),\tfrac{\diam S_1}{\diam S}\bfalpha(.);S_1,V_1\rp\in L^1(\Omega), 
  \end{aligned} 
\nonumber   \end{equation}
%
thus   \eqref{estimintcap0}    results  from \eqref{VnV} and the dominated convergence theorem.  \qed

   \noindent (ii)  Assume   that $1<p<2$ and let $\bfvarphi_n$ be  a 
   solution to 
    $ \calP^f_K\lp \bfa,\tfrac{\diam S'_n}{\diam S}\bfalpha;S'_n,V\rp$. 
   By \eqref{croissance} and \eqref{hypSnprim}, we have
  \begin{equation}\begin{aligned}
&
    \int_V \hskip-0,1cmf(\bfe_y(\bfvarphi_n)) dy
=  \capsca^f\lp \bfa,\tfrac{\diam S'_n}{\diam S}\bfalpha;S'_n,V\rp \leq  \capsca^f(\bfa,\bfalpha;S,V),
\end{aligned}
\label{S'2}
  \end{equation}
thus, by  \eqref{growthp} and \eqref{minKorn},     $(\bfvarphi_n)_{n\in \NN^{\star}}$  is bounded in $K^p_0(V;\RR^3)$, hence   weakly converges in $K^p_0(V;\RR^3)$, 
and   a.e.   converges  in $V$, 
up to a subsequence, 
to  some $\bfvarphi$.  Since $\bfvarphi_n\in  \K^p\lp \bfa,\tfrac{\diam S'_n}{\diam S}\bfalpha;S'_n,V\rp$, we deduce from  \eqref{hypSnprim} that   $\bfvarphi =\bfa+\tfrac{1}{\diam S} \bfalpha\wedge \bfy\,$ a. e.  in $S$, thus  $\bfvarphi\in  \K^p(\bfa,\bfalpha;S,V)$ and, by   Lemma \ref{reg} (i),   \eqref{minK} and   \eqref{S'2}, 
  \begin{equation}\begin{aligned}
\capsca^f(\bfa,\bfalpha;S,V)&\le\int_V f(\bfe_y(\bfvarphi)) dy\leq \liminf_{n\to+\infty}  \int_V f(\bfe_y(\bfvarphi_n)) dy
\\& \le\limsup_{n\to+\infty} \capsca^f\lp \bfa,\tfrac{\diam S'_n}{\diam S}\bfalpha;S'_n,V\rp 
  \leq  \capsca^f(\bfa,\bfalpha;S,V).
  \end{aligned}
\nonumber
  \end{equation}
The proof of  \eqref{lbS'}   when  $V$ is  bounded in one direction is similar. \qed    

\noindent 
The next lemma is  crucial
 for the proof of the upper bound.
We  set: 
     \begin{equation}
  \begin{aligned} &\bfgotW^p_b(\bfa,\bfalpha;  S, V )  := 
\\&  \qquad \qquad
\la  \begin{aligned}  \bfeta \in  C^1(\ov\Omega;  \D(V;\RR^3)) \lb
 \begin{aligned}&  \bfeta(x, .)  \in \W^p(\bfa(x),\bfalpha(x);  S, V) \  \forall x  \in  \Omega,
 \\& \bfeta=0 \hbox{ on } (\Omega'\times  \{0\})\times V
 \end{aligned}  \right.
\end{aligned}\ra.
\end{aligned}
\label{defgotWaalpha}
    \end{equation}
 \begin{lemma}\label{lemapproxcap}  
For every $(\bfa,\bfalpha)\in (C^\infty(\ov\Omega;\RR^3))^2$,  
  \begin{equation}
 \begin{aligned}
  \int_\Omega \hskip-0,1cm \capsca^f (\bfa(x),\bfalpha(x)&;S, V) dx = 
\inf_{\bfeta\in \bfgotW^p_b(\bfa,\bfalpha;S,V)  }  \int_{\Omega\times V} \hskip-0,1cm f(\bfe_{y}(\bfeta(x,y))) dxdy. 
\end{aligned}
 \label{estimintcap}
    \end{equation}
   \end{lemma}

\noindent
 {\bf Proof.} \  The direct proof of this lemma is considerably shortened as follows:
we  set 
 $$
 \calH:=\la g(x)=\int_Vf(\bfe_y(\bfeta(x,y)))dy, \   \bfeta\in \bfgotW^p_b(\bfa,\bfalpha; S,V) \ra.  
 $$
We check below that 
\begin{equation}
 \begin{aligned}
 \forall (g_i)_{i\in I} \subset \calH \hbox{ ($I$ finite)}, \quad   \forall (\varphi_i)_{i\in I}\subset C^1(\ov\Omega;[0,1]), \quad \exists g\in \calH, \quad  g\le \sum_I \varphi_i g_i, 
\end{aligned}
 \label{propH}
    \end{equation}
therefore,  by \cite[Lemma  4.3]{BoDa} 
(a corollary  from  \cite[Th.  1]{BoVa}),  
$$
 \int_\Omega  \hbox{ess} \inf_{g\in \calH} g(x)  dx=\inf_{g\in \calH} \int_\Omega g(x) dx,
$$
which, by Lemma \ref{reg} (ii), is equivalent to   \eqref{estimintcap}.

\noindent{\it Verification of \eqref{propH}.}  For each $i\in I$, let $\bfeta_i\in \bfgotW^p_b(\bfa,\bfalpha; S,V)$  be  such that   $g_i(x)=\int_Vf(\bfe_y(\bfeta_i(x,y)))dy$. Then, 
by \eqref{defgotWaalpha}, $\sum_{i\in I} \varphi_i (x)\bfeta_i(x,y)) \in  \bfgotW^p_b(\bfa,\bfalpha; S,V) $,
hence   $g(x):=\int_Vf(\bfe_y( \sum_{i\in I} \varphi_i (x)\bfeta_i(x,y)))dy\in \calH$ and, 
by  the convexity of $f$,  $g\le \sum_I \varphi_i g_i$.
      \qed

\begin{remark}\label{remghom} If $0<k<+\infty$, by   the same  argument,  one can check  that
   \begin{equation}
\begin{aligned}
& 
\inf_{\bfq\in C^1_c(\Omega; C^\infty(\ov S;\RR^3))}\hskip-0,1cm k\int_\Omega  \hskip-0,1cm  \intb_{S} \hskip-0,1cmg\lp  \hskip-0,1cm \bfgotL\lp \bfq, \frac{\partial \psi_3}{\partial x_3}, \frac{\partial \zeta}{ \partial x_3}\rp \hskip-0,1cm \rp\hskip-0,1cm dx dy  
= \int_\Omega\hskip-0,1cm g^{hom}  \lp  \bfgotD(\bfpsi^{tuple})\rp dx, 
\end{aligned} 
\nonumber  \end{equation}
where $g^{hom}$ is given by \eqref{defghom1}.
 An analogous statement  holds in the case $1<\kappa<+\infty$.
The first equation in \eqref{lirem} and the existence of $\bfvarphi_0$ satisfying \eqref{intint} in the proof of Theorem \ref{thsoft} 
 can be justified in the same manner.
\end{remark}
  \noindent
 We  particularize below  the case of a  positively homogeneous function.
 %
\begin{lemma}\label{lemhomogeneous}      
If $f$ is positively homogeneous of degree $p\in (1,+\infty)$,  then for every $ (\bfa,\bfalpha)\in (\RR^3)^2$, the following holds: 
\begin{equation}\begin{aligned}
&
\capsca^f( \lambda\bfa,\lambda\bfalpha;S,V) =  \lambda^p \capsca^f\lp \bfa,\bfalpha; S, V\rp \qquad \hskip0,4cm&&\forall \lambda\ge0,
\\ &
\capsca^f( \bfa,\bfalpha;\lambda S,\lambda V) = \lambda^{2-p} \capsca^f\lp \bfa,\bfalpha; S, V\rp \quad &&\forall \lambda>0 .\end{aligned}
\label{rel3}
  \end{equation}
\end{lemma}

\noindent
\noindent {\bf Proof.} \  The first line of \eqref{rel3} is straightforward.  Fix  $(t,\lambda )\in (0,+\infty)^2$,   $(\bfa,\bfalpha)\in (\RR^3)^2$,    $\bfvarphi \in  \W^p(\bfa,\bfalpha;\lambda S, \lambda V) $        such that  
 $\capsca^f(\bfa,\bfalpha;\lambda S, \lambda V)  +t\ge   \int_{\lambda V } f(   \bfe_y(\bfvarphi ) ) dy$ and set   $\tilde\bfvarphi(y):=  \bfvarphi (\lambda  y )$. Then  $\tilde\bfvarphi\in  \W^p(\bfa,\bfalpha;  S,  V  )$,   $  \bfe_y(\tilde\bfvarphi)(y) = \lambda   \bfe_y(\bfvarphi)(\lambda y)$, and,  by  the 
   $p$-positive homogeneousness   of $f$ and   the change of variables formula,  
 \begin{equation}\begin{aligned}
\capsca^f(\bfa, \bfalpha;\lambda S, \lambda V)  +t &   \ge  \int_{\lambda V }f(   \bfe_y(\bfvarphi )  ) dy  
  = \lambda^2  \int_V f(   \bfe_y(\bfvarphi )(\lambda y)) dy
   \\ &=\lambda^{2-p} \int_V f(  \bfe_y(\tilde\bfvarphi )(y)) dy
\ge  \lambda^{2-p} \capsca^f\lp \bfa,\bfalpha;S,V\rp.
  \end{aligned}
\nonumber
   \end{equation}
  \qed

\noindent  The main purpose of  the  following  three propositions is to prove that $c^f$ is well defined by \eqref{defcfintro} and  satisfies 
 the  estimates \eqref{estimcfp<2} and \eqref{estimcfp=2} and  the formula \eqref{cfp<2}  synthetized in 
\eqref{cf=} (see Corollary \ref{corcf}).
To that aim, we investigate the asymptotic behavior as $r\to 0$ of the sequence   $\lp\capsca^f(\bfa,\bfalpha;r S,R_r D)\rp_{r>0}$  when   $(R_r)_{r>0}$   satisfies  
\begin{equation}
\begin{aligned}
&0< r\ll R_r \ll1,   \quad  R_r^2\ll r^{2-p},
 \quad 
  R_r\ll \frac{1}{\sqrt{|\log r|}} \  \hbox{ if } \ p=2  .
 \end{aligned}
\label{hypRr}
\end{equation}
%
 %
 \begin{proposition}\label{propcapreRep<2} 
 Under 
   \eqref{finftyfg0g} and  \eqref{hypRr}, the following holds 
for $p\in (1,+\infty)$: 
 \begin{equation}
 \begin{aligned} & \lim_{r\to0} \frac{ \capsca^f(\bfa,\bfalpha;r S,R_r  D)}{r^{2-p}}=  \capsca^{f^{\infty,p}}( \bfa,\bfalpha;S,\RR^2) \quad \forall (\bfa,\bfalpha)\in \RR^3\times\RR^3
.
\end{aligned}
 \label{capp<2}
  \end{equation}
\end{proposition}

 \noindent{\bf Proof.}         
Since  $f^{\infty,p}$ is  $p$-positively homogeneous, by   \eqref{VnV} and   \eqref{rel3} we have 
 \begin{equation}
 \begin{aligned}
  \lim_{r\to0} \frac{ \capsca^{f^{\infty,p}}(\bfa,\bfalpha;r S,R_r  D)}{r^{2-p}}&=
  \lim_{r\to0} \capsca^{f^{\infty,p}}( \bfa,\bfalpha;S,R_r/r D)
  \\&=
   \capsca^{f^{\infty,p}}( \bfa,\bfalpha;S,\RR^2) .
\end{aligned}
 \label{capp<2finfty}
  \end{equation}
By  \eqref{growthp} and  \eqref{deffinfty},  
 $f \le C  f^{\infty,p}$
for some  $C>0$,  hence
$   \capsca^f \le C \capsca^{f^{\infty,p}} $.
Denoting by  $\bfvarphi_r$     a solution to 
  $\calP_K^{f^{\infty,p}}  (\bfa,\bfalpha;r S,R_r D)$ if 
  $\capsca^f  (\bfa,\bfalpha;r S,R_r D)\!\ge \!\capsca^{f^{\infty,p}}  (\bfa,\bfalpha;r S,R_r D)$,  
  to   $\calP_K^{f}  (\bfa,\bfalpha;r S,R_r D)$ otherwise, we deduce 
  \begin{equation}
\begin{aligned}
    \int_{R_r D}   
  |\bfe (\bfvarphi_r)|^p dx & \le C  \capsca^{f^{\infty,p}}  (\bfa,\bfalpha;r S,R_r D),
  \end{aligned}
\nonumber
\end{equation}
and infer from  \eqref{finftyfg0g},  \eqref{hypRr}  and H\"older's inequality
  \begin{equation}
\begin{aligned}
&\lb(\capsca^f - \capsca^{f^{\infty,p}} ) (\bfa,\bfalpha;r S,R_r D)\rb
  \le   \int_{R_r D}  
\lb f(\bfe(\bfvarphi_r)) -  f^{\infty, p} (\bfe (\bfvarphi_r))\rb dx 
\\&  \leq C  \int_{R_r D}   
 1+  |\bfe (\bfvarphi_r)|^{p-\varsigma} dx  
\le  C R_r^2+C  \lp R_r^2\rp^{\frac{\varsigma}{p}}\lp  \capsca^{f^{\infty,p}}  (\bfa,\bfalpha;r S,R_r D) \rp^{\frac{p-\varsigma}{p}}.
\end{aligned}
\label{capfcapfinfty}
\end{equation}
Observing that, by  \eqref{croissance} and \eqref{rel3}, $\capsca^{f^{\infty,p}}  (\bfa,\bfalpha;r S,R_r D) \le r^{2-p}  \capsca^{f^{\infty,p}}  (\bfa,\bfalpha; S,V)$
as soon as   $(R_r/r) D\supset V$,   we deduce from \eqref{boundp} and \eqref{hypRr}   that 
  \begin{equation}
\begin{aligned}
&\lb(\capsca^f - \capsca^{f^{\infty,p}} ) (\bfa,\bfalpha;r S,R_r D)\rb
   \le  C R_r^2+C  \lp R_r^2\rp^{\frac{\varsigma}{p}}\lp r^{2-p}(|\bfa|^p+|\bfalpha|^p)\rp^{\frac{p-\varsigma}{p}}
=o(r^{2-p}),
\end{aligned}
\nonumber
\end{equation}
which, combined with \eqref{capp<2finfty}, yields \eqref{capp<2}. 
 \qed 
 
\noindent
When $p=2$, 
explicit computations 
   developed in   \cite[pp. 73-75]{BeGr} and  \cite[p. 147]{BeArma}
 in the setting of  linear isotropic elasticity,   show  that (see \eqref{fisot}) 
    \begin{equation}
\capsca^{f_{\lambda_0,\mu_0}}(\bfa,\zeta\bfe_3;rD,R_r D) =\begin{pmatrix} \bfa\\  \zeta \bfe_3\end{pmatrix}\cdot \bfCap^{f_{\lambda_0,\mu_0}}(rD,R_r D)\begin{pmatrix} \bfa\\  \zeta \bfe_3\end{pmatrix},
\label{capquad}\end{equation}
  where  $ \bfCap^{f_{\lambda_0,\mu_0}}(r D, R_r D)\in \SSym^4$   is  diagonal   and satisfies, as $r\to0$,
 \begin{equation}
 \hskip-0,2cm  \begin{aligned}
 &  \Capa^{f_{\lambda_0,\mu_0}}(r D, R_r D)_{\alpha\alpha}= 4\pi\mu_0  \frac{\lambda_0+2\mu_0}{\lambda_0+3\mu_0}\frac{1}{|\log r|}(1+o(1)) \quad \forall  \alpha\in \{1,2\},
 \\&  \Capa^{f_{\lambda_0,\mu_0}}(r D, R_r D)_{33}= 
  \frac{2\pi\mu_0}{\log R_r-\log r},
  \\& \Capa^{f_{\lambda_0,\mu_0}}(r D, R_r D)_{44}= 
 \frac{ 4\pi\mu_0 R_r^2}{R_r^2-r^2}.
 \end{aligned}
\label{dddii}
  \end{equation}
Hence, by  \eqref{croissance} and \eqref{cf1lecf2},  if $f$ is an arbitrary  convex function 
  satisfying  \eqref{growthp},
 \begin{equation}
\label{bdcfSrR}
\begin{aligned}
& c|\bfa|^2  \leq  |\log r|  \capsca^{f}(\bfa, 0;r S,R_r D )  \leq  C|\bfa|^2 \qquad & &\forall \bfa\in \RR^3,
\\& c|\zeta|^2  \leq   \capsca^{f}(0,\zeta\bfe_3;r S,R_r D ) \leq  C|\zeta|^2 \qquad & & \forall \zeta\in \RR,
\end{aligned}
  \end{equation}
 for some  $c,C\in(0,+\infty)$. In what follows, we denote by $\T$  the Fr\'echet  topology on $C(\RR^3)$  of the uniform 
   convergence on the compact subsets of $\RR^3$.  
  We set
 \begin{equation}
\label{defcfSrR}
\begin{aligned}
&  c^{f,S}_{r, R_r}(\bfa) :=  |\log r |  \capsca^{f}(\bfa,0;r S,R_r D ) & &\forall \bfa\in \RR^3.
\end{aligned}
  \end{equation}

\begin{proposition}\label{propcapreRep=2} 
Assume that  $p=2$. There exists a convex function  $c^f_0$
  independent of    $S$,    verifing   
 \begin{equation}
\label{growthpcf0}
\begin{aligned}
&c_0^f=c_0^{f^{\infty,2}}, \qquad  c|\bfa|^2  \leq c_0^{f}(\bfa)  \leq  C|\bfa|^2 \qquad & &\forall \bfa\in \RR^3,
\end{aligned}
  \end{equation}
and such that, for every sequence   $(R_r)_{r>0}$   satisfying  \eqref{hypRr}, 
\begin{equation}\begin{aligned}
\hbox{ $c^{f,S}_{r,R_r}$ $\T$-converges to $c^f_0$  as $r\to0$.}
  \end{aligned}
  \label{Tcv0}
  \end{equation}
\end{proposition}

\noindent
 The demonstration   of this   result, postponed to Section \ref{secprooflemA}, 
is   delicate:    the convergence 
 \eqref{Tcv0} 
 is   obtained 
 by capitalizing on its consequence,  Theorem \ref{th}. 
\begin{remark} 
  \label{remMf}   
 By   \eqref{growthp},  \eqref{VnV},  \eqref{rel3},   and    \eqref{dddii},   the following holds
 \begin{equation}\begin{aligned}
   &  \capsca^{f_{\lambda_0,\mu_0}}(\bfa,\zeta\bfe_3; D,\RR^2) =4\pi\mu_0 \zeta^2 & &  \forall (\bfa,\zeta) \in \RR^3\times \RR,  
\\& C_1|\zeta|^2\leq \capsca^f(\bfa, \zeta \bfe_3; S,\RR^2) =  \capsca^f(0,\zeta\bfe_3; S,\RR^2)\leq C_2\zeta^2 & &  \forall (\bfa,\zeta) \in \RR^3\times \RR .
  \end{aligned}
\label{explicit3}
  \end{equation} 
\end{remark}

\noindent We turn to   the case  $p>2$. 

\begin{proposition}\label{propcapreRep>2} 
If $2<p<+\infty$,  
 \begin{equation}
 \capsca^f(\bfa,\zeta\bfe_3;r S,R_r D)\ge \frac{C}{R_r^{p-2}}\lp |\bfa|^p+|\zeta|^p
  \rp \qquad \forall (\bfa, \zeta)\in \RR^3\times  \RR .
 \label{capp>2}
\end{equation}
\end{proposition}

\noindent
 \noindent{\bf Proof.}      Let $\bfvarphi$  be  
a solution to   $\calP^f\hskip-0,05cm(\bfa,\hskip-0,05cm\tfrac{\diam  D}{\diam S}\zeta\bfe_3;  \hskip-0,05cmr D\hskip-0,05cm,\hskip-0,05cm R_r D)$.
By  \eqref{DsubsetS},    \eqref{growthp},    \eqref{LpfirstKorn}, and  \eqref{croissance}, we have
 \begin{equation}
\begin{aligned}
 &  \capsca^f\lp \bfa,\zeta\bfe_3  ;r S,R_r D\rp \ge  \capsca^f\lp \bfa, \tfrac{\diam D}{\diam S}\zeta \bfe_3 ;r D,R_r D\rp
 \\&\qquad\qquad= 
   \int_{R_r D} f(\bfe_y(\bfvarphi)) dy
 \ge C   \int_{R_r D} |\bfe_y(\bfvarphi)|^p dy
 \ge C    \int_{R_r D} |\bfnabla \bfvarphi|^p dy.
 \end{aligned}
  \label{partcr4}
  \end{equation}  
  \noindent  Applying  to each component of $\bfvarphi$  the  estimate  (see  \cite[Lemma  A3]{BeBo})
     \begin{equation}
\begin{aligned}
    \int_{R_r D\setminus r D} |\bfnabla \eta|^p dy & \ge    \frac{C}{R_r^{p-2}} \lb\intb_{\partial R_r D} \eta d\calH^1- \intb_{\partial r D}\eta d\calH^1\rb^p
\quad \forall \eta\in W^{1,p}(R_r D),
 \end{aligned}
 \nonumber
  \end{equation}  
noticing that,  by  \eqref{infW},
$ {\int\!\!\!\!\!-}_{\partial R_r D} \bfvarphi d\calH^1=0$,  $ {\int\!\!\!\!\!-}_{\partial r D} \bfvarphi d\calH^1=\bfa$,
we  deduce 
   \begin{equation}
\begin{aligned}
  \capsca^f\lp \bfa,\zeta\bfe_3  ;r S,R_r D\rp  \ge \frac{C}{R_r^{p-2}} |\bfa|^p.
 \end{aligned}
  \label{partcr51}
  \end{equation}  
  \noindent 
Since $p>2$, by H\"older's inequality we have     $\bfvarphi\in\W^2(\bfa,\tfrac{\diam  D}{\diam S}\zeta\bfe_3;  r D, R_r D)$,
hence, by \eqref{growthp},  \eqref{croissance},   \eqref{bdcfSrR}, \eqref{partcr4}, and   the quadratic nature of the mapping $(\bfa,\zeta)\to\capsca^{|.|^2}(\bfa,\zeta; $ $ r D, R_r D)$ (see \eqref{capquad}),  the following holds:
 \begin{equation}
 \begin{aligned}
 \zeta^2  \hskip-0,1cm &\leq \hskip-0,1cmC  \capsca^{|.|^2}\hskip-0,1cm\lp0,\tfrac{\diam D}{\diam S}\zeta\bfe_3; r D, R_r D\rp
\hskip-0,1cm \leq  C \capsca^{|.|^2} \hskip-0,1cm\lp  \bfa, \tfrac{\diam D}{\diam S}\zeta\bfe_3 ; r D,  R_r D \rp \hskip-0,1cm +\hskip-0,1cm \frac{C|\bfa|^2}{|\log r|} 
    \\& \leq  C  \int_{R_r D} |\bfe_y(\bfvarphi)|^2 dy   +  \frac{C}{|\log r|} \leq C R_r^{\frac{2(p-2)}{p}}  \lp \int_{R_r D} |\bfe_y(\bfvarphi)|^p dy \rp^{\frac{2}{p}} + \frac{C}{|\log r|} 
    \\& \leq C  R_r^{\frac{2(p-2)}{p}}\hskip-0,1cm  \lp  \capsca^f\lp \bfa, \frac{\diam D}{\diam S}\zeta \bfe_3 ;r D,R_r D\rp\rp^{\frac{2}{p}}+\frac{C}{|\log r|} 
\\& \leq C  R_r^{\frac{2(p-2)}{p}}\hskip-0,1cm  \lp  \capsca^f\lp \bfa, \zeta \bfe_3 ;r S,R_r D\rp\rp^{\frac{2}{p}}+\frac{C}{|\log r|} . 
  \end{aligned}
\nonumber 
  \end{equation}  
We infer that  $ \capsca^f\hskip-0,1cm\lp \bfa,\zeta\bfe_3  ;r S,R_r D\rp  \hskip-0,1cm \ge\hskip-0,1cm  \frac{C}{ R_r^{p-2}} \lb \zeta^2- \frac{C}{|\log r|}\rb^{\tfrac{p}{2}}$ which, combined with 
  \eqref{partcr51}, yields   \eqref{capp>2}.  \qed

\begin{remark}\label{rem0} 
One can check that, if   $p\not=2$ and  $0<R_1<R_2$,  
  \begin{equation}
 \begin{aligned}
  \inf_{\eta\in H^1(0,R_2)}\hskip-0,1cm\la \int_{R_1}^{R_2} |\eta'(\rho)|^p \rho d\rho,\  \ \eta\hskip-0,1cm=\hskip-0,1cm1  \hbox{ in } (0,R_1), \ \eta(R_2)\hskip-0,1cm=\hskip-0,1cm0 \ra= \hskip-0,1cm\lp \frac{s}{R_2^s-R_1^s}\rp^{p-1}\hskip-0,1cm\hskip-0,1cm, 
  \end{aligned}
\nonumber 
  \end{equation} 
  where $s:=\tfrac{p-2}{p-1}$.
The field  defined in polar coordinates by 
  $\bfvarphi(\rho,\theta):= \eta(\rho) \bfa$, where $ \eta$  is     the solution to the above problem,
belongs to   $\W^p(\bfa, 0;R_1 D, $ $R_2 D)$, therefore 
$
 \capsca^{|.|^p}(\bfa,0;R_1 D,R_2 D) \leq  C \int_{R_2 D} |\bfnabla\bfvarphi|^p dx \leq  \frac{C |\bfa|^p }{(R_2^s-R_1^s)^{p-1}}.
$
Setting $R_1=1$, letting  $R_2$ tend  to $+\infty$,   we  deduce from \eqref{growthp},  \eqref{cf1lecf2} and  \eqref{VnV}   that $\capsca^{|.|^p}(\bfa,0;  D,\RR^2)=\capsca^f(\bfa,0;  S,\RR^2)=0$.
\end{remark}

\noindent The following corollary   results from 
   \eqref{boundp},   \eqref{boundp2} and 
   Proposition \ref{propcapreRep<2}    if  $p<2$,   from \eqref{bdcfSrR} and  Proposition \ref{propcapreRep=2} if $p=2$, and  from Proposition \ref{propcapreRep>2} if $p>2$. 
%
    \begin{corollary} \label{corcf}  
For every sequence $(R_\e)$ verifying  \eqref{Re}, the  convergence  \eqref{defcfintro} holds, where 
 $c^f$ is  convex,  
 independent of  $(R_\e)$, satisfies \eqref{cf=}.  If $p=2$, $c^f$  is independent of $S$  and 
  \begin{equation}
 \begin{aligned}
 c^f(\bfa,0) = c^f_0(\bfa), \quad c^f(\bfa,\theta)=+\infty \quad \hbox{ if } \theta\not=0,
  \end{aligned}
\label{cf=p=2}
  \end{equation} 
where $c^f_0$ is given by \eqref{Tcv0}.
   \end{corollary}

\begin{remark} 
  \label{remMf2}  By 
   \eqref{capquad},   \eqref{dddii}, \eqref{defcfSrR} and \eqref{cf=p=2},  Formula \eqref{explicit2} holds.
\end{remark}

\section{Proofs of the main results}\label{secproof}

\subsection{Proof of Theorem \ref{th}}

Following the $\Gamma$-convergence method \cite{Da}, we establish  
a lower bound  in  Section   \ref{seclowerbound}, an upper bound
in  Section  \ref{secupperbound},
 and conclude 
 the proof of  Theorem \ref{th} 
  in Section \ref{secproofth}.
 %
\subsubsection{Lower bound}\label{seclowerbound}
  
$ $

    \begin{proposition}
     \label{proplowerbound}
  Let  $(\bfu_\e)$ be a sequence in $W^{1,p}_b(\Omega;\RR^3)$ verifying (\ref{supFeuefini}) and 
  $(\bfu_{\e_k})_{k\in \NN}$   a subsequence satisfying  the convergences  
stated in   Propositions \ref{prop1}, \ref{prop2}  and  \ref{prop3}. 
Then  
\begin{equation}\begin{aligned}
&
\liminf_{k\to +\infty} F_{\e_k}(\bfu_{\e_k}) \ge \Phi(\bfu,\bfv^{tuple}).
 \end{aligned}\label{linfF}  \end{equation}
   \end{proposition}

\noindent 
{\bf Proof.} \quad To simplify, the subsequence will still be denoted by $(\bfu_\e)$.
By Lemma \ref{lemslicing}, there exists  a sequence $(R_\e)$
verifying \eqref{Re} and   
\begin{equation}
\begin{aligned}
\limsup_{\e\to 0} \int_{D_{  R_\e}\setminus D_{\frac{1}{2}R_\e} \times (0,L)} |\bfnabla   \bfu_\e |^p dx =0.
\end{aligned} 
\label{slicing} \end{equation} 
We set
\begin{equation}
\begin{aligned}
&  F_\e( \bfu_\e) = I_{\e1} + I_{\e2} +I_{\e3} ; 
\qquad I_{\e1} := \int_\Omega f\lp\bfe(\bfu_\e)\mathds{1}_{\Omega\setminus  (D_{ R_\e}\times(0,L))}\rp dx ;
\\& I_{\e2} :=  \int_{ ( D_{ R_\e} \times(0,L))\setminus {T_{r_\e}} }  f(\bfe(\bfu_\e)) dx;
\qquad I_{\e3}:=  k_\e  \int  g(\bfe(\bfu_\e))  dm_\e.
  \end{aligned} 
\label{splitlowerbound}    \end{equation} 
Since   $|D_{R_\e}|\ll 1$,  by  \eqref{cvbasic}  the sequence    $\lp\bfe(\bfu_\e)  \mathds{1}_{\Omega\setminus  (D_{ R_\e}\times(0,L))}\rp$  weakly converges  to  $\bfe(\bfu)$  in $L^p(\Omega;\SSym^3)$. Therefore, by 
  the weak lower-semicontinuity  on $L^p(\Omega;\SSym^3)$ of  the functional 
  $\bfPsi \to \int_\Omega f(\bfPsi) dx $,    resulting from \eqref{growthp} and    the convexity of $f$,  we have  
  \begin{equation}
\begin{aligned}
 \liminf_{\e\to0}  I_{\e1} \ge \int_\Omega f(\bfe(\bfu )) dx.
 \end{aligned} 
 \label{linf1}    \end{equation} 
 Let us check that 
  \begin{equation}
\begin{aligned}
  \liminf_{\e\to0}  I_{\e3}\ge   \int_\Omega g^{hom}(\bfgotD\,\bfv^{tuple})dx.  
 \end{aligned} 
 \label{linf3}    \end{equation} 
If $(k,\kappa)\in \{(+\infty,0), (+\infty,+\infty)\}$, then by \eqref{defghom1}, \eqref{defghom2},  \eqref{defD} and \eqref{vtupleinD}, we have 
$\bfgotD\,\bfv^{tuple}=0$, thus  there is nothing to prove.
  If $0<k<+\infty$, by   
   (\ref{cvbasice})  we have   $\bfe (\bfu_\e)\mdto $ $ \bfgotL\lp \bfq, \frac{\partial v_3}{\partial x_3}, \frac{\partial \theta}{\partial x_3}
\rp$  for some   $\bfq\in L^p(\Omega; W^{1,p}(S;\RR^3))$. Applying  (\ref{linfcarre}), taking   \eqref{kkappa} and   \eqref{defghom1} into account, we obtain 
\begin{equation}
\begin{aligned}
&   \liminf_{\e\to0}   I_{\e3}\ge 
  \frac{ k}{|S|} \int_{\Omega\times  S}  
g\lp   \bfgotL\lp \bfq, \frac{\partial v_3}{\partial x_3}, \frac{\partial \theta}{\partial x_3}
\rp   \rp  dx dy 
   \ge 
    \int_\Omega g^{hom}(\bfgotD\,\bfv^{tuple})dx, 
  \end{aligned} 
\nonumber   \end{equation} 
yielding \eqref{linf3}.
If  $0<\kappa<+\infty$,  by     \eqref{kkappa},  
    \eqref{finftyfg0g} and \eqref{estimbasic1}, 
  %
  \begin{equation}
\begin{aligned}
\hskip-0,1cm\lb I_{\e3}\hskip-0,1cm-r_\e^p  k_\e \int \hskip-0,1cm g^{0,p}\lp\frac{\bfe(\bfu_\e)}{r_\e}\rp dm_\e\rb&
\leq  r_\e^p   k_\e \int   \lb\frac{g\lp r_\e\frac{\bfe(\bfu_\e)}{r_\e}\rp}{r_\e^p}  -   g^{0,p}\lp\frac{\bfe(\bfu_\e)}{r_\e}\rp\rb  dm_\e
\\&\leq C \varpi(r_\e) \int  \lb\frac{\bfe(\bfu_\e)}{r_\e}\rb^p dm_\e\leq C  \varpi(r_\e),
 \end{aligned} 
\label{I31}
  \end{equation} 
and by 
   \eqref{defghom2}, \eqref{linfcarre} and \eqref{cvkappa},
 \begin{equation}
\begin{aligned}
    \liminf_{\e\to0} r_\e^p  k_\e \int \hskip-0,1cm g^{0,p}\lp\frac{\bfe(\bfu_\e)}{r_\e}\rp dm_\e
    &\ge     \frac{\kappa}{ |S|}  \hskip-0,05cm \int_{\Omega\times S}   \hskip-0,6cm
g^{0,p}\lp\bfgotJ\lp  \hskip-0,05cm\bfl,  \hskip-0,05cm\frac{\partial^2 v_1  }{ \partial x_3^2},  \hskip-0,05cm\frac{\partial^2 v_2  }{ \partial x_3^2},
 \hskip-0,05cm\frac{\partial w}{ \partial x_3},  \hskip-0,05cm\frac{\partial \delta }{ \partial x_3}\rp \hskip-0,1cm \rp  \hskip-0,1cmdxdy \hskip-0,05cm 
 \\&\ge \int_\Omega g^{hom}(\bfgotD\,\bfv^{tuple})dx,
  \end{aligned} 
\nonumber  \end{equation} 
for some $\bfl\in L^p(\Omega; W^{1,p}(S;\RR^3))$.
Taking \eqref{condvarpi}  into account,   \eqref{linf3} is proved.
By \eqref{defPhi}, \eqref{vtupleinD}, \eqref{splitlowerbound}, \eqref{linf1} and \eqref{linf3},    the demonstration of  Proposition \ref{proplowerbound}  is achieved 
provided we show that
 \begin{equation}
\begin{aligned}
  \liminf_{\e\to0}  I_{\e2}  \ge \int_\Omega c^f(\bfv-\bfu, \theta) dx   .
 \end{aligned} 
 \label{linf2}    \end{equation} 
    If $\gamma^{(p)}=0$,   by  \eqref{cf=}  we have $\int_\Omega c^f(\bfv-\bfu, \theta) dx =0$, hence there is nothing to prove.
Likewise, if    $\gamma^{(p)}=+\infty$,  
by \eqref{regthetav}  there holds   $\bfu=\bfv$ and $\theta=0$,  hence, by  \eqref{cf=},  $\int_\Omega c^f(\bfv-\bfu, \theta) dx =0$. 
From now on, we assume  that $0< \gamma^{(p)}<+\infty$,  thus  $p\in(1,2]$.
We fix     $c_0\in (0,1)$ and  a Lipschitz  domain $S'$  such that 
\begin{equation}
\begin{aligned}
&c_0S\subset S'\subset\subset S, \qquad \bfy_{S'}=0.
\end{aligned}
\label{S'}
\end{equation}  
\noindent By  Lemmas \ref{lemtriangle} and \ref{lemparen}
there exists 
    sequences
 $(\widetriangle\bfu_\e)$, $(\wideparen\bfu_\e)$ in $L^p(0,  L;W^{1,p}(\Omega' ;\RR^3) )$ 
satisfying
\begin{equation}
\begin{aligned}
&   \ov\bfgotu_\e(\widetriangle\bfu_\e) \rightharpoonup \bfu, \quad \ov\bfgotv_\e^S(\widetriangle\bfu_\e) \rightharpoonup\bfv,  \quad  \ov\theta_\e^S(\widetriangle \bfu_\e)\bfe_3 \rightharpoonup \theta\bfe_3 \quad  \hbox{weakly in } L^p(\Omega;\RR^3),\hskip-0,1cm 
 \end{aligned} 
\label{cvtriangle}    \end{equation}
\begin{equation}
\lpt\begin{aligned}
&  \liminf_{\e\to0}  I_{\e2} \ge  \liminf_{\e\to0}    \int_{(D_{R_\e}\setminus S'_{r_\e})  \times (0,L)}f^{\infty, p}(\bfe_{x'}(\wideparen  \bfu_\e)) dx,
\\& \wideparen  \bfu_\e=\ov\bfgotu_\e(\widetriangle\bfu_\e) \    \text{ on } \partial D_{R_\e}\times (0,L);
\quad 
   \wideparen \bfu_\e = \ov\bfgotr_\e^S(\widetriangle\bfu_\e) \ \hbox{ on } \ \partial S'_{r_\e}\times (0,L)  
 \end{aligned}\ra \quad\hbox{if } p<2,
\label{basedonI2ep<2}    \end{equation}
\begin{equation}
\lpt\begin{aligned}
&  \liminf_{\e\to0}  I_{\e2} \ge  \liminf_{\e\to0}    \int_{(D_{R_\e}\setminus S'_{r_\e})  \times (0,L)}\hskip-1cm f^{\infty, 2}(\bfe_{x'}(\wideparen  \bfu_\e)) dx -C\frac{|S\setminus S'|}{|S|},
\\& \wideparen  \bfu_\e=\ov\bfgotu_\e(\widetriangle\bfu_\e) \    \text{ on } \partial D_{R_\e}\times (0,L); \quad
   \wideparen \bfu_\e = \ov\bfgotv_\e^S(\widetriangle\bfu_\e) \ \hbox{ on } \ \partial S'_{r_\e}\times (0,L)  \end{aligned}\ra \ \hbox{if } p=2.
\label{basedonI2ep=2}    \end{equation}
  We set
%
\begin{equation}
\begin{aligned}
&
\bfa_\e :=\ov\bfgotv_\e^S(\widetriangle\bfu_\e) - \ov\bfgotu_\e(\widetriangle\bfu_\e), \quad 
 \bfalpha_\e :=\begin{cases} \tfrac{\diam S'}{\diam S} \ov\theta_\e^S(\widetriangle\bfu_\e) \bfe_3 & \hbox{ if } p<2,
 \quad
 \\   0 & \hbox{ if } p=2,
 \end{cases}
 \end{aligned} 
\label{defaalpha}
    \end{equation}
The following weak convergences in $L^p(\Omega)$ hold:
\begin{equation}
\begin{aligned}
&
\bfa_\e \rightharpoonup \bfv-\bfu; \quad 
 \bfalpha_\e  \rightharpoonup    \tfrac{\diam S'}{\diam S} \theta \bfe_3  \hbox{ if } p<2.
 \end{aligned} 
\label{cvaalpha}
    \end{equation}
We  respectively  denote  by   
$\bfa_\e^{i}(x_3)$,  $\bfalpha_\e^{i}(x_3)$,  and   $\ov\bfgotu_\e(\widetriangle\bfu_\e)^{i}(x_3)$  the constant values   taken by   $\bfa_\e$,  $\bfalpha_\e$,   and  $\ov\bfgotu_\e(\widetriangle\bfu_\e)$     in  $Y_\e^{i}\times\{x_3\}$.
 For  each  $i\in I_\e$ and a. e. $x_3\in (0,L)$,   $ \wideparen \bfu_\e(.,x_3)$ belongs to $W^{1,p}(D_{R_\e}^{i};$ $\RR^3)$, 
   $ \bfe_{x'}\lp \wideparen \bfu_\e(.,x_3)\rp = \bfe_{x'}\lp \wideparen \bfu_\e-\ov\bfgotu_\e(\widetriangle\bfu_\e)^{i}\rp(.,x_3)$  in $D_{R_\e}^{i}$
  and, 
by \eqref{infW},  \eqref{defovgotA},  \eqref{basedonI2ep<2},   \eqref{basedonI2ep=2},  \eqref{defaalpha},
\begin{equation}
\begin{aligned}
& \wideparen  \bfu_\e(.,x_3)- \ov\bfgotu_\e(\widetriangle\bfu_\e)^{i}(x_3) \in \W^p\lp \bfa_\e^{i}(x_3), \bfalpha_\e^{i}(x_3) ;(S')_{r_\e}^{i}  , D_{R_\e}^{i}\rp,
 \\&   \int_{D_{R_\e}^{i} \setminus (S')_{r_\e}^{i}}\hskip-0,6cm f^{\infty,p}( \bfe_{x'}(\wideparen \bfu_\e) (x' ,x_3))dx' 
\ge  
 \capsca^{f^{\infty,p}}\Big( \bfa_\e^{i}(x_3), \bfalpha_\e^{i}(x_3) ;(S')_{r_\e}^{i}  , D_{R_\e}^{i} \Big).
      \end{aligned} \hskip-0,2cm 
 \label{listin} 
   \end{equation}
$\bullet${ \it Case   $1<p<2$ and   $0<\gamma^{(p)}<+\infty$.}   
By      \eqref{croissance},  \eqref{cftranslation},    \eqref{rel3},   \eqref{S'} and \eqref{listin}, 
\begin{equation}
\begin{aligned}  
   \int_{D_{R_\e}^{i} \setminus (S')_{r_\e}^{i}}\hskip-0,6cm f^{\infty,p}( \bfe_{x'}(\wideparen \bfu_\e) (x' ,x_3))dx' 
 & = \frac{r_\e^{2-p}}{\e^2}\int_{Y_\e^{i}}  \capsca^{f^{\infty,p}}\Big( \bfa_\e^{i}(x_3),\bfalpha_\e^{i}(x_3) ;S'  , \frac{R_\e}{r_\e}D\Big)dx'
 \\&\ge \gamma_\e^{(p)}(r_\e) \hskip-0,1cm\int_{Y_\e^{i}}  \hskip-0,1cm \capsca^{f^{\infty,p}}\Big( \bfa_\e^{i}(x_3),\bfalpha_\e^{i}(x_3) ;S'  , \RR^2\Big)dx'.
     \end{aligned} \hskip-0,2cm 
 \nonumber
   \end{equation}
 After summation over $i\in I_\e$ and integration w.r.t. $x_3$, we obtain 
\begin{equation}
\begin{aligned}
     &  \hskip-0,1cm \int_{(D_{R_\e}\setminus S'_{r_\e}) \times (0,L)} \hskip-0,5cm f^{\infty,p}( \bfe_{x'}(\wideparen \bfu_\e)) dx 
     \ge\gamma_\e^{(p)}(r_\e)  \int_\Omega  \capsca^{f^{\infty,p}}\lp\bfa_\e ,\bfalpha_\e;S', \RR^2\rp dx, 
        \end{aligned} 
\label{hatfcap0}    \end{equation}
and deduce from 
  \eqref{basedonI2ep<2},  \eqref{cvaalpha} and    Lemma \ref{lemconvex} that
  \begin{equation}
\begin{aligned}
& \liminf_{\e\to0} I_{\e2} 
\ge    \gamma^{(p)}   \int_\Omega \capsca^{f^{\infty,p}}\lp\bfv-\bfu,\tfrac{\diam  S'}{\diam S} \theta \bfe_3;S',\RR^2\rp dx.
\end{aligned} 
\nonumber
    \end{equation} 
Substituting $S'_n$ for $S'$,
 where 
  $(S'_n)_{n\in \NN}$ is an nondecreasing sequence of Lipschitz domains such that 
   $\bigcup_{n\in \NN} \uparrow S'_n=S$, $S'_n\subset\subset S$, $\bfy_{S'_n}=0$, we infer    
   \begin{equation}
\begin{aligned}
& \liminf_{\e\to 0}I_{\e2}
 \ge    \limsup_{n\to+\infty}    \int_\Omega \hskip-0,1cm c_n dx;\quad 
   c_n\hskip-0,1cm := \gamma^{(p)} \capsca^{f^{\infty,p}}\hskip-0,1cm\lp \bfv-\bfu ,  \tfrac{\diam  S'_n}{\diam S} \theta \bfe_3; S'_n,\RR^2\hskip-0,1cm\rp 
.
\end{aligned} 
\nonumber
   \end{equation} 
By  \eqref{cf=},  \eqref{croissance},   \eqref{lbS'},    the  sequence   $(c_n)_{n\in \NN}$     is positive,   nondecreasing, and   pointwise converges to $\gamma^{(p)}\capsca^{f^{\infty,p}}\lp\bfv-\bfu ,   \theta \bfe_3; S ,\RR^2\rp=c^f(\bfv-\bfu,\theta)$. 
Applying the     
monotone convergence theorem,    we obtain  \eqref{linf2}. 
\\
$\bullet${ \it Case    $p=2$, $0<\gamma^{(2)}<+\infty$}
By    \eqref{defgammae},   \eqref{defcfSrR}  and \eqref{listin}, we have
\begin{equation}
  \begin{aligned}   \int_{D_{R_\e}^{i} \setminus (S')_{r_\e}^{i}}   f^{\infty,2}( \bfe_{x'}(\wideparen \bfu_\e)(x',x_3)) dx'
& \ge  
    \gamma_\e^{(2)}(r_\e)  \int_{Y_\e^{i}}  c^{f^{\infty,2},S'}_{r_\e,R_\e}(\bfa_\e^{i}(x_3) ) dx'.
      \end{aligned}
 \nonumber   \end{equation}
Summing w.r.t. $i$ over $I_\e$ and integrating w.r.t. $x_3$ over $(0,L)$, 
 taking   \eqref{S'}  and   \eqref{propcf} into account, we infer
\begin{equation}
\begin{aligned}
     &  \hskip-0,2cm \int_{(D_{R_\e}\setminus S'_{r_\e}) \times (0,L)} f^{\infty,2}( \bfe_{x'}(\wideparen \bfu_\e)) dx 
      \ge   \gamma_\e^{(2)}(r_\e)    \int_\Omega      c^{f^{\infty,2},S'}_{r_\e,R_\e}(  \bfa_\e) dx
      \\&\hskip1cm \ge   \gamma_\e^{(2)}(r_\e)  \hskip-0,1cm  \int_\Omega    \hskip-0,1cm c^{f^{\infty,2},c_0S}_{r_\e,R_\e}( \bfa_\e )dx=   \gamma_\e^{(2)}(r_\e)    \hskip-0,1cm \int_\Omega           \hskip-0,1cm c^{f^{\infty,2},S}_{r_\e,R_\e/c_0}(\bfa_\e) dx.
     \end{aligned} 
\label{hatfcap1}   \end{equation}
After possibly extracting a subsequence,  by \eqref{defgammae} we can assume that 
\begin{equation}
\begin{aligned}
     &  \liminf_{\e\to0}   \gamma_\e^{(2)}(r_\e)    \hskip-0,1cm \int_\Omega           \hskip-0,1cm c^{f^{\infty,2},S}_{r_\e,R_\e/c_0}(\bfa_\e) dx
     =\lime \gamma_\e^{(2)}(r_\e)    \hskip-0,1cm \int_\Omega           \hskip-0,1cm c^{f^{\infty,2},S}_{r_\e,R_\e/c_0}(\bfa_\e) dx.
     \end{aligned} 
\label{liminf=lime}   \end{equation}
%
%
For each  $N\in \NN$, we set 
\begin{equation}
\begin{aligned}
     & \bfa_\e^N:= \bfa_\e \quad \hbox{ if }  |\bfa_\e|\leq N; \quad  \bfa_\e^N:=  N \frac{\bfa_\e}{|\bfa_\e|} \quad \hbox{ if }  |\bfa_\e|>N.
     \end{aligned} 
\label{defaN}   \end{equation}
By  \eqref{rel3} and \eqref{defcfSrR},   $c^{f^{\infty,2},S}_{r_\e,R_\e/c_0}(t\bfxi)=t^2 c^{f^{\infty,2},S}_{r_\e,R_\e/c_0}(\bfxi)$,   hence
\begin{equation}
\begin{aligned}
     &  \int_\Omega           c^{f^{\infty,2},S}_{r_\e,R_\e/c_0}(\bfa_\e) dx
     \ge     \int_\Omega           c^{f^{\infty,2},S}_{r_\e,R_\e/c_0}(\bfa^N_\e) dx.
     \end{aligned} 
\label{limegelimN}   \end{equation}
Since $|\bfa_\e^N|_{L^2(\Omega;\RR^3)}\leq |\bfa_\e|_{L^2(\Omega;\RR^3)}\le C$,
 by the metrizability  of the weak topology on   bounded subsets of $L^2(\Omega;\RR^3)$,
   there exists a subsequence  of  $(\bfa_\e)_{\e>0}$ (denoted the same) such that, for every $N\in\NN$, 
   $(\bfa_\e^N)_{\e>0}$ weakly converges in $L^2(\Omega;\RR^3)$  to some $\bfa^N$.
As   $(\bfa^N)_{N\in\NN}$  is bounded in $L^2(\Omega;\RR^3)$, a subsequence still denoted   $(\bfa^N)_{N\in\NN}$  
 weakly converges  as $N\to+\infty$  to some $\bf\tilde\bfa$. 
We fix $N\in \NN$ and $\alpha>0$. By Proposition \ref{propcapreRep=2},   $(c^{f^{\infty,2},S}_{r_\e,R_\e/c_0})_{\e>0}$ uniformly converges 
on $B_{\RR^3}(0,N)$ to $c_0^f$, thus
%
\begin{equation}
\begin{aligned}
     &  \exists \e_0>0, \ \ 
      c^{f^{\infty,2},S}_{r_\e,R_\e/c_0}(\bfxi) \ge c_0^f(\bfxi)- \frac{\alpha}{\gamma^{(2)}|\Omega| } \quad \forall \bfxi\in B_{\RR^3}(0,N), \ \forall \e<\e_0.
     \end{aligned} 
\nonumber  \end{equation}
We infer from   the weak lower semicontinuity on $L^2(\Omega;\RR^3)$ of
  $\bfpsi\to  \int_\Omega         c_0^f(\bfpsi) dx$ 
\begin{equation}
\begin{aligned}
    \liminf_{\e\to0}    \gamma_\e^{(2)}(r_\e)     \int_\Omega           \hskip-0,1cm c^{f^{\infty,2},S}_{r_\e,R_\e/c_0}(\bfa_\e^{N}) dx 
    & \ge    \liminf_{\e\to0}    \gamma^{(2)}     \int_\Omega         c_0^f(\bfa_\e^N) dx  -\alpha  
    \\&
     \ge      \gamma^{(2)}    \int_\Omega         c_0^f(\bfa^{N}) dx  -\alpha,  
     \end{aligned} 
\nonumber  \end{equation}
and then from  \eqref{basedonI2ep=2}, \eqref{hatfcap1}, \eqref{liminf=lime}, \eqref{limegelimN}
     and the arbitrariness of $\alpha$,  
\begin{equation}
\begin{aligned}
    \liminf_{\e\to0}   I_{\e2} 
  \ge     \gamma^{(2)}    \int_\Omega         c_0^f(\bfa^{N}) dx   -C \frac{|S\setminus S'|}{|S|}\qquad \forall N \in \NN.  
     \end{aligned} 
\nonumber  \end{equation}
By passing to the limit inferior as $N\to+\infty$, we obtain 
\begin{equation}
\begin{aligned}
    \liminf_{\e\to0}   I_{\e2} 
  \ge     \gamma^{(2)}    \int_\Omega         c_0^f({\bf\tilde\bfa}) dx   -C \frac{|S\setminus S'|}{|S|}.  
     \end{aligned} 
\label{liminf1}   \end{equation}
%
Let us fix $\bfvarphi\in \D(\Omega;\RR^3)$. By \eqref{defaN},  
\begin{equation}
\begin{aligned}
   \lb   \int_\Omega      (\bfa_\e-\bfa_\e^N)\cdot\bfvarphi dx \rb = \lb \int_{|\bfa_\e|>N} \bfa_\e\lp 1-\frac{N}{|\bfa_\e|}\rp\cdot \bfvarphi dx\rb\leq   \int_{|\bfa_\e|>N} | \bfa_\e| dx .
     \end{aligned} 
\nonumber  \end{equation}
By passing to the limit as $\e\to0$, taking \eqref{cvaalpha} into account, we infer
\begin{equation}
\begin{aligned}
   \lb   \int_\Omega      (\bfv-\bfu -\bfa^N)\cdot\bfvarphi dx \rb  \leq C  \sup_{\e>0}  \int_{|\bfa_\e|>N} | \bfa_\e| dx .
     \end{aligned} 
\nonumber   \end{equation}
The sequence  $(\bfa_\e)_{\e>0}$ is bounded in $L^2(\Omega;\RR^3)$, 
hence equiintegrable, 
therefore  (see   \cite[Prop. 1.27]{AmFuPa})
$  \lim_{N\to+\infty}   \sup_{\e>0}         \int_{|\bfa_\e|>N} | \bfa_\e| dx =0$.
By passing to the limit as $N\to+\infty$, we infer  
$  \int_\Omega      (\bfv-\bfu -{\bf\tilde\bfa})\cdot\bfvarphi dx=0$
and 
  deduce
 ${\bf\tilde\bfa}=\bfv-\bfu$. By      \eqref{liminf1}   
and    the arbitrariness of $S'$,    \eqref{linf2} holds. 
     Proposition \ref{proplowerbound} is proved.  \qed

     %
  %
\subsubsection{Upper  bound}\label{secupperbound}
$ $

 \begin{proposition}
 \label{propupperbound}
Let   $(\bfvarphi, \bfpsi, \zeta, a,b)$
 be such that
\begin{equation}\begin{aligned}
 &
(\bfvarphi, \bfpsi, \zeta, a,b)\in  (C^\infty(\ov\Omega;\RR^3)^2\times (C^\infty(\ov\Omega))^3),\quad 
\\& (\bfvarphi, \bfpsi^{tuple}) \in  W^{1,p}_b(\Omega;\RR^3)\times \D, \qquad \Phi(\bfvarphi, \bfpsi^{tuple})<+\infty,
\end{aligned}
\label{psitupleinD}
  \end{equation} 
  where  $\D$ is defined by (\ref{defD}) and 
    \begin{equation}
   \begin{aligned}
& 
 \bfpsi^{tuple}:=  ( \bfpsi, \zeta)\   \text{if }  \kappa=0; \qquad  
\bfpsi^{tuple}:=  ( \bfpsi, \zeta, a,b) \  \text{if } 0<\kappa\le+\infty,
\end{aligned} 
 \label{defpsituple}
  \end{equation}
  and let $\alpha>0$.
There exists  a sequence $(\bfvarphi_\e)$   in $W^{1,p}_b(\Omega; \RR^3)$   such that 
\begin{equation}\begin{aligned}
&
\bfvarphi_\e   \rightharpoonup \bfvarphi \quad   \hbox{   weakly  in } \  W^{1,p} (\Omega;\RR^3) , 
\\&  \ov\bfgotv^S_\e(\bfvarphi_\e)  \rightharpoonup \bfpsi   \quad     \hbox{ weakly in    }     L^p(\Omega;\RR^3),  
\quad 
\ov\theta^S_\e(\bfvarphi_\e)   \rightharpoonup \zeta   \   \hbox{ weakly in    }     L^p(\Omega),  
\\& \hskip-0,1cm   \frac{1}{r_\e}\lp\ov \gotv^S_{\e3}(\bfvarphi_\e)  ,  \ov\theta^S_\e(\bfvarphi_\e)  \rp  \rightharpoonup (a,b)   \     \hbox{ weakly in    }     L^p(\Omega)^2 \  \text{ if } \quad 0<\kappa<+\infty,  
\\& \tau^{-1}  \frac{\ov\gots^{S}(\bfvarphi_\e)}{r_\e}  m_\e
 \buildrel \star \over \rightharpoonup  
b \quad   \hbox{weakly$^\star$  in }   \M (\ov\Omega) \  \text{ if } \quad 0<\kappa<+\infty, 
\end{aligned}\hskip-0,2cm 
\label{ub}
  \end{equation} 
 and  
  \begin{equation}\begin{aligned}
&  \limsup_{\e\to 0}  F_\e(\bfvarphi_\e)  
\leq \Phi(\bfvarphi, \bfpsi^{tuple})
+\alpha.
\end{aligned} 
\label{lsupF}
  \end{equation} 
     \end{proposition}
 
   \noindent
     \begin{remark}\label{remFGamma}
Proposition  \ref{propupperbound} implies,  by  density and diagonalisation arguments,   that for every $\bfu\in W^{1,p}_b(\Omega;\RR^3)$, there exists a sequence $(\bfu_\e)$ weakly converging to $\bfu$ in $W^{1,p}_b(\Omega;\RR^3)$
such that $\limsup_{\e\to0} F_\e(\bfu_\e)$ $ \leq F(\bfu)$ (see \eqref{defF}).  By
  Proposition  \ref{proplowerbound}, for every $\bfu\in W^{1,p}_b(\Omega;\RR^3)$
and every  sequence $(\bfu_\e)$ weakly converging   in $W^{1,p}_b(\Omega;\RR^3)$ to  $\bfu$,  
 $\liminf_{\e\to0} F_\e(\bfu_\e)\ge F(\bfu)$. 
Thus, since   $(F_\e)_{\e>0}$ is equicoercive on $W^{1,p}_b(\Omega;\RR^3)$, 
 it   $\Gamma$-converges to $F$ in the weak topology of $W^{1,p}_b(\Omega;\RR^3)$ (see \cite{Da}).
 \end{remark}
 %
\noindent {\bf Proof.} \ 
Let us fix $R>0$, a sequence $(R_\e)$  satisfying \eqref{Re},  and  a Lipschitz  domain $S'$  such that 
\begin{equation}
\begin{aligned}
& S \subset\subset S'\subset\subset RD, \qquad \bfy_{S'}=0.
\end{aligned}
\label{S'sup}
\end{equation}  
The  approximating sequence will be  defined by 
\begin{equation}
\begin{aligned}
&
\bfvarphi_\e(x):=  \lp1-\xi_\e(x) \rp  (\bfvarphi(x)+  \bfeta_\e(x)) + \xi_\e(x) \bfchi_\e(x),
  \end{aligned} 
\label{defvarphie}  
  \end{equation} 
where   $\xi_\e$ is such that
\begin{equation}
\begin{aligned}
  &\xi _\e  \in \D(\RR^2),
\qquad     0 \le \xi_\e \le 1, \qquad
 | \bfnabla \xi_\e| \le  \frac{C}{r_\e},
 \\&  \xi_\e = 1\quad  \hbox{ in } S_{r_\e}, \qquad \xi_\e = 0 \quad \hbox{ in }  \RR^2\setminus S'_{r_\e}, 
  \end{aligned} 
\label{defxi}    \end{equation} 
  $\bfchi_\e$  is defined  either  by (\ref{chik}),  (\ref{chikappa}),  or      \eqref{defchiothercases}, depending on  the magnitude of   $k$ and $\kappa$, and $\bfeta_\e$ is given by \eqref{defetaep<2} 
  if  $\gamma^{(p)}<+\infty$ and  $p<2$, 
 by \eqref{defetaep=2}   if  $\gamma^{(p)}<+\infty$ and  $p=2$, and  otherwise  $\bfeta_\e=0$.
 These choices ensure  that $\bfvarphi_\e$ belongs to $W^{1,p}_b(\Omega;\RR^3)$ and satisfies  (\ref{ub}). 
 The support of the  function $\bfeta_\e$ is included in $D_{\widetilde R_\e}\times (0,L)$, where 
 \begin{equation}
\begin{aligned}
  &
  \widetilde R_\e = \begin{cases} R r_\e & \hbox{if } \ p<2, \\ R_\e & \hbox{if } \ p=2, \end{cases}
  \end{aligned} 
\label{Relimsup}    \end{equation} 
therefore, the following decomposition holds
\begin{equation}
\begin{aligned}
&F_\e( \bfvarphi_\e)= I_{\e1}+I_{\e2}+I_{\e3}+I_{\e4} ;
\quad  I_{\e1}:= \int_{ \Omega\setminus  (D_{\widetilde R_\e}\times(0,L))} f(\bfe(\bfvarphi )) dx;
\\& I_{\e2}:= \int_{(D_{\widetilde R_\e}\setminus S'_{r_\e} )\times(0,L)}   f\lp \bfe\lp \bfvarphi+  \bfeta_\e  \rp \rp dx;\quad  I_{\e3}  :=k_\e \int g \lp   \bfe\lp \bfchi_\e \rp\rp dm_\e;
\\&I_{\e4} :=   \int_{S'_{r_\e} \setminus  S_{r_\e}    \times(0,L) }f(\bfe(\bfvarphi_\e)) dx.
  \end{aligned} 
\label{splitsup0}     \end{equation} 
Clearly,  
\begin{equation}
\begin{aligned}
&\lim_{\e\to0}   I_{\e1}= \int_{ \Omega}f(\bfe(\bfvarphi )) dx.
  \end{aligned} 
\label{limI1}    \end{equation} 
The assertion \eqref{lsupF} will result from 
\begin{eqnarray}
 & & \limsup_{\e\to0}  I_{\e3} \leq \int_\Omega g^{hom}(\bfgotD(\bfpsi^{tuple}))+ \frac{\alpha}{2},
\label{limI3} 
\\
& & \lime I_{\e4} =0, \qquad \label{limI4} 
\\
& &\limsup_{\e\to 0}  I_{\e2}\le \int_\Omega c^f(\bfpsi-\bfvarphi, \zeta) dx + \frac{\alpha}{2}.
  \label{limI2}
\end{eqnarray} 
 %
%
 {\bf Proof of \eqref{limI3}.} If  $0<k<+\infty$,  
we  fix    $\bfq\in C^1_c(\Omega; C^\infty(\ov S;\RR^3))$ such that  (see \eqref{defghom1} and  Remark \ref{remghom})
   \begin{equation}
\begin{aligned}
& k\int_\Omega    \intb_{S} g\lp   \bfgotL\lp \bfq, \frac{\partial \psi_3}{\partial x_3}, \frac{\partial \zeta}{ \partial x_3}\rp  \rp dx dy  
 \leq \int_\Omega g^{hom}  \lp  \bfgotD(\bfpsi^{tuple})\rp dx + \frac{\alpha}{2}, 
\end{aligned} 
\label{condh}   \end{equation}
and  set 
(see  \eqref{defovzeta})   
  \begin{equation}
\begin{aligned}
\hskip-0,1cm \bfchi_\e\hskip-0,1cm:=    \ov  \bfpsi_\e\hskip-0,1cm  +\hskip-0,1cm\frac{2}{\diam S }\ov \zeta_\e  \bfe_3  \wedge  \frac{ \bfy_\e}{r_\e}  
\hskip-0,05cm - r_\e \hskip-0,05cm\lp \frac{\partial \ov\psi_{\e1} }{ \partial x_3}\frac{y_{\e1}}{r_\e}\hskip-0,1cm+ \hskip-0,1cm \frac{\partial \ov\psi_{\e2} }{ \partial x_3}\frac{y_{\e2}}{r_\e}\rp\bfe_3
 + 
  r_\e   \bfq\lp x,\frac {y_\e(x')}{r_\e} \rp\hskip-0,1cm  . 
  \end{aligned} 
\label{chik}     \end{equation} 
  By \eqref{defghom1}, \eqref{zetaovzeta} and \eqref{chik}, the following estimate holds 
 \begin{equation}
\begin{aligned}
&
\lb
\bfe\lp \bfchi_\e\rp - \bfgotL \lp \bfq,  \frac{\partial \psi_3 }{ \partial x_3} , \frac{\partial \zeta }{ \partial x_3} \rp\lp x,\frac{y_\e(x')}{r_\e}\rp
\rb
  \leq C r_\e\qquad \hbox{in } T_{r_\e}.
   \end{aligned} 
\nonumber
 \end{equation} 
Applying  \eqref{inconv},  we deduce 
\begin{equation}
\begin{aligned}
&
\lb
g\lp \bfe\lp \bfchi_\e\rp\rp - g\lp  \bfgotL \lp \bfq,  \frac{\partial \psi_3 }{ \partial x_3} , \frac{\partial \zeta }{ \partial x_3} \rp\lp x,\frac{y_\e(x')}{r_\e}\rp\rp
\rb
  \leq C r_\e\qquad \hbox{in } T_{r_\e}.
   \end{aligned} 
\nonumber \end{equation} 
Taking  (\ref{twoscalecontinuous}), \eqref{splitsup0}, 
 and    (\ref{condh})   into account, we  infer
\begin{equation}
\begin{aligned}
\lime I_{\e3} 
  &=\lime  k_\e \int g  \Bigg (\bfgotL \lp \bfq,  \frac{\partial \psi_3 }{ \partial x_3} , \frac{\partial \zeta }{ \partial x_3} \rp\lp x,\frac{y_\e(x')}{r_\e}\rp \Bigg)
dm_\e   
\\&=  k \int_\Omega    \intb_{S} g \lp   \bfgotL\lp \bfq, \frac{\partial \psi_3}{\partial x_3}, \frac{\partial \zeta}{ \partial x_3}\rp   \rp dx dy  
 \leq    \int_\Omega g^{hom} \lp  \bfgotD(\bfpsi^{tuple})
  \rp dx + \frac{\alpha}{2}. 
   \end{aligned} 
\nonumber   \end{equation} 
 %
\\  If  $0<\kappa<+\infty$, we fix 
 $\bfq\in C^1_c(\Omega; C^\infty(\ov S;\RR^3))$ 
 such that
   (see (\ref{defghom2}))
   \begin{equation}
\begin{aligned}
& \kappa\int_\Omega    \intb_{S} g^{0,p}\lp   \bfgotJ\lp  \bfq, \frac{\partial^2 \psi_1}{ \partial x_3^2},\frac{\partial^2 \psi_2}{ \partial x_3^2}, \frac{\partial a}{ \partial x_3},\frac{\partial b}{ \partial x_3}\rp  \rp dx dy  
\\&\hskip5cm  \leq    \int_\Omega g^{hom} \lp  \bfgotD(\bfpsi^{tuple})\rp dx +\frac{\alpha}{2},
\end{aligned} 
\label{condhkappa}   \end{equation}
and   set 
  \begin{equation}
\begin{aligned}
  \bfchi_\e(x)\hskip-0,1cm:=    \ov  \bfpsi_\e   (x)
+ 
    r_\e&\hskip-0,1cm \begin{pmatrix}-2\frac{2}{\diam S}\ov a_\e  \frac{y_{\e2} }{ r_\e} 
 \\ 2\frac{2}{\diam S}\ov a_\e  \frac{y_{\e1} }{ r_\e}  
 \\ \ov b_\e(x)-  \frac{y_{\e1}}{r_\e}  \frac{\partial}{\partial x_3} \ov \psi_{\e1}-  \frac{y_{\e2}}{r_\e} \frac{\partial}{\partial x_3} \ov \psi_{\e2}
 \end{pmatrix}\hskip-0,1cm
  +   r_\e^2 \bfq\lp x, \frac{y_\e }{ r_\e} \rp
.   \hskip-0,1cm
  \end{aligned} 
\label{chikappa}     \end{equation} 
Noticing that by      \eqref{defD} and   \eqref{psitupleinD}, $\psi_3= 0$,  a straightforward computation yields
 \begin{equation}
\begin{aligned}
 \hskip-0,2cm  \bfe( \bfchi_\e)  = 
 r_\e &\frac{4}{\diam S}\frac{\partial \ov a_\e}{ \partial x_3}\lp -\frac{y_{\e2} }{ r_\e} \bfe_1\odot\bfe_3+ \frac{y_{\e1} }{ r_\e} \bfe_2\odot\bfe_3\rp  + r_\e \bfe_{y}(\bfq)\lp x, \frac{y_\e}{ r_\e}\rp 
 \\ &+r_\e\hskip-0,1cm\lp \frac{\partial \ov b_\e}{ \partial x_3}\hskip-0,1cm -\hskip-0,1cm\frac{\partial^2 \ov \psi_{\e1}}{ \partial x_3^2}\frac{y_{\e1} }{ r_\e} 
\hskip-0,1cm -\hskip-0,1cm\frac{\partial^2 \ov \psi_{\e2}}{ \partial x_3^2}\frac{y_{\e2} }{ r_\e}\rp \bfe_3\hskip-0,05cm\odot\hskip-0,05cm\bfe_3
 + r_\e^2  \bfe_x(\bfq)\lp x, \frac{y_\e}{ r_\e}\rp   .   \hskip-0,1cm
  \end{aligned} 
\label{edchi}   \end{equation} 
Substituting $\bfchi_\e$ for $\bfu_\e$ in \eqref{I31}, we infer
  \begin{equation}
\hskip-0,1cm \begin{aligned}
\lb I_{\e3} -r_\e^p  k_\e \int  g^{0,p}\lp\frac{\bfe(\bfchi_\e)}{r_\e}\rp dm_\e\rb&
\leq C  \varpi(r_\e). \hskip-0,2cm 
 \end{aligned} 
\label{li31}
  \end{equation} 
By    (\ref{defghom2}),  \eqref{zetaovzeta} and \eqref{edchi}, we have 
\begin{equation}
\begin{aligned}
&
\Bigg| \frac{\bfe\lp\bfchi_\e \rp}{ r_\e}
-\bfgotJ\lp  \bfq, \frac{\partial^2 \psi_1}{ \partial x_3^2},\frac{\partial^2 \psi_2}{ \partial x_3^2}, \frac{\partial a}{ \partial x_3},\frac{\partial b}{ \partial x_3}\rp\lp x,\frac{y_\e(x')}{r_\e}\rp   \Bigg|
  \leq C r_\e\qquad \text{ in } T_{r_\e},
   \end{aligned} 
\nonumber   \end{equation} 
hence, by   \eqref{inconv}, 
\begin{equation}
\begin{aligned}
\hskip-0,2cm &
\Bigg| g^{0,p}\lp\hskip-0,1cm\frac{\bfe\lp\bfchi_\e \rp}{r_\e}\hskip-0,1cm\rp
-g^{0,p}\hskip-0,1cm\lp\bfgotJ\lp  \bfq, \frac{\partial^2 \psi_1}{ \partial x_3^2},\frac{\partial^2 \psi_2}{ \partial x_3^2}, \frac{\partial a}{ \partial x_3},\frac{\partial b}{ \partial x_3}\rp\hskip-0,1cm\lp x,\frac{y_\e(x')}{r_\e}\rp\hskip-0,1cm\rp  \hskip-0,1cm  \Bigg|\hskip-0,1cm
  \leq C r_\e \text{ in } T_{r_\e}.\hskip-0,3cm
   \end{aligned} 
\nonumber\end{equation} 
 It follows from  \eqref{condvarpi},  (\ref{twoscalecontinuous}),     \eqref{li31} and  the abobe inequality,    that 
\begin{equation}
\begin{aligned} \lime I_{\e3} 
 & =\lime   r_\e^p  k_\e\int g^{0,p}\!  \Bigg (   \bfgotJ\lp  \bfq, \frac{\partial^2 \psi_1}{ \partial x_3^2},\frac{\partial^2 \psi_2}{ \partial x_3^2}, \frac{\partial a}{ \partial x_3},\frac{\partial b}{ \partial x_3}\rp\lp x,\frac {y_\e(x')}{r_\e}\rp\Bigg)
dm_\e   
\\&=    \kappa\int_\Omega    \intb_{S} g^{0,p} \lp \bfgotJ\lp  \bfq, \frac{\partial^2 \psi_1}{ \partial x_3^2},\frac{\partial^2 \psi_2}{ \partial x_3^2},
\frac{\partial a}{ \partial x_3},\frac{\partial b}{ \partial x_3}\rp\rp dx dy,  
 \end{aligned} 
\nonumber  \end{equation} 
which, combined with  \eqref{condhkappa},  yields   \eqref{limI3}.
   In the remaining  cases, we set 
\begin{equation}
\begin{aligned} 
&\bfchi_\e= \ov \bfpsi_\e -r_\e \lp  \frac{\partial \ov\psi_{\e1}}{ \partial x_3}\frac{y_{\e1}}{r_\e}
 +\frac{\partial \ov\psi_{\e2}}{ \partial x_3}\frac{y_{\e2}}{r_\e}\rp\bfe_3\qquad & & \hbox{if  $k=+\infty$ and $\kappa=0$,} 
 \\&\bfchi_\e=0  & & \hbox{if $k=\kappa=+\infty$.}
 \end{aligned} 
\label{defchiothercases}
 \end{equation} 
If $k=+\infty$ and $\kappa=0$, by    (\ref{defD}) and   (\ref{psitupleinD}),
  $\zeta=\psi_3= 0$, $\bfpsi^{tuple}=0$, hence $|\bfe(\bfchi_\e)|\leq C r_\e$  in $T_{r_\e} $ and, by  \eqref{kkappa},  
  $ I_{\e3}\le C r_\e^p k_\e \to 0=\int_\Omega g^{hom}(\bfgotD(\bfpsi^{tuple})) dx$. 
The case $k=\kappa=+\infty$ is straightforward. 
 The assertion \eqref{limI3} is proved.  \qed 

$ $
\\
{\bf Proofs of  \eqref{limI4} and \eqref{limI2}.} 
If $\gamma^{(p)}=+\infty$,   these assertions   straightforwardly  results  from   $\bfeta_\e=0$.  Otherwise, 
we distinguish two cases:
\\
\\
$\bullet${ \it Case   $p<2$ and   $\gamma^{(p)}<+\infty$.}  
By      \eqref{cf=} and  \eqref{estimintcap0},   we can assume that $R$ and $S'$ satisfy, besides \eqref{S'sup}, 
\begin{equation}
\begin{aligned}
&  
 \gamma^{(p)} \int_\Omega \hskip-0,1cm  \capsca^{f^{\infty,p}}\lp \bfpsi-\bfvarphi,  \tfrac{\diam  S'}{\diam S} \zeta\bfe_3 ;S',R D\rp  dx
\leq \int_\Omega\hskip-0,1cm c^f (\bfpsi-\bfvarphi, \zeta) dx
 +\frac{\alpha}{4}.
  \end{aligned} 
\nonumber
  \end{equation} 
  Lemma \ref{lemapproxcap} ensures the existence of      $\bfeta$
such that
\begin{equation}
\begin{aligned}
&\bfeta \in  \bfgotW^p_b\lp\bfpsi-\bfvarphi, \tfrac{\diam  S'}{\diam S}  \zeta\bfe_3;S',RD\rp, 
\\&
  \gamma^{(p)}  \int_{\Omega\times RD}   f^{\infty,p}(\bfe_{y}(\bfeta(x,y)) dxdy 
 \leq   \int_\Omega c^f (\bfpsi-\bfvarphi, \zeta) dx+\frac{\alpha}{2}.
  \end{aligned} 
\label{etaPhicap}    \end{equation} 
Extending  $\bfeta$ by $0$ to $\ov\Omega\times \RR^2$, we set
\begin{equation}
\begin{aligned}
&
 \bfeta_\e(x):= \bfeta\lp x, \frac{y_\e(x')}{r_\e}\rp,
   \end{aligned} 
\label{defetaep<2}    \end{equation} 
By \eqref{etaPhicap} and  \eqref{defetaep<2}, we have
 \begin{equation}
\begin{aligned}
& \bfeta_\e =0 \quad \hbox{ in }\Omega \setminus D_{R r_\e}\times (0,L), \quad
\\& \bfeta_\e= \bfpsi-\bfvarphi  +\frac{2}{\diam S } \zeta  \bfe_3  \wedge  \frac{\bfy_\e(x')}{r_\e} \  \hbox{ in } S'_{r_\e}\times(0,L),
 \end{aligned} 
\label{eta=}    \end{equation} 
and, by  \eqref{splitsup0},  
 %
\begin{equation}
\begin{aligned}
& I_{\e4} :=   \int_{S'_{r_\e} \setminus  S_{r_\e}    \times(0,L) } f\lp \bfe\lp \lp1-\xi_\e \rp \lp  \bfpsi   +\frac{2}{\diam S } \zeta  \bfe_3  \wedge  \frac{ \bfy_\e}{ r_\e} \rp 
+ \xi_\e \bfchi_\e\rp\rp dx.
  \end{aligned} 
\nonumber  \end{equation} 
\noindent   The assertion   \eqref{limI4}    results from the estimates 
 $\Big|  \bfe\Big( \lp1-\xi_\e \rp \lp  \bfpsi   +\frac{2}{\diam S } \zeta  \bfe_3  \wedge  \frac{ \bfy_\e}{ r_\e} \Big)
+ \xi_\e \bfchi_\e\rp\Big| $ $\leq C$   and $|S'_{r_\e}\setminus S_{r_\e}\times(0,L)|\le C \frac{r_\e^2}{\e^2}$, 
     deduced from  \eqref{defvarphie}, \eqref{chik},  \eqref{chikappa} and  \eqref{defchiothercases}.
Applying  \eqref{inconv} with 
  \begin{equation}
\begin{aligned}
h=f, \quad    \bfM_\e(x)=\bfe\lp \bfvarphi +  \bfeta_\e\rp (x),\quad \bfM'_\e(x)= \frac{1}{r_\e} \bfe_{y}(\bfeta)\lp x, \frac{y_\e(x')}{r_\e}\rp,
  \end{aligned} 
\label{MM'}    
 \end{equation} 
   integrating   over $D_{Rr_\e}\setminus S'_{r_\e} \times(0,L)$,  noticing that $|\bfM_\e-\bfM_\e'|\leq C$,
$|\bfM_\e|+ |\bfM_\e'|\leq \frac{C}{r_\e}$, and 
 $ |D_{Rr_\e\setminus S_{r_\e}'\times(0,L)}| \leq C \frac{r_\e^2}{\e^2}$,    we infer  
 \begin{equation}
\begin{aligned}
    \lb  I_{\e2} -\int_{(D_{Rr_\e}\setminus S'_{r_\e} )\times(0,L)}  \hskip-1cm  f\lp\bfM_\e'(x)\rp  dx \rb
&\leq C\frac{r_\e^2}{\e^2}\lp 1 + \frac{1}{r_\e^{p-1}} \rp
  \leq C r_\e\gamma_\e(r_\e) \leq C r_\e.
  \end{aligned} 
\nonumber
 \end{equation} 
By  \eqref{finftyfg0g}, we have
 \begin{equation}
\begin{aligned}
  \lb  \int_{(D_{Rr_\e}\setminus S'_{r_\e} )\times(0,L)}  \hskip-1cm ( f-   f^{\infty,p})\lp \bfM_\e'(x)\rp  dx \rb
&\leq C  \int_{(D_{Rr_\e}\setminus S'_{r_\e} )\times(0,L)}  1 +|\bfM_\e'(x)|^{p-\varsigma} dx
\\& \leq  C\frac{r_\e^2}{\e^2}\lp 1 + \frac{1}{r_\e^{p-\varsigma}} \rp
   \leq C r_\e^\varsigma \gamma_\e(r_\e) \leq C r_\e^\varsigma.
  \end{aligned} 
\nonumber
 \end{equation} 
 By  \eqref{etaPhicap},      $\bfe_{y}( \bfeta)\lp x, \frac{y_\e(x')}{r_\e}\rp=0$ in $S'_{r_\e}\times(0,L)$, 
hence, by    \eqref{MM'},
 %
  \begin{equation}
\begin{aligned}
   \int_{(D_{Rr_\e}\setminus S'_{r_\e} )\times(0,L)}  \hskip-1cm  & f^{\infty,p} \lp  \bfM_\e'(x)\rp  dx 
 =  \gamma_\e^{(p)}(r_\e) |RD| \hskip-0,1cm \int f^{\infty,p}  \lp  \bfe_{y}( \bfeta)\lp x, \frac{y_\e(x')}{r_\e}\rp \rp  d\nu_\e,
 \end{aligned} 
\nonumber
 \end{equation} 
where    $\nu_\e :=\frac{\e^2}{r_\e^2 |RD|}   \mathds{1}_{(RD)_{r_\e} \times(0,L)} \calL^3_{\lfloor \Omega}$. Notice that  $\nu_\e$ is deduced from  $m_\e$  by substituting $RD$ for $S$ in \eqref{defme}.
By applying 
 \eqref{twoscalecontinuous}  to  $(\nu_\e,f^{\infty,p}(\bfe_{y}(\bfeta)))$  in place of $(m_\e,\bfpsi)$, taking the above estimates into account,
 we obtain
  \begin{equation}
\begin{aligned}
\lime I_{\e2}&=\lime
   \gamma_\e^{(p)}(r_\e) |RD| \int f^{\infty,p}   \lp  \bfe_{y}( \bfeta)\lp x, \frac{y_\e(x')}{r_\e}\rp \rp  d\nu_\e
 \\& =   \gamma^{(p)}  \int_{\Omega\times RD}  f^{\infty,p}  \lp  \bfe_{y}( \bfeta)\lp x, y \rp \rp  dxdy.
 \end{aligned} 
\nonumber
 \end{equation} 
 %
The assertion \eqref{limI2} is proved.  
 The proof of Proposition \ref{propupperbound} in the case $p<2$, $\gamma^{(p)}<+\infty$ is achieved.  \qed

$ $

\noindent  $\bullet$  {\it Case $p=2$,   $\gamma^{(2)}<+\infty$.} 
By Lemma \ref{reg}, we can fix  for each $\bfa\in \RR^3$, $\e>0$  a field   $\bfeta_\e^{\bfa}\in \W^2(\bfa,0;r_\e S',R_\e D)$ such that  
\begin{equation}
\begin{aligned}
  &\bfeta_\e^{\bfa} \ \hbox{ is a solution to } \  \calP^{ f} (\bfa,0; r_\e S', R_\e D) \quad \forall \bfa\in \RR^3, \ \forall \e>0.
  \end{aligned} 
\label{defetaa} \end{equation} 
We    extend each $\bfeta_\e^{\bfa}$    by $0$ to $\RR^2$. We
  split $(0,L)$ into a  suitable family of  intervals $(J_\e^k)_{k\in \{1,..,n_\e\}}$  and set
%
 %
\begin{equation}
\begin{aligned}
& \bfphi:= \bfpsi-\bfvarphi,  \qquad\bfphi_\e^{ik} :=  \intb_{Y_\e^{i}\times J_\e^k}
 \bfphi (s)d\calL^3(s), 
 \\ &\widecheck \bfphi_\e (x ):=   \sum_{(i,k)\in I_\e\times\{1,..,n_\e\}}\bfphi_\e^{ik}   \mathds{1}_{Y_\e^{i}\times J_\e^k}(x).
\end{aligned}
\label {defcheckphi}
\end{equation} 
The field $\bfeta_\e$  in \eqref{defvarphie}  will  be defined by 
\begin{equation}
\begin{aligned}
&\bfeta_\e(x):= \rho_\e(x_3)  \tilde\bfeta_\e(x), 
\quad  \tilde\bfeta_\e(x):=  \sum_{(i,k)\in I_\e\times\{1,..,n_\e\}}   \bfeta_\e^{\bfphi_\e^{ik}}(y_\e(x')) \mathds{1}_{Y_\e^{i}\times J_\e^k}(x) , 
\end{aligned}
\label{defetaep=2}
\end{equation} 
where $\bfeta_\e^{\bfphi_\e^{ik}}$ is given by  \eqref{defetaa}. The choice of the  function $\rho_\e$ (see below)
 ensures    that  $\bfeta_\e$ approximates to $\tilde\bfeta_\e$ and 
 belongs to $  H^1(\Omega;\RR^3)$.  
   By \eqref{defetaa} and \eqref{defcheckphi}, we have
 \begin{equation}
\begin{aligned}
& \spt \bfeta_\e\subset D_{R_\e}\times (0,L); \quad \bfeta_\e(x)= \rho_\e(x_3) \widecheck \bfphi_\e (x)  \quad \hbox{ in } S'_{r_\e}\times(0,L),
\\&\tilde\bfeta_\e(x)=   \widecheck \bfphi_\e (x) \quad   \hbox{ and }  \quad \bfe_{x'}(\tilde\bfeta_\e) =0 \quad \hbox{ in }  S_{r_\e}'\times(0,L).
 \end{aligned} 
\label{eta=2}    \end{equation} 
The intervals $J_\e^k$  will be defined  by fixing 
  sequences    $(a_\e)$ and $(b_\e)$   such that 
   \begin{equation}
\begin{aligned}
& 0< b_\e\ll a_\e\ll 1, \qquad \frac{r_\e^2}{\e^2} + R_\e^2 \ll a_\e b_\e,
  \end{aligned} 
\label{hypab1}\end{equation}
and setting  ($[s]$:    integer part of   $s$)
  \begin{equation}
\begin{aligned}
&  J_\e^k:=  (l_{k,\e}, l_{k+1,\e}); \qquad 
  l_{k,\e}:= k a_\e; \qquad n_\e:= \lc \frac{L}{a_\e}-1\rc.
  \end{aligned} 
\nonumber \end{equation}
The function $\rho_\e$ introduced in \eqref{defetaep=2} is chosen such that 
\begin{equation}
\begin{aligned}
& \rho_\e\in C^\infty([0,L]);\quad    \rho_\e=1  \  \text{ in } \  (0,L) \setminus \bigcup_{k=1}^{ n_\e}     \left(  l_{k,\e}-\frac{1}{2}b_\e;  l_{k,\e}+\frac{1}{2}b_\e \right) ,  
 \\&  \rho_\e=0  \quad \text{ on }  \quad \bigcup_{k=1}^{ n_\e}   \{ l_{k,\e} \},
\qquad 0\leq \rho_\e \leq 1,\qquad 
|\rho_\e'|<\frac{C}{b_\e}.
  \end{aligned} 
\label{defrhoe}\end{equation}
By \eqref{defvarphie},  \eqref{splitsup0} and  \eqref{eta=2}, we  have 
\begin{equation}
\begin{aligned}
I_{\e4} :=   \int_{S'_{r_\e} \setminus  S_{r_\e}    \times(0,L) } f\lp \bfe\lp \lp1-\xi_\e \rp ( \rho_\e(x_3) \widecheck \bfphi_\e (x) +\bfvarphi)
+ \xi_\e \bfchi_\e\rp\rp dx.
  \end{aligned} 
\label{splitsup2}     \end{equation} 
The  estimates 
$\lb  \bfe\lp \lp1-\xi_\e \rp  (\rho_\e(x_3) \widecheck \bfphi_\e (x)  +\bfvarphi)
+ \xi_\e \bfchi_\e\rp\rb^2 \leq  C (1+ |\rho'(x_3)|^2)$,    deduced from  \eqref{defvarphie}, \eqref{chik},  \eqref{chikappa} and  \eqref{defchiothercases}, 
 and   $\int_0^L  |\rho_\e'|^2 dx_3\leq  C  \frac{1}{a_\e b_\e}$,
infered  from \eqref{defrhoe}, yield $I_{\e4}\leq  
C  \frac{r_\e^2}{\e^2}\lp 1+ \frac{1}{a_\e b_\e}\rp$,
which by  \eqref{hypab1}  implies    \eqref{limI4}. 
 By  \eqref{defgammae},  \eqref{defcfSrR}  and \eqref{defetaa}    we have, for every $(i,\bfa)\in I_\e\times \RR^3$, 
    \begin{equation}
\begin{aligned}
  &\int_{D_{R_\e}^{i} }   f\lp  \bfe_y (\bfeta_\e^{\bfa})(y_\e(x')) \rp dx' 
=  \capsca^{ f}(\bfa, 0; r_\e S', R_\e D)
\\&\hskip3cm = \frac{1}{|\log r_\e|}  c^{ f,S'}_{r_\e, R_\e}(\bfa)
  = \gamma_\e^{(2)}(r_\e) \int_{Y_\e^{i}}  c^{ f,S'}_{r_\e, R_\e}(\bfa) dx,
  \end{aligned} 
\label{fetaa}\end{equation} 
hence, by  \eqref {defetaep=2} and  \eqref{eta=2},   
   \begin{equation}
\begin{aligned}
 \hskip-0,1cm & \int_{(D_{R_\e}\setminus S'_{r_\e})\times(0,L)} \hskip-0,1cm   f  \lp  \bfe_{x'}(\tilde\bfeta_\e) \rp  dx  
   = \hskip-0,1cm\hskip-0,1cm\hskip-0,1cm\hskip-0,1cm \sum_{(i,k)\in I_\e\times\{1,..,n_\e\}}\int_{J_\e^k} \hskip-0,1cm\hskip-0,1cmdx_3   \int_{D_{R_\e}^{i} } \hskip-0,1cm\hskip-0,1cm f\lp  \bfe_{y}(\bfeta_\e^{\bfphi_\e^{ik}})(y_\e(x')) \rp dx' 
\\& \quad=  \gamma_\e^{(2)}(r_\e)  \int_0^L dx_3 \sum_{i\in I_\e} \int_{Y_\e^{i}} c^{ f,S'}_{r_\e, R_\e} (\widecheck \bfphi_\e) dx
  =  \gamma_\e^{(2)}(r_\e) \int_\Omega c^{ f,S'}_{r_\e, R_\e} (\widecheck \bfphi_\e) dx
. 
  \end{aligned} 
\label{lsI23} \end{equation} 
Since $(\widecheck \bfphi_\e )$ is bounded in $L^\infty$ and uniformly converges   
 on each compact subset of $\Omega$   to $\bfpsi-\bfvarphi$, we deduce from 
 Proposition \ref{propcapreRep=2}  and \eqref{defgammae} that
   \begin{equation}
\begin{aligned}
& 
\lim_{\e\to0}  \int_{(D_{R_\e}\setminus S'_{r_\e})\times(0,L)} \hskip-0,1cm   f  \lp  \bfe_{x'}(\tilde\bfeta_\e) \rp  dx  
 =  \gamma^{(2)}\int_\Omega c^f_0(\bfpsi-\bfvarphi) dx
. 
  \end{aligned} 
\label{lsI231} \end{equation} 
 We prove below that
\begin{equation}
\begin{aligned}
&\lb  I_{\e2}- \int_{(D_{R_\e}\setminus S'_{r_\e})\times(0,L)}    f(\bfe_{x'}(\tilde\bfeta_\e))dx \rb= o(1).
  \end{aligned} 
\label{lsI21}\end{equation} 
Taking \eqref{cf=} and    \eqref{lsI231}  into account, we obtain \eqref{limI2}. 
  We turn to the proof of    \eqref{lsI21}.
By \eqref{defcheckphi},  
 \eqref{defetaep=2} and \eqref{eta=2},   we have 
\begin{equation}
\begin{aligned}
\bfe(\bfvarphi+ \bfeta_\e)(x)  -   \bfe_{x'}(\tilde\bfeta_\e) =
  \bfe(\bfvarphi)+  &(\rho_\e(x_3)-1)\bfe_y (\tilde\bfeta_\e)(y_\e(x'))
  + \rho'(x_3) \tilde\bfeta_\e \odot \bfe_3.
  \end{aligned} 
\label{Ee}\end{equation}
By  \eqref{growthp}, \eqref{defetaep=2} and  \eqref{fetaa}, 
   \begin{equation}
\begin{aligned}
 \hskip-0,1cm  \int_{D_{R_\e}} \hskip-0,1cm     \lb  \bfe_{x'}(\tilde\bfeta_\e) \rb^2 (x',x_3) dx'  
\le C\int_{\Omega'} c^{ f^{\infty,2},S'}_{r_\e, R_\e} (\widecheck \bfphi_\e) (x',x_3) dx'
   \le C \quad \forall x_3\in(0,L)
. 
  \end{aligned} 
\label{ex'bd} \end{equation} 
Taking \eqref{hypab1} and  \eqref{defrhoe} into account, we infer 
\begin{equation}
\begin{aligned}
 \int_{D_{R_\e}\times (0,L)} \hskip-0,5cm  | (\rho_\e(x_3)-1  )\bfe_{x'}  (\tilde\bfeta_\e)|^2  dx
& \leq   \sum_{k=1}^{n_\e} \int_{l_{k,\e}-b_\e/2}^{l_{k,\e}+ b_\e/2} dx_3  \int_{D_{R_\e}}|\bfe_{x'}(\tilde\bfeta_\e) |^2dx'  
 \\& \leq C  \frac{b_\e}{a_\e} =o(1) .
  \end{aligned} 
\label{due}\end{equation}
 By making  suitable changes  of variables
in the inequality $ \int_D\lb  \bfeta
\rb^2 dy  \leq  C 
 \int_D   |\bfe_y (\bfeta)|^2dy$ $ \forall \bfeta\in H^1_0(D;\RR^3)$ (see \eqref{minKornVbounded}), we deduce from   \eqref{ex'bd}
  that 
   \begin{equation}
\begin{aligned}
 \int_{D_{R_\e}} |\tilde\bfeta_\e(y_\e(x'))|^2dx'
\le 
 \hskip-0,1cm  C R_\e^2 \int_{D_{R_\e}} \hskip-0,1cm     \lb  \bfe_{x'}(\tilde\bfeta_\e) \rb^2 (x',x_3) dx'  
\le   C R_\e^2 \quad \forall  x_3\in(0,L)
,
  \end{aligned} 
\nonumber \end{equation} 
yielding, by \eqref{hypab1}  and  \eqref{defrhoe},  
\begin{equation}
\begin{aligned}
\int_{D_{R_\e}\times (0,L)} \hskip-0,5cm |  \rho'(x_3) \tilde\bfeta_\e  \odot \bfe_3 |^2 dx
 & \leq \frac{C}{b_\e^2}  \sum_{k=1}^{n_\e} \int_{l_{k,\e}-b_\e/2}^{l_{k,\e}+ b_\e/2} dx_3 \int_{D_{R_\e}} |\tilde\bfeta_\e(y_\e(x'))|^2dx'
 \\&     \leq  C   \frac{R_\e^2}{a_\e b_\e} =o(1) .
  \end{aligned} 
\label{tre}\end{equation}
Combining   \eqref{inconv},  \eqref{Relimsup}, \eqref{splitsup0},  \eqref{Ee}, 
  \eqref{due},  and \eqref{tre}, we obtain  \eqref{lsI21}.
   \qed

\noindent 
\subsubsection{Conclusion}
\label{secproofth}
\noindent  
The solution  $\bfu_\e$   to \eqref{Pe} 
 satisfies  
 $ F_{\e}(\bfu_\e)-\int_\Omega \bff. \bfu_\e dx \leq  F_{\e}(0) =0$, thus, since $W^{1,p}_b(\Omega;\RR^3)\cap \R=\{0\}$, 
 by \eqref{LpKorn}  the following holds
  \begin{equation}\begin{aligned}
 F_{\e}(\bfu_\e) \leq C |\bfu_\e|_{L^p(\Omega;\RR^3)}  \leq C | \bfu_\e|_{W^{1,p}_b(\Omega;\RR^3)}  \leq (F_\e(\bfu_\e))^{\frac{1}{p}}.\end{aligned}
\nonumber
  \end{equation}
We deduce  
   (\ref{supFeuefini}), and then   the  existence of   a subsequence
    $(\bfu_{\e_k})_{k\in\NN}$  satisfying     the convergences   stated in   Propositions \ref{prop1}, \ref{prop2}  and  \ref{prop3} for some $(\bfu,\bfv^{tuple})\in W^{1,p}_b(\Omega;\RR^3)\times\D$. 
Applying 
 Proposition \ref{proplowerbound}, we infer
\begin{equation}\hskip-0,1cm\begin{aligned}
 &
\liminf_{k\to+\infty}F_{\e_k}(\bfu_{\e_k})-\hskip-0,1cm \int_\Omega\hskip-0,1cm \bff\cdot  \bfu_{\e_k} dx \ge \Phi(\bfu,\bfv^{tuple})-\hskip-0,1cm\int_\Omega\hskip-0,1cm \bff\cdot \bfu  dx\ge \inf (\calP^{hom}) .\hskip-0,1cm
\end{aligned}
\label{linf}
  \end{equation}
  Let us   fix $\alpha  >0$.  
The functional $\Phi$ is  convex  and bounded on the unit ball of the Banach space $W^{1,p}_b(\Omega;\RR^3)\times \D$, hence strongly continuous.
By  density, there exists 
  a  couple  $(\bfvarphi, \bfpsi^{tuple})$ satisfying \eqref{psitupleinD},  \eqref{defpsituple}, \eqref{ub}
   and 
  such that 
    \begin{equation}\begin{aligned}
&
  \Phi(\bfvarphi, \bfpsi^{tuple})-\int_\Omega \bff\cdot  \bfvarphi dx  \leq \inf (\calP^{hom}) +\alpha. 
\end{aligned}
\nonumber
  \end{equation}
Applying  Proposition \ref{propupperbound}, we fix  $(\bfvarphi_\e)\subset W^{1,p}_b(\Omega; \RR^3)$  verifying \eqref{ub} and \eqref{lsupF}. 
 Since $\bfu_\e$ is the solution to  (\ref{Pe}),  we  have
 \begin{equation}\begin{aligned}
  \limsup_{\e\to0} F_\e(\bfu_\e)-\int_\Omega& \bff\cdot   \bfu_\e dx  
\leq \limsup_{\e\to 0}   F_\e(\bfvarphi_\e)\hskip-0,1cm-\hskip-0,1cm\int_\Omega\hskip-0,1cm  \bff\hskip-0,1cm\cdot\hskip-0,05cm\bfvarphi_\e dx 
\\& \leq \Phi(\bfvarphi, \bfpsi^{tuple})
+\alpha  -\hskip-0,1cm\int_\Omega\hskip-0,1cm \bff\hskip-0,1cm\cdot\hskip-0,05cm\bfvarphi  dx
 \leq \inf (\calP^{hom}) +2\alpha.
 \end{aligned}
\label{lsup}
  \end{equation}
By \eqref{linf},  \eqref{lsup}, and  the arbitrariness  of $\alpha$,  
   $(\bfu,\bfv^{tuple})$ is a   solution to  \eqref{Phom}, hence $\bfu$ is a solution to \eqref{defF}. 
By   the  strict convexity of  $F$, this solution is unique, hence the entire sequence $(\bfu_\e)$ weakly  converges to $\bfu$ in $W^{1,p}(\Omega;\RR^3)$. Theorem \ref{th} is proved.
    \qed 
 
   \subsection{Proof of Theorem \ref{thsoft}}\label{secproofthsoft}
 The proof follows the same pattern as that of Theorem \ref{th}. 
We  establish 
the analog of Proposition \ref{proplowerbound},
   assuming   \eqref{condsoft} and 
     \begin{equation} 
\begin{aligned}
   \sup_{\e>0}\  F_\e^{soft}(\bfu_\e) <+\infty.
   \end{aligned}
\label{supFesoftuefini}
  \end{equation}
By   \eqref{resoft}, \eqref{condsoft} and  \eqref{supFesoftuefini},  whatever the operator $\bfgotp_\e \in \{\bfgotv_\e,...\}$
    introduced in Section \ref{secaux}, 
    the sequence  $(\bfgotp_\e(\bfu_\e))$ is bounded in $L^p$ and 
satisfies \eqref{hypfie}.
Repeating their proofs, we  find  that  every assertion     stated  in   Propositions \ref{prop1}, \ref{prop2} and 
\ref{prop3}  holds, except the first line of \eqref{cvbasic}  where the convergence of $(\bfu_\e)$ only holds in  the weak  topology of $L^p(\Omega;\RR^3)$ 
and that of $(\ov\bfgotu_\e( \bfu_\e))$ is irrelevant. 
Setting 
\begin{equation}
\begin{aligned}
&  F_\e^{soft}( \bfu_\e) = I_{\e1} +  I_{\e3} ; 
\quad I_{\e1} := \int_{\Omega\setminus  T_\e} \e^p f(\bfe(\bfu_\e)) dx ;
\quad   I_{\e3}:=  k_\e  \int  g(\bfe(\bfu_\e))  dm_\e,
  \end{aligned} 
\nonumber  \end{equation} 
we infer   \eqref{linf3} (same proof). 
We deduce from 
\cite[Lemma  3.2]{BeBo2}  that, up to a subsequence, the following holds for some $\bfu_0 \in L^p(\Omega; W^{1,p}_\sharp(Y;\RR^3))$:
   \begin{equation}\begin{aligned}
& \bfu_\e \, \dto \, \bfu_0   , \qquad  \e\bfe (\bfu_\e) \, \dto \,  \bfe_y(\bfu_0),
\qquad  \int_Y \bfu_0(.,y)dy=\bfu,
     \end{aligned}
\label{Ldto}
  \end{equation}
  where   $"{\dto}"$  denotes  the  "usual"  two-scale convergence, defined  by substituting
$\calL^3$ for $m_\e$ in  \eqref{defmdto}. In other words,  $\bfu_\e \, \dto \, \bfu_0$ means that    for all $\bfeta\in 
C(\ov{\Omega\times Y};\RR^3)$, 
\begin{equation}\begin{aligned}
\lim_{\e \to 0}\int_\Omega  \bfu_\e(x ) \cdot \bfeta\lp x,  \frac{y_\e(x')}{r_\e}\rp   dx  =   \int_ {\Omega\times Y}  \bfu_0(x,y) \cdot& \bfeta(x,y)  dx dy .
\end{aligned}
\nonumber
  \end{equation}
On the other hand, by  \eqref{defmdto}  and   the third line of \eqref{cvbasic}, 
   \begin{equation}\begin{aligned}
\lim_{\e \to 0} \int  \bfu_\e(x ) \cdot \bfeta\lp x,  \frac{y_\e(x')}{r_\e}\rp   dm_\e(x)   = \frac{1}{|S|} \int_ {\Omega\times S} ( \bfv(x)+\frac{2}{\diam S}\theta(x)\bfe_3 \wedge \bfy)  \cdot \bfeta(x,y)  dx dy .
\end{aligned}
\nonumber
  \end{equation}
By varying  $\bfeta$ in $\D(\Omega\times S;\RR^3)$, taking  \eqref{resoft} and  \eqref{defme}   into account,  we infer that $\bfu_0(x,y)= \bfv(x)+\frac{2}{\diam S}\theta(x)\bfe_3 \wedge \bfy \ \hbox{ in } \Omega\times S$,
yielding, by \eqref{defcfsoft} and \eqref{Ldto},   
\begin{equation}\begin{aligned} 
&  \bfu_0(x,.)-\bfu(x)\in \W(\bfv(x)-\bfu(x), \theta(x)) \quad \hbox{ a.e. } \ x\in \Omega,
  \\&
  \int_{\Omega\times Y\setminus S} f^{\infty,p}(\bfe_y(\bfu_0)) dxdy \ge \int_\Omega c^f_{soft}(\bfv-\bfu, \theta) dx \qquad  (\bfe_y(\bfu(x))=0).
\end{aligned}
\label{u0inW}
  \end{equation}
By  \eqref{Ldto}  we have  $\e\bfe (\bfu_\e) \mathds1_{\Omega\setminus T_\e} \, \dto \,  \bfe_y(\bfu_0)\mathds1_{Y\setminus S}$
 (see \cite[Lemma  1]{BeSiam}).
In view of  \eqref{finftyfg0g} and 
  \eqref{linfcarre}  applied to "$\dto$", we infer 
  $ \liminf_{\e\to0} I_{\e1} \ge \int_{\Omega\times Y\setminus S} f^{\infty,p}(\bfe_y(\bfu_0)) dxdy$  and deduce from  \eqref{Phomsoft} and  \eqref{linf3} and \eqref{u0inW} that  
$$
\liminf_{\e\to 0} F_{\e}(\bfu_{\e}) \ge \Phi^{soft}(\bfu,\bfv^{tuple}).
$$

\noindent 
To prove the couterpart of Proposition \ref{propupperbound}, we fix  $\alpha>0$ and   $(\bfvarphi, \bfpsi, \zeta, a,b)$  satisfying  \eqref{psitupleinD},  \eqref{defpsituple}. Mimicking  the proof of Lemma \ref{lemapproxcap}, we obtain  (see \eqref{defcfsoft})
 \begin{equation}
 \begin{aligned}
&  \int_\Omega c^f_{soft}(\bfpsi- \bfvarphi, \zeta) dx= 
\inf_{\bfeta\in \bfgotW(\bfpsi-\bfvarphi,\zeta)}  \int_{\Omega\times (Y\setminus S)} \hskip-0,1cm f^{\infty,p}(\bfe_{y}(\bfeta(x,y))) dxdy,
\\& \bfgotW(\bfa,\zeta):= \la \bfeta\in C^1(\ov\Omega; C^\infty_\sharp(Y;\RR^3)), \ \bfeta(x,.)\in \W(\bfa(x),\zeta(x)) \ \forall x\in \Omega 
\ra, 
\end{aligned}
\nonumber
    \end{equation}
hence there exists    $\bfvarphi_0$ such that $\bfvarphi_0-\bfvarphi \in\bfgotW(\bfpsi-\bfvarphi,\zeta)$   and 
 \begin{equation}\begin{aligned}
&  \int_{\Omega\times (Y\setminus S)} f^{\infty,p}(\bfe_y(\bfvarphi_0)) dxdy \le \int_\Omega c^f_{soft}(\bfpsi- \bfvarphi, \zeta) dx
+\frac{\alpha}{2}. 
     \end{aligned}
\label{intint}
  \end{equation}
Setting $T_\e^\e:= \{ x\in \Omega,   \dist ( x, T_\e ) <\e^2 \}$,  we fix
 $(\xi_\e )\subset C^\infty(\ov\Omega)$ satisfying
\begin{equation}
\begin{aligned}
&
0\le \xi_\e\le 1, \quad \xi_\e =1 \hbox{ in } T_\e, \quad \xi_\e =0 \hbox{ in }  \Omega\setminus T_\e^\e, \quad |\bfnabla \xi_\e|\le \frac{C}{\e^2},
  \end{aligned} 
  \nonumber
  \end{equation} 
and  put
$$
\bfvarphi_\e(x):= 
\lp1-\xi_\e  \rp   \bfvarphi_0\lp x, \frac{y_\e(x')}{\e}\rp + \xi_\e \bfchi_\e,
$$
where     $\bfchi_\e$   is defined  
by (\ref{chik}),  (\ref{chikappa}) %
  or      \eqref{defchiothercases} 
 (depending on $k$, $\kappa$). 
  Writing 
\begin{equation}
\begin{aligned}
&F^{soft}_\e( \bfvarphi_\e)\le I_{\e1}+I_{\e2}+I_{\e3}+ I_{\e4}; \quad I_{\e2}:= \int_{T_\e^\e\setminus T_\e}  \e^p  f\lp \bfe\lp \bfvarphi_\e  \rp \rp dx;
\\& I_{\e1}:= \int_{ \Omega } \e^p f(\bfe(\bfvarphi_0( x, {y_\e(x')}/{\e}))) \mathds{1}_{Y\setminus S} ({y_\e(x')}/{\e})dx;
 \quad  I_{\e3}  :=k_\e \int g \lp   \bfe\lp \bfchi_\e \rp\rp dm_\e,
\\& I_{\e4}:=  \int_{ \Omega'_\e\times(0,L) } \e^p f(\bfe(\bfvarphi_0( x, {y_\e(x')}/{\e}))) dx; \quad \Omega'_\e:= \{x'\in\Omega', \dist(x',\partial\Omega')<\sqrt2\e\},
  \end{aligned} 
  \nonumber
     \end{equation} 
noticing that  
 $\bfe\lp\bfvarphi_0\lp x, \frac{y_\e(x')}{r_\e}\rp\rp = \frac{1}{\e} \bfe_y(\bfvarphi_0)\lp x, \frac{y_\e(x')}{r_\e}\rp  + \hbox{O}(1)$ and 
   $|\bfe (\bfvarphi_\e )| \le \frac{C}{\e}$ in $T^\e_\e\setminus T_\e$,  
taking  \eqref{finftyfg0g} and   \eqref{twoscalecontinuous}  applied to  $"\dto"$ 
into account, 
we deduce that 
\begin{equation}
\begin{aligned}
&\limsup_{\e\to0}   I_{\e1}\le  \int_{ \Omega\times (Y\setminus S)}f^{\infty,p}(\bfe_y(\bfvarphi_0)) dxdy ;\ 
   |I_{\e2}|\le  C |T_\e^\e\setminus T_\e|\le C\e;\  |I_{\e4}|\le C |\Omega_\e'|\le C\e,
  \end{aligned} 
  \nonumber
     \end{equation} 
yielding, along  with  \eqref{Phomsoft},   \eqref{limI3} and  \eqref{intint},  
$$
 \limsup_{\e\to 0}  F^{soft}_\e(\bfvarphi_\e)  
\leq \Phi^{soft}(\bfvarphi, \bfpsi^{tuple})
+\alpha.
$$
 We conclude  like   in Section  \ref{secproofth}. 
\qed
   \subsection{Proof of Theorem \ref{thL}}\label{secproofthL}
   \noindent  
{\it Problem \ref{Pefe}, case   $0<\kappa<+\infty$}.  The solution  $\bfu_\e$  to \eqref{Pefe} satisfies 
      \begin{equation}\begin{aligned}
       F_{\e}(\bfu_\e) &\le \int_\Omega \bff_\e\cdot \bfu_\e dx.
          \end{aligned} 
  \label{Ffu} 
    \end{equation}   
By    (\ref{fe2}), \eqref{defovgotA}, \eqref{defovgots}, we have
      \begin{equation}\begin{aligned}
 &    \int_\Omega \bff_\e\cdot \bfu_\e dx
      =\int_\Omega \bff\cdot  \bfu_\e dx 
\\&\qquad+ 
\int \bfg_0 \lp x,\frac{y_\e}{r_\e}\rp\cdot  \bfu_\e    +    \frac{u_{\e3}}{ r_\e}   a_0\lp x,\frac{y_{\e }}{r_\e } \rp  +     \frac{\ov\gots^S(\bfu_\e)}{r_\e}   \beta_0\lp x,\frac{y_{\e }}{r_\e } \rp   dm_\e,
          \end{aligned} 
  \label{hu} 
    \end{equation}   
thus,  by   \eqref{kkappa},    \eqref{estimDirichletb},   \eqref{votheta3} and the first inequality in \eqref{estimbasic1}, the following holds  
      \begin{equation}\begin{aligned}
   \hskip-0,2cm \lb \int_\Omega   \bff_\e\cdot \bfu_\e dx \rb^p & \le 
     C|\bfu_\e|^p_{L^p(\Omega;\RR^3)}+ C  \int  \hskip-0,1cm   |\bfu_\e|^p+ \lb\frac{u_{\e3}}{r_\e}\rb^p+ \lb\frac{\ov\gots^S(\bfu_\e)}{r_\e}\rb^p dm_\e 
   \\& \leq C   F_\e(\bfu_\e). 
       \end{aligned} 
  \label{fuF} 
    \end{equation}   
Combining  \eqref{Ffu} and \eqref{fuF}, we infer  \eqref{supFeuefini}, hence  the convergences stated in Propositions \ref{prop1}, \ref{prop2}, \ref{prop3} hold
for some subsequence $(\bfu_{\e_k})$. We deduce from  \eqref{L2},  \eqref{cvbasic}, \eqref{cvbasickappa>0},  \eqref{regthetav},   \eqref{cvkappa},  and  \eqref{hu},    that 
   \begin{equation}\begin{aligned} 
    \lim_{k\to+\infty}      \int    \bff_{\e_k} \cdot \bfu_{\e_k} dx =
    \int_\Omega\bff\cdot \bfu dx +L( \bfv^{tuple}), 
  \end{aligned} 
\nonumber
    \end{equation}   
 and then from 
 Proposition \ref{proplowerbound} that 
 \begin{equation}\hskip-0,1cm\begin{aligned}
\liminf_{k\to+\infty}F_{\e_k}(\bfu_{\e_k})-\hskip-0,1cm \int_\Omega\hskip-0,1cm \bff\cdot  \bfu_{\e_k} dx 
&\ge \Phi(\bfu,\bfv^{tuple})-\hskip-0,1cm\int_\Omega\hskip-0,1cm \bff\cdot \bfu  dx -L(\bfv^{tuple})
\\&\ge \inf \eqref{PhomL} .\hskip-0,1cm
\end{aligned}
\label{linfL}
  \end{equation}
By density,  there exists 
   $(\bfvarphi, \bfpsi^{tuple})$ satisfying \eqref{psitupleinD},  \eqref{defpsituple} 
     and  
    \begin{equation}\begin{aligned}
&
  \Phi(\bfvarphi, \bfpsi^{tuple})-\int_\Omega \bff\cdot  \bfvarphi dx -L(  \bfpsi^{tuple})  \leq \inf \eqref{PhomL}+\alpha,
\end{aligned}
\nonumber
  \end{equation}
and, by  Proposition \ref{propupperbound}    a sequence   $(\bfvarphi_\e)$  verifying  \eqref{ub} and  \eqref{lsupF}. By \eqref{L2},   \eqref{ub} and \eqref{hu}, 
  $\lim_{\e\to0} \int_\Omega \bfvarphi_\e\cdot\bff_\e dx = \int_\Omega \bfvarphi\cdot\bff dx + L(\bfpsi^{tuple})$, hence, by \eqref{lsupF} , 
   \begin{equation}\begin{aligned}
  \limsup_{\e\to0} F_\e(\bfu_\e)-\int_\Omega& \bff_\e\cdot   \bfu_\e dx  
\leq \limsup_{\e\to 0}   F_\e(\bfvarphi_\e)\hskip-0,1cm-\hskip-0,1cm\int_\Omega\hskip-0,1cm  \bff_\e\hskip-0,1cm\cdot\hskip-0,05cm\bfvarphi_\e dx 
\\&\hskip-2cm 
 \leq \Phi(\bfvarphi, \bfpsi^{tuple})
+\alpha  -\hskip-0,1cm\int_\Omega\hskip-0,1cm \bff\hskip-0,1cm\cdot\hskip-0,05cm\bfvarphi  dx -L(\bfpsi^{tuple})
 \leq \inf \eqref{PhomL} +2\alpha.
 \end{aligned}
\nonumber
  \end{equation}
This, along with  \eqref{linfL}, shows, by the arbitrariness of $\alpha$, that $(\bfu,\bfv^{tuple})$ is a solution to Problem \ref{PhomL}.
The case $\kappa=0$ can be  proved in a similar manner,  by using   the weak$^\star$ convergence   in $\M(\ov\Omega)$,    
stated in  \eqref{cvbasic},  of $ \lp -\frac{y_{\e2}}{r_\e} u_{\e1} +\frac{y_{\e1}}{r_\e} u_{\e2}\rp m_\e$  to $\tau\theta$.
  \qed
  
     \noindent  
{\it Problem \ref{Pesoftfe}}. 
Noticing that the solution $\bfu_\e$ to \eqref{Pesoftfe}   satisfies 
\eqref{Ffu}, \eqref{fuF}  with $F_\e^{soft}$ in place of $F_\e$  (same proof), hence, by 
  \eqref{condsoft},  is bounded in $L^p(\Omega;\RR^3)$, we  
 adapt   the proof of Theorem \ref{thsoft}  in the same way as  above.\qed

\section{Appendix}\label{secappendix}
%
 \subsection{Technical lemmas related to  the lower bound}\label{technicallemmas}
In this section, given a sequence $(\bfu_\e)$ satisfying \eqref{supFeuefini}, 
we establish the existence of 
   sequences $(\widetriangle\bfu_\e)$ and $(\wideparen\bfu_\e)$  verifying 
\eqref{cvtriangle},    \eqref{basedonI2ep<2}   and 
    \eqref{basedonI2ep=2}.
\noindent\begin{lemma}\label{lemtriangle} Assume  that   $p\le2$,   $(R_\e)$       verifies  \eqref{Re}, and 
   $(\bfu_\e)$ 
 satisfies (\ref{supFeuefini}) and \eqref{cvbasic}.
Let    $(I_{\e2})$ be  defined by \eqref{splitlowerbound}. There exists  a  sequence  $ (\widetriangle\bfu_\e)$  in $L^p(0,L; W^{1,p}(\Omega';$ $\RR^3))$  such that   
\begin{equation}
\begin{aligned}
& \ov\bfgotu_\e(\widetriangle\bfu_\e) \rightharpoonup \bfu, \quad \ov\bfgotv_\e^S(\widetriangle\bfu_\e) \rightharpoonup\bfv,  \quad  \ov\theta_\e^S(\widetriangle \bfu_\e)\bfe_3 \rightharpoonup \theta\bfe_3 \quad \hbox{weakly in } L^p(\Omega;\RR^3), 
\\&
\liminf_{\e\to0}\  I_{\e2}
  \geq \liminf_{\e\to0}\int_{{(D_{R_\e}\setminus S_{r_\e})\times (0,L)}}  f^{\infty, p} (\bfe_{x'}( \widetriangle \bfu_{\e  }))\ dx,
 \\& \int_{{E_\e \times (0,L)}}\hskip-0,2cm | \bfe_{x'}( \widetriangle \bfu_{\e  })|^p\ dx \leq  \int_{{E_\e \times (0,L)}} \hskip-0,2cm| \bfe (   \bfu_{\e  })|^p\ dx +o(1)   \quad \forall E_\e \subset \B(D_{R_\e}).
\end{aligned}
\label{triangle}
\end{equation} 
 \end{lemma}

\noindent
\begin{proof}   To construct $\widetriangle\bfu_\e$,  we fix     two sequences     $(a_\e)$ and $(b_\e)$   such that 
\begin{equation}
\begin{aligned}
&   1 \gg  a_\e \gg b_\e,     \qquad a_\e b_\e^2 \gg \frac{R_\e^2}{ \e^2},
  \end{aligned} 
\label{ab}    \end{equation} 
and  choose
a family  $(l_{k,\e})_{\e>0, \ k\in \{1,..., n_\e\}}$      verifying 
%
\begin{equation}
\begin{aligned}
&n_\e:=\lc\frac{L}{a_\e}\rc-1, \quad \frac{3}{4} a_\e \leq   l_{1,\e}\leq   \frac{5}{4}  a_\e,
  \quad 
L-\frac{9}{4} a_\e < l_{n_\e,\e}< L- \frac{3}{4}a_\e,
 \\& 
 l_{k,\e}=l_{1,\e}+ (k-1) a_\e
\qquad\forall k\in \{2,..,n_\e\},
  \end{aligned} 
\label{condl}    \end{equation}
 and such that, setting 
\begin{equation}
\begin{aligned}
&   H_\e :=   D_{R_\e} \times   
  \bigcup_{k=1}^{ n_\e  }    \left( l_{k,\e}-\frac{b_\e}{2};  l_{k,\e}+\frac{b_\e}{2}\right)   ,
  \end{aligned} 
\label{defHe}    \end{equation}
the following holds:
\begin{equation}
\begin{aligned}
&\int_{H_\e} |\bfe(\bfu_\e)|^p dx+  \int_{H_\e} |\ov\bfgotv_\e^S(\bfu_\e)|^p +  |\ov\bfgotu_\e(\bfu_\e)|^p + |\ov\theta_\e^S(\bfu_\e)|^p d m_\e  
\\&\quad \leq C\frac{b_\e}{a_\e}\lp  \int_{\Omega} |\bfe( \bfu_\e)|^p dx+ \int |\ov\bfgotv_\e^S(\bfu_\e)|^p +  |\ov\bfgotu_\e(\bfu_\e)|^p + |\ov\theta_\e^S(\bfu_\e)|^p d m_\e \rp.
  \end{aligned} 
\label{De}    \end{equation}
To prove  that this choice is possible, 
 we set, for every
  $m\in \la 0,... ,\lc \frac{a_\e}{2 b_\e}\rc \ra$  
\begin{equation}
\begin{aligned}
 & H_\e^m:=   D_{R_\e}\times 
 \lc \bigcup_{k=1}^{n_\e} \lp l_{k,\e}^m-\frac{1}{2} b_\e;  l_{k,\e}^m+\frac{1}{2} b_\e\rp \rc;
\quad  l_{k,\e}^m:= \lp k-\frac{1}{4}\rp a_\e + m b_\e.
  \end{aligned} 
\nonumber   \end{equation}
Since $H_\e^m\cap H_\e^n=\emptyset$ if $m\not= n$,  we have
\begin{equation}
\begin{aligned}
  \sum_{m=0}^{\lc \frac{a_\e}{2 b_\e}\rc}  \int_{H^m_\e} &|\bfe(\bfu_\e)|^p dx +  \int_{H^m_\e} |\ov\bfgotv_\e^S(\bfu_\e)|^p +  |\ov\bfgotu_\e(\bfu_\e)|^p + |\ov\theta_\e^S(\bfu_\e)|^p d m_\e  
\\&\quad \leq   \int_{\Omega} |\bfe( \bfu_\e)|^p dx+ \int |\ov\bfgotv_\e^S(\bfu_\e)|^p +  |\ov\bfgotu_\e(\bfu_\e)|^p + |\ov\theta_\e^S(\bfu_\e)|^p d m_\e.
  \end{aligned} 
\nonumber   \end{equation}
 Accordingly, for each $\e>0$,  there exists $\check m_\e\in \la 0,... ,\lc \frac{a_\e}{2 b_\e}\rc \ra$ such that 
\begin{equation}
\begin{aligned}
&\int_{H^{\check m_\e}_\e}  |\bfe(\bfu_\e)|^p dx +  \int_{H^{\check m_\e}_\e} |\ov\bfgotv_\e^S(\bfu_\e)|^p +  |\ov\bfgotu_\e(\bfu_\e)|^p + |\ov\theta_\e^S(\bfu_\e)|^p d m_\e  
\\&  \leq  \frac{1}{\lc \frac{a_\e}{2 b_\e}\rc+1} \lp\int_{\Omega} |\bfe( \bfu_\e)|^p dx+ \int |\ov\bfgotv_\e^S(\bfu_\e)|^p +  |\ov\bfgotu_\e(\bfu_\e)|^p + |\ov\theta_\e^S(\bfu_\e)|^p d m_\e\rp.
  \end{aligned} 
\label{satme}   \end{equation}
Therefore
$(l_{k,\e})_{k\in \{1,..., n_\e \}}$  defined by $l_{k,\e}:= l_{k,\e}^{\check m_\e}$ 
satisfies \eqref{condl}  and 
  \eqref{De}.
Let us  fix   $(\rho_\e)$ verifying \eqref{defrhoe}.  
To each  sequence $(\bfpsi_\e)$ in $L^p(\Omega;\RR^3)$, we  associate   the sequence $(\widetriangle\bfpsi_\e)$  defined by 
\begin{equation}
 \hskip-0,2cm\begin{aligned}
& \widetriangle{\bfpsi}_\e (x):=    \hskip-0,1cm \sum_{k=2}^{ n_\e}      \left(   \intb_{ \hskip-0,1cmJ_\e^k}   \rho_\e(s_3)   \bfpsi_\e  (x',s_3)  ds_3 \right)   \mathds{1}_{ J_\e^k}(x_3); \quad  J_\e^k:= (l_{k-1 ,\e}; l_{k,\e}). 
  \end{aligned} 
\label{deftriangleue}    \end{equation}
Since $0\le \rho_\e\le1$, for every Borel subset $E$ of $\Omega'$,  we have 
\begin{equation}
\begin{aligned}
\int_{E\times(0,L)}  | \widetriangle \bfpsi_\e |^p dx   &= 
\sum_{k=2}^{ n_\e} \int_E  dx'  \int_{J_\e^k} dx_3\lb \intb_{  J_\e^k}   \rho_\e(s_3)   \bfpsi_\e  (x',s_3)  ds_3\rb^p 
\\&\leq \int_{E\times(0,L)} \lb  \bfpsi_\e \rb^p dx. 
  \end{aligned} 
\label{inttrileint}    \end{equation}
By  \eqref{defovgotue} and \eqref{deftriangleue}, the following holds
\begin{equation}
\begin{aligned}
 &\hskip-0,1cm \ov\bfgotu_\e  \lp\widetriangle{\bfpsi}_\e\rp (x)   = 
\sum_{i\in I_\e} \left (
\intb_{D_{R_\e}^{i}\setminus D_{R_\e\!/2}^{i}}\!\!\widetriangle{\bfpsi}_\e(s', x_3  ) ds' \right )  \mathds{1}_{Y_\e^{i}}(x' )  
\\&\hskip-0,1cm=  \sum_{i\in I_\e}  \sum_{k=2}^{ n_\e}    
  \left(  \intb_{D_{R_\e}^{i}\setminus D_{R_\e\!/2}^{i}}\!\! \lp\intb_{J_\e^k}   \rho_\e(s_3)  \bfpsi_\e  (s',s_3)  ds_3\rp ds' \right)   \mathds{1}_{Y_\e^{i}\times  J_\e^k}(x)
  \\&\hskip-0,1cm=   \hskip-0,1cm \sum_{k=2}^{ n_\e} \hskip-0,1cm      \left( \hskip-0,1cm  \intb_{\hskip-0,1cm J_\e^k} \hskip-0,2cm  \rho_\e(s_3) \hskip-0,1cm \lp
   \sum_{i\in I_\e} \hskip-0,1cm \intb_{\hskip-0,1cmD_{R_\e}^{i}\setminus D_{R_\e\!/2}^{i}}\!\! \hskip-0,2cm{\bfpsi}_\e(s', s_3  ) ds'    \mathds{1}_{Y_\e^{i}}(x')\hskip-0,1cm \rp \hskip-0,1cm ds_3\hskip-0,1cm
    \right)\hskip-0,1cm \mathds{1}_{ J_\e^k}(x_3) \hskip-0,1cm 
 = \widetriangle{ \ov\bfgotu_\e(\bfpsi_\e)}(x). 
  \end{aligned} 
  \nonumber
   \end{equation}
%
 Likewise,
 \begin{equation}
\begin{aligned}
 & \ov\bfgotu_\e \lp\widetriangle{\bfpsi}_\e\rp   
 = \widetriangle{ \ov\bfgotu_\e(\bfpsi_\e)}; 
 \qquad \ov\bfgotv_\e^S \lp\widetriangle{\bfpsi}_\e\rp   
 = \widetriangle{ \ov\bfgotv_\e^S(\bfpsi_\e)}; \qquad  \ov\theta_\e^S \lp\widetriangle{\bfpsi}_\e\rp    
 = \widetriangle{ \ov\theta_\e^S(\bfpsi_\e)},
   \end{aligned} 
\label{gottri=trigot}    \end{equation}
yielding,  by  \eqref{mdx} and  \eqref{inttrileint},  
\begin{equation}
 \hskip-0,1cm \begin{aligned}
&  \int  | \ov\bfgotv_\e^S(\widetriangle\bfpsi_\e)|^p dm_\e  \hskip-0,1cm = \hskip-0,1cm \int  \lb\widetriangle{ \ov\bfgotv_\e^S(\bfpsi_\e)}\rb^p dm_\e  \hskip-0,1cm\leq  \hskip-0,1cm \int  | \ov\bfgotv_\e^S(\bfpsi_\e)|^p dm_\e \hskip-0,1cm \leq \hskip-0,1cm \int  |  \bfpsi_\e|^p dm_\e,
\\&
\int |\ov\bfgotu_\e\lp\widetriangle{\bfpsi}_\e\rp|^p dm_\e\le
C \int  |  \bfpsi_\e|^p dm_\e
;
\quad \int  |\ov\theta_\e^S\lp\widetriangle{\bfpsi}_\e\rp|^p dm_\e  
\leq C \int  |  \bfpsi_\e|^p dm_\e.\hskip-0,2cm 
\end{aligned}
\label{sectri1}
\end{equation}
 Since   $\int  |  \bfu_\e|^p dm_\e \leq C$ (see \eqref{uemebounded}),
  by \eqref{sectri1} and  Lemma \ref{lemtwoscale}  
 the following   holds,  up to a subsequence,  for some    $(\bfg,\bfh,\eta)\in (L^p(\Omega;\RR^3))^2\times L^p(\Omega)$:
  \begin{equation}
\begin{aligned}
&\ov\bfgotu_\e(\widetriangle \bfu_\e) m_\e    \buildrel \star \over \rightharpoonup \bfg {\calL}^3_{\lfloor \Omega} , \quad 
\ov\bfgotv_\e^S(\widetriangle  \bfu_\e) m_\e    \buildrel \star \over \rightharpoonup \bfh {\calL}^3_{\lfloor \Omega} \quad &&\text{ weakly$^\star$ in }  \M(\overline\Omega;\RR^3),
\\& \ov\theta_\e^S(\widetriangle\bfu_\e) m_\e  \buildrel \star \over \rightharpoonup \eta {\calL}^3_{\lfloor \Omega} \quad &&\text{ weakly$^\star$ in }  \M(\overline\Omega).
\end{aligned}
\label{triuemuetown}
\end{equation}
 Taking \eqref{mdx2}  into account,  the convergences stated  in \eqref{triangle} are proved provided we establish that 
 $(\bfg,\bfh,\eta)=(\bfu,\bfv,\theta)$.  
To that aim, we  fix   $\bfxi\in \D(\Omega;\RR^3)$ and set 
$   \bfxi_\e(x):=    \sum_{k=2}^{ n_\e}      \left(   \intb_{J_\e^k}    \bfxi  (x',s_3)  ds_3 \right)   \mathds{1}_{ J_\e^k}(x_3)$. 
 By   \eqref{defovgotA} and    \eqref{deftriangleue}, we have
 \begin{equation}
\begin{aligned}
&  \int    \bfxi_\e \cdot \rho_\e \ov\bfgotv_\e^S(\bfu_\e)d m_\e   
\\&=  \frac{\e^2}{r_\e^2|S|} \sum_{i\in  {I}_\e} \sum_{k=2}^{n_\e}\int_{S_{r_\e}^{i}  \times J_\e^k} \hskip-0,8cm
 \rho_\e (x_3) \lp \intb_{J_\e^k} \bfxi(x',s_3) ds_3 \rp\cdot
 \lp \intb_{S_{r_\e}^{i} } \bfu_\e(s',x_3) ds'\rp dx
 \\&= \frac{\e^2}{r_\e^2|S|}  \sum_{i\in  {I}_\e}\hskip-0,05cm  \sum_{k=2}^{n_\e}\hskip-0,05cm\int_{S_{r_\e}^{i}  \times J_\e^k} \hskip-0,8cm \bfxi(x',s_3) \cdot  \lp \intb_{S_{r_\e}^{i} } \hskip-0,2cm
 \lp \intb_{J_\e^k}  \rho_\e (x_3)  \bfu_\e(s',x_3)  dx_3\hskip-0,1cm\rp 
\hskip-0,1cm ds'\hskip-0,1cm\rp dx' ds_3
\\ &=  \int   \bfxi\cdot\ov\bfgotv_\e^S(\widetriangle \bfu_\e)   d m_\e. 
  \end{aligned} 
\label{fivtrianglev} \end{equation}
By  \eqref{mdx2} and \eqref{cvbasic},   $\lp\ov\bfgotv_\e^S(\bfu_\e) m_\e  \rp$  weakly$^\star$ converges  to $\bfv$ in $\M(\ov\Omega;\RR^3)$ and, 
by  \eqref{estimbasic3},   \eqref{defrhoe},   \eqref{ab}, \eqref{defHe},  and    \eqref{De},   
 \begin{equation}
\begin{aligned}
\hskip-0,1cm  \int     |(1-\rho_\e )\ov\bfgotv_\e^S(\bfu_\e)|^p   d m_\e  
\leq  \int_{H_\e}    |\ov\bfgotv_\e^S(\bfu_\e)|^p   d m_\e 
  \leq C \frac{b_\e}{a_\e} =o(1),
\end{aligned}
\nonumber
\end{equation}
therefore  $\lp\rho_\e  \ov\bfgotv_\e^S(\bfu_\e) m_\e  \rp$  weakly$^\star$ converges  to $\bfv$ in $\M(\ov\Omega;\RR^3)$.
Observing that  $|\bfxi-  \bfxi_\e|_{L^\infty(\Omega;\RR^3)}\leq C a_\e \ll 1$,
we deduce from \eqref{triuemuetown} and \eqref{fivtrianglev} that
 \begin{equation}
\begin{aligned}
\int_\Omega \bfxi\cdot\bfv dx= \lime \int \bfxi_\e\cdot \rho_\e  \ov\bfgotv_\e^S(\bfu_\e) d m_\e = \lime  \int   \bfxi\cdot \ov\bfgotv_\e^S(\widetriangle \bfu_\e)   d m_\e = \int_\Omega \bfxi\cdot\bfh dx,
\end{aligned}
\nonumber 
\end{equation}
\noindent  and infer from the arbitrariness of $\bfxi$ that  $\bfh=\bfv$.  We  likewise  obtain    $\eta=\theta$ and $\bfg=\bfu$. The convergences stated in    (\ref{triangle}) are  proved. 
Let   $E_\e$ be a Borel subset of $D_{R_\e}$.
By  (\ref{finftyfg0g}),   (\ref{estimbasic1}),   and the estimate  $|D_{R_\e}|\leq C\frac{R_\e^2}{\e^2}$, we have
\begin{equation}\label{e316}
\begin{aligned}
&   \int_{E_\e\times (0,L)}  
\lb f^{\infty, p} (\bfe(\bfu_\e)) -  f(\bfe (\bfu_\e))\rb\ dx  \leq C  \int_{E_\e\times (0,L)}  
\lp1+  |\bfe (\bfu_\e)|^{p-\varsigma}\rp\ dx  
 \\& \hskip3cm \leq  C\lp \frac{R_\e^2}{\e^2}+ \lp \frac{R_\e^2}{\e^2}\rp^{\frac{\varsigma}{p}}
\vert\bfe(\bfu_\e)\vert ^{\frac{p-\varsigma}{p}}_{L^p(\Omega;\SSym^3)}\rp
 \leq C \lp \frac{R_\e^2}{\e^2}\rp^{\frac{\varsigma}{p}}. 
\end{aligned}
\end{equation}
By     \eqref{defrhoe} and  \eqref{defHe},     $\rho_\e (x_3)   \bfu_\e(x)=    \bfu_\e(x)$ in $(E_\e\times (0,L))\setminus H_\e$ and    $|\rho_\e'|<\frac{C}{b_\e}$,
hence by 
  (\ref{growthp}),   \eqref{estimbasic1},   \eqref{ab}, and  (\ref{De}),  
\begin{equation}
 \begin{aligned}  
   \int_{E_\e\times (0,L)}
    & 
    |f^{\infty, p}(\bfe(\bfu_\e)) -f^{\infty, p}(\bfe (\rho_\e (x_3)   \bfu_\e(x))) |dx
  \\&\le\int_{H_\e} f^{\infty, p}(\bfe(\bfu_\e) ) + f^{\infty, p}(\bfe(\rho_\e (x_3)  \bfu_\e(x)))  dx
\\&\leq C \int_{H_\e}
|\bfe(\bfu_\e) |^p+  \left|     \frac{\bfu_\e(x)}{b_\e}\right|^p    dx
 \leq  C  \int_{H_\e} \left|     \frac{\bfu_\e(x)}{b_\e}\right|^p    dx +o(1).
 \end{aligned}    
\label{hatue0}  
 \end{equation}
It follows  from \eqref{estimbasic1},  H\"older's inequality, \eqref{ab},  \eqref{defHe},  and  the continuous embedding of $W^{1,p}(\Omega;\RR^3)$    into   $L^{\frac{3p}{3-p}}(\Omega;\RR^3)$ 
that
 \begin{equation}
\begin{aligned}  
     \int_{H_\e}
  \left|     \frac{\bfu_\e(x)}{b_\e}\right|^p    dx 
&\leq    \frac{1}{ b_\e^p} \left(  \int_{H_\e}   \left|     \bfu_\e(x) \right|^{\frac{3p}{3-p}}    dx \right) ^{\frac{3-p}{3}} {\calL}^3(H_\e)^{ \frac{p}{3}} 
\\& \hskip-0,5cm\leq    \frac{C}{b_\e^p}  |  \bfu_\e|_{W^{1,p}(\Omega;\RR^3)}^p   \lp\frac{R_\e^2}{\e^2}\frac{b_\e}{a_\e  } \rp^{^{  \frac{p}{3} }} 
 \leq   C\left( \frac{R_\e^2}{ \e^2}\frac{1}{ a_\e b_\e^{2}}\right)^{  \frac{p}{3} } =o(1).
 \end{aligned}    
\label{hatue3}  
 \end{equation}
 Since   $\rho_\e =0$ on $ \bigcup_{k=1}^{n_\e}  \{ l_{k,\e}\}$ (see \eqref{defrhoe}), we have $   \intb_{J_\e^k} \frac{\partial}{\partial s_3}( \rho_\e (s_3)  \bfu_\e(x',s_3) )\,  ds_3 =0$ a. e. $x'\in \Omega'$,
$ \forall k\in \{2,..., n_\e\}$. 
Hence,    for a. e. $ x\in \Omega'\times J_\e^k$, there holds   (see  \eqref{defexprim}, \eqref{deftriangleue}) 
\begin{equation}
\begin{aligned}
\!\intb_{J_\e^k}  \! \bfe ( \rho_\e (s_3)  \bfu_\e(x',s_3)) \,  ds_3& = \!    \intb_{J_\e^k}  \!\bfe_{x'}\lp  \rho_\e (s_3)   \bfu_\e(x',s_3))  \rp  ds_3 
\\&=
  \!   \bfe_{x'}\lp  \intb_{J_\e^k}  \! \rho_\e (s_3)   \bfu_\e(x',s_3))    ds_3\rp
 = 	 \bfe_{x'}( \widetriangle\bfu_\e )(x) .  
  \end{aligned} 
\nonumber
 \end{equation}
Applying Jensen's inequality, we  infer 
\begin{equation}
\begin{aligned}
 &  \int_{ E_\e\times (0,L)} \hskip-1cm  f^{\infty, p}(\bfe(\rho_\e (s_3)  \bfu_\e(x',s_3))) \, dx'ds_3
   =\sum_{k=2}^{n_\e}\int_{E_\e} \hskip-0,3cm  dx'   \hskip-0,15cm    \int_{J_\e^k}   \hskip-0,1cm\hskip-0,1cm  f^{\infty, p}(\bfe(\rho_\e (s_3)   \bfu_\e(x',s_3))) \,  ds_3 
\\&\qquad \ge  \sum_{k=2}^{n_\e} \int_{E_\e} dx'     (l_{k  ,\e}- l_{k-1,\e}) f^{\infty, p}\left(  \intb_{J_\e^k} \bfe(\rho_\e (x_3)   \bfu_\e(x',s_3)) \,  ds_3 \right)  
\\&\qquad=   \sum_{k=2}^{n_\e} \int_{E_\e}   dx'     
\int_{J_\e^k}  f^{\infty, p}(\bfe_{x'}(\widetriangle\bfu_\e )(x) )  dx_3
 =  \int_{E_\e\times (0,L)}  f^{\infty, p}(\bfe_{x'}(\widetriangle\bfu_\e  )) \, dx, 
  \end{aligned} 
\nonumber
 \end{equation}
which, combined with (\ref{e316}),  \eqref{hatue0},  and (\ref{hatue3}), yields
\begin{equation}
\begin{aligned}
  \int_{E_\e\times (0,L)}  f^{\infty, p}(\bfe_{x'}(\widetriangle\bfu_\e  )) \, dx
 \leq    \int_{ E_\e\times (0,L)}  f(\bfe(   \bfu_\e )) \, dx +o(1).
  \end{aligned} 
\nonumber   \end{equation}
Choosing first $E_\e= D_{R_\e}\setminus S_{r_\e}$, taking \eqref{splitlowerbound} into account, 
we deduce 
the lower bound  stated  in the second line of  (\ref{triangle}). 
Next, substituting $|.|^p$ for $f$ in the above argument, we obtain   
   the    third line of \eqref{triangle}. \qed
\end{proof}
 \noindent 
We check below the  existence of a sequence  $(R_\e)$  
verifying (\ref{Re})   and  \eqref{slicing}.
\begin{lemma}\label{lemslicing}  Assume that 
$(\bfu_\e)$ verifies  \eqref{supFeuefini}. 
For any sequence  $(R'_\e)$   satisfying \eqref{Re},  there exists    $(R_\e)$  
verifying (\ref{Re})   and  \eqref{slicing}  such that $0<R_\e \leq R_\e'$ \ $\forall\e>0$.
  \end{lemma}

\noindent
{\bf Proof.}   Let  $(R'_\e)$ be   satisfying \eqref{Re}.
There exists  a sequence   of positive integers    $(n_\e)$  such that $\lime n_\e=+\infty$ and  $(2^{-n_\e} R_\e')$ verifies \eqref{Re}.
By \eqref{supFeuefini}, we have 
\begin{equation}
\begin{aligned}
 \sum_{m=0}^{n_\e-1}  \int_{D_{  2^{-m} R_\e' }\setminus D_{2^{-m-1}R_\e'} \times (0,L)} |\bfnabla   \bfu_\e  |^p dx & \leq  \int_{\Omega} |\bfnabla   \bfu_\e  |^p dx
\leq C ,   \end{aligned} 
\nonumber \end{equation} 
 hence  there exists     $ m_\e \in\{0,..,n_\e-1\}$ such that
$$
  \int_{D_{  2^{-m_\e} R_\e'  }\setminus D_{2^{-m_\e-1}R_\e'} \times (0,L)} |\bfnabla   \bfu_\e  |^p dx   \leq   \frac{ C }{ n_\e} .
  $$ 
We set  $R_\e  := 2^{-m_\e}R_\e'$.  
\qed
  %
%
 \begin{lemma}  \label{lemparen}  Assume that $p\le2$,  
   $(\bfu_\e)$ 
 satisfies (\ref{supFeuefini}),  $(R_\e)$       verifies  \eqref{Re} and \eqref{slicing}, and   
 the convergences  \eqref{cvbasic}  hold. Let 
 $(\widetriangle \bfu_\e)$  given by Lemma \ref{lemtriangle} and 
 $S'$     a Lipschitz  domain verifying \eqref{S'}.
Then, there exists   $(\wideparen\bfu_\e)\subset L^p(0,L; W^{1,p}(\Omega';\RR^3))$   satisfying 
  \eqref{basedonI2ep<2} and \eqref{basedonI2ep=2}.
 \end{lemma}

\noindent 
{\bf Proof.}     We  fix  
   $\psi_\e, \zeta_\e \in \D(\Omega')$,  
such that 
 \begin{equation}
\hskip-0,2cm \begin{aligned}
&\spt \psi_\e  \subset D_{2R_\e}\hskip-0,1cm\setminus D_{\frac{1}{2}R_\e}  ;\ 
    \psi_\e = 1 \   \  \text{on} \ \partial D_{R_\e}; \   \  0\leq \psi_\e\leq 1; 
\   \   |\nabla \psi_\e |\leq \frac{C}{R_\e }, 
\\&\spt \zeta_\e \subset S_{r_\e};
\qquad \zeta_\e = 1 \quad \text{ in  }  S_{r_\e}' ;\qquad 0\leq \zeta_\e\leq 1;
\qquad  |\nabla \zeta_\e |\leq \frac{C}{r_\e},
  \end{aligned} 
\label{filower}    \end{equation} 
 and set   (see (\ref{defovgotA}), \eqref{deftriangleue})
 \begin{equation}
\begin{aligned}
&   \wideparen \bfu_\e := \widetriangle\bfu_\e + \psi_\e (\ov\bfgotu_\e(\widetriangle\bfu_\e)-\widetriangle\bfu_\e) + \zeta_\e 
 \lp\ov\bfgotr^S_\e(\widetriangle\bfu_\e)   -\widetriangle\bfu_\e\rp & & \hbox{if } \ p<2,
 \\&    \wideparen \bfu_\e := \widetriangle\bfu_\e + \psi_\e (\ov\bfgotu_\e(\widetriangle\bfu_\e)-\widetriangle\bfu_\e) + \zeta_\e 
 \lp\ov\bfgotv^S_\e(\widetriangle\bfu_\e)   -\widetriangle\bfu_\e\rp \qquad & & \hbox{if } \ p =2.
  \end{aligned} 
\label{wideparrenue=}\end{equation} 
The boundary conditions stated in  (\ref{basedonI2ep<2}),  \eqref{basedonI2ep=2} are satisfied.  By  \eqref{inconv},  we have
 \begin{equation}
\begin{aligned}
&\Bigg| \int_{{(D_{R_\e}\setminus S_{r_\e})\times (0,L)}}   f^{\infty, p} (\bfe_{x'}( \widetriangle \bfu_{\e  }))\ dx
 -\int_{(D_{R_\e}\setminus S_{r_\e}')\times(0,L)}  f^{\infty, p}(\bfe_{x'}(\wideparen  \bfu_\e)) dx \Bigg|
\\&\leq C \int_{(D_{R_\e}\setminus S_{r_\e})\times (0,L)}  |\bfe_{x'}(\widetriangle\bfu_\e-\wideparen \bfu_\e)|\lp 1+  |\bfe_{x'}(\widetriangle\bfu_\e )|^{p-1}+ |\bfe_{x'}( \wideparen \bfu_\e)|^{p-1}\rp dx 
\\& \hskip5cm+  \int_{(S_{r_\e}\setminus S_{r_\e}')\times (0,L)  } f^{\infty, p} (\bfe_{x'}(\wideparen  \bfu_\e)) dx.
  \end{aligned} 
\nonumber  \end{equation} 
We prove below that 
\begin{eqnarray}
& &  \hskip-0,3cm \int_{(D_{R_\e}\setminus S_{r_\e})\times (0,L)} \hskip-0,9cm  |\bfe_{x'}(\widetriangle\bfu_\e-\wideparen \bfu_\e)|\lp 1+  |\bfe_{x'}(\widetriangle\bfu_\e )|^{p-1}\hskip-0,1cm + |\bfe_{x'}( \wideparen \bfu_\e)|^{p-1}\rp dx 
=o(1),\hskip-0,2cm 
\label{bound1}    
\\
& &  \hskip-0,3cm \int_{(S_{r_\e}\setminus S_{r_\e}')\times (0,L)  } f^{\infty, p} (\bfe_{x'}(\wideparen  \bfu_\e)) dx=o(1) \hskip2,25cm\hbox{if } 1<p<2,
\label{bound2}  
\\
& &  \hskip-0,3cm \int_{(S_{r_\e}\setminus S_{r_\e}')\times (0,L)  } f^{\infty, p} (\bfe_{x'}(\wideparen  \bfu_\e)) dx\leq C   \frac{|S\setminus S'|}{|S|}  \hskip1,5cm\hbox{if } p=2.
\label{bound3}  
\end{eqnarray}
Taking  \eqref{triangle} into account,  Lemma \ref{lemparen} is proved. \qed

\noindent {\it Proof of \eqref{bound1}.} 
By  \eqref{estimbasic1} and the   last line of \eqref{triangle},   the following holds
 \begin{equation}
\begin{aligned}
 \int_{(D_{R_\e}\setminus S_{r_\e})\times (0,L)}  |\bfe_{x'}(\widetriangle\bfu_\e)|^pdx\leq C \int_{(D_{R_\e}\setminus S_{r_\e})\times (0,L)}  |\bfe ( \bfu_\e)|^pdx+o(1)\leq C.
 \end{aligned} 
\label{applig00}    \end{equation} 
We deduce
 \begin{equation}
\begin{aligned}
  & \int_{(D_{R_\e}\setminus S_{r_\e})\times (0,L) }   |\bfe_{x'}(\widetriangle\bfu_\e -\wideparen \bfu_\e)| \lp 1+  |\bfe_{x'}(\widetriangle\bfu_\e )|^{p-1}+ |\bfe_{x'}( \wideparen \bfu_\e)|^{p-1}\rp dx  
\\& \le C  \int_{(D_{R_\e}\setminus S_{r_\e})\times (0,L) } \hskip-1cm |\bfe_{x'}(\widetriangle\bfu_\e-\wideparen \bfu_\e)| \lp 1+  |\bfe_{x'}(\widetriangle\bfu_\e )|^{p-1}+ |\bfe_{x'}(\wideparen \bfu_\e-\widetriangle\bfu_\e)|^{p-1}\rp dx  
\\&
\leq C\lp\int_{(D_{R_\e}\setminus S_{r_\e})\times (0,L) } \hskip-1cm |\bfe_{x'}(\widetriangle\bfu_\e-\wideparen \bfu_\e)|^p dx\rp^{\frac{1}{p}} +C    \int_{(D_{R_\e}\setminus S_{r_\e})\times (0,L) } 
\hskip-1cm |\bfe_{x'}(\widetriangle\bfu_\e-\wideparen \bfu_\e)|^p dx. 
  \end{aligned} 
\label{applig22}    \end{equation} 
By \eqref{defexprim} and  (\ref{defovgotue})  we have    $\bfe_{x'}(\ov\bfgotu_\e(\widetriangle\bfu_\e))=0$ in $D_{R_\e}\times (0,L)$ and, by  \eqref{gottri=trigot},    (\ref{filower})  and  (\ref{wideparrenue=}), 
$\widetriangle\bfu_\e- \wideparen\bfu_\e=  \psi_\e (\widetriangle\bfu_\e-\ov\bfgotu_\e(\widetriangle\bfu_\e))=  \psi_\e \lp\widetriangle{\bfu_\e-\ov\bfgotu_\e(\bfu_\e)}\rp$ in 
$(D_{R_\e}\setminus S_{r_\e})\times (0,L)$.
Applying the last line of   \eqref{triangle} to $E_\e:= D_{R_\e}\setminus D_{R_\e/2}$, taking    \eqref{slicing}, \eqref{inttrileint}, 
 (\ref{filower}) and \eqref{applig00}  into account, we infer 
\begin{equation}
\begin{aligned}
   \int_{(D_{R_\e}\setminus S_{r_\e})\times (0,L)}\hskip-0,6cm &  |\bfe_{x'}(\widetriangle\bfu_\e- \wideparen\bfu_\e)|^p  dx
   = 
     \int_{ ( D_{R_\e}\setminus D_{R_\e/2} )  \times (0,L) } \hskip-1cm
    \lb\bfe_{x'}\lp \psi_\e (\widetriangle\bfu_\e-\ov\bfgotu_\e(\widetriangle\bfu_\e))\rp\rb^p dx 
    \\& \leq C  \int_{ ( D_{R_\e}\setminus D_{R_\e/2} )  \times (0,L) }  
    |\bfe_{x'}( \widetriangle\bfu_\e)|^p+  \frac{1}{ R_\e^p} \lb\widetriangle{\ov\bfgotu_\e(\bfu_\e)-\bfu_\e}\rb^p
     dx
      \\&
      \leq     \frac{C}{R_\e^p}  \int_{ \lp D_{R_\e}\setminus D_{R_\e/2} \rp \times (0,L) }   
  | \ov\bfgotu_\e(\bfu_\e)- \bfu_\e  |^p dx+o(1).
  \end{aligned} 
\label{c2}    \end{equation} 
Let $E$ be a bounded Lipschitz  domain of $\RR^2$. 
One can check that  
$
\int_E \big| \psi- {\int\hskip-0,3cm-}_E\psi d\calL^2 \big|^p d\calL^2 \leq C \int_E \lb \nabla \psi  \rb^p  d\calL^2$ \ $\forall \psi \in W^{1,p}(E)$.
   By making suitable changes of variables,  
 we infer  
 $$
 \int_{E_{R_\e}^{i}} \Big| \psi- \intb_{E_{R_\e}^{i}} \psi d\calL^2 \Big|^p d\calL^2\leq C R_\e^p \int_{E_{R_\e}^{i}} \lb \nabla \psi  \rb^p  d\calL^2 \quad  \forall i\in I_\e,\ \forall \psi \in W^{1,p}\Big( E_{R_\e}^{i}\Big).
 $$
Applying this to  $(E, \psi(.))= (D\setminus \frac{1}{2} D,  \bfu_\e(.,x_3))$,    summing w.r.t. $i$ over $ I_\e$ and integrating w.r.t. $x_3$ over $(0,L)$, 
 taking  \eqref{defovgotue}  and \eqref{slicing}  into account, 
 we deduce 
\begin{equation}
\begin{aligned}
     \int_{ \lp D_{R_\e}\setminus D_{R_\e/2} \rp \times (0,L) }  \hskip-0,5cm | \bfu_\e-\ov\bfgotu_\e(\bfu_\e)|^p dx
     \leq C R_\e^p   \int_{ \lp D_{R_\e}\setminus D_{R_\e/2} \rp \times (0,L) }  \hskip-0,5cm | \bfnabla \bfu_\e  |^p dx =  o(R_\e^p). 
  \end{aligned} 
\nonumber\end{equation} 
Combining this with  \eqref{applig22}  and  \eqref{c2}   the proof of \eqref{bound1} is achieved. \qed
\noindent  {\it Proof of \eqref{bound2}.}  If $1<p<2$, by \eqref{growthp},   (\ref{filower}) and \eqref{wideparrenue=}, we have 
\begin{equation}
\begin{aligned}
 \int_{(S_{r_\e} \setminus  S'_{r_\e}) \times (0,L)  }\hskip-1,5cm &  \hskip1,5cm f^{\infty, p}  (\bfe_{x'} (\wideparen  \bfu_\e)) dx
  \leq  C \int_{T_{r_\e }}  
   \lb  \bfe_{x'}\lp \widetriangle\bfu_\e + \zeta_\e \lp  \ov\bfgotr^S_\e(\widetriangle\bfu_\e)    -\widetriangle\bfu_\e\rp\rp\rb^p dx.
  \end{aligned} 
\nonumber
 \end{equation} 
Noticing that, by  \eqref{defovgotA},  $\bfe_{x'}(\ov\bfgotr^S_\e(\widetriangle\bfu_\e))=0$ in $T_{r_\e}$,  
we deduce from   \eqref{defme},   (\ref{estimrigidprim}),  (\ref{estimbasic1}), and  the last line of \eqref{triangle}  applied to $E_\e=S_{r_\e}$,   that
\begin{equation}
\begin{aligned}
  \int_{T_{r_\e }}  &
   \lb  \bfe_{x'}\lp \widetriangle\bfu_\e + \zeta_\e \lp  \ov\bfgotr^S_\e(\widetriangle\bfu_\e)    -\widetriangle\bfu_\e\rp\rp\rb^p dx
  \\& \hskip0,5cm   \leq    C  \int_{T_{r_\e }} 
   \lb  \bfe_{x'}\lp \widetriangle\bfu_\e \rp\rb^p dx + C    \int_{T_{r_\e }}
\frac{1}{ r_\e^{p}} \lb  \ov\bfgotr^S_\e(\widetriangle\bfu_\e)   -\widetriangle\bfu_\e \rb^p dx 
    \\&  \hskip0,5cm  \leq    C  \int_{T_{r_\e }} 
   \lb \bfe_{x'}\lp \widetriangle\bfu_\e \rp\rb^p dx
   \leq      C \frac{r_\e^2  }{ \e^2}  \int   \lb\bfe(\bfu_\e)  \rb^p dm_\e +o(1)=  o(1).
  \end{aligned} 
\label{c51}
 \end{equation} 
 The estimate \eqref{bound2} is proved.\qed

%
\noindent  {\it Proof of \eqref{bound3}.}  If $p=2$, setting $G_{r_\e}:= (S_{r_\e} \setminus  S'_{r_\e})  \times (0,L)$, we infer from \eqref{defovgotA},    \eqref{filower},  \eqref{wideparrenue=}, and \eqref{c51} (which also holds for $p=2$)  that  %
  \begin{equation}\begin{aligned}  
 & \int_{(S_{r_\e} \setminus  S'_{r_\e})\times (0,L)  } f^{\infty, 2}  (\bfe_{x'}(\wideparen  \bfu_\e)) dx
\leq C  \int_{G_{r_\e}}\lb \bfe_{x'}\lp \widetriangle\bfu_\e + \zeta_\e \lp  \ov\bfgotv^S_\e(\widetriangle\bfu_\e)  -\widetriangle\bfu_\e\rp\rp\rb^2dx
   \\&  \leq C  \int_{G_{r_\e}} \hskip-0,1cm\lb  \bfe_{x'}\lp \widetriangle\bfu_\e + \zeta_\e \lp  \ov\bfgotr^S_\e(\widetriangle\bfu_\e)    -\widetriangle\bfu_\e\rp \hskip-0,1cm\rp \hskip-0,05cm\rb^2  \hskip-0,1cm \hskip-0,1cm
   +   \lb  \bfe_{x'} \hskip-0,1cm\lp   \hskip-0,1cm \zeta_\e \frac{2}{\diam S}  \ov\theta_\e^{S} (\widetriangle\bfu_\e) \bfe_3   \hskip-0,1cm\wedge   \hskip-0,1cm \frac{\bfy_\e(x')}{r_\e}  \hskip-0,1cm
        \rp\rb^2 \hskip-0,1cmdx
    \\&    \leq C \int_{G_{r_\e }}  
   \lb  \bfe_{x'}\lp   \zeta_\e \frac{2}{\diam S}  \ov\theta_\e^{S} (\widetriangle\bfu_\e) \bfe_3  \wedge   \frac{\bfy_\e(x')}{r_\e} 
        \rp\rb^2 dx
+o(1).
    \end{aligned}  \nonumber
    \end{equation}
 Noticing 
 that $\bfe_{x'}\lp  \ov\theta_\e^{S} (\widetriangle\bfu_\e) \bfe_3  \wedge   \frac{\bfy_\e(x')}{r_\e} \rp=0$ in $G_{r_\e}$, $ \ov\theta^S_\e(\widetriangle\bfu_\e)$ takes constant values in each set $S_{r_\e}^{i}\times \{x_3\}$, 
 and  $\frac{|G_{r_\e}|}{|T_{r_\e}|}=\frac{|S\setminus S'|}{|S|}$, taking \eqref{defme} and  \eqref{filower} into account, we  deduce
\begin{equation}\begin{aligned}    &\int_{G_{r_\e }}  
 \lb  \bfe_{x'}\lp   \zeta_\e \frac{2}{\diam S}  \ov\theta_\e^{S} (\widetriangle\bfu_\e) \bfe_3  \wedge   \frac{\bfy_\e(x')}{r_\e} 
        \rp\rb^2 dx 
   \leq   C\frac{1}{r_\e^{2}}  \int_{G_{r_\e }}  
   \lb     \ov\theta^S_\e(\widetriangle\bfu_\e)  \rb^2dx
 \\&  \hskip 1cm \leq    C\frac{1}{r_\e^{2}}  \frac{|G_{r_\e}|}{|T_{r_\e}|} \int_{T_{r_\e }}   \lb     \ov\theta^S_\e(\widetriangle\bfu_\e)  \rb^2dx
  \leq   C   \frac{1}{\e^2}    \frac{|S\setminus S'|}{|S|}\int  \lb     \ov\theta^S_\e(\widetriangle\bfu_\e)  \rb^2dm_\e.
 \end{aligned}   \nonumber
    \end{equation}
 By  
  \eqref{estimgotw},  \eqref{estimovA}, 
  \eqref{estimbasic1}, \eqref{inttrileint},  and  \eqref{gottri=trigot},  we have, recalling  that  $\theta_\e(\bfu_\e)=\gotw_{\e3}(\bfu_\e)$ 
  and $c_\e=\e^2$ if $p=2$,  
 \begin{equation}\begin{aligned}          \int  \lb     \ov\theta^S_\e(\widetriangle\bfu_\e)  \rb^2 dm_\e
  &=         \int  \lb     \widetriangle{\ov\theta^S_\e(\bfu_\e)}  \rb^2 dm_\e\leq C    \int  \lb      \ov\theta^S_\e(\bfu_\e)   \rb^2 dm_\e
\\&  \leq   C        \int  \lb      \theta_\e(\bfu_\e)  \rb^2 dm_\e 
 + C         r_\e^2 \int  \lb     \bfe(\bfu_\e)  \rb^2 dm_\e 
  \\&  \leq C          r_\e^2 \int  |\bfe ( \bfu_\e)|^2dm_\e 
+ {C
\e^2} \int_{{T_{r_\e}}}   |\bfnabla \bfu_\e|^2dx   \leq C   \e^2   .
  \end{aligned} 
  \nonumber
    \end{equation}
Combining  the above estimates, we obtain \eqref{bound3}.   \qed

\subsection{Proof of Proposition \ref{propcapreRep=2}}\label{secprooflemA}
When $p=2$,
 by  \eqref{hypRr}, \eqref{capfcapfinfty}, \eqref{bdcfSrR}  and \eqref{defcfSrR}, we have 
 \begin{equation}
\begin{aligned}
&\lb c^{f,S}_{r, R_r}(\bfa)-c^{f^{\infty,2},S}_{r, R_r}(\bfa)\rb
  \le   C|\log r|  R_r^2+ C\lp |\log r|R_r^2\rp^{\frac{\varsigma}{2}}|\bfa|^{2-\varsigma}  = o(|\bfa|),
\end{aligned}
\label{c0fc0finfty}
\end{equation}
hence $c_0^f=c_0^{f^{\infty,2}}$ and  it suffices to prove  Proposition \ref{propcapreRep=2}   when $f$ is $2$-positively homogeneous, which is  assumed in the sequel. 
We  first  establish    the $\T$-relative compactness of $(c^{f,S}_{r,R_r})_{r>0}$. 
\begin{lemma}\label{lemascoli}
Under \eqref{hypRr}, the sequence $(c^{f,S}_{r,R_r})_{r>0}$ defined by \eqref{defcfSrR} is uniformly equicontinuous  on the compact  subsets of $\RR^3$, 
   hence   $\T$-relatively compact.
   Its cluster points are 
 convex  and satisfy \eqref{growthpcf0}. 
\end{lemma}

\noindent
{\bf Proof.} By  \eqref{bdcfSrR} and Lemma \ref{lemconvex},  any cluster point of $(c^{f,S}_{r,R_r})$ is  convex and satisfies \eqref{growthpcf0}.
We prove below that 
   \begin{equation}\begin{aligned} 
 |c^{f,S}_{r, R_r}(\bfa+ \bfh)-c^{f,S}_{r, R_r}(\bfa)|
  \leq C |\bfh|     \lp 1+|\bfa| \rp,
   \end{aligned}
\label{estascoli}
  \end{equation} 
hence  $(c^{f,S}_{r,R_r})_{r>0}$ is uniformly equicontinuous  on the compact  subsets of $\RR^3$, thus, by 
  Ascoli's Theorem  and Cantor's diagonal process,
$\T$-relatively compact. We turn to the proof of \eqref{estascoli}:
by 
\eqref{croissance},    \eqref{rel3}  and \eqref{defcfSrR}, for all $r,R,R',t>0$  and 
  all  bounded   Lipschitz  domains $S_1$,  $S_2$  of $\RR^2$ such that  $R\le R'$,  $\ov S_1 \subset  \ov S_2 \subset RD $, 
  $\bfy_{S_1}\hskip-0,05cm=\hskip-0,05cm\bfy_{S_2}\hskip-0,05cm=\hskip-0,05cm0$,
  $tS\subset RD$, and  $S\subset RD$,  we   have
   \begin{equation}\begin{aligned} 
  &  c^{f,S}_{r,R'} \le c^{f,S}_{r,R}  ; \quad 
   c^{f,S_1}_{r,R}\leq c^{f,S_2}_{r,R};\quad  
  c^{f,tS}_{r,R} = c^{f,S}_{r,R/t}  
  ;\quad  
c^{f,S}_{tr,tR} = \tfrac{\log(rt)}{\log r}c^{f,S}_{r,R}.
 \end{aligned}
\label{propcf}
  \end{equation} 
Let   $\varphi_r$   be the radial function defined   on $R_r D$  in polar coordinates by 
  $\varphi_r(\rho,\theta):= \psi_r(\rho) $, where $\psi_r$ is the solution to
  \begin{equation}
 \begin{aligned}
  \inf_{\psi\in H^1(0,R_r)}\la \int_{r}^{R_r} |\psi'(\rho)|^2 \rho d\rho,\  \ \psi =1\ \hbox{ in } (0,r), \ \psi(R_r)=0 \ra\lp = \frac{1}{\log\frac{R_r}{r}}\rp. 
  \end{aligned}
\nonumber 
  \end{equation} 
 One can check that 
 \begin{equation}
\begin{aligned}
&\varphi_r\in  H^1_0(R_r D), \quad \varphi_r=1 \ \hbox{ in } \ rD, \quad \int_{R_r D} |\bfnabla\varphi_r|^2 dx \leq \frac{C}{|\log r|},
\\&\int_{R_r D} |\bfe(\bfh \varphi_r)| dx \leq C \frac{|\bfh|R_r}{\sqrt{|\log r|}}, \quad \int_{R_r D} |\bfe(\bfh \varphi_r)|^2 dx  \leq C \frac{|\bfh|^2}{|\log r|}\quad \forall \bfh\in\RR^3.
\end{aligned}
\label{fir}
\end{equation}
Let $\bfpsi_{\bfa,r}$ be  a solution  to 
 $\calP^f(\bfa,0; rS, R_r D)$.
Then  $\bfpsi_{\bfa,r} + \bfh \varphi_r \in \W^2(\bfa+\bfh,0, rS, R_r D)$ and, by  \eqref{growthp}, \eqref{inconv}, \eqref{hypRr},  \eqref{bdcfSrR}, and \eqref{fir}, 
 if $|\bfh|\leq C$, 
 \begin{equation}
\begin{aligned}
&  c^{f,S}_{r, R_r}(\bfa+\bfh)-c^{f,S}_{r, R_r}(\bfa)  \leq |\log r| \int_{R_r D} \hskip-0,1cm \hskip-0,1cm f(\bfe(\bfpsi_{\bfa,r} + \bfh \varphi_r ))- f(\bfe(\bfpsi_{\bfa,r})) dy
    \\& \leq C|\log r| \int_{R_r D} |\bfe(\bfh\varphi_r)|  +  |\bfe(\bfh\varphi_r)| |\bfe(\bfpsi_{\bfa,r})|+ |\bfe(\bfh\varphi_r)|^2 dy 
   \\& \leq C|\log r| \lp \frac{|\bfh|  }{\sqrt{|\log r|}} \lp R_r  +   \sqrt{ \int_{R_r D} f(\bfe(\bfpsi_{\bfa,r})) dy}\rp + \frac{|\bfh|^2}{|\log r|} \rp
       \\& \leq C  |\bfh|     \lp 1 + \sqrt{|\log r|\capsca^f(\bfa,0;rS,R_r D)} +|\bfh| \rp
 \leq C |\bfh|     \lp 1+|\bfa| \rp.
  \end{aligned}
\nonumber
\end{equation}
Applying the same argument to $(\bfa', \bfh')\hskip-0,05cm := \hskip-0,05cm(\bfa+\bfh, -\bfh)$, we deduce \eqref{estascoli}.
  \qed 
\\
\noindent By Lemma \ref{propcapreRep=2}, Proposition \ref{propcapreRep=2} is proved if we show that 
$(c^{f,S}_{r, R_r})$ has a unique $\T$-cluster point as $r\to0$,  independent of $(R_r)_{r>0}$ and $S$.
  We achieve this
task  by making use of 
 the proof of Theorem \ref{th}. Setting 
\begin{equation}
  \e:= \frac{1}{\sqrt{|\log r|}},\quad
r_\e=\exp\lp-\frac{1}{\e^2} \rp(=r), \quad l_\e:= \frac{\e^2}{r_\e^5},
\label{re=}
\end{equation}
we consider the functional  $F_\e$   given  by 
\eqref{Pe}. By    \eqref{kkappa}, \eqref{defgammae} and \eqref{re=}, we have
  \begin{equation}
\kappa=+\infty \quad \hbox{ and } \quad \gamma_\e^{(2)}(r_\e) =1 \quad \forall \e>0.
\label{kgamma}
\end{equation}
Let   $\I^f$  be  the weak $\Gamma$-lower limit  of $(F_\e)$
 in  $W^{1,2}_b(\Omega;\RR^3)$, that is the functional defined by (see \cite{Da})
  \begin{equation}\begin{aligned}
&
\I^f(\bfvarphi):= \inf_{(\bfvarphi_\e)\in \C_{\bfvarphi} } \liminf_{\e\to0} F_\e(\bfvarphi_\e),
\\&  \C_{\bfvarphi}:= \la
(\bfvarphi_\e)\subset W^{1,p}_b(\Omega;\RR^3) \lb
\begin{aligned} &\bfvarphi_\e\rightharpoonup  \bfvarphi  \ \hbox{weakly  in } W^{1,p}_b(\Omega;\RR^3) 
\\&\sup_{\e>0} F_\e (\bfvarphi_\e)<+\infty
\end{aligned} \right. 
\ra.
\end{aligned}
 \label{defI}
  \end{equation}
We set 
   \begin{equation}\begin{aligned}
&
\C:=\la  c\in C(\RR^3), \ 
 \int_\Omega c(\bfvarphi)dx=
   \I^f(\bfvarphi)- \int_\Omega  f(\bfe(\bfvarphi)) dx   \quad \forall \bfvarphi \in \D(\Omega;\RR^3)\ra.
\end{aligned}
 \label{C}
  \end{equation}
We will prove that every cluster point of $(c^{f,S}_{r,R_r})$  as $r\to0$ belongs to $\C$. Hence, by  Lemma \ref{lemascoli}, $\C \not=\emptyset$. The next  implication, easily proved by contradiction, 
  \begin{equation}\begin{aligned}
&
 \lc  c_1,c_2\in C(\RR^3), \ 
 \int_\Omega c_1(\bfvarphi)dx= \int_\Omega c_1(\bfvarphi)dx   \quad \forall \bfvarphi \in \D(\Omega;\RR^3)\rc
 \Longrightarrow c_1=c_2,
 \end{aligned}
\nonumber
  \end{equation}
shows  that  $\C$ has a unique  element 
$c_0^f$. It follows that  the whole sequence $(c^{f,S}_{r,R_r})$ $\T$-converges  
 to $c_0^f$   as $r\to 0$, for every  sequence $(R_r)_{r>0}$ verifying  \eqref{hypRr}. 
 To prove that    $c_0^f$  is 
  independent of $S$, let us  fix $\lambda>1$   such that $D\subset S\subset \lambda D$. 
By what precedes,
 $ (c^{f,D}_{r,R_r/\lambda})$ and $(c^{f,D}_{r,R_r})$  
 both  $\T$-converge  to  some  $c_0^{f,D}\in C(\RR^3)$.
 By passing to the limit  as $r\to 0$  in the   inequalities
   $ c^{f,D}_{r,R_r}\leq c^{f,S}_{r,R_r}\leq  c^{f,\lambda D}_{r,R_r}=  c^{f,D}_{r,R_r/\lambda} $, deduced from   \eqref{propcf},     we infer   $ c_0^f= c_0^{f,D}$,
  thus $c_0^f$  is 
  independent of $S$.
The proof of Proposition \ref{propcapreRep=2} is achieved  provided we establish the following  lemma:
\begin{lemma}\label{lemA} Every $\T$-cluster point of $(c^{f,S}_{r,R_r})$  as $r\to0$
  belongs to $\C$.
  \end{lemma}
  
\noindent{\bf  Proof.}   
Let   $c$ be   such  a cluster point and  $(r_k)_{k\in \NN}$   a  decreasing sequence   converging to $0$   such that  $r_1<1$ and 
\begin{equation}\begin{aligned}
\hbox{ $c^{f,S}_{r_k,R_{r_k}}$ $\T$-converges to $c$.}
  \end{aligned}
  \label{Tcv}
  \end{equation}
Then, the sequence  $(\e_k)_{k\in \NN}$ defined by 
 \begin{equation}\begin{aligned}
  \e_k:= \frac{1}{\sqrt{|\log r_k|}},
  \end{aligned}
  \label{ek}
  \end{equation}
is decreasing and, by  \eqref{re=}, 
  \begin{equation}\begin{aligned}
  r_k= r_{\e_k} \qquad \forall k\in \NN.
  \end{aligned}
  \label{rkek}
  \end{equation}
  We  prove  below that for all $ \bfvarphi \in \D(\Omega;\RR^3)$, there exists   
$(\bfvarphi_\e)_{\e>0}\in \C_{\bfvarphi }$ such that
  \begin{equation}\begin{aligned}
 \limsup_{k\to+\infty} F_{\e_k}(\bfvarphi_{\e_k})
   \le\hskip-0,1cm \int_\Omega \hskip-0,1cm f(\bfe(\bfvarphi))  
  +   c(\bfvarphi) dx.
\end{aligned}
\label{lsupek}
  \end{equation}
Since $ \liminf_{\e\to0} F_\e(\bfvarphi_\e)\le  \limsup_{k\to+\infty} F_{\e_k}(\bfvarphi_{\e_k})$, we deduce from   \eqref{defI}    that
   \begin{equation}\begin{aligned}
 \I^f(\bfvarphi)
\leq 
 \int_\Omega   f(\bfe(\bfvarphi))  
  +c(\bfvarphi) dx \qquad \forall  \bfvarphi \in \D(\Omega;\RR^3).
\end{aligned}
\label{linfek0}
  \end{equation}
Next, we  establish   that  
   \begin{equation}\begin{aligned}
\hskip-0,2cm   \I^f(\bfu)
 \ge \int_\Omega   f(\bfe(\bfu))  
  +c(\bfu) dx   \qquad  \forall \bfu \in \D(\Omega;\RR^3).
\end{aligned}
\label{linfek}
  \end{equation}
\noindent  
We infer  from  \eqref{linfek0}, \eqref{linfek}, and 
  the arbitrariness of  $\bfvarphi$, $\bfu$  that  $c\in \C$.
It remains to exhibit   $(\bfvarphi_\e)\in \C_\bfvarphi$  satisfying   \eqref{lsupek}
and to prove\eqref{linfek}. 

 \noindent  $\bullet$  Given $ \bfvarphi \in \D(\Omega;\RR^3)$, the sequence   $(\bfvarphi_\e)\in \C_\bfvarphi$  verifying   \eqref{lsupek} 
 will be  deduced from \eqref{defvarphie} (upperbound) 
 by specifying the choice of  
 $R_\e$ in \eqref{defetaa},  namely  by setting 
 \begin{equation}
\begin{aligned}
  &
  R_\e :=   R_\e^{(2)},
  \end{aligned} 
\label{Relimsup2}    \end{equation} 
where,  
fixing     $\lambda>1$  such that     $ S'\subset\lambda S$ (see \eqref{S'sup}),
\begin{equation}\begin{aligned}
&R_\e^{(2)}:= \lambda R_\e^{(1)}, \qquad  R_\e^{(1)}:= \frac{r_\e}{r_k} R_{r_k} \qquad \forall   k\in \NN, \ \forall \e\in ]{\e_{k+1}},{\e_k}].
\end{aligned}
\label{defRe1}
  \end{equation}
 By  \eqref{S'sup},  \eqref{propcf} and  \eqref{rkek}, we have 
%
   \begin{equation}\begin{aligned}
   &c^{f,S'}_{r_k ,R_{\e_k}^{(2)}}\leq c^{f,\lambda S}_{r_k ,R_{\e_k}^{(2)}}
  = c^{f,S}_{r_k ,R_{\e_k}^{(2)}/\lambda}
 = c^{f,S}_{r_k ,R_{\e_k}^{(1)}} = c^{f,S}_{r_k ,R_{r_k}}  \qquad \forall   k\in \NN.  
      \end{aligned}
\label{nnn}
  \end{equation} 
We check below that 
 \begin{equation}\begin{aligned}
 &
 (\e, r_\e,R_\e^{(1)}) \hbox{ verifies  \eqref{Re}};\quad 
  c^{f,S}_{r_k,R_{r_k}}\leq c^{f,S}_{r_\e,R_\e^{(1)}}\     \forall \e\in ]{\e_{k+1}},{\e_k}], \ \forall   k\in \NN,
    \end{aligned}
 \label{propRe1} 
  \end{equation}
  thus  $R_\e$,  defined by  \eqref{Relimsup2}, satisfies  \eqref{Re}.
We consider   the decomposition \eqref{splitsup0}.   
By  \eqref{defD}, \eqref{psitupleinD}, \eqref{defchiothercases},  and \eqref{kgamma}, we have
  \begin{equation}
\begin{aligned}
  &\bfpsi^{tuple}=\bfchi_\e=0.
  \end{aligned}
\label{psichi=0}
  \end{equation}
By  \eqref{splitsup0},  \eqref{limI1},  \eqref{limI4} and \eqref{psichi=0},
     \begin{equation}\begin{aligned}
 &\lim_{\e\to0}  I_{\e1}+  I_{\e3}  +  I_{\e4}=  \int_\Omega f(\bfe(\bfvarphi )) dx.
   \end{aligned}
\label{number9}
  \end{equation}
 By  \eqref{defcheckphi} and  \eqref{psichi=0}, 
   $\widecheck \bfphi_\e=-\widecheck \bfvarphi_\e$  and  by    \eqref{lsI23},   \eqref{lsI21},   \eqref{kgamma},  and \eqref{nnn}
%
        \begin{equation}\begin{aligned}
  I_{\e_k 2}& =    
  \int_\Omega c^{ f ,S'}_{r_{\e_k}, R_{\e_k}^{(2)}} (\widecheck \bfvarphi_\e) dx + o(1)
  \le \int_\Omega c^{ f ,S}_{r_k, R_{r_k}} (\widecheck \bfvarphi_\e) dx + o(1).
   \end{aligned}
\label{nnnn}
  \end{equation}
By  \eqref{defcheckphi},  $(\widecheck \bfvarphi_\e )$ is bounded in $L^\infty$ and uniformly converges   
 on each compact subset of $\Omega$   to $\bfvarphi$, 
therefore, by \eqref{Tcv}, 
        \begin{equation}\begin{aligned}
 \limsup_{k\to+\infty}  I_{\e_k 2}& \leq     \lim_{k\to+\infty} 
  \int_\Omega c^{ f ,S}_{r_k, R_{r_k} } (\widecheck \bfvarphi_{\e_k}) dx = \int_\Omega c(\bfvarphi)dx.
   \end{aligned}
\nonumber
  \end{equation}
Taking   \eqref{splitsup0}  and    \eqref{number9}    into account, the assertion  \eqref{lsupek} is proved.

\noindent{\it Verification of   \eqref{propRe1}.}  
 By \eqref{hypRr},  $r_k\ll R_{r_k}$, hence, by \eqref{defRe1},
   $r_\e\ll R_\e^{(1)} $.
By \eqref{re=},   
the mapping  $\e\to \frac{\e}{r_\e}$   is decreasing on $\lp0, \sqrt{2}\rp$   hence
  $R_\e^{(1)}=\frac{r_\e}{r_k} R_{r_k}   \leq    \frac{\e}{\e_k} R_{r_k} $ for all $ \e\in ]\e_{k+1}, \e_k]$
if  $\e_k <\sqrt2$.
Since
 $R_{r_k} \ll \frac{1}{\sqrt{|\log r_k|}}={\e_k}$,
we deduce  that  $R_\e^{(1)} \ll\e
   =\frac{1}{\sqrt{|\log r_\e|}}$ (see \eqref{re=}), thus
   $(\e,r_\e,R_\e^{(1)})$ verifies   \eqref{Re}. 
By  \eqref{propcf}, \eqref{rkek},    \eqref{defRe1},  and the
decrease  of       $\e\in (0,+\infty)\to |\log r_\e|$,  
   %
   \begin{equation}\begin{aligned}
c^{f,S}_{r_\e,R_\e^{(1)}}= c^{f,S}_{\frac{r_\e}{r_k}  r_k ,\frac{r_\e}{r_k} R_{r_k} } 
 =    \frac{|\log r_\e|}{|\log r_k|} c^{f,S}_{r_k,R_{r_k} }  =   & \frac{|\log r_\e|}{|\log r_{\e_k}|} c^{f,S}_{r_k,R_{r_k} } \ge \    c^{f,S}_{r_k,R_{r_k} } 
 \\&\quad \forall k\in \NN, \  \forall \e\in ]\e_{k+1},\e_k].
    \end{aligned}
\nonumber
  \end{equation}
The assertion  \eqref{propRe1} is proved.  

  \noindent $\bullet$ {\it Proof of \eqref{linfek}.}    %
  We fix  $\bfu\in \D(\Omega;\RR^3)$,
  $(\bfu_\e)_{\e>0}\in \C_{\bfu}$, and a subsequence $(\bfu_{\e_l})_{l\in\NN}$    such that 
    \begin{equation}\begin{aligned}
   \liminf_{\e\to0} F_\e(\bfu_\e)
=\lim_{l\to+\infty} F_{\e_l}(\bfu_{\e_l}).
\end{aligned}
\label{el}
  \end{equation}
By \eqref{defI}, $(\bfu_{\e_l})$  verifies  \eqref{supFeuefini},
 hence, up to a subsequence,  every assertions stated in  Propositions \ref{prop1},\ref{prop2},\ref{prop3}.
 Let us fix $S'$ and $c_0$ as in   \eqref{S'}.
By   \eqref{propRe1} and Lemma  \ref{lemslicing}, there exists  $(R_\e) $ such that 
 \begin{equation}\begin{aligned}
(R_\e)  \ \hbox{ verifies \eqref{Re},  \eqref{slicing} and} \quad
\frac{R_\e}{c_0}\leq R_\e^{(1)} \quad \forall \e>0.
\end{aligned}
\label{R3/c0leR1}
  \end{equation}
By    
\eqref{propcf},  \eqref{propRe1} and \eqref{R3/c0leR1},
 \begin{equation}\begin{aligned}
c^{f,S}_{r_k ,R_{r_k} }\leq c^{f,S}_{r_\e,R_\e^{(1)}}\leq c^{f,S}_{r_\e,R_\e/c_0}
\quad \forall k \in \NN, \ \forall \e\in ]{\e_{k+1}},{\e_k}].
  \end{aligned}
 \label{MegeAk2} 
  \end{equation}
By Lemma \ref{lemascoli},  
\begin{equation}\begin{aligned}
\hbox{ $ c^{f,S}_{r_{\e_l},R_{\e_l}/c_0}$ $\T$-converges to $\tilde c$,}
  \end{aligned}
  \label{Tcvtilde}
  \end{equation}
  %
 up to a further  subsequence,
for some $\widetilde c\in C(\RR^3)$ satisfying,  by \eqref{Tcv} and   \eqref{MegeAk2},  $ \widetilde c \ge c.$
Using \eqref{Tcvtilde} in place of \eqref{Tcv0} in the argument of the lower bound,  
noticing that by (\ref{defD}), \eqref{vtupleinD} and (\ref{kgamma}),  $\bfv^{tuple}=0$, we obtain 
   \begin{equation}\begin{aligned}
\lim_{l\to+\infty} F_{\e_l}(\bfu_{\e_l}) \ge  \int_\Omega f(\bfe(\bfu))  + \tilde c(\bfu) dx 
  \ge  \int_\Omega f(\bfe(\bfu))  +   c(\bfu) dx.
      \end{aligned}
\nonumber
  \end{equation}
The assertion  \eqref{linfek} is proved.
 \qed
\end{document}